\documentclass[11pt,twoside]{book}
\usepackage{a4wide}
\usepackage{amssymb}
\usepackage{amsmath}
\usepackage{hhline}
\makeindex
\usepackage[english]{babel}
\usepackage[applemac]{inputenc}
\usepackage{latexsym}
\usepackage{graphicx}
\usepackage{epsfig}
\usepackage{makeidx}
\usepackage{color}
\usepackage{theorem}
\theorembodyfont{\itshape} 
\definecolor{grey}{rgb}{0.6,0.6,0.6}

\newcommand{\grey}[1]{\textcolor{grey}{#1}}
\newcommand{\drawat}[3]{\makebox[0pt][l]{\raisebox{#2}{\hspace*{#1}#3}}}

\newtheorem{Prp}{Proposition}[section]
\newtheorem{Def}[Prp]{\newline Definition}
\newtheorem{Not}[Prp]{\newline Notation}
\newtheorem{Exa}[Prp]{\newline Example}
\newtheorem{Lem}[Prp]{\newline Lemma}
\newtheorem{Cor}[Prp]{\newline Corollary}
\newtheorem{proof}{\newline Proof}
	
\newtheorem{remark}{\newline Remark}
	
\newtheorem{remarks}{\newline Remarks}
	
\newtheorem{remarque}{\newline Remarque}

\newcommand{\proofend}{\hfill $\square$}

\usepackage{amsmath}
\usepackage{amsfonts}
\usepackage{amssymb}
\usepackage{verbatim}
\usepackage[all]{xy}
\xyoption{matrix}

\newcommand{\im}{{\bf i}}
\newcommand{\refd}[1]{\ensuremath{\ref{#1}}}

\newcommand{\styleoftitlesofpages}{\sf}
\newcommand{\markboths}[2]{\markboth{\styleoftitlesofpages #1}{\styleoftitlesofpages #2}}
\newcommand{\pagetitless}[2]{\markboth{\styleoftitlesofpages Chapter \thechapter. #1}{\styleoftitlesofpages \thesection. #2}}

\newcommand{\titrepageresume}[1]{\markboth{\styleoftitlesofpages Résumé de la thèse: Chapitre #1}{\styleoftitlesofpages Résumé de la thèse: Chapitre #1}}
\newcommand{\titrepagesresume}[1]{\markboth{\styleoftitlesofpages Résumé de la thèse: Chapitres #1}{\styleoftitlesofpages Résumé de la thèse: Chapitres #1}}
\newcommand{\hts}{\vphantom{\Big(}}

\newcommand{\h}[1]{\hspace{-#1mm}}
\newcommand{\titreProp}[1]{{\sc : #1} }
\newcommand{\titrePropRef}[2]{{\sc : #1} {\upshape (#2)} }
\newcommand{\PropRef}[1]{{\upshape (#1)}\\ }
\newcommand{\fauxtitre}{{\sc: }\\ }
\newcommand{\fauxtitred}{{\sc: }}
\newcommand{\K}{\mathbb{C}}
\newcommand{\Pn}{\mathbb{P}^2}
\newcommand{\PGLn}[1]{\PGL(#1,\K)}
\newcommand{\PGL}{\mathrm{PGL}}
\newcommand{\GLZ}[1]{GL_{#1}(\mathbb{Z})}
\newcommand{\CrP}{\Cr{\mathbb{P}^2}}
\newcommand{\CrPp}{\Crc{\mathbb{P}^1\times\mathbb{P}^1}{\pi_1}}
\newcommand{\Cr}[1]{\mathit{Cr}(#1)}
\newcommand{\Crc}[2]{\mathit{Cr}(#1,#2)}
\newcommand{\Diag}[3]{[#1:#2:#3]}
\newcommand{\DiaG}[4]{[#1:#2:#3:#4]}
\newcommand{\Sym}{\mathrm{Sym}}
\newcommand{\Z}[1]{{\mathbb{Z}/#1\mathbb{Z}}}

\newcommand{\Pic}[1]{\mathrm{Pic}(#1)}

\newcommand{\rkPic}[1]{\mathrm{rk\ Pic}(#1)}
\newcommand{\num}[1]{\ensuremath{\begin{array}{|c|}\hline {\mathbf{#1}} \\ \hline\end{array}}}
\newcommand{\nump}[1]{\ensuremath{\begin{array}{|c|}\hline {\mathrm{#1}} \\ \hline\end{array}}}
\newcommand{\cadre}[1]{\begin{tabular}{|c|}\hline {{#1}} \\ \hline\end{tabular}}
\newcommand{\defn}[1]{\emph{#1}}
\newcommand{\ChapterBegin}[1]{\emph{#1}}
\newcommand{\ChapterNotation}[1]{\begin{itemize}\item \emph{#1}\end{itemize}}
\newcommand{\Fn}{\mathbb{F}_n}
\newcommand{\Aut}{\mathrm{Aut}}

\newcommand{\btriangle}{\includegraphics[height=0.2cm]{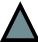}}
\newcommand{\bsquare}{\includegraphics[height=0.2cm]{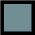}}

\newcommand{\refThmAbelGenusCurve}{5}
\newcommand{\refThmExistenceNonLinear}{1}
\newcommand{\refThmRootsLinearAutomorphisms}{2}
\newcommand{\refThmCastelnuovoGeneralis}{3}
\newcommand{\refThmCyclicGenusCurve}{4}
\newcommand{\refThmNonCyclicGroup}{B}
\newcommand{\refThmIsomorphyclasses}{6}
\newcommand{\refThmCyclicGroup}{A}

\newcommand{\Citation}[2]{\begin{tabular}{p{10 cm}}{\small\tt "#1" \hfill {\it #2}} \end{tabular}}

\newcommand{\AutSp}{\Aut(S,\pi)}
\begin{document}
\frontmatter
\begin{titlepage}
\begin{center}

{\scshape Université de Genève \hfill Faculté des Sciences} \\
{Section de Math\'ematiques \hfill Professeur F. Ronga} \\[-2mm]
\hrulefill

\vskip 35mm

{\huge\scshape Finite abelian subgroups of the  \\[2mm]
Cremona group of the plane \\}

\vskip 25mm

{\scshape Thèse} \\[2mm]
{présentée à la Faculté des sciences
de l'Université de Genève \\
pour l'obtention du grade de Docteur ès sciences,
mention mathématiques \\}

\vskip 7mm

par \\[7mm]
Jérémy BLANC \\
de \\
Pregny-Chambésy (GE) \\

\vskip 30mm

{Thèse N° 3777}

\vfill

{\footnotesize
Genève\\
Atelier de reproduction de la Section de physique \\
Septembre 2006 \\[-2mm]}
\end{center}
\end{titlepage}
\markboths{Finite abelian subgroups of the Cremona group of the plane}{Finite abelian subgroups of the Cremona group of the plane}

\noindent {\sf La Faculté des sciences, sur le préavis de Messieurs F. RONGA, professeur ordinaire et directeur de thèse (Section de mathématiques), T. VUST, docteur et co-directeur de thèse (Section de mathématiques), A. BEAUVILLE, professeur (Université de Nice - Laboratoire J. A. Dieudonné - Nice, France), et I. PAN, professeur adjoint (Universidade Federal do Rio Grande do Sul - Instituto de Matematica - Porto Alegre, Brasil), autorise l'impression de la présente thèse, sans exprimer d'opinion sur les propositions qui y sont énoncées.}

\vspace{ 3 cm}

\noindent {\small \sf Genève, le 14 septembre 2006}

\tableofcontents
\chapter*{Remerciements}
Le travail de th\`ese est long et dur et c'est gr\^ace \`a un climat agréable et un entourage bienveillant qu'on arrive à surmonter les moments difficiles. Il est temps de remercier tous ceux qui m'ont aidé durant ces 4 années à l'université de Genève.

\begin{itemize}
\item
Tout d'abord, le plus important: le directeur de thèse, \emph{Thierry Vust}. C'est lui qui en été 2002 m'a fait découvrir la géométrie algébrique en me proposant un livre de Miles Reid. Puis, il a passé de longues heures à m'expliquer les techniques de ce sujet à premier abord très rude, m'en faisant découvrir la beauté cachée. Il a également sué en lisant mes textes et surtout il a su, à travers de longues discussions mathématiques avec moi, me donner la passion de la recherche mathématique. Pour tout cela, un tout grand merci!
\item
J'ai eu la chance de pouvoir également parler avec d'autres mathématiciens sur le sujet, et non des moindres. Je souhaiterais donc les remercier chacun pour leur aide précieuse.
\begin{itemize}
\item
\emph{Arnaud Beauville}, que j'ai rencontré par e-mail d'abord puis "en chair et en os" à Nice en janvier 2006 et à qui j'ai pu poser de nombreuses questions.
\item
\emph{Ivan Pan}, qui faisait ses vacances d'été en février à Genève, sous la neige et qui m'a fait découvrir la géométrie algébrique "à l'ancienne".
\item
\emph{Felice Ronga}, le directeur de thèse officiel, dont la porte au bout du couloir du 2ème étage était toujours ouverte pout toutes questions de topologie, géométrie, voire même d'informatique.
\item
\emph{Igor Dolgachev}, que j'ai rencontré à Turin en septembre 2005 et qui a travaillé sur le même sujet en parrallèle. Nous avons pu converser ensuite longuement par e-mail, se posant quelques questions précises au coeur du sujet.
\item
\emph{Pierre de la Harpe}, toujours disponible pour des discussions très intéressantes sur la théorie des groupes et les représentations.
\end{itemize}
\item
Ayant fait la folie de ne pas écrire ma thèse dans ma langue maternelle, j'ai du avoir recours à l'aide d'experts anglophones. Je tiens ici à remercier \emph{John Steinig} pour ses multitudes de remarques pertinentes sur la langue de Shakespeare et ses longues relectures de la totalité du texte.
\item
Les mathématiques n'étant pas l'unique ingrédient d'un travail si long, je suis heureux de pouvoir remercier tous ceux qui m'ont suivi pendant ces années.
\begin{itemize}
\item
Merci à \emph{Sophie}, mon épouse qui m'a soutenu pendant tous ces moments qui furent aussi durs pour elle que pour moi.
\item
Merci à \emph{David}, qui même depuis l'Allemagne a montré qu'il était toujours présent et indispensable.
\item
Merci à \emph{Nicolas}, le collègue de bureau inépuisable qui m'a longuement écouté lui racouter du charabia de géométrie algébrique.
\item
Merci à \emph{Luc, Eugenio, Ghislain, Shaula}, toujours présents dans cette bonne vieille section de mathématiques.
\item
Merci à tous les autres (la liste est longue...) pour leur soutien moral. Ma famille, celle de Sophie, mes amis, mes collègues...
\end{itemize}
Enfin, j'aimerais encore dire un mot à \emph{Bernard} qui m'a toujours aidé à trouver les articles les plus farfelus au travers de la bibliothèque et \emph{Martin}, le maître de l'informatique, véritable Hotline-Maple-Matlab dans le bureau d'à côté.
\end{itemize}

\vspace{3 cm}
\begin{flushright}\textit{Many thanks.}\end{flushright}
\begin{flushright}\textit{Jérémy}\end{flushright}
\mainmatter
\chapter{Introduction}
\label{Chap:Introduction}
\ChapterBegin{The aim of this work is to classify the finite abelian subgroups of the Cremona group, up to conjugation. This classification gives rise to many  new results about the Cremona group, and contributes to understanding the properties of this classical group, which has been studied for more than a hundred years.}
\section{An informal introduction to some basic notions}
\label{Sec:InformalIntro}
\pagetitless{Introduction}{An informal introduction to some basic notions}
\ChapterBegin{We introduce our subject for the non-specialist, as well as for the reader who is not familiar with the study of mathematics. The language used will therefore be as simple as possible and we will introduce only a few mathematical symbols.}

A recurrent subject of study in mathematics is to define some spaces, equipped with some structure, and to study the transformations of these spaces that preserve the structure.

\subsection{Isometries of the plane}
\label{Subsec:Isometries}
Take for example the Euclidean plane equipped with the distance between points. A transformation  of the plane that preserves this distance is called an \defn{isometry}. The set of such transformations is called  \defn{the group of isometries}. 

This group contains for example rotations, reflections in a line, and translations (see the figures below):

\begin{center} \includegraphics[width=3.75cm,height=2.92cm]{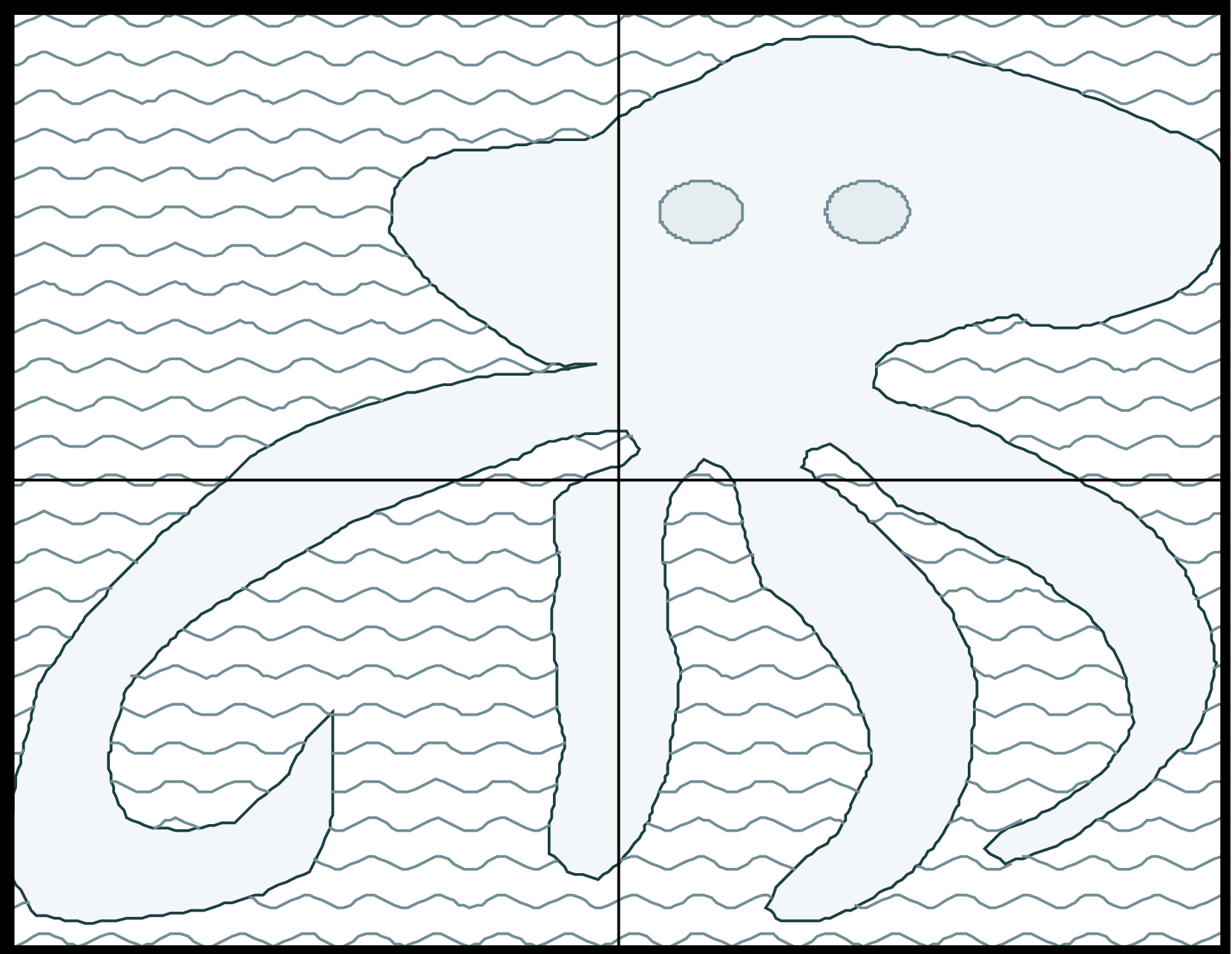} \includegraphics[width=3.75cm,height=2.92cm]{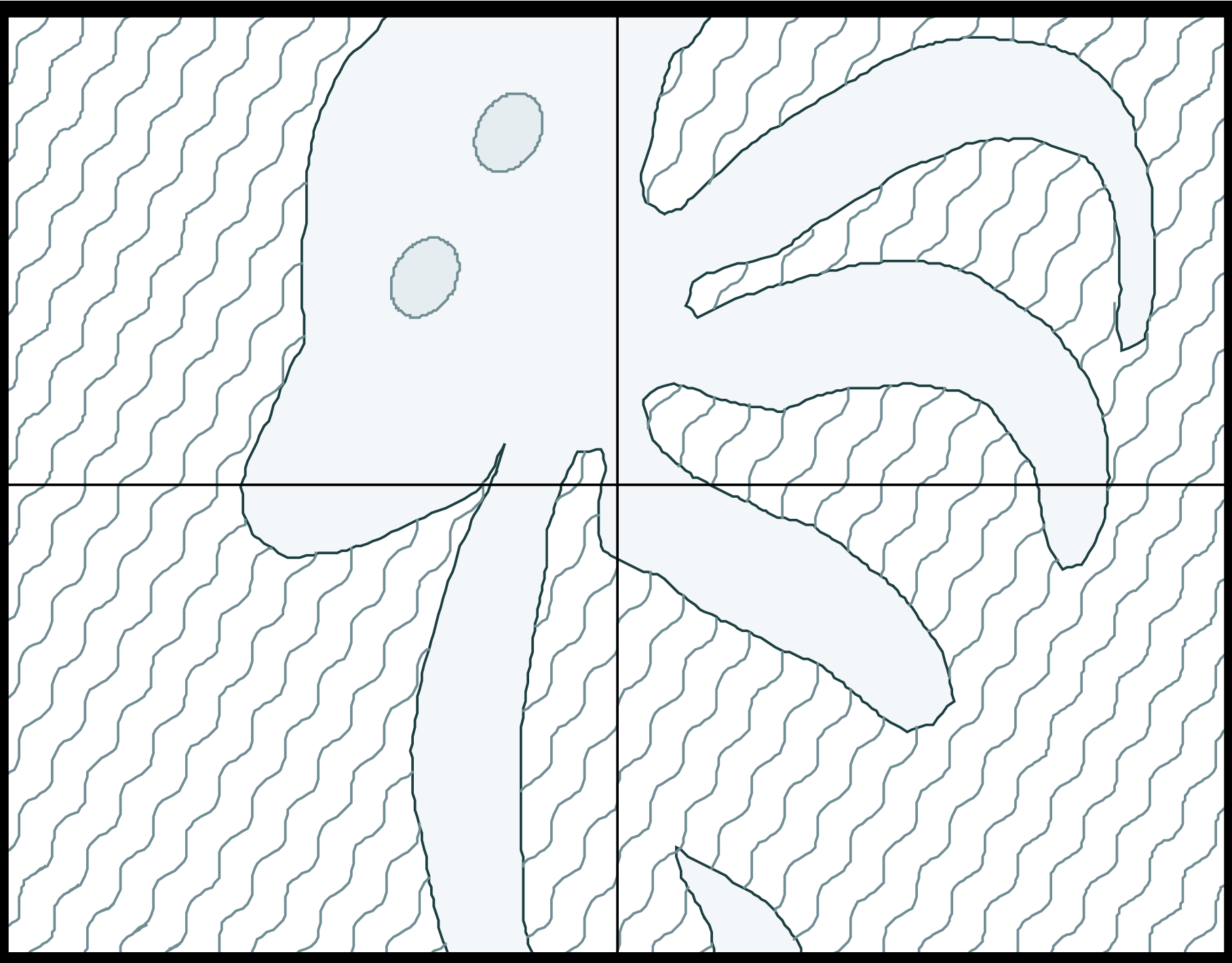} \includegraphics[width=3.75cm,height=2.92cm]{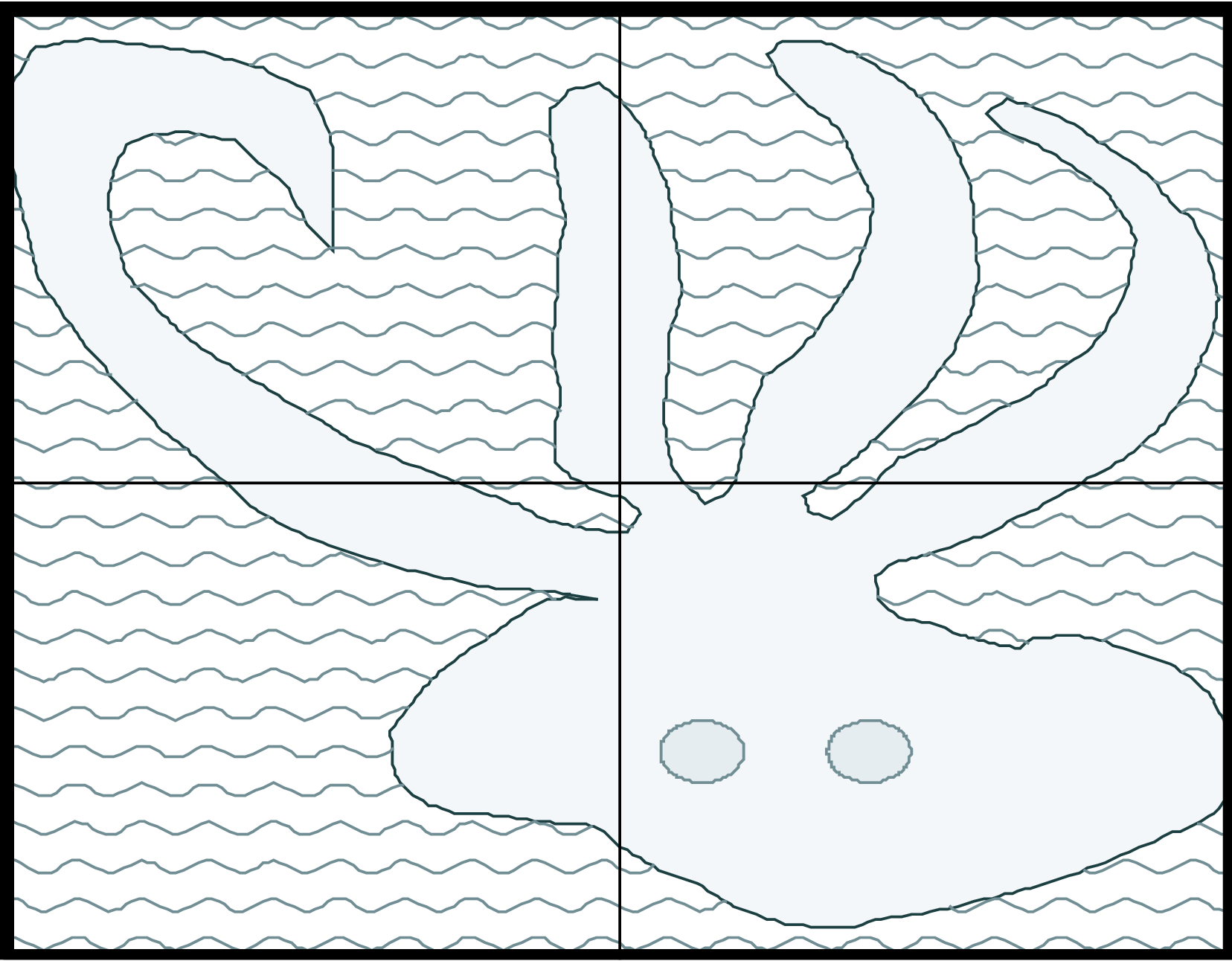}\hspace{0.1 cm}\includegraphics[width=3.75cm,height=2.92cm]{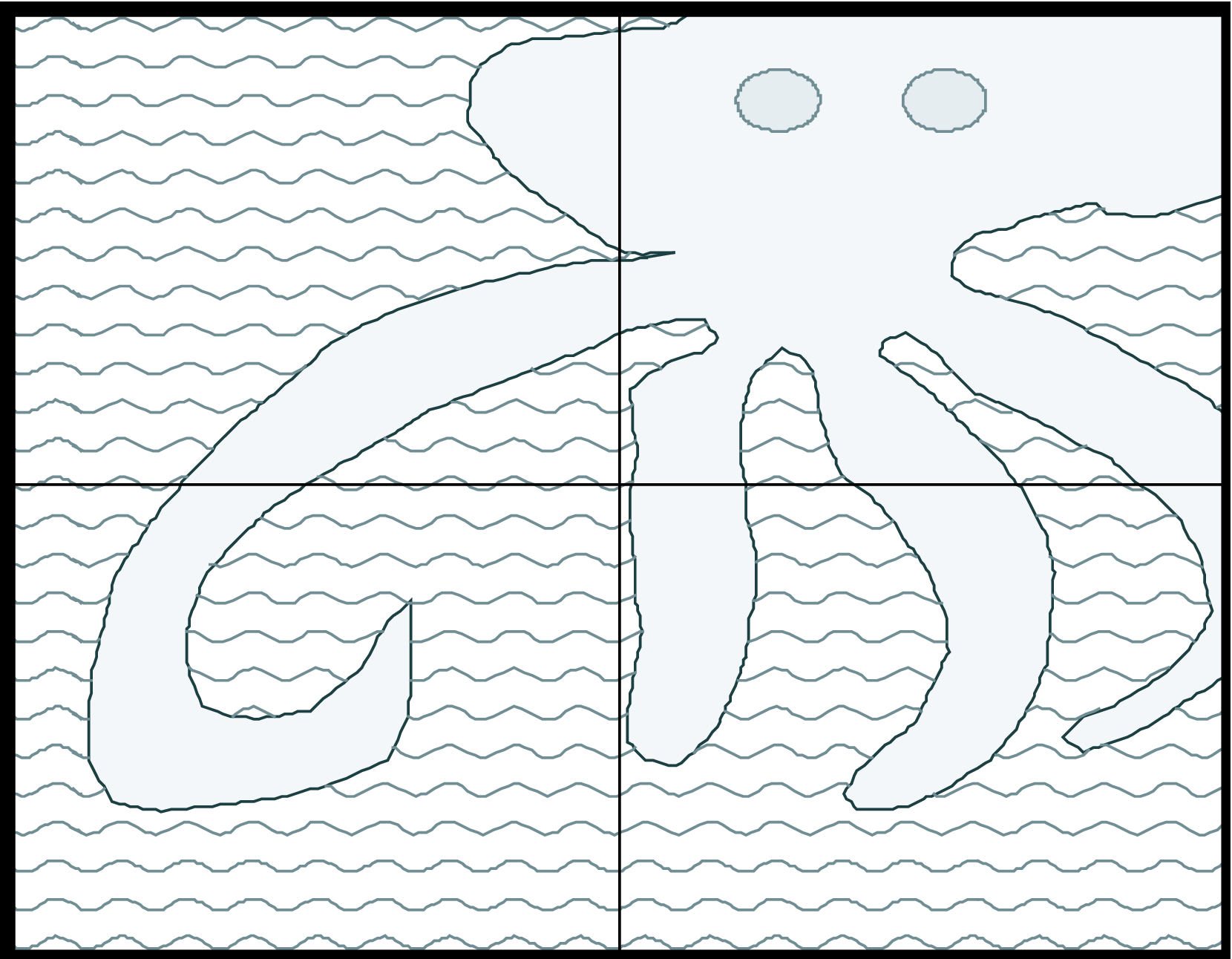}\\
{\it \small On the left, some part of the plane, in which we put an octopus and some waves. Then, its image by respectively a $60$° rotation, a reflection with horizontal axis and a translation.}\end{center}

One of the isometries is the \defn{identity}. This is the unique isometry that sends every point onto itself. 
Take an isometry which is not the identity. If its composition with itself equals the identity, we say that the isometry is an \defn{involution}. For example, all reflections are involutions;  rotations are involutions if and only if their angle is $180$°; no translation is an involution. 

\bigskip

We say that two isometries are \defn{conjugate} if there exists a third isometry which intertwines the first with the second. This amounts to saying that two isometries are conjugate if they coincide after a change of coordinates, which is itself an isometry. 

Consider for example two rotations by the same angle, but centred at two different points. The translation that sends the first of these points to the second conjugates the two rotations. We illustrate this by a diagram:

\begin{center}
\hspace{0 cm}\xymatrix{ \includegraphics[height=3cm]{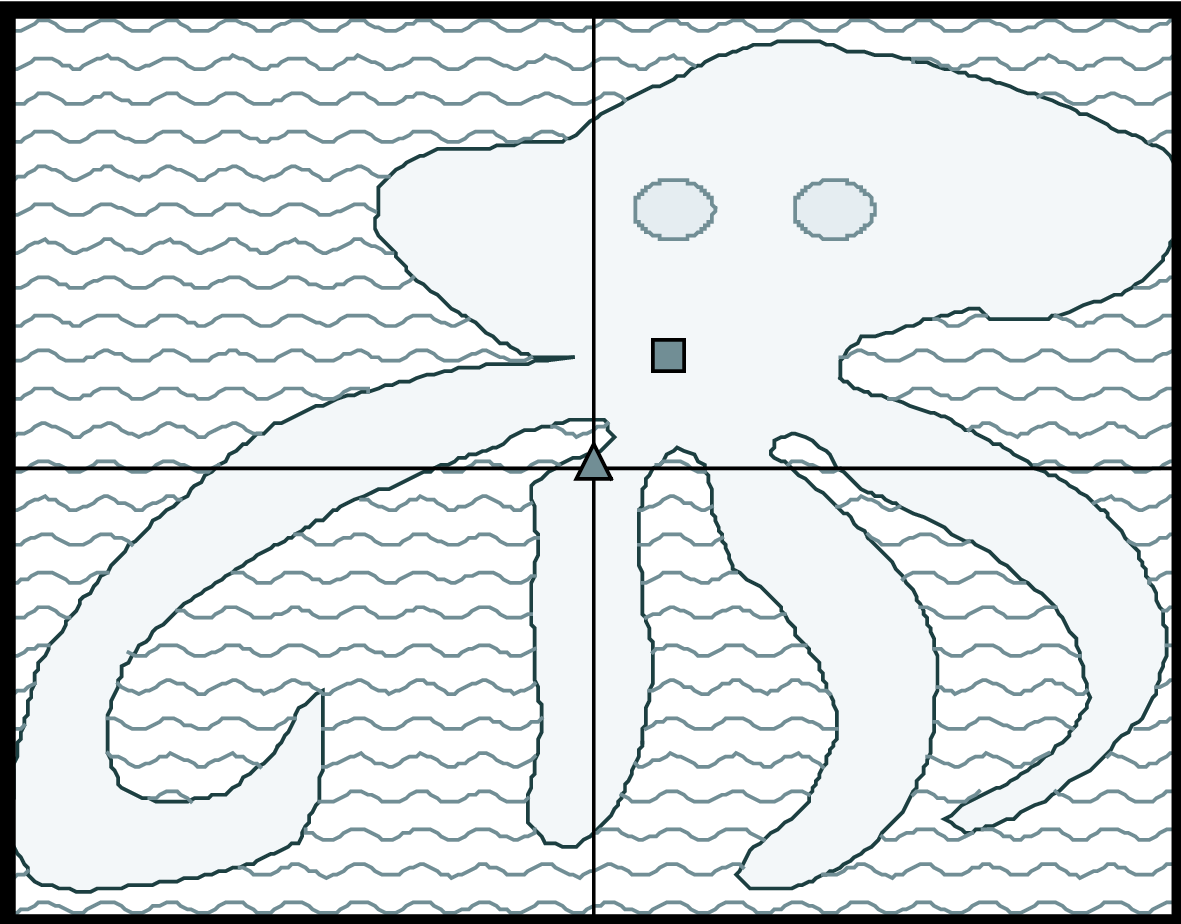}  \ar^{\mbox{\it \small $60$° rotation about }\btriangle}[rrrr] \ar_{\mbox{\it \small translation from }\btriangle \mbox{\it \small\ to }\bsquare}[d] & & & &\includegraphics[height=3cm]{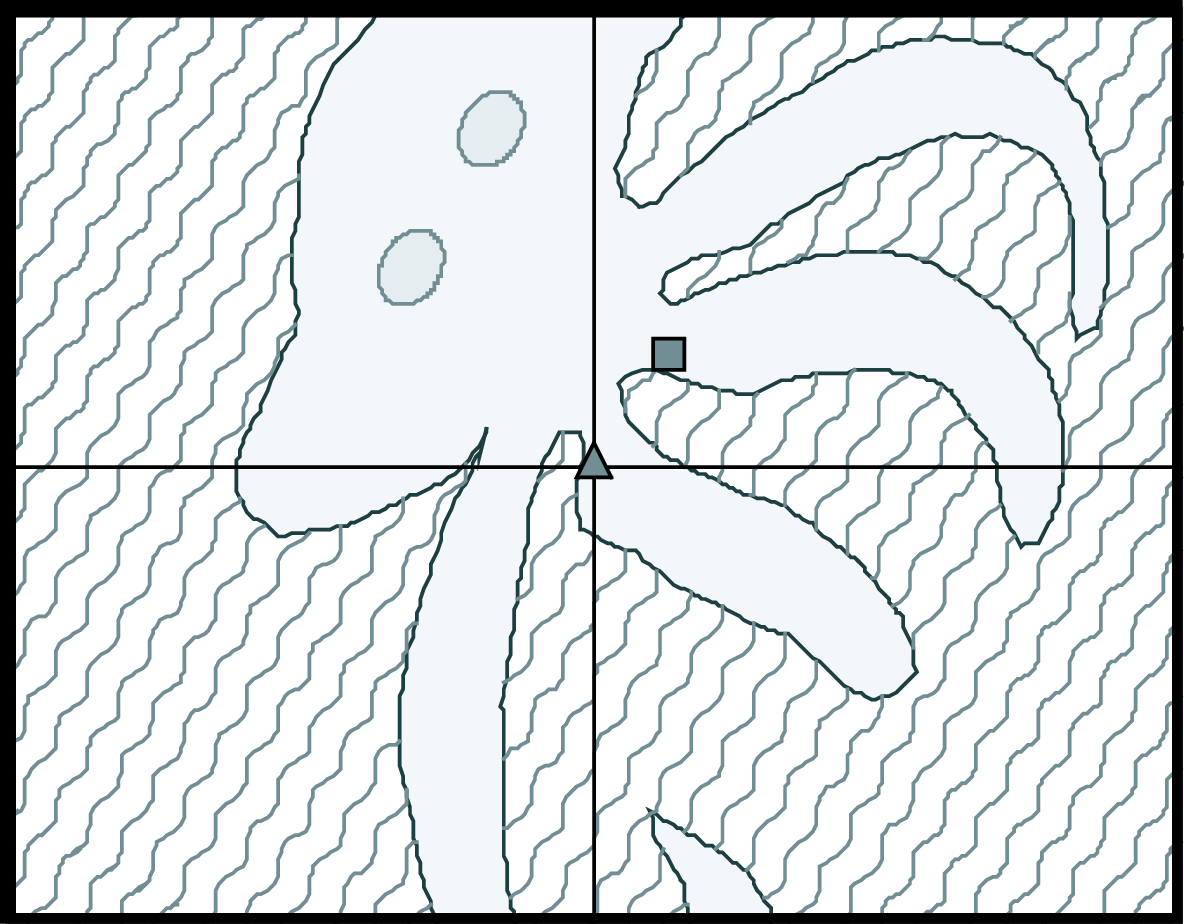} \ar_{\mbox{\it \small translation from }\btriangle \mbox{\it \small\ to }\bsquare}[d] \\
\includegraphics[height=3cm]{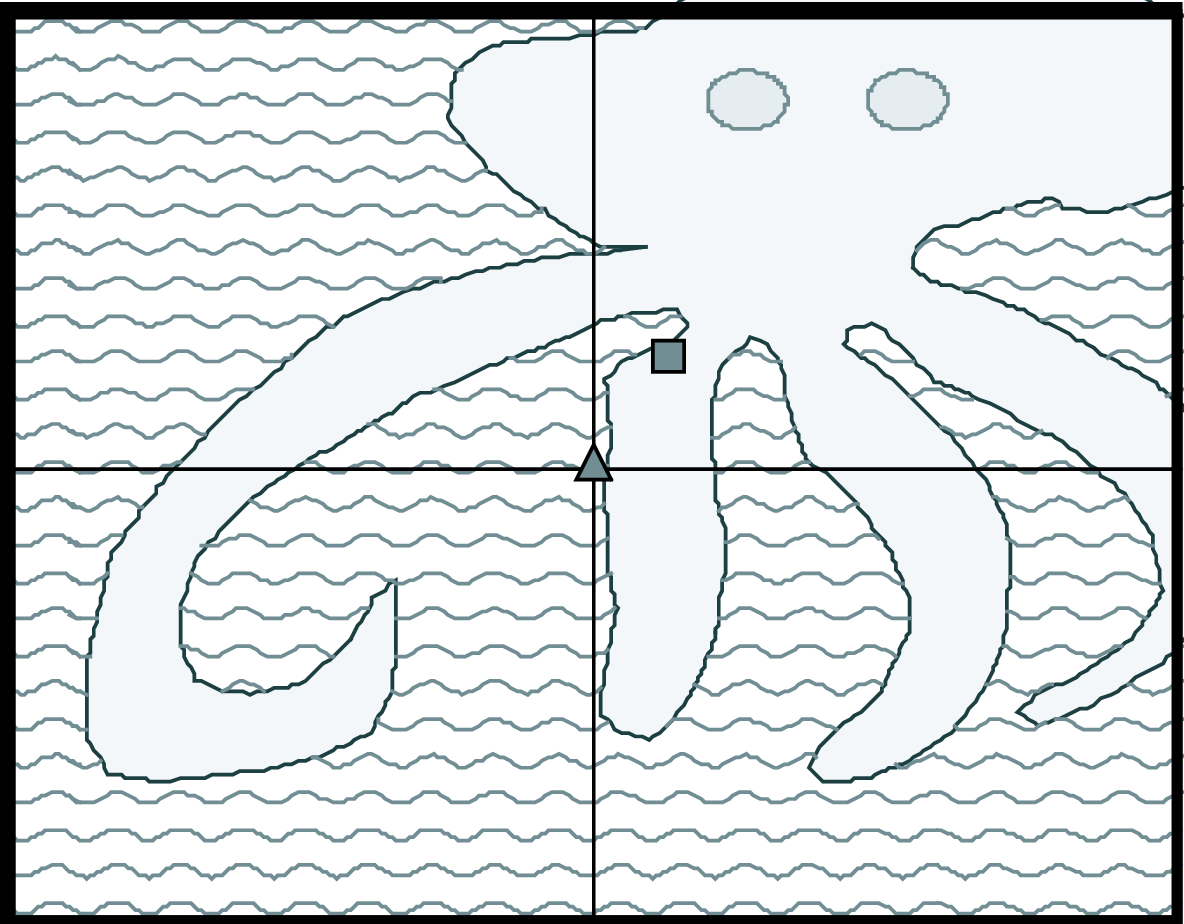} \ar^{\mbox{\it \small $60$° rotation about } \bsquare}[rrrr] & & & & \includegraphics[height=3cm]{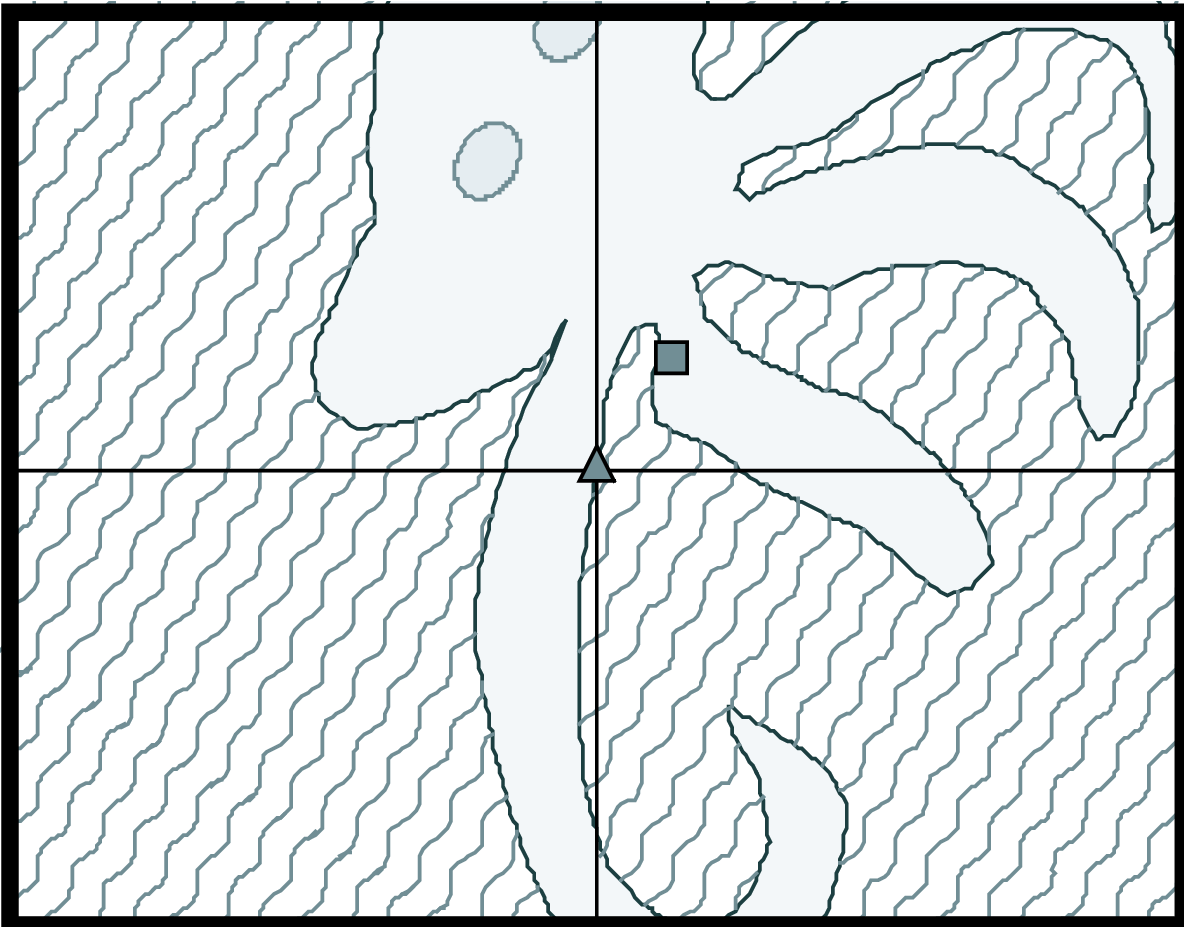}}
\end{center}

The sets obtained by grouping transformations which are conjugate to each other are called  \defn{conjugacy classes}.

In the group of isometries, there are exactly two conjugacy classes of involutions. The first consists of the involutions which fix exactly one point of the plane, i.e.\ the $180$° rotations about one point. The second consists of the involutions which fix a line in the plane, i.e.\ the reflections.

\bigskip

\emph{We now study the conjugacy classes of involutions in some larger groups...}

\subsection{Automorphisms of the plane}
We put a new structure on the Euclidean plane, by using the polynomial functions defined on it.
We say that a curve of the plane is an \defn{algebraic plane curve}  if it is defined by a polynomial equation. We have for example:
\begin{itemize}
\item
\begin{tabular}{l}
Lines (curves of degree $1$)\\
\end{tabular}
  \begin{flushright}\vspace{-1.5cm}\includegraphics[height=1.6cm]{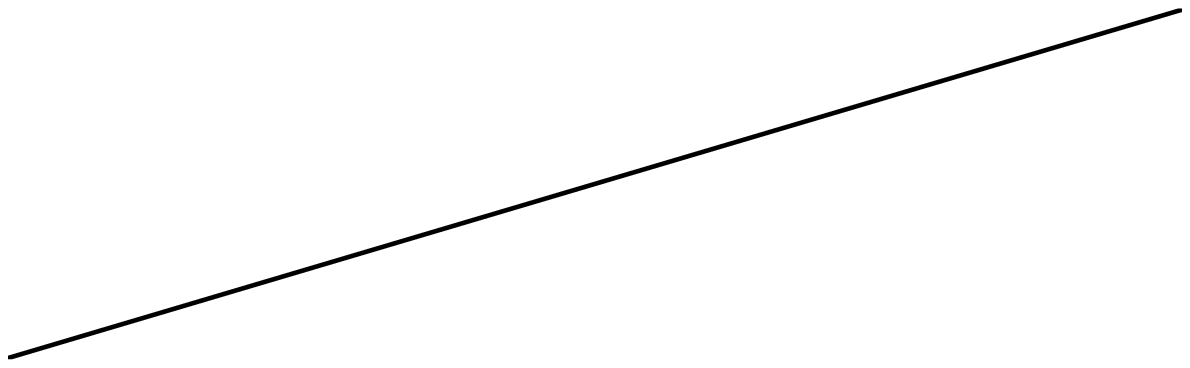}\end{flushright}
\item
\begin{tabular}{l}
Conics (curves of degree $2$)\\
{\it ellipses, parabolas, hyperbolas:}\end{tabular}
  \begin{flushright}\vspace{-1.5cm}\includegraphics[height=2.2cm]{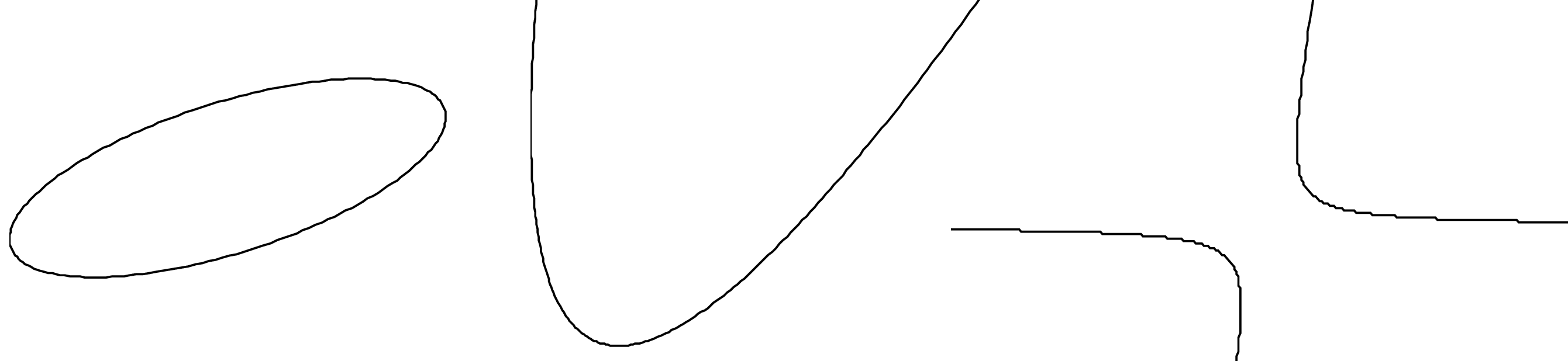}\end{flushright}
\end{itemize}

\break

\begin{itemize}
\item
\begin{tabular}{l}
Curves of higher degree\\
{\it these two have respectively degree $4$ and $3$:}\end{tabular}
 \begin{flushright}\vspace{-1.5cm}\includegraphics[height=2.3cm]{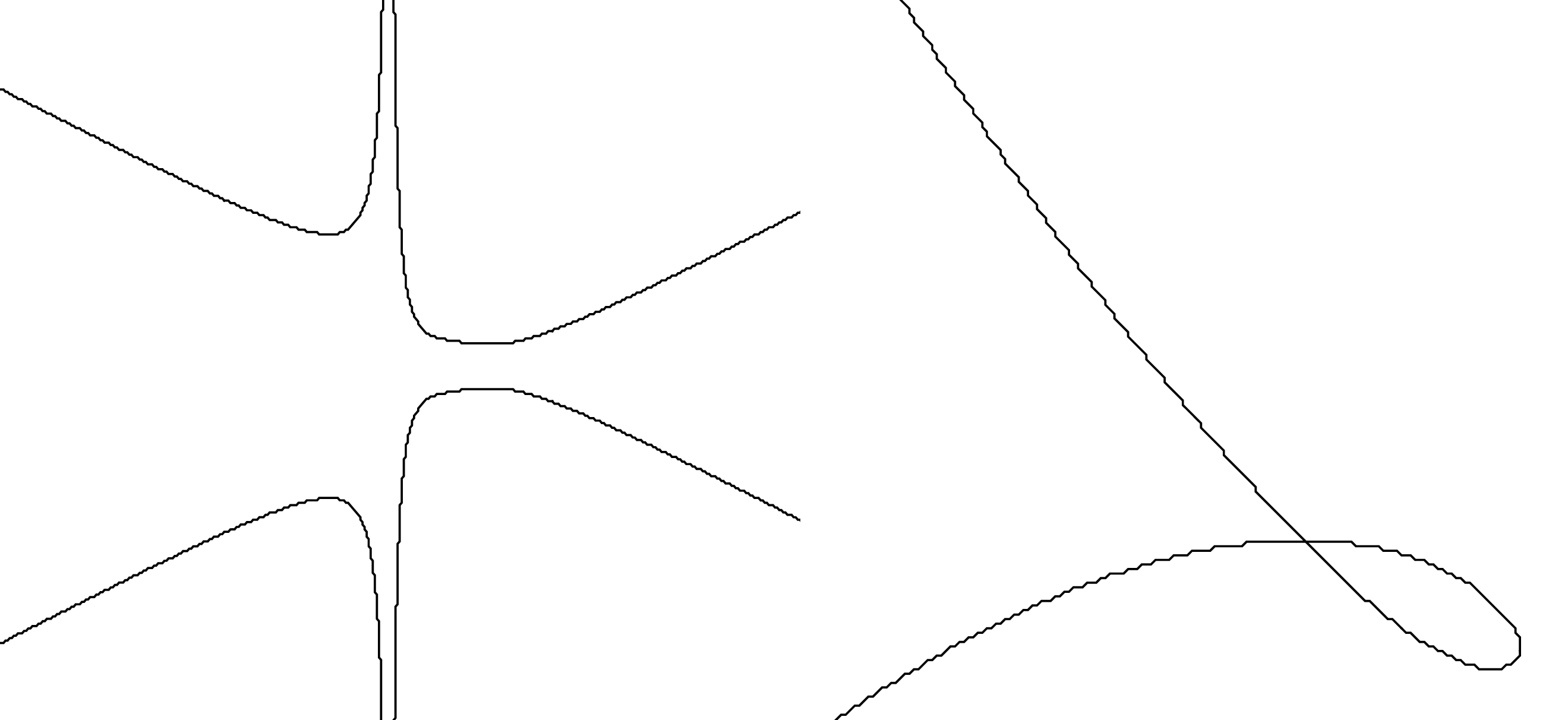}\end{flushright}
\end{itemize}
\emph{Note that not all curves are algebraic. For example, these two are not:}\begin{center}
\includegraphics[height=1cm]{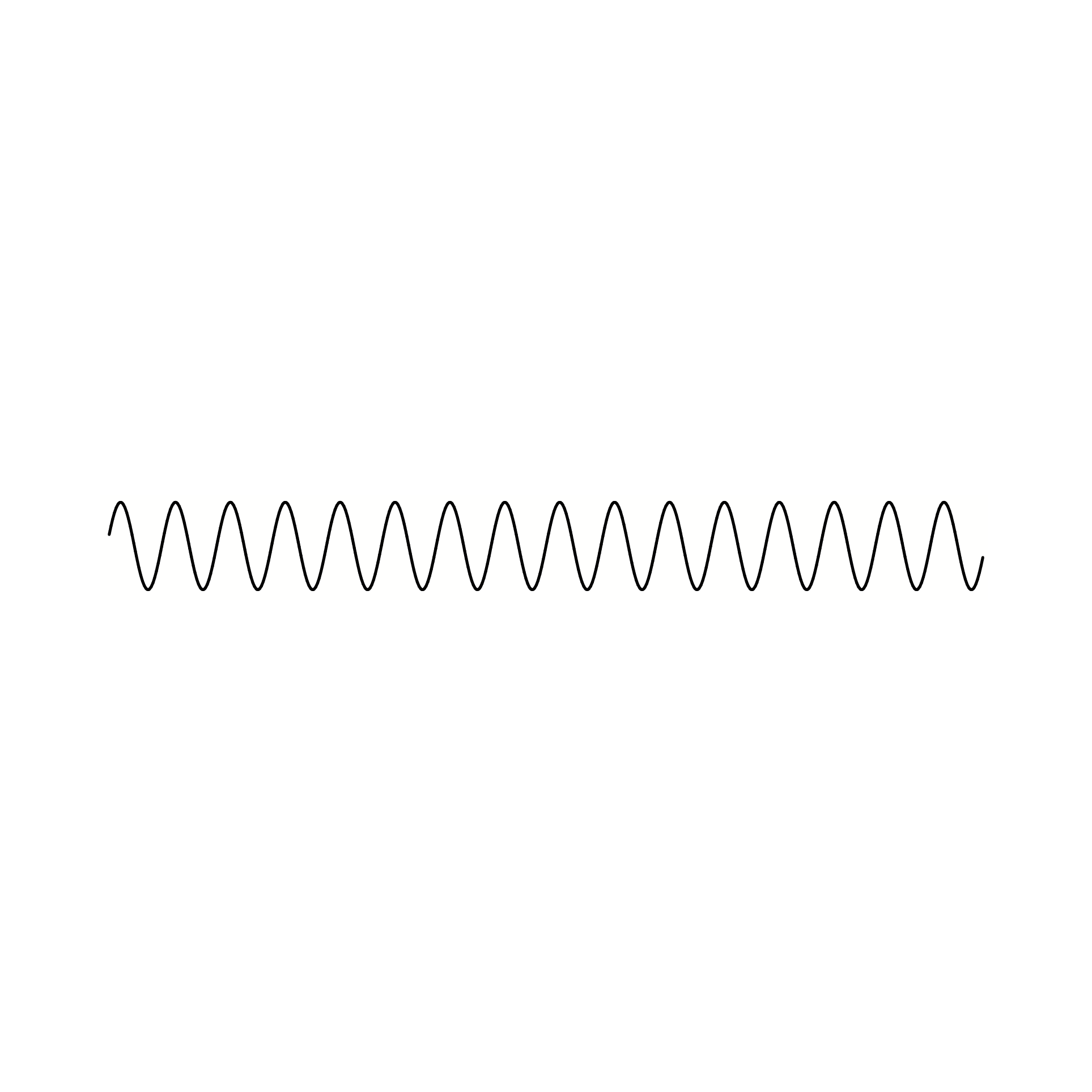}
\includegraphics[height=1cm]{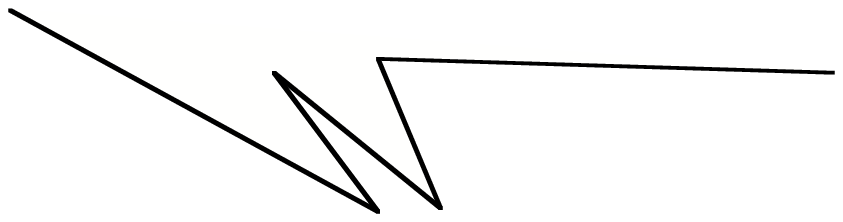}
\end{center}
\bigskip

We say that a transformation of the plane that preserves the set of algebraic curves, i.e.\ that sends any algebraic curve to an algebraic curve, is an \defn{automorphism of the plane}. Equivalently, this is a transformation of the plane defined by polynomials. For example, the isometries considered in Section \refd{Subsec:Isometries} are automorphisms of the plane defined by polynomials of degree $1$. There exist other automorphisms. Consider the following example:

\begin{center}  \includegraphics[height=3cm]{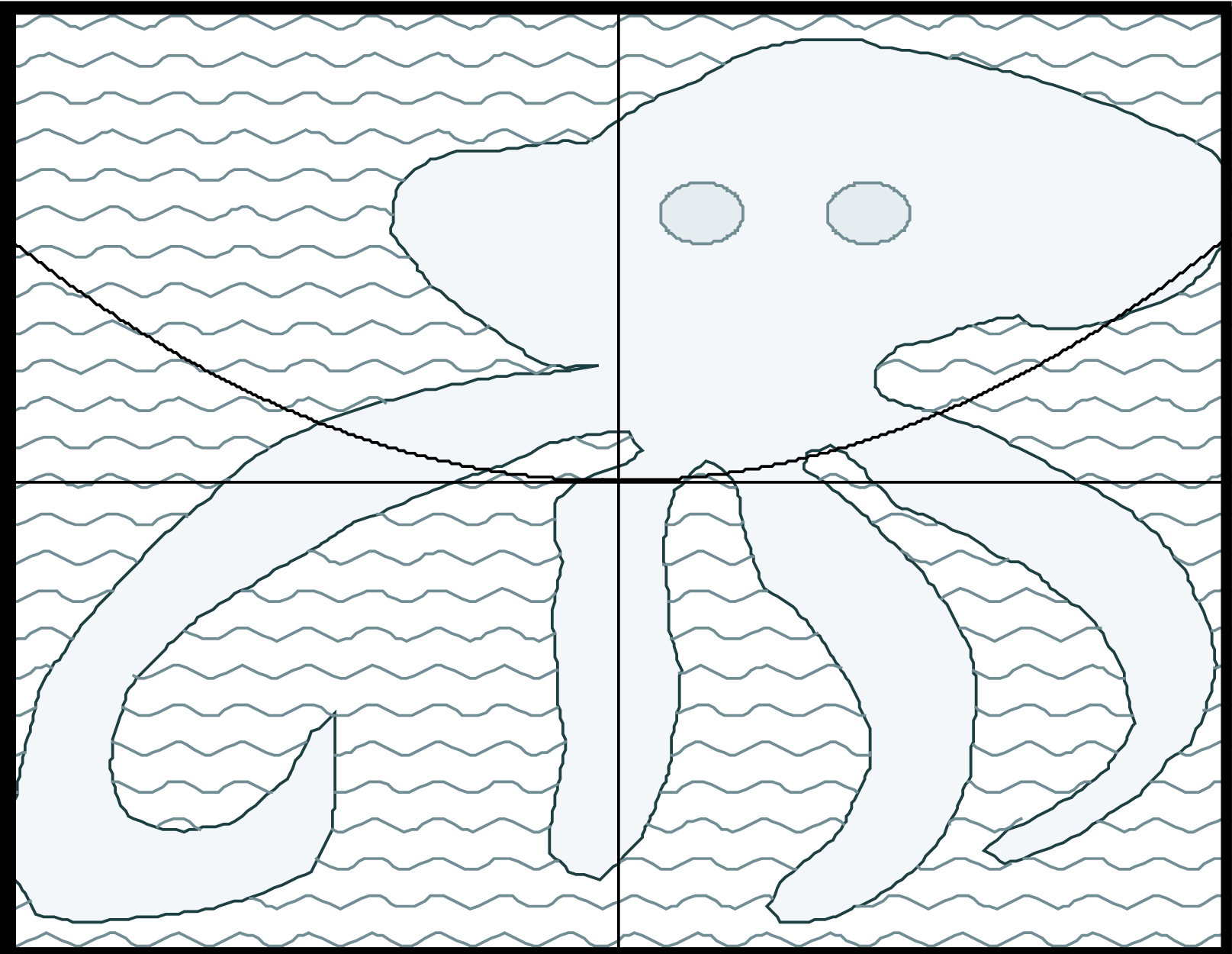} \includegraphics[height=3cm]{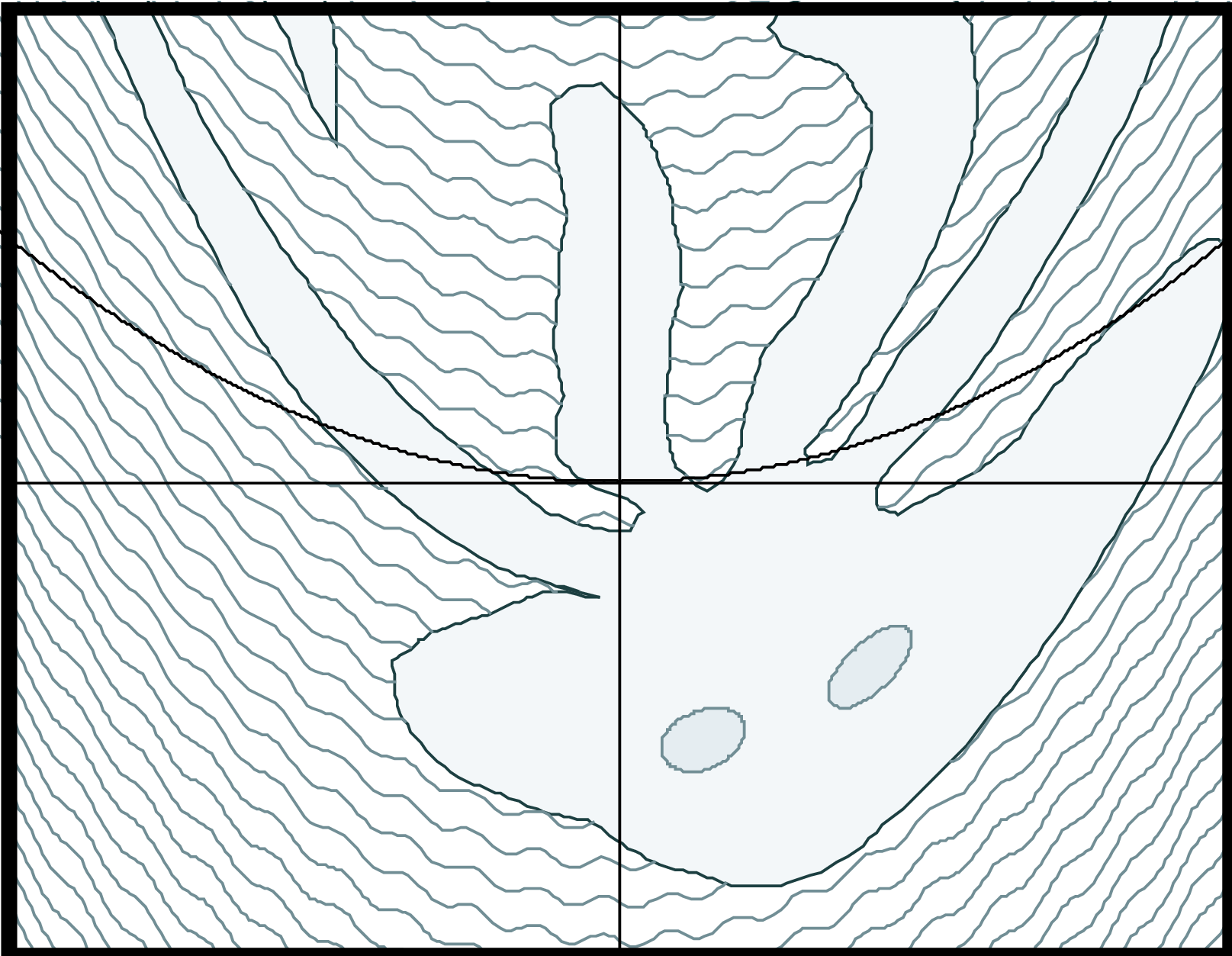}\\
{\it \small An automorphism of the plane defined by $(x,y) \mapsto (x,-y+2x^2)$.\\
 Its set of fixed points is a parabola, indicated on the figures.}\end{center}

The automorphism in this example is an involution which is not an isometry. As in Section \refd{Subsec:Isometries}, we say that two automorphisms are conjugate if some third automorphism intertwines the first with the second.
For example, the automorphism defined above is conjugate to a reflection:

\begin{center}
\hspace{0 cm}\xymatrix{ \includegraphics[height=3cm]{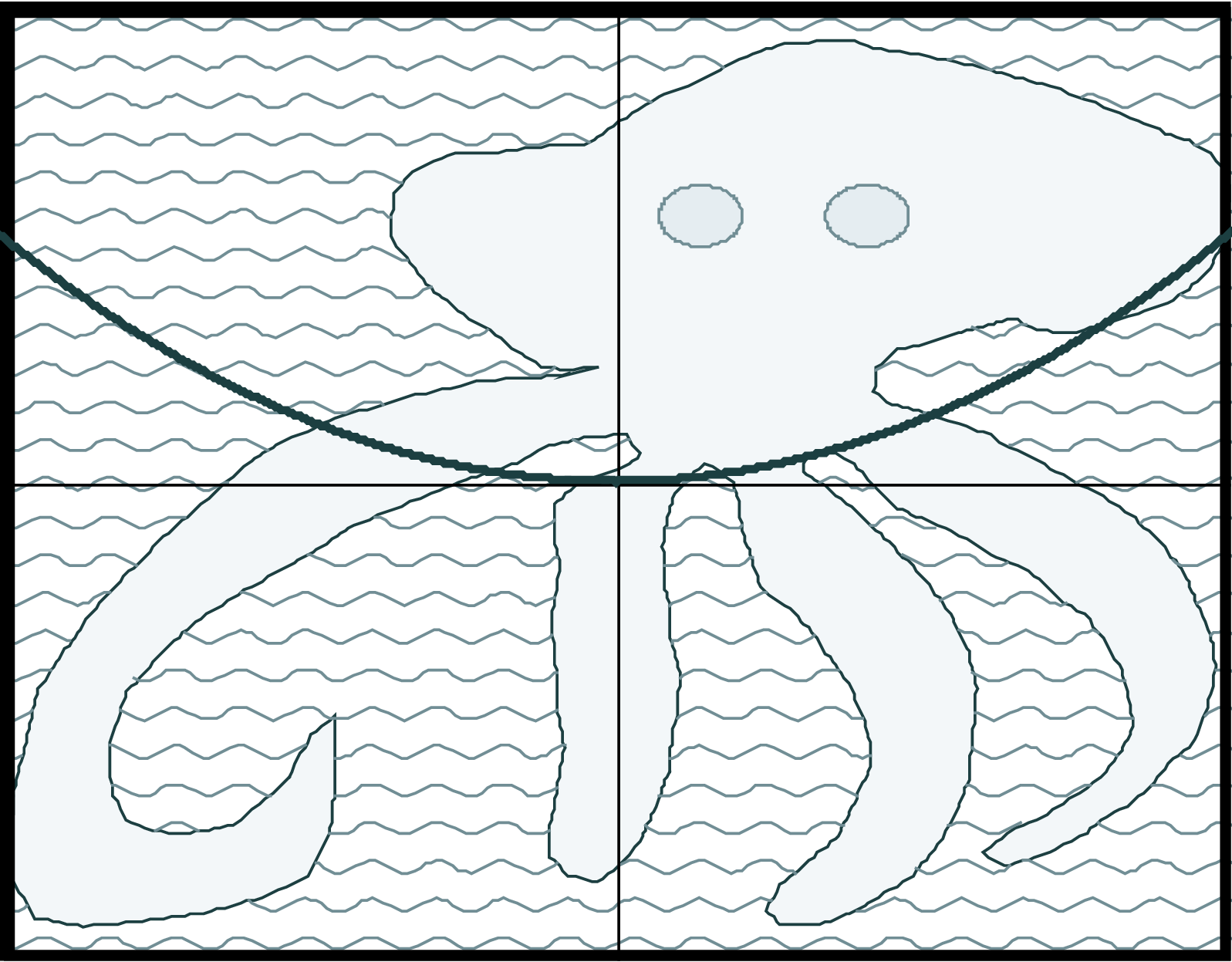}  \ar^{(x,y)\mapsto (x,-y+2x^2)}[rrrrr] \ar_{(x,y)\mapsto (x,y-x^2)}[d] & &  & & &\includegraphics[height=3cm]{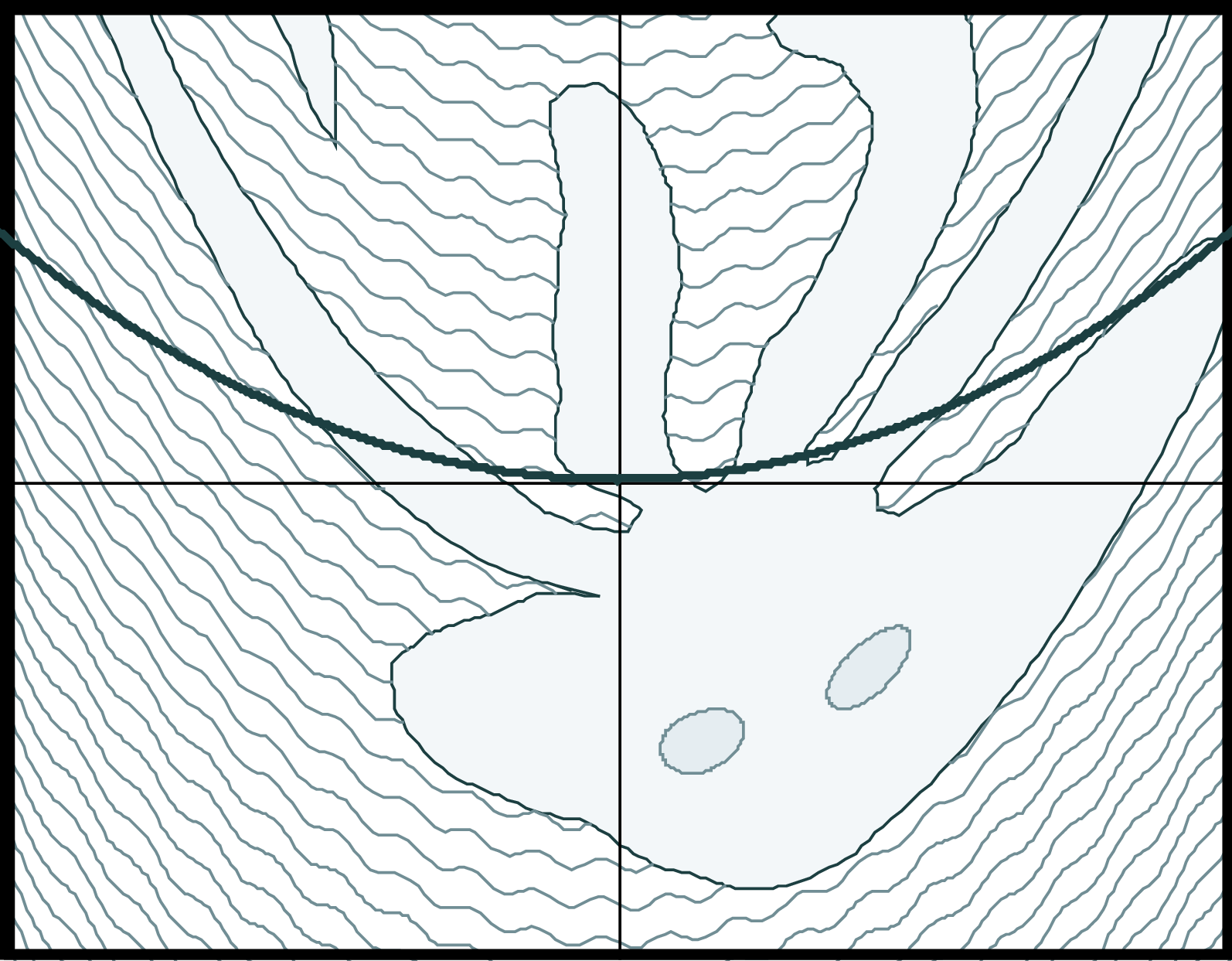} \ar_{(x,y)\mapsto (x,y-x^2)}[d] \\
\includegraphics[height=3cm]{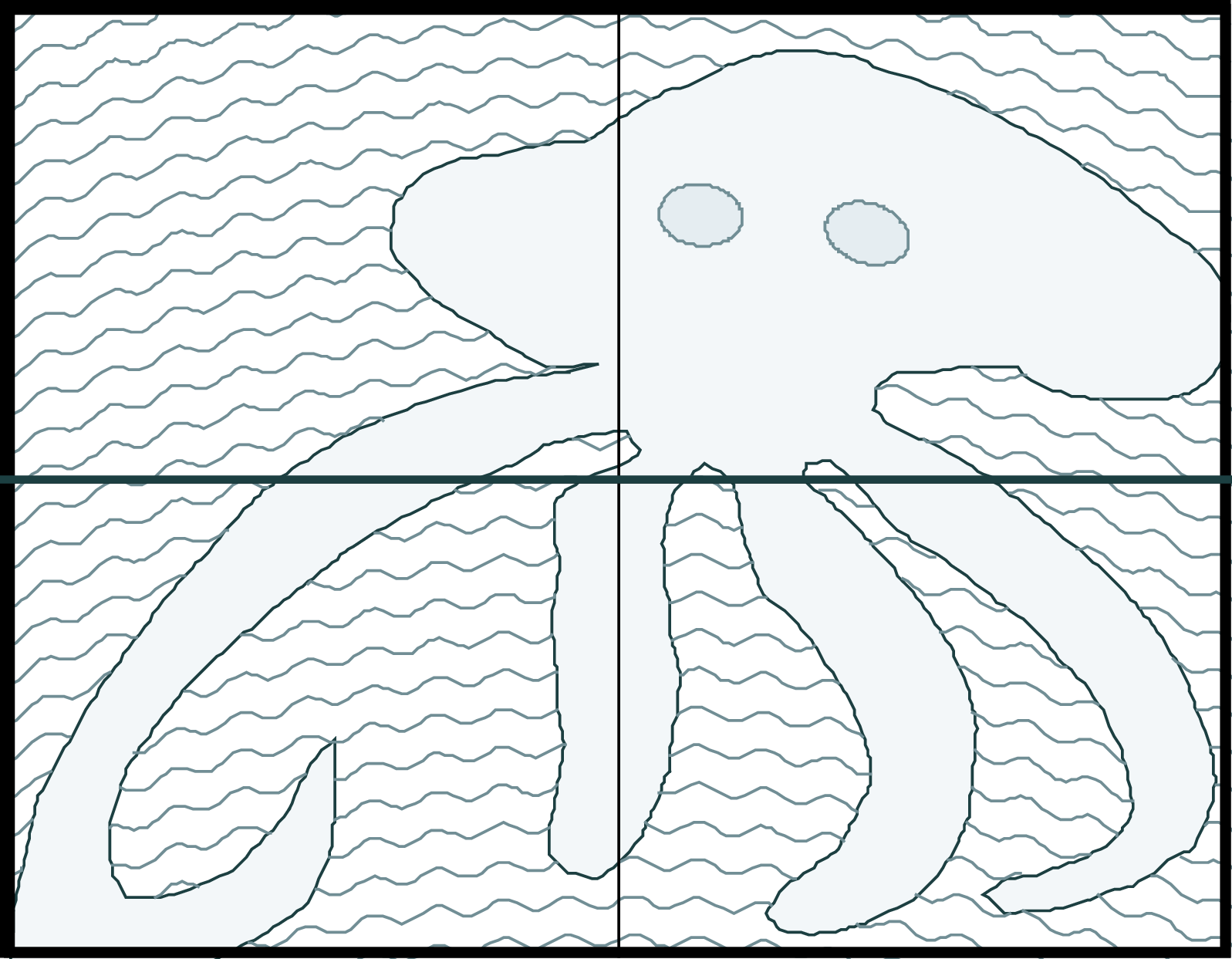} \ar^{\mbox{\it \small reflection with horizontal axis}}[rrrrr] & & & &  & \includegraphics[height=3cm]{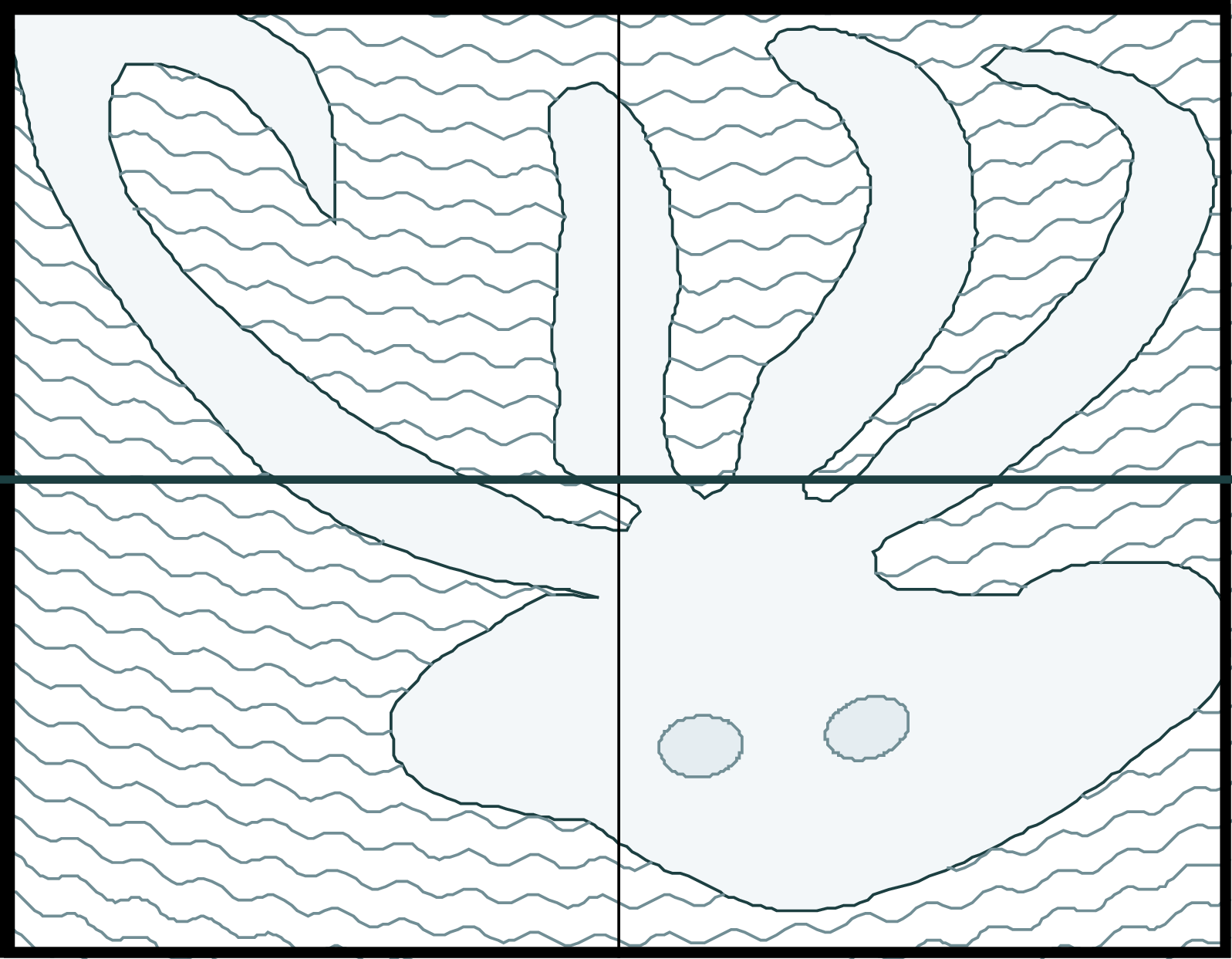}}
\end{center}

In fact, there are also exactly two conjugacy classes of involutions in the group of automorphisms of the plane. Every element of the first class is conjugate to a rotation and every element of the second one to a reflection. This follows from a famous theorem of Jung (see \cite{bib:Jun}).

\subsection{The Cremona group}
Finally, we extend our group. A \defn{birational map} of the plane is a transformation defined by quotients of polynomials. Such a map always sends most algebraic curves to algebraic curves, except for a finite number of curves, which might be mapped to points. Moreover, some points might have no image. For this reason, we denote the map by the symbol "$\dasharrow$".
For example, the map $(x,y) \dasharrow (xy,y)$ has inverse $(x,y) \dasharrow (x/y,y)$:
\begin{center}
\includegraphics[height=3cm]{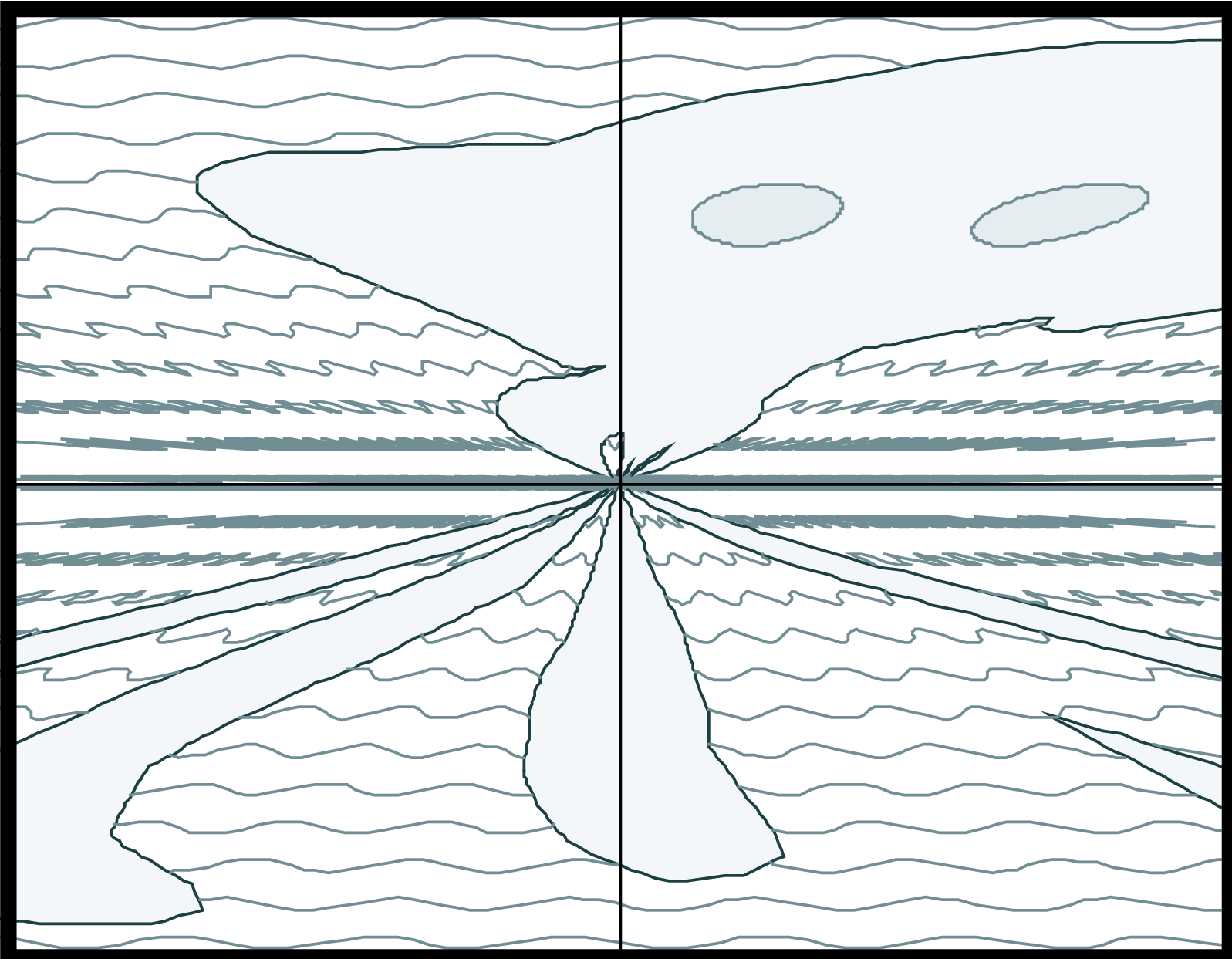}
\includegraphics[height=3cm]{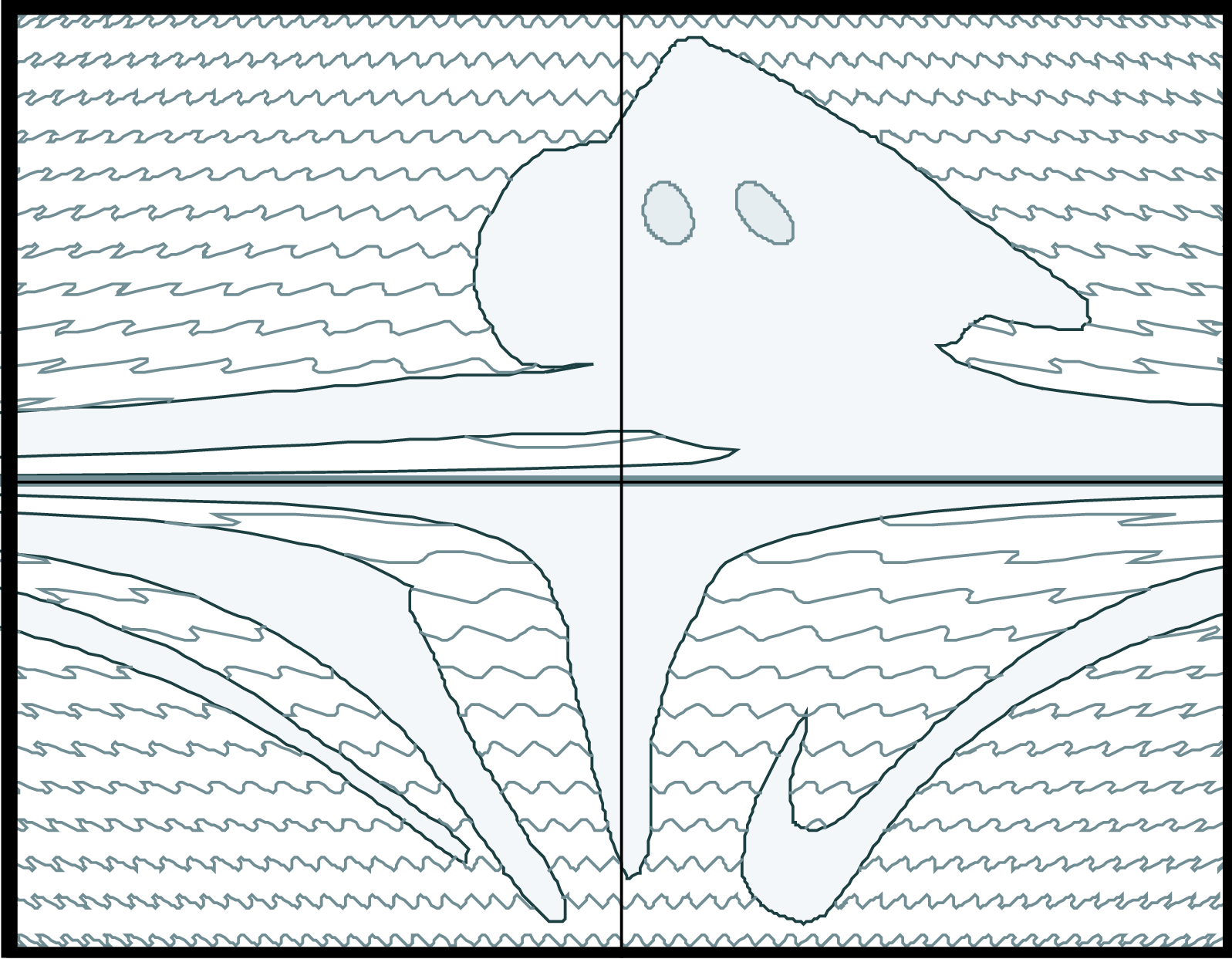}\\
{\it \small On the left, the image of the octopus by $(x,y) \dasharrow (xy,y)$, which sends the line $y=0$ to a point.\\ On the right, the image of the octopus by $(x,y) \dasharrow (x/y,y)$; the point $(0,0)$ has no image.}
\end{center}

The set of all birational maps of the plane is called \defn{the Cremona group}.

Note that the reflection and the $180$° rotation, which are not conjugate in the automorphism group, are conjugate in the Cremona group (in fact, this result is true in any dimension, as proved in \cite{bib:JB}, Theorem 1).
We can show this explicitly:
\begin{center}
\hspace{0 cm}\xymatrix{ \includegraphics[height=3cm]{Id.ps}  \ar^{\mbox{\it \small reflection with horizontal axis}}[rrrrr] \ar_{(x,y)\mapsto (xy,y)}[d] & & & & &\includegraphics[height=3cm]{Sym.ps} \ar_{(x,y)\mapsto (xy,y)}[d] \\
\includegraphics[height=3cm]{xy-y.ps} \ar^{\mbox{\it \small  $180$° rotation}}[rrrrr] &  & & & & \includegraphics[height=3cm]{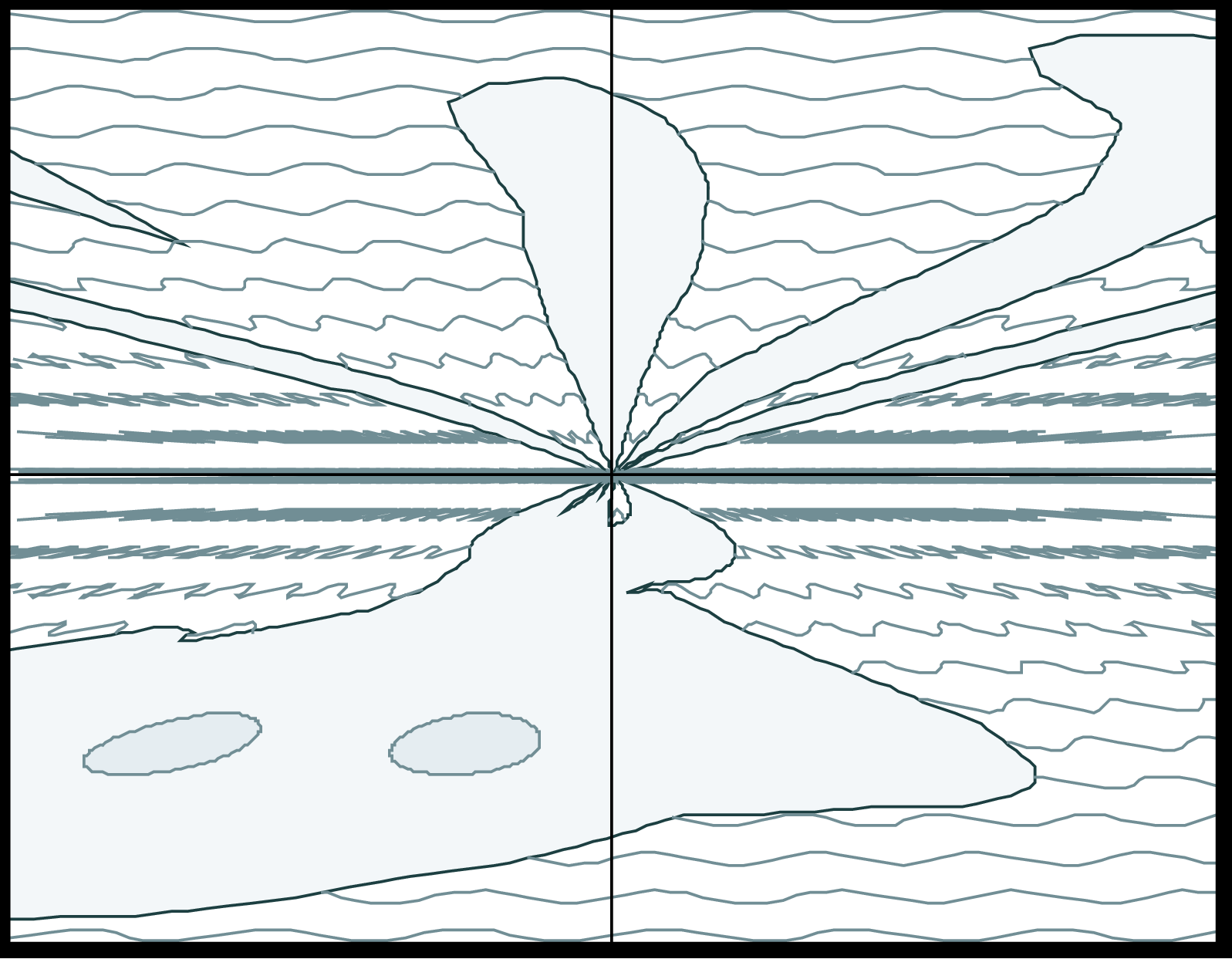}}
\end{center}
Thus, the two conjugacy classes of involutive isometries belong to the same conjugacy class of the Cremona group.
 
However, there are infinitely many conjugacy classes of involutions in the Cremona group; each class belongs to one of three families, called de Jonqui\`eres, Geiser and Bertini involutions.\index{Involutions!Bertini|idi}\index{Involutions!Geiser|idi}
\begin{center}  \includegraphics[height=3cm]{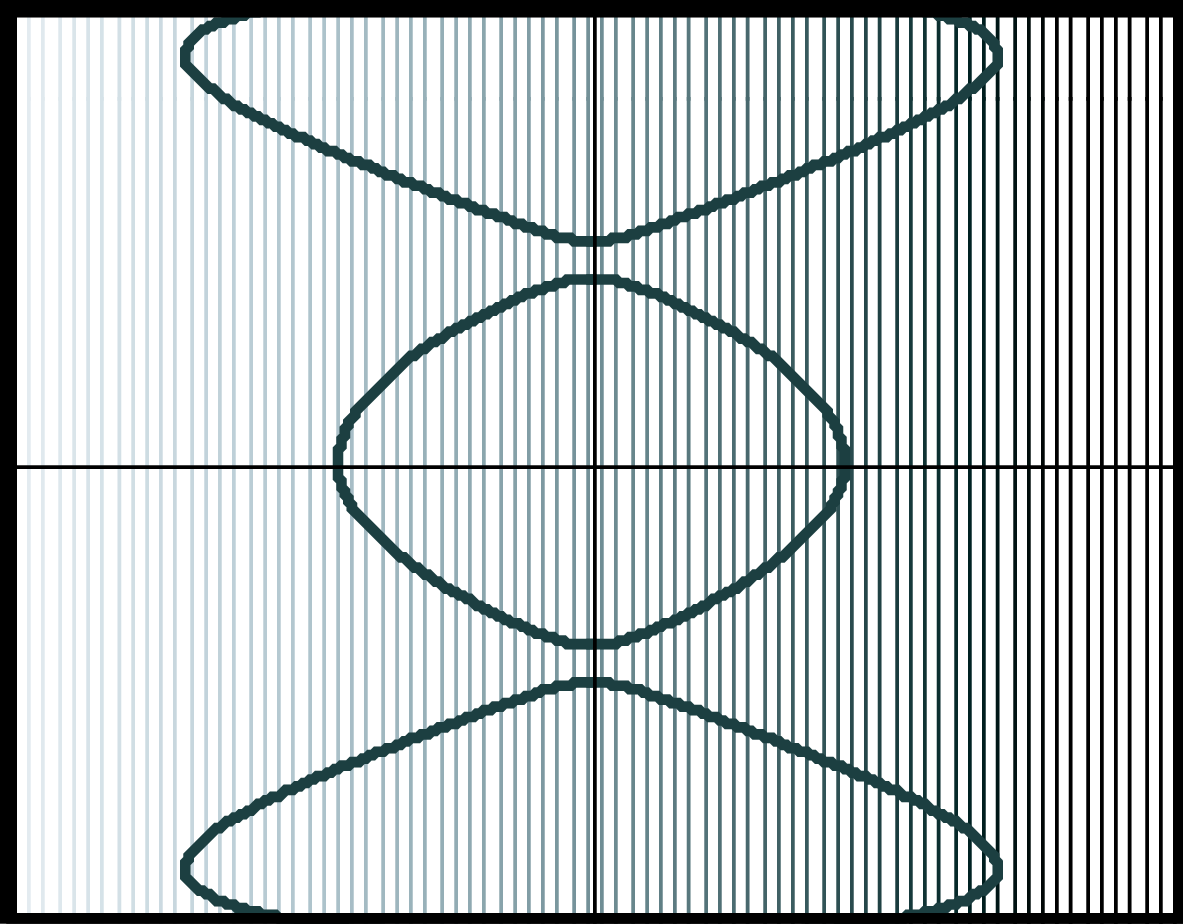} \includegraphics[height=3cm]{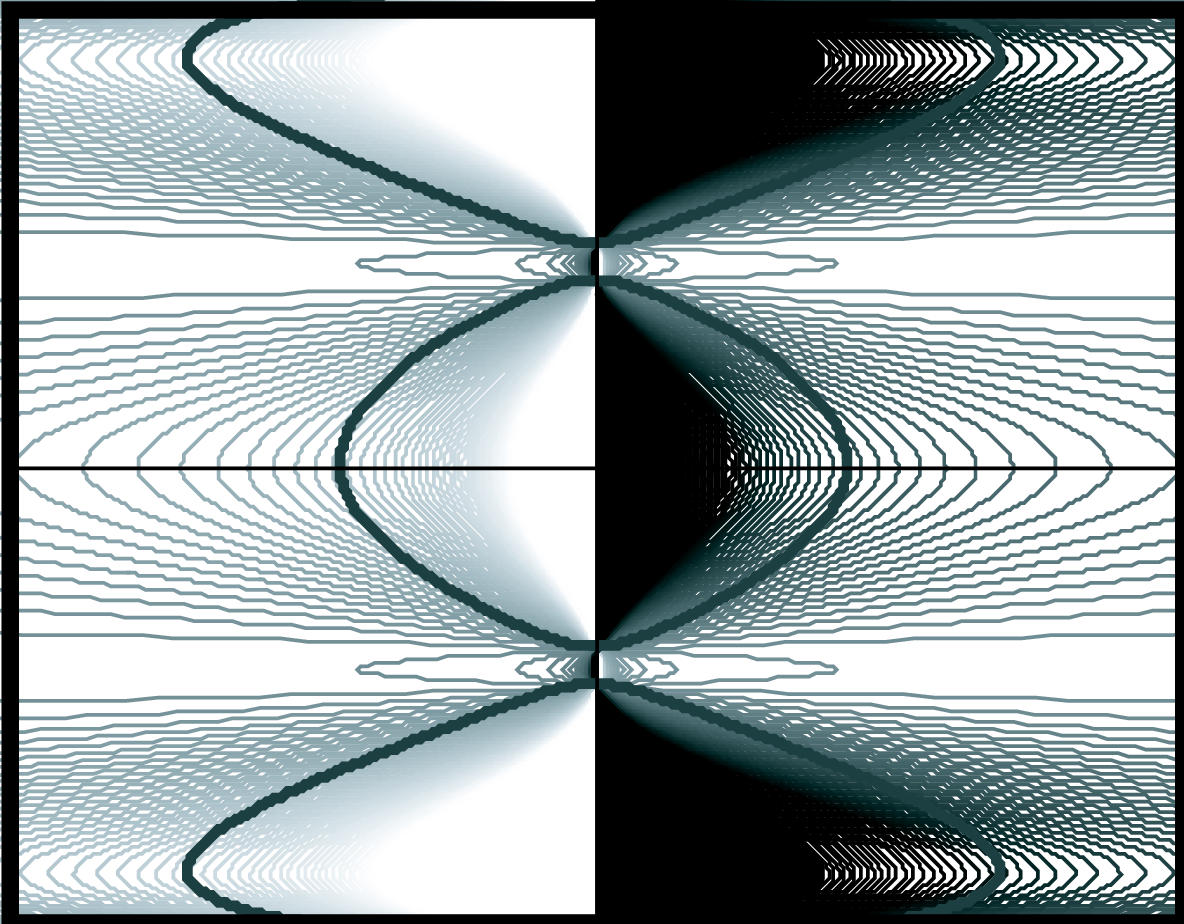}\\
{\it \small The image of some vertical lines by a de Jonqui\`eres involution. The set of fixed points is in grey.}\end{center}

Our aim is to obtain such a classification for involutions, but also for some other finite groups of birational maps. In fact we shall work with the complex plane rather than the real plane. This simplifies many of the results (but not the figures).

\emph{So, let the story begin....}

\section{A long history of results}
\label{Sec:LongHistory}
\pagetitless{Introduction}{A long history of results}
As we explained above, our subject is the Cremona group of the plane, which is the group of birational maps of the surface $\mathbb{P}^2(\K)$ (or $\K^2$), or equivalently, the group of $\K$-automorphisms of the field $\K(X,Y)$.

This very large group has been a subject of research for many years. We refer to \cite{bib:AC} for a modern survey about its elements. Some presentations of the group by generators and relations are available (see \cite{bib:Gi2} and \cite{bib:Isk4}), but these results does not provide substantial insight into the algebraic properties of the group. For example, given an abstract group, it is not possible to say whether it is isomorphic to a subgroup of the Cremona group. Moreover, the results of \cite{bib:Gi2} and \cite{bib:Isk4} do not allow one to decide whether the Cremona group is isomorphic to a linear group.

The study of the finite subgroups of the Cremona group is therefore fundamental to the understanding of this group.

\bigskip

\emph{Our aim is to describe the conjugacy classes of finite abelian subgroups of the Cremona group.}\\
The study of finite subgroups of the Cremona group was begun over one hundred years ago and it appears unlikely that it will be completed in the near future. Let us give a historical review of the main results: 

\begin{itemize}
\item
The first results are attributed to Bertini, for his work on involutions in 1877 (see \cite{bib:EB}). He identified three types of conjugacy classes, which are now known as de Jonqui\`eres, Geiser and Bertini involutions.\index{Involutions!Bertini|idi}\index{Involutions!Geiser|idi} However, his proof of the classification of conjugacy classes in each type is generally considered incomplete (see \cite{bib:BaB}).
\item
In 1895, S. Kantor \cite{bib:SK} and A. Wiman \cite{bib:Wim} gave a description of finite subgroups. The list is exhaustive, but not precise in two respects:
\begin{itemize}
\item
Given some finite group, it is not possible using their list to say whether this group is isomorphic to a subgroup of the Cremona group.
\item
The possible conjugation between the groups of the list is not considered.
\end{itemize}
\item
A great deal of work was done by the Russian school and in particular by M.K. Gizatullin, V.A. Iskovskikh and Yu. Manin.
\begin{itemize}
\item
They obtained many results on $G$-surfaces (rational surfaces with a biregular action of some group $G$, see Section \refd{Sec:GSurfaces}). Our main interest is in the classification of minimal $G$-surfaces into two types (see Proposition \refd{Prp:TwoCases}). The description of decompositions of birational maps into elementary links (see \cite{bib:Isk5}) is also a very useful tool. 
\item
We refer to the articles
\cite{bib:Gi1}, \cite{bib:Isk1},\cite{bib:Isk2}, \cite{bib:Isk3} and \cite{bib:Isk5} for more information.
\end{itemize}
\item
The modern approach started with the work of L. Bayle and A. Beauville on involutions  (see \cite{bib:BaB}). They used the classification of minimal $G$-surfaces to classify the subgroups of order $2$ of the Cremona group. This is the first example of a precise description of conjugacy classes:
\begin{itemize}
\item
The number of conjugacy classes and their descriptions are precise and clear, parametrised by isomorphism classes of curves.
\item
One can decide directly whether two involutions are conjugate or not.\end{itemize}
\item
The techniques of \cite{bib:BaB} were generalised by T. de Fernex (see \cite{bib:deF}) to cyclic groups of prime order. The list is as precise as one can wish, except for two classes of groups of order $5$, for which the question of their conjugacy is not answered.\\
The complete classification was obtained in \cite{bib:BeB}, by A. Beauville and the author. For example, the following result was proved: a cyclic group of prime order is conjugate to a linear automorphism of the plane if and only if it does not fix a curve of positive genus.\index{Curves!of positive genus!birational maps that do not fix a curve of positive genus}\\
{\it We generalise this result to finite abelian groups in Theorem \refThmAbelGenusCurve\ (see Section \refd{Sec:Results}).}
\item
A. Beauville has further classified the $p$-elementary maximal groups up to conjugation (see \cite{bib:Be2}). He obtains for example the following results:
\begin{itemize}
\item No group $(\Z{p})^3$ belongs to the Cremona group if $p$ is a prime $\not=2,3$.
\item There exist infinitely many conjugacy classes of groups isomorphic to $(\Z{2})^4$.
\item He also gave some useful algebraic tools, using group cohomology, which we use in the sequel (see Lemma \refd{Lem:splitsequence}).
\end{itemize}
Note that the conjugacy classes of subgroups $(\Z{2})^4$ of the de Jonqui\`eres group are well described (\cite{bib:Be2}, Proposition 2.6). However it is not clear whether two groups non-conjugate in the de Jonqui\`eres group are conjugate in the Cremona group.\item
More recently, I.V. Dolgachev and V.A. Iskovskikh had updated the list of S. Kantor and A. Wiman, using the modern theory of $G$-surfaces, the theory of elementary links of V.A. Iskovskikh, and the conjugacy classes of Weyl groups (see \cite{bib:Dol}). This manuscript in preparation contains many new results and is currently the most precise classification of conjugacy classes of finite subgroups. \\
However, the following questions remain open:
\begin{itemize}
\item
Given two subgroups of automorphisms of the same rational surface, not conjugate by an automorphism of the surface, are they birationally conjugate?
\item
\cite{bib:Dol} gives a list of elements of the Cremona group which are not conjugate to linear automorphisms. All these elements have order $\leq 30$. The complete list of possible orders is not given.\\
{\it  This complete list is given in our Theorem \refThmExistenceNonLinear.}
\item
Given some finite group, it is still not possible using \cite{bib:Dol} to say whether the group is isomorphic to a subgroup of the Cremona group.\\
{\it  The answer to this question, in the abelian case, is given in Theorem \refThmIsomorphyclasses\ below.}
\end{itemize}
Moreover, the conjugacy classes of automorphisms of conic bundles are only partially described.
\end{itemize}

We give in Chapter \refd{Chap:List}  a new list which is the most complete one to date in the case of finite abelian groups, but some questions remain open and some new ones arise. 
For example, we discover infinitely many conjugacy classes of cyclic groups of order $2n$, for any integer $n$, but the precise conjugation between them is not established.

Even though the classification is not yet complete, it implies several results, and answers some questions that were previously open (see Section \refd{Sec:Results}).

\section{The techniques of this paper}
\pagetitless{Introduction}{The techniques of this paper}
The main technique is the same as in \cite{bib:BaB}, \cite{bib:deF}, \cite{bib:Be2} and \cite{bib:Dol}: 
\begin{itemize}
\item
We consider finite subgroups of the Cremona group as biregular automorphisms of some complete rational smooth surfaces (or equivalently as $G$-surfaces).
\item
We use the classification of minimal $G$-surfaces (Proposition \refd{Prp:TwoCases}), which comprises two cases: 
\begin{itemize}
\item
groups of automorphisms of conic bundles (see Chapters \refd{Chap:ConicBundles}, \refd{Chap:LargeDegreeSurfaces}
and \refd{Chap:FiniteAbConicBundle}),
\item
groups of automorphisms of surfaces where the canonical class is, up to a multiple,  the only invariant class of divisors  (see Chapters \refd{Chap:LargeDegreeSurfaces} and \refd{Chap:RankOne}).
\end{itemize}
This dichotomy organises our investigations. We describe the finite abelian groups in both cases.
\item
We use some conjugacy invariants (see Chapter \refd{Chap:ConjugacyInvariants}) to say that the groups obtained are non-conjugate, or we conjugate them directly if we can.
\end{itemize}
In addition, we also use many other tools, for example:
\begin{itemize}
\item
invariant pencils of rational curves (see Sections \refd{Sec:PencilRatCurves} and \refd{Sec:PairInvariantPencils});
\item
actions of a group on the set of points fixed by one of its elements (see Section \refd{Sec:ActionFixedCurves});
\item
conic bundle structures on Del Pezzo surfaces\index{Del Pezzo surfaces!others|idi} (see Section \refd{Sec:ConicBundlonDP});
\item
the existence of fixed points (see Section \refd{Sec:ExistenceFixedPoints});
\item
actions of the automorphisms on the exceptional divisors which induce the actions on the Picard group (see Propositions \refd{Prp:rootsDeJ}, \refd{Prp:GeoQuad} and \refd{Prp:GeoCub}).
\end{itemize}

\section{The results}
\pagetitless{Introduction}{The results}
\label{Sec:Results}
{\it In Chapter \refd{Chap:List}, Theorems \refThmCyclicGroup\ and \refThmNonCyclicGroup, we summarise our investigations and give an exhaustive list of the finite abelian subgroups of the Cremona group.} \\
As we mentioned above, this is the most precise list to date.

Using this list, we prove the following new theorems:

\bigskip

\noindent
{\bf THEOREM 1: \sc Non-linear birational maps of large order}
\begin{itemize}\itshape
\item
For any integer $n\geq 1$, there are infinitely many conjugacy classes of birational maps of the plane of order $2n$, that are non-conjugate to a linear automorphism. 
\item
If $n>15$, a birational map of order $2n$ is an $n$-th root of a de Jonqui\`eres involution and preserves a pencil of rational curves.
\item
If a birational map is of finite odd order and is not conjugate to a linear automorphism of the plane, then its order is $3,5,9$ or $15$. In particular,  any birational map of the plane of odd order $>15$ is conjugate to a linear automorphism of the plane.
\end{itemize}

\begin{remark}\fauxtitred
Prior to this work, the largest known order of  a non-linear birational map was $30$, given in $\cite{bib:Dol}$.\end{remark} 

\bigskip

\noindent{\bf THEOREM 2: \sc Roots of linear automorphisms}\\ \itshape
Any birational map which is a root of a non-trivial linear automorphism of finite order of the plane  is conjugate to a linear automorphism of the plane. \upshape

\bigskip

\begin{flushleft}{\bf THEOREM 3: \sc Groups which fix a curve of positive genus}\\ \itshape\index{Curves!of positive genus!birational maps that fix a curve of positive genus|idb}
Let $G$ be a finite abelian group which fixes some curve of positive genus. Then $G$ is cyclic, of order $2$, $3$, $4$, $5$ or $6$, and all these cases occur. If the curve has genus $>1$, the order is $2$ or $3$.\end{flushleft} \upshape
\begin{remark}\fauxtitred
This generalises a theorem of Castelnuovo  \cite{bib:Cas}, which states that an element of finite order which fixes a curve of genus $\geq 2$ has order $2$, $3$, $4$ or $6$.\end{remark}
\bigskip

\bigskip

\begin{flushleft}{\bf THEOREM 4: \sc Cyclic groups whose non-trivial elements do not fix a curve of positive genus}\\ \itshape\index{Curves!of positive genus!birational maps that do not fix a curve of positive genus|idb}
Let $G$ be a finite cyclic subgroup of the Cremona group. The following conditions are equivalent:
\begin{itemize}
\item
If $g \in G$, $g\not=1$, then $g$ does not fix a curve of positive genus.
\item
$G$ is birationally conjugate to a subgroup of $\Aut(\Pn)$.
\item
$G$ is birationally conjugate to a subgroup of $\Aut(\mathbb{P}^1\times\mathbb{P}^1)$.
\end{itemize}
\end{flushleft}
\upshape
\begin{remark}\fauxtitred
This generalises  Corollary $1.2$ of \cite{bib:BeB}, which states the same result for cyclic subgroups of the Cremona group of prime order.\end{remark}

\break

\begin{flushleft}{\bf THEOREM 5: \sc Abelian groups whose non-trivial elements do not fix a curve of positive genus}\\ \itshape\index{Curves!of positive genus!birational maps that do not fix a curve of positive genus|idb}
Let $G$ be a finite abelian subgroup of the Cremona group. The following conditions are equivalent:
\begin{itemize}
\item
If $g \in G$, $g\not=1$, then $g$ does not fix a curve of positive genus.
\item
$G$ is birationally conjugate to a subgroup of $\Aut(\Pn)$, or to a subgroup of $\Aut(\mathbb{P}^1\times\mathbb{P}^1)$ or to the group isomorphic to $ \Z{2}\times\Z{4}$, generated by the two elements $\begin{array}{lll}(x:y:z)&\dasharrow&(yz:xy:-xz),\\
(x:y:z)&\dasharrow &( yz(y-z):xz(y+z):xy(y+z)).\end{array}$
\end{itemize}
Moreover, this last group is conjugate neither to a subgroup of $\Aut(\Pn)$, nor to a subgroup of $\Aut(\mathbb{P}^1\times\mathbb{P}^1)$.\end{flushleft}\upshape
\begin{remark}\fauxtitred
This generalises Theorems $2$ and $\refThmCyclicGenusCurve$.\end{remark}

\bigskip

\begin{flushleft}{\bf THEOREM 6: \sc Isomorphy classes of finite abelian groups}\\ \itshape
The \emph{isomorphism} classes of finite abelian subgroups of the Cremona group are the following:
\begin{itemize}
\item
$\Z{m}\times\Z{n}$, for any integers $m,n \geq 1;$
\item
$\Z{2n}\times(\Z{2})^2$, for any integer $n\geq 1;$
\item
$(\Z{4})^2\times\Z{2};$
\item
$(\Z{3})^3;$
\item
$(\Z{2})^4$.
\end{itemize}
\end{flushleft} \upshape
\begin{remark}\fauxtitred
This generalises a result of \cite{bib:Be2}, which gives the isomorphism classes of $p$-elementary subgroups of the Cremona group.\end{remark}

\bigskip

Note that, apart from these theorems, we have proved more precise results of interest. See for example \begin{itemize}
\item
Proposition \refd{Prp:NbConicDP}: a description of conic bundle structures on Del Pezzo surfaces\index{Del Pezzo surfaces!others|idi}.
\item
Proposition \refd{Prp:LargeDegree}: any finite abelian group of automorphisms of some projective smooth rational surface, whose square of the canonical divisor is $\geq 5$, is birationally conjugate to a subgroup of $\Aut(\Pn)$ or $\Aut(\mathbb{P}^1\times\mathbb{P}^1)$.
\item
Proposition \refd{Prp:CBdistCases}: four types of finite abelian groups of automorphisms of conic bundles.
\item
Proposition \refd{Prp:ConjBetweenRankPic}: birational conjugation between the minimal pairs $(G,S)$.
\item
Proposition \refd{Prp:HardChgExacrSeq}: changes in the presentations of automorphism groups of conic bundles.
\end{itemize}

\chapter{Finite groups of automorphisms of rational surfaces}
\label{Chap:GrAutRat}
\ChapterBegin{In this chapter, we recall some results on rational surfaces and their groups of automorphisms and birational maps. We also give the representation of finite subgroups of the Cremona group as groups of automorphisms of rational surfaces.}
\section{Some facts about rational and minimal surfaces}
\label{Sec:FactrationalMinimalSurf}
\pagetitless{Finite groups of automorphisms of rational surfaces}{Some facts about rational and minimal surfaces}
\begin{flushright}\Citation{The most basic algebraic varieties are the projective spaces, and rational varieties are their closest relatives. In many applications where algebraic varieties appear in mathematics and the sciences, we see rational ones emerging as the most interesting examples...}{Preface of \cite{bib:KoSC}}\end{flushright}

We begin by recalling some classical definitions and results on surfaces, without any reference to group action. 

\begin{Def} \fauxtitre
Let $S$ be a surface. 
\begin{itemize}
\item
We say that $S$ is \defn{rational} if there exists a birational map $\varphi:S\dasharrow \Pn$. 
\item
We say that $S$ is \defn{minimal} if any birational morphism $S\rightarrow S'$ is an isomorphism (i.e.\ is biregular).
\item
An irreducible curve $\Gamma \subset S$ is called an \defn{exceptional curve} if there exists some birational morphism $\varphi:S\rightarrow S'$ for which $\varphi(\Gamma)$ is a point. The divisor associated to $\Gamma$ is called an \defn{exceptional divisor}.
\end{itemize}\end{Def}

\bigskip

\noindent\emph{Note that all surfaces considered in this paper are complete, rational and smooth. As usually the surfaces considered are smooth and complete, we will often recall only the rationality.}

\bigskip

Here is a classical result, a recent proof of which can be found in \cite{bib:Be1}, Theorem II.17:
\noindent\begin{Prp}\titreProp{Castelnuovo's contractibility criterion}\\ 
\label{Prp:Castelnuovocontractibilitycriterion}
Let $E \subset S$ be an irreducible curve. Then $E$ is an exceptional curve if and only if $E^2=-1$ and $E\cong \mathbb{P}^1$.
\end{Prp}
Let $\mathbb{F}_n=\mathbb{P}_{\mathbb{P}^1}(\mathcal{O}_{\mathbb{P}^1} \oplus \mathcal{O}_{\mathbb{P}^1}(n))$ be a Hirzebruch surface, for some integer $n\geq 0$ (see \cite{bib:Be1}, Proposition IV.1 for more details on these surfaces). We have the following result:
\begin{Prp}\titrePropRef{Minimal rational surfaces}{\cite{bib:Be1}, Theorem V.10}\\
\label{Prp:MinimalRatSurfaces}
Let $S$ be a minimal rational surface. Then $S$ is isomorphic to $\Pn$ or to one of the surfaces $\mathbb{F}_n$, with $n\not=1$.\end{Prp}

\begin{Def} \fauxtitre
The \emph{degree} of a surface $S$ is the integer $K_S^2$, where $K_S$ denotes the canonical divisor of $S$.
\end{Def}

We recall that the \defn{Picard group $\Pic{S}$} of $S$ is the group of numerical equivalence classes of divisors of $S$. In our case, this is equal to the N\'eron-Severi group $NS(S)$.

\section{Some facts about groups acting biregularly on surfaces}
\pagetitless{Finite groups of automorphisms of rational surfaces}{Some facts about groups acting biregularly on surfaces}
\label{Sec:GSurfaces}
We now consider groups acting on a surface. Note that although our work concerns finite abelian groups, the groups in this chapter  are  finite or abelian only if this is stated explicitly.

\begin{Not}\fauxtitre
Let $S$ be a rational surface. We denote its group of automorphisms by $\Aut(S)$  and its group of birational maps by $\Cr{S}$.\end{Not}

Here is a simple but important observation: a birational map $\varphi:S\dasharrow \Pn$ yields an isomorphism between $\Cr{S}$ and $\Cr{\Pn}$. The Cremona group, which is the group of birational maps of\hspace{0.2 cm}$\Pn$, is then (isomorphic to) the group of birational maps of any rational surface. We say that two groups $G \subset \Cr{S}$ and $G'\subset \Cr{S'}$ are \defn{birationally conjugate} if there exists a birational map $\varphi: S \dasharrow S'$ such that $\varphi G\varphi^{-1}=G'$. This means that the two groups represent the same conjugacy class in the Cremona group.

\begin{Def} \fauxtitred
\label{Def:Gsurfaces}
\begin{itemize}
\item
We denote by \defn{$(G,S)$} a pair in which $S$ is a (rational) surface and $G\subset \Aut(S)$ is a group acting biregularly on $S$. A pair $(G,S)$ is also classically called a \defn{$G$-surface}.
\item
Let  $(G,S)$ be a $G$-surface. We say that a birational map $\varphi: S \rightarrow S'$ is \defn{$G$-equivariant} if the $G$-action on $S'$ induced by $\varphi$ is biregular. The birational map $\varphi$ is called \defn{a birational map of $G$-surfaces.}
\item
We say that a pair $(G,S)$ is \defn{minimal} (or equivalently that $G$ \defn{acts minimally on $S$}) if any $G$-equivariant birational morphism  $\varphi:S\rightarrow S'$ is an isomorphism.
\end{itemize}\end{Def}

Note that a pair $(G,S)$ in fact represents a conjugacy class of subgroups of $\CrP$. Indeed, by taking any birational map $\varphi:S\dasharrow \Pn$, we get a subgroup $\varphi G \varphi^{-1}$ of $\CrP$. The choice of $\varphi$ does not change the conjugacy class of the group and any element of the conjugacy class can be obtained in this manner.

The converse is true, if $G$ is finite (see for example a proof in  \cite{bib:DFE}, Theorem 1.4):
\begin{Prp}\fauxtitre
\label{Prp:GS}
If $G \subset \CrP$ is a finite group, then there exist a rational surface $S$ and a birational map $\varphi: \Pn\dasharrow S$ such that $\varphi G \varphi^{-1}$ acts biregularly on $S$.
\end{Prp}

This proposition is false if $G$ is infinite, as is shown for example by the following result:

\begin{Prp}\fauxtitre
\label{Prp:InfiniteNo}
Let $G \subset \CrP$ be some group, such that $\PGLn{3}\subsetneqq G$.
There do not exist a surface $S$ and a birational map $\varphi:S\dasharrow \Pn$ such that $\varphi^{-1}G\varphi\subset \Aut(S)$.
\end{Prp}
\begin{proof}\fauxtitred\upshape
We shall use the fact that the group $\PGL(3,\K)$ is simple (see for example \cite{bib:Gri}, Theorem 3.6.7).

Suppose that such $\varphi$ and $S$ exist. In particular, $\varphi^{-1}\PGL(3,\K)\varphi\subset \Aut(S)$. The natural homomorphism $\rho: \Aut(S) \rightarrow GL(\Pic{S})$ gives a homomorphism $\tilde{\rho}:\PGLn{3}\rightarrow GL(\Pic{S})$ which we prove to be trivial. As $\PGL(3,\K)$ is simple, either $\ker\tilde{\rho}$ or $Im\tilde{\rho}$ is trivial. The homomorphism $\tilde{\rho}$ is not injective since $GL(\Pic{S})$ is countable but $\PGL(3,\K)$ is not, so $Im \tilde{\rho}$ is trivial.

Let us prove that $S=\Pn$. Indeed, the group $\PGL(3,\K)$ acts trivially on $\Pic{S}$ and therefore on the set of exceptional divisors of $S$. If $S\not=\Pn$, we may blow-down some exceptional divisors, and get a $\PGL(3,\K)$-equivariant birational morphism  $\pi:S\rightarrow \mathbb{F}_n$ for some integer 
$n\geq 0$. If $n=0$, the action of $\PGL(3,\K)$ on the two conic fibrations gives a homomorphism $\PGL(3,\K) \rightarrow \Z{2}$ whose image is trivial. Thus, $\PGL(3,\K)$ preserves the two conic bundle fibrations.
If $n>0$, there is a unique conic bundle structure, so $\PGL(3,\K)$ acts on it. In each case, we get a homomorphism $\PGL(3,\K) \rightarrow \PGL(2,\K)$, which is not injective ($\PGL(3,\K)$ contains groups isomorphic to $(\Z{n})^2$ for any positive integer $n$ but $\PGL(2,\K)$ does not) and hence is trivial. Hence, $\PGL(3,\K)$ embeds in the group of automorphisms of the generic fibre with a trivial action on the basis, which is $\PGL(2,\K(x))$. The same argument as immediately above shows that this is not possible.

We have proved that $S=\Pn$. Since $\varphi^{-1}$ is $\PGL(3,\K)$-equivariant, its base points are orbits of $\PGL(3,\K)$. As the number of base points is finite, we see that this number is $0$, so $\varphi$ is an isomorphism. This implies that $G \subset \Aut(\Pn)=\PGL(3,\K)$, which is false by hypothesis.\proofend
\end{proof}

\begin{Def} \fauxtitre
We say that two pairs $(G,S)$ and $(G',S')$ are \defn{birationally conjugate} if they represent the same conjugacy class in the Cremona group.
\end{Def}

\section{Two distinct cases}
\pagetitless{Finite groups of automorphisms of rational surfaces}{Two distinct cases}
\label{Sec:TwoDistinctCases}

The modern approach (used in $\cite{bib:BaB}$, $\cite{bib:deF}$, $\cite{bib:Be2}$ and $\cite{bib:Dol}$) is to consider any finite subgroup $G \subset \CrP$ as a group of automorphisms of  a surface $S$ (using Proposition \refd{Prp:GS}) and to suppose that the pair $(G,S)$ is minimal.

Then, we apply the following proposition, due to Yu.\ Manin in the abelian case (see \cite{bib:Man}) and V.A. Iskovskikh in the general case (see \cite{bib:Isk3}). 
\begin{Prp}\fauxtitre
\label{Prp:TwoCases}\index{Conic bundles!automorphisms}
Let $S$ be a (rational smooth) surface and $G \subset \Aut(S)$ be a finite subgroup of automorphisms of $S$.
If the pair $(G,S)$ is minimal then one \emph{and only one} of the following holds:
\begin{itemize}
\item[\upshape 1.]
The surface $S$ has a conic bundle structure invariant by $G$, and $\rkPic{S}^G=2$, i.e.\ the fixed part of the Picard group is generated by the canonical divisor and the divisor class of a fibre.
\item[\upshape 2.]
$\rkPic{S}^{G}=1$, i.e.\ the fixed part of the Picard group is generated by the canonical divisor.\proofend
\end{itemize}
\end{Prp}

Let us comment on this proposition:
\begin{itemize}
\item
In the first case, there exists a morphism $\pi:S \rightarrow \mathbb{P}^1$ with general fibres rational and irreducible and such that every singular fibre is the union of two rational curves $F$, $F'$ with $F^{2}=F'^{2}=-1$ and $FF'=1$. The group $G$ embeds in the group of automorphisms of the generic fibre $\mathbb{P}^1_{\K(x)}$ of $\pi$.
This group is abstractly $\PGL_{2}(\K(x))\rtimes \PGLn{2}$ and can be viewed as the group of birational maps of $\mathbb{P}^1\times \mathbb{P}^1$ that preserve the first projection (see Section \refd{Sec:CBTwoRepr} and in particular Proposition \refd{Prp:BeauTsen}).

The \defn{de Jonqui\`eres involutions} are examples of this case: they are given in this context as birational maps of $\mathbb{P}^1\times\mathbb{P}^1$ of the form
\begin{center}$(x_1:x_2) \times (y_1:y_2) \dasharrow (x_1:x_2) \times (y_2\prod_{i=1}^n(x_1-a_i x_2):y_1\prod_{i=1}^n(x_1-b_i x_2))$,\end{center}
for some $a_1,...,a_n,b_1,...,b_n \in \K$ all distinct.

We explicit all possible finite abelian groups of automorphisms of conic bundles in Chapter $\refd{Chap:ConicBundles}$.

\item

In the second case, $S$ is a \defn{Del Pezzo surface}\index{Del Pezzo surfaces!definition|idb} (see Lemma \refd{Lem:rkDelPezzo}). We recall that a Del Pezzo surface is a rational surface whose anti-canonical divisor $-K_S$ is ample. In fact, it is either
$\Pn$, $\mathbb{P}^1 \times \mathbb{P}^1$ or the blow-up of $r$ points 
$A_1,A_2,...,A_r \in \Pn$, $1 \leq r \leq 8$, in general position. (See \cite{bib:Dem} and \cite{bib:Be1} for more details.) We describe in Chapter $\refd{Chap:RankOne}$ the finite abelian subgroups $G \subset \Aut(S)$ when the group of invariant divisors is of rank $1$ (case $2$ of Proposition \refd{Prp:TwoCases}).

Here is an example: the $3$-torsion of the diagonal torus of $\PGL(4,\K)$, isomorphic to $(\Z{3})^3$, acting on the cubic surface in $\mathbb{P}^3$ with equation $w^3+x^3+y^3+z^3=0$.\index{Del Pezzo surfaces!automorphisms} This surface is a Del Pezzo surface of degree~3\index{Del Pezzo surfaces!of degree $3$|idi} (see Section \refd{Sec:Cubic}). Other famous examples are Geiser and Bertini involutions, acting minimally respectively on Del Pezzo surfaces of degree $2$ and $1$ (see Sections \refd{Sec:DP2} and \refd{Sec:DP1}).\index{Del Pezzo surfaces!of degree $2$|idi}\index{Del Pezzo surfaces!of degree $1$|idi}\index{Involutions!Bertini|idi}\index{Involutions!Geiser|idi}
\end{itemize}

\bigskip

\begin{remark}\fauxtitred
\begin{itemize}
\item
In the second case, no conic bundle structure can be invariant (for this would imply that both the canonical divisor and the fibre are invariant, so $\rkPic{S}^G>1$). However, in the first case, the underlying surface may be a Del Pezzo surface (see Section \refd{Sec:ConicBundlonDP}).
\item
Although the two cases are distinct, they are not birationally distinct (see Section \refd{Sec:PencilRatCurves}).
\end{itemize}
\end{remark}

\bigskip

After enumerating these two cases and their subcases, we have to decide whether two groups viewed as groups of automorphisms of different surfaces are birationally conjugate or not. This question is considered in Chapters \refd{Chap:ConjBetweenTheCases} and \refd{Chap:List}.

\chapter{Conjugacy invariants}
\label{Chap:ConjugacyInvariants}
\ChapterBegin{A conjugacy class of finite subgroups of the Cremona group has infinitely many representations as a group of birational maps of some rational surface.
In this chapter, we define some invariants of these representations which allow us in certain cases to prove that two groups acting birationally on some rational surface are not birationally conjugate (i.e.\ they represent distinct conjugacy classes of subgroups of the Cremona group).}

\ChapterNotation{In this chapter, $G\subset \Cr{S}$ denotes a finite abelian group of birational maps of the rational surface $S$.}

\section{Invariant pencils of rational curves}
\pagetitless{Conjugacy invariants}{Invariant pencils of rational curves}
\label{Sec:PencilRatCurves}
\emph{The first conjugacy invariant is the existence of a pencil of rational curves of the surface $S$ which is invariant by the group $G$.}
\index{Conic bundles!automorphisms|idi}
Note that this invariant is equivalent to the existence of a birational map $\varphi:S\dasharrow S'$ to some conic bundle $(S',\pi)$ such that $\varphi G \varphi^{-1}$ acts biregularly on the conic bundle (case $1$ of Proposition \refd{Prp:TwoCases}). This follows from Proposition \refd{Prp:BeauTsen}.

Note also that there exist some groups $G$ acting biregularly on a Del Pezzo surface $S$\index{Del Pezzo surfaces!automorphisms} with $\rkPic{S}^G=1$ (case $2$ of Proposition \refd{Prp:TwoCases}), which leave invariant a pencil of rational curves. The two cases of Proposition \refd{Prp:TwoCases} are therefore not birationally distinct.
\begin{itemize}
\item
For example, let $G$ be a finite diagonal group of automorphisms of $\mathbb{P}^2$. This group acts biregularly on $\mathbb{P}^2$ and $\rkPic{\Pn}^G=1$. The pair $(G,\Pn)$ is minimal and corresponds to case $2$ of Proposition \refd{Prp:TwoCases}. However, $G$ leaves invariant the pencil of rational curves passing through one  fixed point.
By blowing-up two fixed points and blowing-down the strict pull-back of the line through these two points, the group is conjugate to a group $G'$ of automorphisms of $\mathbb{P}^1\times\mathbb{P}^1$, leaving invariant the two rulings of the surface. Hence, $(G',\mathbb{P}^1\times\mathbb{P}^1)$ is a minimal pair and corresponds to case $1$ of Proposition \refd{Prp:TwoCases}.\end{itemize}

However, we give in Chapters \refd{Chap:ConicBundles} and \refd{Chap:FiniteAbConicBundle} some information on groups that preserve a conic bundle structure, which allows us to say that some groups are not in this case. Proposition \refd{Prp:ConjBetweenRankPic} explains the distinction between cases.
\section{Pairs of invariant pencils}
\pagetitless{Conjugacy invariants}{Pairs of invariant pencils}
\label{Sec:PairInvariantPencils}
\emph{A second invariant is the existence of two $G$-invariant pencils of rational curves on the surface $S$, with some given free intersection.}

Note that this is equivalent to the existence of a rational map $\varphi:S \dasharrow \mathbb{P}^1\times\mathbb{P}^1$ of some given degree, which is $G$-equivariant, where $G$ acts on the target as a subgroup of $\PGL(2,\K)\times \PGL(2,\K)$.

For example, the existence of two pencils of rational curves on the surface $S$ with free intersection $1$ is equivalent to being conjugate to a subgroup of $\PGLn{2}\times\PGLn{2} \subset \Aut(\mathbb{P}^1\times\mathbb{P}^1)$. In Proposition \refd{Prp:AutP1Bir} we use this invariant to show that the following two groups $G_1,G_2 \subset \Aut(\mathbb{P}^1\times\mathbb{P}^1)$ are not birationally conjugate:
\begin{center}$\begin{array}{lll}
G_1=<\alpha_1,\beta_1>\cong \Z{4}\times{\Z{2}},&\mbox{ where }&\alpha_1:(x,y) \mapsto (-x,\im y),\ \ \beta_1:(x,y) \mapsto (x^{-1},y),\vspace{0.1 cm}\\ 
G_2=<\alpha_2,\beta_2>\cong \Z{4}\times{\Z{2}},&\mbox{ where }&\alpha_2:(x,y) \mapsto (-y,x),\ \ \beta_2:(x,y) \mapsto (x^{-1},y^{-1}).\end{array}$
\end{center}

\section{The existence of fixed points}
\label{Sec:ExistenceFixedPoints}
\pagetitless{Conjugacy invariants}{The existence of fixed points}
\emph{A third invariant is the existence of points of $S$ fixed by $G$.}

Indeed, let $G \subset \Aut(S)$ and $G'\subset \Aut(S')$ be two finite abelian groups. Suppose that $G$ fixes at least one point of $S$ and that $G$ and $G'$ are conjugate by some birational map $\varphi:S\dasharrow S'$.
Then, by \cite{bib:KoS}, Proposition A.2, the group $G'$ fixes at least one point of $S'$.

For example, we use this invariant to see that the two groups of automorphisms of $\mathbb{P}^1\times\mathbb{P}^1$ \begin{center}$(x,y) \mapsto (\pm x^{\pm 1},y)$\hspace{1 cm}and\hspace{1 cm}$(x,y) \mapsto (\pm x,\pm y)$\end{center} are not birationally conjugate (see Proposition \refd{Prp:AutP1Bir}).

In Lemma \refd{Lem:Cs24notconjugate}, we also use this invariant (on the groups and their subgroups) to see that the following group $G_3$ is not conjugate to the groups $G_1$ and $G_2$ of Section \refd{Sec:PairInvariantPencils}:
\begin{center}$\begin{array}{l}
G_3=<\alpha_3,\beta_3>\cong \Z{4}\times{\Z{2}},  \mbox{ where}\\
\alpha_3:(x:y:z)\dasharrow (yz:xy:-xz),\\
\beta_3:(x:y:z) \dasharrow (x(y+z):z(y-z):-y(y-z)).\end{array}$
\end{center}
\section{The non-rational curves of fixed points}\index{Curves!of positive genus!birational maps that fix a curve of positive genus}
\label{Sec:CurveFixedPoints}
\pagetitless{Conjugacy invariants}{The non-rational curves of fixed points}
We use the notion of \textit{normalized fixed curve} defined in \cite{bib:BaB}:\index{Curves!of positive genus!n\vphantom{n}ormalized fixed curve (NFC)|idb}

Let $\varphi \in \Cr{S}$ be a birational map that fixes a curve $C$ (i.e.\ that fixes every point of $C$). 
Then, for any birational map $\alpha:S\dasharrow S'$, the birational map $\alpha \varphi\alpha^{-1}\in \Cr{S'}$ fixes the strict transform of $C$ by $\alpha$, which is a curve isomorphic to $C$, 
except possibly if $C$ is rational - a case in which the curve may be contracted to a point.
We define the \defn{normalized fixed curve} of $\varphi$ to be the union of the normalizations of the non-rational curves fixed by $\varphi$. This is an invariant of the conjugacy class of $\varphi$ in the Cremona group.
\begin{remark}\fauxtitred
In fact, using the classification of plane Cremona maps of prime order (see 
\cite{bib:BaB} and \cite{bib:deF}), we see that any  such map may fix at most one non-rational curve. Therefore, the same property holds for any element of finite order. 
\end{remark}

We extend this definition slightly to groups and elements with no non-rational fixed curve.

\break

\begin{Def} \fauxtitred
\label{Def:IsoFixedPoints}
\begin{itemize}
\item
Let $\mathcal{C}$ denote the set of isomorphism classes of non-rational curves.
We define the map \begin{center}$NFC:\Cr{S}\rightarrow \mathcal{C} \cup \{\emptyset\}$\end{center}
by associating to a birational map its normalized fixed curve (or the empty set $\emptyset$ if the birational map does not fix any non-rational curve).

\item
We say that two groups $G \subset \Cr{S}$ and $G' \subset \Cr{S'}$ \defn{have isomorphic sets of fixed points} if there is an isomorphism $\gamma:G\mapsto G'$ such that $NFC(g)=NFC(\gamma(g))$ for every element $g \in G$.
\end{itemize}
\end{Def}

The map $NFC$ induces a new invariant:
\begin{Prp}\fauxtitre
\label{Prp:NFC}
Let $G \subset \Cr{S}$ and $G' \subset \Cr{S'}$ be two finite subgroups of birational maps. If $G$ and $G'$ are birationally conjugate, they have isomorphic sets of fixed points. The converse is true for cyclic groups of prime order, but not in general.
\end{Prp}\begin{proof}\fauxtitred\upshape
Suppose that $G$ and $G'$ are conjugate by a birational map $\varphi:S\dasharrow S'$. It is clear that $\varphi$ gives rise to an isomorphism from $G$ to $G'$ defined by $g \mapsto \varphi g \varphi^{-1}$. We have $NFC(\varphi g \varphi^{-1})=\varphi(NFC(g))$ for any $g\in G$; the groups $G$ and $G'$ therefore have isomorphic sets of fixed points.

The converse is not true in general, for example:
\begin{itemize}
\item
The groups $\nump{P1s.222}$ and $\nump{P1.222}$ of Proposition \refd{Prp:AutP1Bir} are both isomorphic to $(\Z{2})^3$ and have no element which fixes a non-rational curve, but these two groups are not birationally conjugate.\index{Curves!of positive genus!birational maps that fix a curve of positive genus!involutions|idb}
\end{itemize}
However, it is true for cyclic groups of prime order:
\begin{itemize}
\item
Two elements of the Cremona group of the same prime order $p$ that fix the same non-rational curve are conjugate. (See \cite{bib:BaB}, Proposition 2.7 and \cite{bib:deF}, Theorem F.)
\item
An element of the Cremona group of prime order $p$ that does not fix a non-rational curve is birationally conjugate to a linear automorphism of\hspace{0.2 cm}$\Pn$ of order $p$. (See \cite{bib:BaB}, Corollary 1.2.)
\item
Two linear automorphisms of\hspace{0.2 cm}$\Pn$ of the same order are birationally conjugate. (See \cite{bib:BaB}, Proposition 2.1.)\proofend
\end{itemize}
\end{proof}

\section{The action on the non-rational fixed curves}\index{Curves!of positive genus!NFC - the action on it|idb}
\label{Sec:ActionFixedCurves}
\pagetitless{Conjugacy invariants}{The action on the non-rational fixed curves}
Since $G$ is abelian, every fixed curve of an element of $G$ is invariant by the whole group.

Let $g \in G \subset \Cr{S}$. We denote by $Fix(g) \subset S$ the set of points of $S$ fixed by $g$. If $Fix(g)$ contains some non-rational curve, its normalisation will be denoted by $\widetilde{Fix(g)}$. Note that $G$ acts birationally (and therefore biregularly) on $\widetilde{Fix(g)}$.

We now give a new, more precise invariant:
\begin{Def} \fauxtitre
\label{Def:actionFixedPoints}
Let $G \subset \Cr{S}$ and $G' \subset \Cr{S'}$ be two finite subgroups of birational maps. We say that these groups  \defn{have the same action on their sets of fixed points} if the two groups \defn{have isomorphic sets of fixed points} and the actions of $G$ and $\gamma(G)$ on the normalized fixed curves (if these are not empty) are the same.

Explicitly, this means that there exist an isomorphism $\gamma:G\rightarrow G'$ and for any $g\in G$,  an isomorphism $\gamma^{*}_g:\widetilde{Fix(g)}\rightarrow \widetilde{Fix(\gamma(g))}$ such that the following diagram commutes:

\begin{center}
\hspace{0 cm}\xymatrix{ G \times \widetilde{Fix(g)} \ar^{\alpha}[rr] \ar_{\gamma \times \gamma^{*}_g}[d] & &\widetilde{Fix(g)} \ar_{\gamma^{*}_g}[d] \\
G' \times \widetilde{Fix(\gamma(g))} \ar^{\alpha'}[rr] & & \widetilde{Fix(\gamma(g)),}}
\end{center}
where $\alpha$ (respectively $\alpha'$) denotes the action of $G$ on $\widetilde{Fix(g)}$ (respectively that of $G'$ on $\widetilde{Fix(\gamma(g)))}$.

\end{Def}

The following proposition is a direct consequence of Definition \refd{Def:actionFixedPoints}:
\begin{Prp}\fauxtitre
\label{Prp:NFCAction}
Let $G$ and $G'$ be two finite abelian subgroups of $\CrP$. If $G$ and $G'$ are conjugate in $\CrP$, they have same action on their sets of fixed points. The converse is true for finite cyclic groups of prime order, but not in general.\proofend
\end{Prp}
\bigskip

Note that, in most cases, the invariants defined above determine  the conjugacy classes of finite abelian subgroups of $\CrP$. There are however a few cases of non-conjugate groups having the same invariants. See for example Proposition \refd{Prp:conjisoA1A2}.

\section{More tools}
\label{Sec:MoreTools}
\pagetitless{Conjugacy invariants}{More tools}
To decide that two groups are not conjugate, we can also use the classification of $G$-invariant elementary links given in \cite{bib:Isk5}, as does \cite{bib:Dol}. This tool is often very useful. 

When $\rkPic{S}^G=1$, the classification says that the isomorphism class of the surface $S$ determines the minimal pair $(G,S)$ in many cases. For example, if $(G,S)$ is birationally conjugate to $(G',S')$, if $\rkPic{S}^G=\rkPic{S'}^{G'}=1$ and $S$ is a Del Pezzo surface of degree $1$, $2$ or $3$, then $S'\cong S$.\index{Del Pezzo surfaces!of degree $1$|idi}\index{Del Pezzo surfaces!of degree $2$|idi}\index{Del Pezzo surfaces!of degree $3$|idi}

However this is not true in general. For example, let $S$ be the Del Pezzo surface of degree $6$ \index{Del Pezzo surfaces!automorphisms}\index{Del Pezzo surfaces!of degree $6$|idi}defined by $S=\{ (x:y:z) \times (u:v:w)\ | \ ux=vy=wz\} \subset \Pn\times\Pn$ (see Section \refd{Subsec:DelPezzo6}). Let $S'=\Pn$ and let $(G,S)$ and $(G,S')$ be the two pairs defined thus:
\begin{center}$\begin{array}{rll}
G=<\alpha>\cong \Z{6},&\mbox{ where }&\alpha:(x:y:z)\times(u:v:w) \mapsto (v:w:u)\times(y:z:x),\\
G'=<\alpha'>\cong \Z{6},&\mbox{ where }&\alpha':(x:y:z) \mapsto (-\omega x:y:z).\end{array}$
\end{center}
Then $\rkPic{S}^{G}=1$ (see Section \refd{Subsec:DelPezzo6}) and $\rkPic{S'}^{G'}=1$, as $\rkPic{S'}=1$. Although the two pairs are birationally conjugate (see Proposition \refd{Prp:DP6conjLin}), $S\not\cong S'$.

\chapter{Del Pezzo surfaces}
\label{Chap:DelPezzoSurfaces}\index{Del Pezzo surfaces|idb}
\emph{In this chapter we describe the Del Pezzo surfaces, which often appear in this paper. In particular we describe their Picard groups, exceptional divisors and conic bundle structures.}
\section{The classification}
\label{Sec:DelPezzoExceptDivisors}
\pagetitless{Del Pezzo surfaces}{The classification}\index{Del Pezzo surfaces!classification|idb}
According to the definition given in Chapter \refd{Chap:GrAutRat}, a smooth rational surface $S$ is a Del Pezzo surface if the anti-canonical divisor $-K_S$ is ample. 
We recall the following result: 
\begin{Prp} \titrePropRef{The adjunction formula}{\cite{bib:Be1}, I.15}\\
\label{Prp:AdjunctionFormula}
Let $C \subset S$ be an irreducible curve on a surface $S$. We have $C \cdot (C+K_S)=-2+2\cdot g(C)$, where $g(C)=H^1(C,\mathcal{O}_C)$ is the (arithmetical) genus of $C$.
\end{Prp}

\bigskip

We first prove a simple lemma, which will help us to explicit the classification of Del Pezzo surfaces in Proposition \refd{Prp:ConditionsDelPezzo}.

\begin{Lem}\fauxtitre
\label{Lem:SelfIntersgeqOne}
Let $S$ be a rational surface such that any irreducible curve of $S$ has self-intersection $\geq-1$. 

Then, $S \cong \mathbb{P}^1\times\mathbb{P}^1$ or $S\cong \Pn$ or $S$ is the blow-up of $r\geq 1$ points of\hspace{0.2 cm}$\Pn$ in general position  (i.e.\ no $3$ are collinear, no $6$ are on the same conic, no $8$ lie on a cubic having a double point at one of them).
\end{Lem}
\begin{proof}\fauxtitre\upshape
By blowing-down some curves, we get a birational morphism $S\rightarrow S'$ to some minimal rational surface $S'$ which is isomorphic to $\Pn$ or to some Hirzebruch surface $\mathbb{F}_n$, with $n\not= 1$ (see Proposition \refd{Prp:MinimalRatSurfaces}). Note that no irreducible curve of $S'$ has self-intersection $\leq -2$, since this is the case for $S$. As $\mathbb{F}_n$ contains the section $E_n$ which has self-intersection $-n$, we have $S'=\Pn$ or $S'=\mathbb{P}^1\times\mathbb{P}^1$. If $S=S'$, we are done, otherwise we may choose $S'=\Pn$ (the blow-up of one point of $\mathbb{P}^1\times\mathbb{P}^1$ is isomorphic to the blow-up of two distinct points of\hspace{0.2 cm}$\Pn$). Therefore, $S$ is the blow-up of some (perhaps infinitely near) points on $\Pn$. In fact, all points must lie in $\Pn$, since no irreducible curve of $S$ has self-intersection $\leq -2$.

Furthermore, the points blown-up are in general position, as otherwise a curve of self-intersection $\leq -2$ would lie on $S$. \proofend
\end{proof}

\bigskip

\begin{Prp}\fauxtitre
\label{Prp:ConditionsDelPezzo}
Let $S$ be a rational surface. The following conditions are equivalent:\index{Del Pezzo surfaces!classification}
\begin{itemize}
\item[\upshape 1.]
$S$ is a Del Pezzo surface;
\item[\upshape 2.]
$S \cong \mathbb{P}^1\times\mathbb{P}^1$ or $S\cong \Pn$ or $S$ is the blow-up of $1\leq r \leq 8$ points of\hspace{0.2 cm}$\Pn$ in general position (i.e.\ no $3$ are collinear, no $6$ are on the same conic, no $8$ lie on a cubic having a double point at one of them);
\item[\upshape 3.]
$K_S^2\geq 1$ and any irreducible curve of $S$ has self-intersection $\geq -1$;
\item[\upshape 4.]
$C\cdot (-K_S)>0$ for any effective divisor $C$.
\end{itemize}
\end{Prp}
\begin{remark}\fauxtitred
In general a divisor is ample if and only if it intersects positively the adherence of the cone of effective divisors (Kleiman's ampleness criterion, \cite{bib:KoM}, Theorem I.18).

In our case (when the surface is rational and the divisor is the anti-canonical divisor), the equivalence of assertions $1$ and $4$ shows that the criterion is true, even if we omit the adherence in the statement.
\end{remark}
\begin{proof}\fauxtitred\upshape
\begin{itemize}
\item[$(1 \Rightarrow 4)$]
As $-K_S$ is ample, $-mK_S$ is very ample, for some integer $m>0$, in which case\\ $m\cdot (-K_S)\cdot C$ is the degree of $C$ in the corresponding embedding, which must be positive.
\item[$(4 \Rightarrow 2,3)$]
We first prove that assertion $4$ implies that any irreducible curve of $S$ has self-intersection $\geq -1$. Suppose that some irreducible curve $C$ of $S$ has self-intersection $\leq -2$. The adjunction formula (Proposition \refd{Prp:AdjunctionFormula}) gives $C \cdot (C+K_S)=-2+2\cdot g(C)\geq -2$, whence $C\cdot (-K_S) \leq 2+C^2\leq 0$, which contradicts assertion 4.

Secondly, using Lemma \refd{Lem:SelfIntersgeqOne}, 
$S \cong \mathbb{P}^1\times\mathbb{P}^1$ or $S\cong \Pn$ or $S$ is the blow-up of $r\geq 1$ points of\hspace{0.2 cm}$\Pn$ in general position.
 Furthermore, if the number of blown-up points is at least $9$, there exists a cubic passing through $9$ of the blown-up points, which is irreducible as the points are in general position. The strict transform of this curve intersects the anti-canonical divisor of $S$ non positively. The number of blown-up points is then at most $8$, and so $K_S^2\geq 1$. We get assertions $2$ and $3$.
\item[$(3 \Rightarrow 2)$]
Using once again Lemma \refd{Lem:SelfIntersgeqOne} we see that, $S \cong \mathbb{P}^1\times\mathbb{P}^1$ or $S\cong \Pn$ or $S$ is the blow-up of $r\geq 1$ points of\hspace{0.2 cm}$\Pn$ in general position. Furthermore, $K_S^2\geq 1$ implies that $r\leq 8$.
\item[$(2 \Rightarrow 1)$]
In \cite{bib:Dem}, Theorem 1, it is proved that the blow-up of $1\leq r \leq 8$ points in general position gives a Del Pezzo surface. The cases of\hspace{0.2 cm}$\Pn$ and $\mathbb{P}^1\times\mathbb{P}^1$ are clear.
\proofend
\end{itemize}
\end{proof}

\begin{Cor}\fauxtitre
\label{Prp:DelPezzoContractDP}\index{Del Pezzo surfaces!others}
Let $S$ be a Del Pezzo surface and let $\varphi:S\rightarrow S'$ be a birational morphism. Then $S'$ is a Del Pezzo surface.
\end{Cor}
\begin{proof}\fauxtitred\upshape
This follows directly from Proposition \refd{Prp:ConditionsDelPezzo}, using the equivalence of assertions $1$ and $3$.\proofend
\end{proof}

\break

We can now explain why case $2$ of Proposition \refd{Prp:TwoCases} imposes that the surface is of Del Pezzo type:
\begin{Lem}\fauxtitre
\label{Lem:rkDelPezzo}\index{Del Pezzo surfaces!others}
Let $G \subset \Aut(S)$ be a finite subgroup of automorphisms of a rational smooth surface $S$. If $\rkPic{S}^G=1$, the surface is a Del Pezzo surface and
\begin{center} $\Pic{S}^G=
\left\{\begin{array}{rl}\frac{1}{3}\mathbb{Z}K_S& \mbox{if }S\cong\Pn\vspace{0.1 cm}\\
\frac{1}{2}\mathbb{Z}K_S & \mbox{if }S\cong \mathbb{P}^1\times\mathbb{P}^1\vspace{0.1 cm}\\
\mathbb{Z}K_S & \mbox{otherwise.}\end{array}\right.$\end{center}
\end{Lem}
\begin{proof}\fauxtitred\upshape
To show that $S$ is a Del Pezzo surface, it suffices to show that $-K_S$ intersects positively any irreducible curve of $S$  (see Proposition \refd{Prp:ConditionsDelPezzo}).
As the canonical divisor is $G$-invariant, we have $\Pic{S}^{G}\otimes_{\mathbb{Z}} \mathbb{Q}=\mathbb{Q} K_S$. Then, for any effective divisor $C$, the divisor $\sum_{g \in G} g(C)$ belongs to $\Pic{S}^{G}$, so $\hts\sum_{g \in G} g(C)=-aK_{S}$ for some $a\in \mathbb{Q}$. And $a$ is positive, since $K_S$ is not effective. \emph{This implies that $-K_S$ is effective.}

From the relation $g(C)\cdot (-K_S)=g(C)\cdot g(-K_S)=C\cdot (-K_{S})$ for every $g\in G$, we deduce the formula:
\begin{center}$C\cdot (-K_{S})=\frac{1}{|G|} \sum_{g \in G} g(C)\cdot(-K_S)=a\frac{1}{|G|}K_S^2$.\end{center} 
The sign of $C\cdot (-K_{S})$ is thus the same, for any effective divisor $C$. As $-K_S$ is effective, this must be positive, so $S$ is a Del Pezzo surface (by Proposition \refd{Prp:ConditionsDelPezzo}).

We now prove the last assertion. The relation $\Pic{S}^{G}\otimes_{\mathbb{Z}} \mathbb{Q}=\mathbb{Q} K_S$ given above implies that $\Pic{S}^G=\mathbb{Q} K_S \cap \Pic{S}$. Note that if $S\cong \Pn$, the divisor class of a line is $-\frac{1}{3}K_S$, so $\mathbb{Q} K_S \cap \Pic{S}=\frac{1}{3}\mathbb{Z}K_S$. If $S\cong\mathbb{P}^1\times\mathbb{P}^1$, we get similarly $\mathbb{Q} K_S \cap \Pic{S}=\frac{1}{2}\mathbb{Z}K_S$, since $K_{\mathbb{P}^1\times\mathbb{P}^1}$ is the double of the diagonal divisor in $\Pic{\mathbb{P}^1\times\mathbb{P}^1}$. Otherwise, there exists some birational morphism $\pi:S\rightarrow \Pn$, which is the blow-up of $r\geq 1$ points of\hspace{0.2 cm}$\Pn$. Then, $K_S=-3L+\sum_{i=1}^r E_i$, where $L$ denotes the pull-back by $\pi$ of a line of\hspace{0.2 cm}$\Pn$ and $E_1,...,E_r$ are the exceptional divisors  (see further explanations below). Observe that $K_S$ is not a multiple of some divisor, so $\Pic{S}^G=\mathbb{Z}K_S$.\proofend
\end{proof}

\bigskip

\section{The Picard group and the exceptional divisors}
\label{Sec:DPPicadExc}
\pagetitless{Del Pezzo surfaces}{The Picard group and the exceptional divisors}\index{Del Pezzo surfaces!curves on the surfaces|idb}
We proceed to describe the Picard group and the exceptional curves of Del Pezzo surfaces (see \cite{bib:Dem} and \cite{bib:Be1} for more details).

The surface $\Pn$ (respectively $\mathbb{P}^1\times\mathbb{P}^1$) is minimal, and every effective divisor has positive self-intersection (respectively non-negative).
Now consider the other cases: we denote by $A_1,...,A_r\in \Pn$, $1 \leq r \leq 8$, $r$ points in general position and by $\pi: S \rightarrow \Pn$ their blow-up. Let $L$ denote the pull-back by $\pi$ of a line of\hspace{0.2 cm}$\Pn$;
the Picard group of $S$ is generated by $L$ and $E_1=\pi^{-1}(A_1),..., E_r=\pi^{-1}(A_r)$. 
The intersection form on $\Pic{S}$ is given by $L^2=1$, $L\cdot E_i =0$, $E_i^2=-1$, $i=1,...,r$, $E_i \cdot E_j=0$ for $i \not= j$. The canonical divisor is $K_{S}=-3L+\sum_{i=1}^r E_i$. 
By definition, the degree of the surface is the self-intersection of the canonical divisor, which is $9-r$.

\break

To simplify the following proofs, we gather some information on Del Pezzo surfaces of degree $1$:
\index{Del Pezzo surfaces!of degree $1$}
Let $S$ be a Del Pezzo surface of degree $1$. The linear sytem $|-2K_S|$ induces a degree $2$ morphism onto a quadric cone $Q \subset \mathbb{P}^3$. The involution $\sigma$ of $S$ induced by this $2$-covering is called the \defn{Bertini involution}.\index{Involutions!Bertini|idb}
\begin{Lem}\fauxtitre
\label{Lem:2KSmD}
Let $S$ be a Del Pezzo surface of degree $1$ and let $C \subset S$ be an exceptional curve. Then,
\begin{itemize}
\item
the Bertini involution sends $C$ on the curve $-2K_S-C$;
\item
in particular, $-2K_S-C$ is an exceptional curve.
\end{itemize}
\end{Lem}
\begin{proof}\fauxtitred\upshape
Let $\sigma$ denote the Bertini involution. We have $\rkPic{S}^{\sigma}=1$ (this will be proved in Section \refd{Sec:DP1}, Lemma \refd{Lem:DP1}). The divisor $C+\sigma(C)$ is then a multiple of the canonical divisor. The adjunction formula gives $C \cdot (C+K_S)=-2$, so $C\cdot K_S=-1$. We thus get $K_S \cdot (C + \sigma(C))= -2$, so $C+\sigma(C)=-2K_S$.\proofend
\end{proof}

\bigskip

The exceptional curves of Del Pezzo surfaces are described in the following proposition:
\index{Del Pezzo surfaces!curves on the surfaces|idb}
\begin{Prp}\fauxtitre
\label{Prp:DescExcCurvDelPezzo}
Let $\pi:S\rightarrow \Pn$ be the blow-up of $1\leq r\leq 8$ points $A_1,...,A_r\in \Pn$ in general position. 
\begin{itemize}
\item
The exceptional curves of $S$ are:
\begin{itemize}
\item
The $r$ exceptional curves corresponding to the total pull-backs $\pi^{-1}(A_1)$,...,$\pi^{-1}(A_r)$ of the blown-up points. 
\item
The strict pull-back of the curves of degree $d$ passing through the $A_i$'s with multiplicities given in the following table:
\begin{center}
$\begin{array}{|l|l|l|llllllll|}
\hline
r & \mbox{degree} & \mbox{multiplicities} & \multicolumn{7}{c}{\mbox{number of such curves for}}& \\
& d & \mbox{at the points} & r=1,& 2, & 3, & 4, & 5, & 6, & 7, & 8\\
\hline
\geq 2 & 1 & (1,1) & & 1 & 3 & 6 & 10 & 15 & 21 & 28\\
\geq 5 & 2 & (1,1,1,1,1) & & & & & 1 & 6 & 21 & 56\\
\geq 7 & 3 & (2,1,1,1,1,1,1) & & & & & & & 7 & 56\\
8 & 4 & (2,2,2,1,1,1,1,1) & & & & & & & & 56\\
8 & 5 & (2,2,2,2,2,2,1,1) & & & & & & & & 28\\
8 & 6 & (3,2,2,2,2,2,2,2) & & & & & & & & 8\\
\hline
\end{array}$\end{center}
\end{itemize}
\item
Except when $r=2$, the number of exceptional curves touching one given exceptional curve is constant.
Explicitly we have:
\begin{center}$\begin{array}{|c|cccccccc|}
\hline
\mbox{Number $r$ of blown-up points} & 1 & 2 & 3 & 4 & 5 & 6 & 7 & 8 \\
\mbox{Degree of the surface} & 8 & 7 & 6 & 5 & 4 & 3 & 2 & 1 \\ \hline
\mbox{Number of exceptional curves} & 1 & 3 & 6 & 10 & 16 & 27 & 56 & 240 \\ \hline
\mbox{Number of neighbours with intersection $1$} & & 1/2 & 2 & 3 & 5 & 10 & 27 & 126 \\
\mbox{Number of neighbours with intersection $2$} & & & & & & & 1 & 56 \\
\mbox{Number of neighbours with intersection $3$} & & & & & & & & 1 \\ \hline \end{array}$\end{center}
\end{itemize}
\end{Prp}
\begin{remark}\fauxtitred
The number of exceptional curves of a Del Pezzo surface is well known, see for example \cite{bib:Dem}, Table 3, page 35.
\end{remark}
\begin{proof}\fauxtitred\upshape
As before we denote by $E_i=\pi^{-1}(A_i)$ the exceptional curve corresponding to the blow-up of $A_i$ for $i=1,...,r$. Apart from these, the divisor class $D$
of an exceptional curve of $S$ is equal to $mL-\sum a_i E_i$ for some non-negative integers $m,a_1,...,a_r$, with $m>0$. The self-intersection $D^2=-1$ and the adjunction formula $D (K_S+D)=-2$ (see Proposition \refd{Prp:AdjunctionFormula}) give the following relations:
\begin{equation}
\label{equation:ExcCurves}
\begin{array}{rcl}\sum_{i=1}^r a_i^2&=&m^2+1,\vspace{0.1cm}\\
\sum_{i=1}^r a_i & = & 3m-1,\end{array}
\end{equation}
which imply that $\sum_{i=1}^r a_i(a_i-1)=(m-1)(m-2)$. 
\begin{itemize}
\item
Replacing $m$ by $1$ we find the $r\cdot (r-1)/2$ lines passing through two points. 
\item
Replacing $m$ by $2$, we find the {\small $ \left(\begin{array}{c}r \\ 5\end{array}\right)$} conics passing through $5$ of the points. 
\item
Replacing $m$ by $3$, we find the {\small$ r\cdot \left(\begin{array}{c}r-1 \\ 6\end{array}\right)$} cubics passing through one of the $A_i$'s with multiplicity $2$ and through $6$ others with multiplicity $1$.
\end{itemize}
Assume now that $r=8$ (if $r<8$, we blow-up some general points). In this case, $K_S^2=1$. Observe that if $D$ is an exceptional curve, then $-2K_S-D$ is also exceptional (see Lemma \refd{Lem:2KSmD}). 
Using this observation, we find the solutions for $m\geq 4$ from those for $6-m$.

It remains to prove the last assertion. 
\begin{itemize}
\item
If $r=1$, there is only one exceptional divisor, so it does not intersect any other.
\item
If $r=2$, there are three exceptional divisors: $E_1,E_2$, and $D_{12}$, the strict pull-back by $\pi$ of the line through $A_1$ and $A_2$. $D_{12}$ touches the other two, since $E_1$ and $E_2$ each have only one neighbour. Note that by blowing-down $E_1$ or $E_2$, we get $\mathbb{F}_1$, but by blowing-down $D_{12}$, we get a birational morphism to $\mathbb{P}^1\times\mathbb{P}^1$. 
\item
We observe that there exists only one Del Pezzo surface of degree $7$ (this is not the case for the degree $8$, where $\mathbb{F}_0$ and $\mathbb{F}_1$ are possible). This surface has three exceptional divisors, one of them touching the other two, these touching only one neighbour.

The blow-down of one of the $6$ exceptional curves on a Del Pezzo surface of degree $6$ gives the Del Pezzo surface of degree $7$. Any exceptional curve on a Del Pezzo surface of degree $6$ thus intersects two other exceptional curves.

By induction, we prove that  for $r>2$, the number of exceptional curves touching one given exceptional curve is constant and does not depend on the surface. Indeed, the blow-up of any exceptional curve gives a Del Pezzo surface of degree $10-r<8$. The intersection form on this surface does not depend on the chosen curve, so the neighbours of $S$ are the same for any exceptional divisor. 

Let us now count the neighbours of $E_1$. Firstly, $E_1$ does not touch any other $E_i$. It touches the strict pull-back of a curve $C$ of\hspace{0.2 cm}$\Pn$  if and only if this curve passes through $A_1$. The intersection between the divisors is the multiplicity of $C$ at $A_1$. We thus count the neighbours of $E_1$ and the intersection multiplicities using the table above:
\begin{itemize}
\item
$E_1$ intersects the pull-back of $r-1$ lines with multiplicity $1$.
\item
$E_1$ intersects the pull-back of {\small $\left(\begin{array}{c}r-1 \\ 4\end{array}\right)$} conics, with multiplicity $1$. This gives respectively $1,5,15$ and $35$ conics for $r=5,6,7,8$.
\item
If $r=7$, there are $6$ lines, $15$ conics and $6$ cubics passing through $A_1$ with multiplicity $1$. There is also one cubic which passes through $A_1$ with multiplicity $2$.
\item
For $r=8$, there are $7$ lines, $35$ conics, $42$ cubics, $35$ quartics and $7$ quintics passing through $E_1$ with multiplicity $1$. There are $7$ cubics, $21$ quartics, $21$ quintics and $7$ sextics passing through $E_1$ with multiplicity $2$. Finally, there is one sextic passing through $E_1$ with multiplicity $3$. \end{itemize}
Summing up, we get the announced results.\proofend
\end{itemize}
\end{proof}

\bigskip

Using this information, we are able to state an important result on Del Pezzo surfaces and their isomorphism classes:\index{Del Pezzo surfaces!isomorphism classes}
\begin{Prp}\fauxtitre
\label{Prp:IsoDPlinear}
Let $r\geq 2$ be some integer.

The map that associates to each set of $r$ points of\hspace{0.2 cm}$\Pn$ in general position the surface obtained by blowing-up these points is a bijection between the sets of $r$ points of\hspace{0.2 cm}$\Pn$ in general position, up to the action of $\Aut(\Pn)$, and the isomorphism classes of Del Pezzo surfaces of degree $9-r$.
\end{Prp}
\begin{proof}\fauxtitred
\upshape
First of all, it is clear that by taking two sets equivalent by an action of $\Aut(\Pn)$, and blowing-up the points, we get two isomorphic surfaces.

We now prove the converse. We take two isomorphic surfaces, obtained by blowing-up $r$ points of\hspace{0.2 cm}$\Pn$ and prove that the points are equivalent under the action of $\Aut(\Pn)$. 
We denote the blow-ups by $\pi:S\rightarrow \Pn$, and $\pi':S'\rightarrow \Pn$. By hypothesis, there exists an isomorphism $\varphi:S\rightarrow S'$.
Observe that the birational morphisms $\pi'\circ\varphi:S\rightarrow \Pn$ and $\pi:S\rightarrow \Pn$ are the blow-downs of two sets of $r$ skew exceptional curves of $S$. 
It remains then to prove the following assertion:
\vspace{0.2 cm}

{\it Let $R_E=\{E_1,...,E_r\}$ and $R_{E'}=\{E_1',...,E_r'\}$ be two sets of $r$ skew exceptional curves of $S$. The two sets of $r$ points obtained by the contraction of these two sets of curves are equivalent under the action of a linear automorphism of\hspace{0.2 cm}$\Pn$.}
\vspace{0.2 cm}

We prove this by induction on $r$:
\begin{itemize}
\item
 If $r\leq 4$, the assertion follows from the fact that all sets of four non collinear points of\hspace{0.2 cm}$\Pn$ are equivalent under the action of linear automorphisms.
 \item
 If $r\geq 5$ and  one exceptional curve belongs to the two sets contracted, we contract the curve, obtain a Del Pezzo surface of degree $10-r\leq 5$ and the result follows from the induction hypothesis. 
 
 \break
 
 \item
Otherwise, we take one exceptional curve in each set, say $E_1$ and $E_1'$. As $r\geq 5$, there exists an exceptional curve $F$ that intersects neither  $E_1$ nor $E_1'$.

Blowing-down $E_1$ and $F$, we get a Del Pezzo surface $S_F$ of degree $11-r\leq 6$. There exists then a birational morphism from $S_F$ to $\Pn$, so there exist $r-2$ exceptional curves $F_1,...,F_{r-2}$ on $S$ such that $R_F=\{E_1,F,F_1,...,F_{r-2}\}$ is a set of $r$ skew exceptional curves. We proceed similarly with $E_1'$ and $F$ and obtain a set of $r$ skew exceptional curves $R_{F'}=\{E_1',F,F_1',...,F_{r-2}'\}$. As $E_1$ belongs to $R_E$ and $R_F$, the blow-down of these two sets of curves yields two equivalent sets of $r$ points of\hspace{0.2 cm}$\Pn$. As $F$ belongs to $R_F$ and $R_{F'}$, and $E_1'$ belongs to $R_{F'}$ and  $R_{E'}$, the blow-down of any of the four sets of curves gives an equivalent set of $r$ points of\hspace{0.2 cm}$\Pn$.\proofend
\end{itemize}
\end{proof}

\bigskip

In case $2$ of Proposition \refd{Prp:TwoCases} (when $\rkPic{S}^G=1$), the following lemma (whose proof uses the idea of \cite{bib:deF}, Proposition 4.1.4) restricts the possibilities for the surface $S$ and the group $G$:
\begin{Lem}\titreProp{Size of the orbits}\\
\label{Lem:SizeOrbits}\index{Del Pezzo surfaces!automorphisms}\index{Del Pezzo surfaces!curves on the surfaces}
Let $S$ be a Del Pezzo surface, which is the blow-up of $1\leq r\leq 8$ points of\hspace{0.2 cm}$\mathbb{P}^2$ in general position,
and let $G \subset \Aut(S)$ be a finite subgroup of automorphisms with $\rkPic{S}^{G}=1$. Then:
\begin{itemize}
\item
$G\not=\{1\}$;
\item
the size of any orbit of the action of $G$ 
on the set of exceptional divisors is divisible by the degree of $S$, which is $9-r$;
\item
in particular, the order of $G$ is divisible by the degree of $S$.
\end{itemize}
\end{Lem}\begin{proof}\fauxtitred\upshape
It is clear that $G \not=\{1\}$, since $\rkPic{S}>1$.
Let $D_1,D_2,...,D_k$ be $k$ exceptional divisors of $S$, forming an orbit of $G$. 
The divisor $\sum_{i=1}^{k} D_i$ is fixed by $G$ and thus is a multiple of $K_S$. 
We can write $\sum_{i=1}^{k} D_i=a K_S$, for some rational number $a \in \mathbb{Q}$. In fact, by Lemma 
\refd{Lem:rkDelPezzo}, we have $a \in \mathbb{Z}$ and $a<0$ as $K_S$ is not effective.
Since the $D_i$'s are irreducible and rational, we deduce from the adjunction formula $D_i (K_S +D_i)=-2$ (Proposition \refd{Prp:AdjunctionFormula}) that $D_i\cdot K_S=-1$.
Hence \begin{center}$K_S \cdot \sum_{i=1}^{k} D_i=\sum_{i=1}^{k}K_S \cdot D_i=-k=K_S \cdot a K_S=a (9-r)$.\end{center} Consequently, the degree $9-r$ divides the size $k$ of the orbit. \proofend
\end{proof}

\bigskip
\begin{remark}\fauxtitred
This lemma shows in particular that a pair $(G,S)$ in which $S$ is the blow-up of $r=1,2$ points of\hspace{0.2 cm}$\Pn$ is not minimal, a result which is obvious when $r=1$, and is clear when $r=2$, since the line joining the two blown-up points  is invariant by any automorphism.
\end{remark}

\section{Conic bundle structures on Del Pezzo surfaces}
\label{Sec:ConicBundlonDP}\index{Del Pezzo surfaces!conic bundles on the surfaces|idb}\index{Conic bundles!on Del Pezzo surfaces|idb}
\pagetitless{Del Pezzo surfaces}{Conic bundle structures on Del Pezzo surfaces}
\emph{We suggest that the reader unfamiliar with conic bundles read Section \refd{Sec:DescConicBundles} before this one. We prefer to keep this section here as it follows naturally on Section \refd{Sec:DPPicadExc}}

\bigskip

Note that except for the surface $\Pn$, every Del Pezzo surface has at least one conic bundle structure. For $\mathbb{F}_0=\mathbb{P}^1\times \mathbb{P}^1$ and $\mathbb{F}_1$ (the blow-up of a point in $\Pn$), the fibration is smooth. For other surfaces, we have several choices for blowing-down some exceptional divisors to get $\mathbb{F}_1$, which give different conic bundles, with at most $7$ singular fibres. The number of conic bundle structures depends on the degree (see Proposition \refd{Prp:NbConicDP}). The situation can be described as follows:

\begin{Prp}\fauxtitre
\label{Prp:ConicBundleDPLines}
Let $S$ be a Del Pezzo surface and let $\pi:S\rightarrow \mathbb{P}^1$ be a morphism which induces a conic bundle structure on $S$. Then, one (and only one) of the following occurs:
\begin{itemize}
\item
$S=\mathbb{P}^1\times\mathbb{P}^1$ and $\pi$ is one of the two projections on a factor.
\item
There exists a birational morphism (which is not an isomorphism) from $S$ to $\Pn$ that sends the fibres of $\pi$ on the pencil of lines passing through some point of\hspace{0.2 cm}$\Pn$.
\end{itemize}
\end{Prp}
\begin{proof}\fauxtitred\upshape
First, we suppose that $S=\mathbb{P}^1\times\mathbb{P}^1$. Let $f_1$ (respectively $f_2$) denote the divisor class of the fibre of the projection on the first factor (respectively on the second). As $f_1$ and $f_2$ generate $\Pic{S}$, we can write the divisor class of the fibre of $\pi$ as $f=af_1+bf_2$. Recall that $K_S=-2f_1-2f_2$, $f_1^2=f_2^2=0$ and $f_1\cdot f_2=1$. The conditions $f^2=0$ (two fibres do not intersect) and $f\cdot (f+K_S)=f\cdot K_S=-2$ (adjunction formula, Proposition \refd{Prp:AdjunctionFormula}) imply that $ab=0$ and $a+b=1$. We therefore have $f=f_1$ or $f=f_2$ and are in the first case of the proposition.

Note that every effective divisor of\hspace{0.2 cm}$\Pn$ has positive self-intersection. This surface therefore has no conic bundle structure.

Suppose now that $S\not\cong \mathbb{P}^1\times\mathbb{P}^1$, $S\not\cong \Pn$. 
As the surface has a conic bundle structure, the blow-down of 
one component in each singular fibre produces a minimal ruled surface, which is also a Del Pezzo surface (see Proposition \refd{Prp:DelPezzoContractDP}), 
and that must be $\mathbb{F}_0=\mathbb{P}^1\times\mathbb{P}^1$ or $\mathbb{F}_1$. We may suppose that it is $\mathbb{F}_1$ (if it is $\mathbb{F}_0$, we change the choice of the component in one singular fibre to get $\mathbb{F}_1$). 
Let $p$ denote the
birational morphism $p:S\rightarrow \mathbb{F}_1$ (observe that $p$ is compatible with the conic bundle structure of $S$ and $\mathbb{F}_1$). Blowing-down the unique exceptional curve of $\mathbb{F}_1$, we obtain the announced morphism.\proofend
\end{proof}

\bigskip

This proposition implies that $\Pn$ has no conic bundle structure and that $\mathbb{P}^1\times\mathbb{P}^1$ has two conic bundle structures. For the other surfaces, the birational morphism on $\Pn$ given in the proposition is not unique, so there are several conic bundle structures. We now give the precise description of such conic bundle structures in this case:
\begin{Prp}\fauxtitre
\label{Prp:NbConicDP}
Let $\pi:S\rightarrow \Pn$ be the blow-up of $1\leq r\leq 8$ points $A_1,...,A_r\in \Pn$ in general position. 
\begin{itemize}
\item
The conic bundle structures on $S$ are determined by divisors $f$ corresponding to the pull-back by $\pi$ of pencils of curves of\hspace{0.2 cm}$\Pn$ of degree $d$ passing through the $A_i$'s with multiplicities given in the following table:
\begin{center}
$\begin{array}{|l|l|l|rlllllll|}
\hline
r & \mbox{degree} & \mbox{multiplicities} & \multicolumn{8}{c|}{\mbox{number of such conic bundles for}}\\
& d & \mbox{at the points} & r=1,& 2, & 3, & 4, & 5, & 6, & 7, & 8\\
\hline
\geq 1 & 1 & (1) & 1 & 2 & 3 & 4 & 5 & 6 & 7 & 8\\
\geq 4 & 2 & (1,1,1,1) & & & & 1 & 5 & 15 & 35 & 70\\
\geq 6 & 3 & (2,1,1,1,1,1) & & & & & & 6 & 42 & 168\\
\geq 7 & 4 & (2,2,2,1,1,1,1) & & & & & & & 35 & 280\\
8 & 4 & (3,1,1,1,1,1,1,1) & & & & & & & & 8\\
\geq 7 & 5 & (2,2,2,2,2,2,1) & & & & & & & 7 & 56\\
8 & 5 & (3,2,2,2,1,1,1,1) & & & & & & & & 280\\
8 & 6 & (3,3,2,2,2,2,1,1) & & & & & & & & 420\\
8 & 7 & (3,3,3,3,2,2,2,1) & & & & & & & & 280\\
8 & 7 & (4,3,2,2,2,2,2,2) & & & & & & & & 56\\
8 & 8 & (3,3,3,3,3,3,3,1) & & & & & & & & 8\\
8 & 8 & (4,3,3,3,3,2,2,2) & & & & & & & & 280\\
8 & 9 & (4,4,3,3,3,3,3,2) & & & & & & & & 168\\
8 & 10 & (4,4,4,4,3,3,3,3) & & & & & & & & 70\\
8 & 11 & (4,4,4,4,4,4,4,3) & & & & & & & & 8\\
\hline
\end{array}$\end{center}
\item
The number of conic bundle structures is then:
\begin{center}$\begin{array}{|c|cccccccc|}
\hline
\mbox{Number $r$ of blown-up points} & 1 & 2 & 3 & 4 & 5 & 6 & 7 & 8 \\
\mbox{Degree of the surface} & 8 & 7 & 6 & 5 & 4 & 3 & 2 & 1 \\
\mbox{Number of conic bundle structures} & 1 & 2 & 3 & 5& 10 & 27 & 146 & 2160 \\ \hline \end{array}$\end{center}
\end{itemize}
\end{Prp}\begin{proof}\fauxtitred\upshape
Let $L$ denote the pull-back by $\pi$ of a line of\hspace{0.2 cm}$\Pn$; recall that
the Picard group of $S$ is generated by $L$ and $E_1=\pi^{-1}(A_1),..., E_r=\pi^{-1}(A_r)$ and that the intersection form on $\Pic{S}$ is given by $L^2=1$, $L\cdot E_i =0$, $E_i^2=-1$, $i=1,...,r$, $E_i \cdot E_j=0$ for $i \not= j$.

To any conic bundle structure on $S$ we associate the effective divisor $f$ of its fibre, which satisfies $f^2=0$ (two different fibres are distinct) and $f\cdot (f+K_S)=f\cdot K_S=-2$ (by the adjunction formula, Proposition \refd{Prp:AdjunctionFormula}). Conversely, any such divisor induces a morphism $S\stackrel{|f|}{\rightarrow}\mathbb{P}^1$ that provides a conic bundle structure on $S$. Indeed, since no irreducible curve of $S$ has self-intersection $<-1$ (see Proposition \refd{Prp:ConditionsDelPezzo}), the singular fibres are a transverse union of two exceptional curves (see Lemma \refd{Lem:Negativity} in the Appendix).

The conic bundle structures on $S$ therefore correspond to effective divisors $f$ such that $f^2=0$ and $f\cdot K_S=-2$. We write $f= mL-\sum_{i=1}^r a_iE_i$ for some integers $m,a_1,...,a_r$ and recall that $K_S=-3L+\sum_{i=1}^rE_i$. The two conditions are equivalent to:
\begin{equation}
\label{equation:ConicBundlesDelPezzo}
\begin{array}{rcl}\sum_{i=1}^r a_i^2&=&m^2,\vspace{0.1cm}\\
\sum_{i=1}^r a_i & = & 3m-2,\end{array}
\end{equation}
which imply that $\sum_{i=1}^r a_i(a_i-1)=(m-1)(m-2)$. Note that if $r\leq 7$ and if the $a_i$'s satisfy these conditions, we find by adding some new coefficient with value $1$, a solution to Equation 
\refd{equation:ExcCurves} in the proof of Proposition \refd{Prp:DescExcCurvDelPezzo}. We find in this manner all possibilities with $7$ multiplicities given above.

To consider the cases with $8$ multiplicities, we suppose now that $r=8$. Since $a_1,...,a_8\geq 1$, we have $\sum_{i=1}^8 a_i=3m-2\geq 8$, so $m\geq 4$. 
\begin{itemize}
\item
If $m=4$, then $\sum_{i=1}^8 a_i(a_i-1)=6$, so apart from the multiplicities $1$, we have either three $2$'s or one $3$. We get the solutions $(2,2,2,1,1,1,1)$ (already known) and $(3,1,1,1,1,1,1,1)$ (a new one).
\item
If $m=5$, then $\sum_{i=1}^8 a_i(a_i-1)=12$, so apart from the multiplicities $1$, we have either six $2$'s, or one $3$ and three $2$'s, or two $3$'s, or one $4$. We get the solutions $(2,2,2,2,2,2,1)$ (already known) and $(3,2,2,2,1,1,1,1)$ (a new one). The case with two $3$'s (respectively with one $4$) is impossible since the number of multiplicities $1$ would be $7$ (respectively $9$).
\end{itemize}
To find the cases for $m\geq 6$, we now prove that if $f$ is the divisor of the fibres of a conic bundle structure, then $-4K_S-f$ is the divisor of another one. Indeed, take one singular fibre of $f$, which is $F=\{F_1,F_2\}$. The image by the Bertini involution of these two exceptional divisors is $\{-2K_S-F_1,-2K_S-F_2\}$ (see Lemma \refd{Lem:2KSmD}). As $f=F_1+F_2$, the image of $f$ by the Bertini involution is $(-2K_S-F_1)+(-2K_S-F_2)=-4K_S-f$, which is thus the divisor of the fibre of a conic bundle.

\begin{itemize}
\item
If $m=6$, a conic bundle with solution $(6,a_1,...,a_8)$ yields another one with solution $(6,4-a_1,...,4-a_8)$, so $1 \leq a_i \leq 3$ for $i=1,...,8$. As $\sum_{i=1}^8 a_i=16$, the number of $3$'s is the same as the number of $1$'s. Denoting this number by $d$, we have $36=\sum_{i=1}^8 a_i^2=9d +4(8-2d)+d= 32+2d$, so $d=2$. We find the solution $(3,3,2,2,2,2,1,1)$ given above.
\end{itemize}

We get the solutions for $m>6$ from those for $12-m$.\proofend 
\end{proof}

\bigskip

\index{Del Pezzo surfaces!automorphisms}
A group of automorphisms of a Del Pezzo surface $S$ may be minimal and preserve a fibration (case $1$ of Proposition \refd{Prp:TwoCases}). Let us give the possible cases explicitly:
\begin{Prp}\fauxtitre
\label{Prp:InterDPCB}
Let $S$ be a Del Pezzo surface, $\pi:S \rightarrow \mathbb{P}^1$ a conic bundle structure on the surface and $G \subset \Aut(S,\pi)$ a finite group of automorphisms of the conic bundle, such that the pair $(G,S)$ is minimal.

Then the surface $S$ is $\mathbb{P}^1\times\mathbb{P}^1$ or a Del Pezzo surface of degree $1$, $2$ or $4$.
\end{Prp}
\begin{proof}\fauxtitred\upshape
As the pair $(G,S)$ is minimal, the rank of $\Pic{S}^{G}$ is $2$, by Proposition \refd{Prp:TwoCases}, generated by $K_S$ and $f$, the divisor class of a fibre of $\pi$. 

The cone $\overline{NE(S)}^G$ of curves of $S$ invariant by $G$ also has rank $2$, since it contains $-K_S$ and $f$. It contains two extremal rays (see \cite{bib:KoM}), each  intersecting $K_S$ negatively, since $-K_S$ is ample. One of the extremal rays is generated by $f$; as the pair $(G,S)$ is minimal, the other is generated by the class of a rational curve $f'$, such that $f'=0$ and $|f'|$ gives another $G$-invariant conic bundle structure on $S$ (see \cite{bib:KoM}, page $48$).

We write $f'=-a K_S +bf$, with $a,b \in \mathbb{Q},a\not=0$. We know that $f'^2=0$, whence $f'\cdot K_S=-2$ (by the adjunction formula, Proposition \refd{Prp:AdjunctionFormula}). These two relations imply that

\begin{center}
$\begin{array}{rrr}a^2 K_S^2+4ab=a(aK_S^2+4b)&=&0, \vspace{0.1cm}\\
-a K_S^2-2b&=&-2.\end{array}$
\end{center}
As $a\not=0$, the first equation implies that $aK_S^2=-4b$. This and the second equation imply that $2b=-2$, so $b=-1$. Thus, we find that $f'=-aK_S-f$ and $aK_S^2=4$.

If $S\not=\mathbb{P}^1\times\mathbb{P}^1$ and $S\not=\Pn$, the canonical divisor is not a multiple in $\Pic{S}$, so $a\in \mathbb{Z}$ and the degree $K_S^2$ of the Del Pezzo surface is equal to $1$, $2$ or $4$.
\proofend\end{proof}

\bigskip

We give three examples:
\begin{Exa}\fauxtitred\upshape
\begin{enumerate}
\item\index{Del Pezzo surfaces!$\mathbb{P}^1\times\mathbb{P}^1$}
Let $G \cong (\Z{2})^4$ be the group of automorphisms of $\mathbb{P}^1 \times \mathbb{P}^1$ of the form
\begin{center}$(x,y) \dasharrow (\pm x^{\pm 1},\pm y^{\pm 1})$.\end{center}
It is minimal (as is any automorphism group of $\mathbb{P}^1 \times \mathbb{P}^1$) and preserves the two rulings.
\item
Let $G\cong \Z{2}$ be the group generated by a cubic involution of a Del Pezzo surface of degree $4$\index{Del Pezzo surfaces!of degree $4$} (see Section \refd{Sec:DP4}). It acts minimally on the surface and preserves two conic bundle structures (see Proposition \refd{Prp:GeoCub}, assertion $4$) corresponding to the divisors $L-E_i$ and $-K_S-(L-E_i)$, for some $i \in \{1,...,5\}$.
\item
The group isomorphic to $(\Z{2})^2$ of Example \refd{Exa:DP2ConicBundle} acts minimally on a Del Pezzo surface of degree $2$ and preserves a conic bundle structure.\index{Del Pezzo surfaces!of degree $2$}
\end{enumerate}
\end{Exa}
\chapter{Automorphisms of conic bundles}
\label{Chap:ConicBundles}
\index{Conic bundles!automorphisms|idb}
\ChapterBegin{In this chapter we prepare some tools which will be useful in the sequel to study the automorphisms of conic bundles.}
\section{Description of conic bundles}
\label{Sec:DescConicBundles}
\pagetitless{Automorphisms of conic bundles}{Description of conic bundles}
\index{Conic bundles!description|idb}
In this section, we give some known and some new results on conic bundles, without reference to any group action. We begin with the definition of a conic bundle:
\begin{Def} \fauxtitre
Let $S$ be a rational surface and $\pi:S\rightarrow \mathbb{P}^1$ a morphism. We say that the pair $(S,\pi)$ is a \defn{conic bundle} if
\begin{itemize}
\item
A general fibre of $\pi$ isomorphic to $\mathbb{P}^1$.
\item
There is a finite number of exceptions: these singular fibres are the union of
rational curves $F$ and $F'$ such that $F^2={F'}^2=-1$ and  $FF'=1$.
\end{itemize}
\end{Def}
Note that the condition that $S$ be rational is in fact induced by the others (using for example a Noether-Enriques theorem whose proof can be found in \cite{bib:Be1}, Theorem III.4).

\bigskip

In this context, the natural notion of minimality is defined as follows:
\begin{Def} \fauxtitred
\begin{itemize}
\item
Let $(S,\pi)$ and $(S',\pi')$ be two conic bundles. We say that $\varphi:S \dasharrow S'$ is a \defn{birational map of conic bundles} if $\varphi$ is a birational map which sends a general fibre of $\pi$ on a general fibre of $\pi'$.
\item
We say that a conic bundle $(S,\pi)$ is \defn{minimal} if any birational morphism of conic bundles $(S,\pi)\rightarrow (S',\pi')$ is an isomorphism.
\end{itemize}\end{Def}

\break

The following lemma is well known:
\begin{Lem}\fauxtitre
\label{Lem:SmoothHirz}
\index{Conic bundles!minimality|idb}
Let $(S,\pi)$ be a conic bundle. The following conditions are equivalent:
\begin{itemize}
\item
$(S,\pi)$ is minimal.
\item
The fibration $\pi$ is smooth, i.e.\ no fibre of $\pi$ is singular.
\item
$S$ is a Hirzebruch surface.
\end{itemize}
\end{Lem}
\begin{proof}\fauxtitred\upshape
Let $\varphi:(S,\pi) \rightarrow (S',\pi')$ be some birational morphism of conic bundles. The morphism $\varphi$ is not an isomorphism if and only if it contracts at least one exceptional curve $F$. If this is the case, as the conic bundle structure is preserved, $F$ does not intersect a regular fibre and so $\pi(F)$ is a point. This implies that $F$ is contained in some singular fibre. The second assertion thus implies the first one.

Conversely, if there exists some singular fibre, we may contract one of its components, by Castelnuovo's contractibility criterion (see Proposition \refd{Prp:Castelnuovocontractibilitycriterion}). This gives a birational morphism of conic bundles which is not an isomorphism.

Finally, the equivalence between the second and third assertions is proved in \cite{bib:Be1}, Proposition III.15.\proofend
\end{proof}

\bigskip

Here is a useful result:
\begin{Lem}\fauxtitre
\label{Prp:IskWithoutG}
\index{Conic bundles!automorphisms}
Let $(S,\pi)$ be a conic bundle on a surface $S\not\cong\mathbb{P}^1\times\mathbb{P}^1$. Let $-n$ be the minimal self-intersection of sections of $\pi$ and let $r$ be the number of singular fibres of $\pi$. \\
Then $n\geq 1$ and: 
\begin{itemize}
\item
There exists a birational morphism of conic bundles $p_{-}:S\rightarrow \mathbb{F}_n$ such that:
\begin{itemize}
\item[\upshape 1.]
$p_{-}$ is the blow-up of\hspace{0.05 cm} $r$ points of\hspace{0.05 cm} $\mathbb{F}_n$, none of them lying on the exceptional section $E_n$.
\item[\upshape 2.]
The strict pull-back $\tilde{E_n}$ of $E_n$ by $p_{-}$ is a section of $\pi$ with self-intersection $-n$.
\end{itemize}
\item
If there exist two different sections of $\pi$ with self-intersection $-n$, then $r\geq 2n$. In this case, there exist birational morphisms of conic bundles $p_0:S\rightarrow \mathbb{F}_0$ and $p_1:S\rightarrow \mathbb{F}_1$.
\end{itemize}
\end{Lem}
\begin{proof}\fauxtitred\upshape
We denote by $s$ a section of $\pi$ of minimal self-intersection $-n$, for some integer $n$ (this integer is in fact positive, as we prove later). Note that this curve intersects exactly one irreducible component of each singular fibre.

If $r=0$, the lemma is trivially true: take $p_{-}$ to be identity map. We suppose now that $r\geq1$. We denote by $F_1,...,F_r$ the irreducible components of the singular fibres which do not intersect $s$. Blowing these down, we get a birational morphism of conic bundles $p_{-}:S\rightarrow \mathbb{F}_m$, for some integer $m\geq 0$.

The image of the section $s$ by $p_{-}$ is a section of the conic bundle of $\mathbb{F}_m$ of minimal intersection. Indeed, let $s'$ be some section of $\mathbb{F}_m$. Then the strict pull-back by $p_{-}$ of $s'$ has self-intersection $\leq (s')^2$, and is at least equal to $s^2$ by hypothesis.

As the section of minimal self-intersection of a section of $\mathbb{F}_m$ is $-m$, we get $m=n$, so $n \geq 0$. If $m=0$, then taking some section $s''$ of $\mathbb{P}^1\times\mathbb{P}^1$ of self-intersection $0$ passing through at least one blown-up point, its strict pull-back by $p_{-}$ is a section of self-intersection $<0$, which is not possible as $s^2=-n$.

We find finally that $m=n>0$, and that $p_{-}(s)$ is the unique section $\mathbb{F}_n$ of self-intersection $-n$. This proves the first assertion. 

We now prove the second assertion. The Picard group of $S$ is generated by $s=p_{-}^{*}(E_n)$, the divisor $f$ of a fibre of $\pi$ and $F_1,...,F_r$. Write some section $t=s+bf-\sum_{i=1}^r{a_i}F_i$, for some integers $b,a_1,...,a_r$, with $a_1,...,a_r \geq 0$, such that $t^2=-n$. We have $t\cdot (t+K_S)=-2$ (adjunction formula), where $K_S=p_{-}^{*}(K_{\mathbb{F}_n})+\sum_{i=1}^r F_i=-(n+2)f-2s+\sum_{i=1}^r F_i$. These relations give:
\begin{center}
$\begin{array}{ccccl}
s^2&=&t^2&=&s^2-\sum_{i=1}^r a_i^2+2b,\vspace{0.1 cm}\\
n-2&=&t\cdot K_S&=&-(n+2)+2n-2b+\sum_{i=1}^r a_i,
\end{array}$\end{center}
whence $\sum_{i=1}^r a_i=\sum_{i=1}^r a_i^2=2b$, so every $a_i$ is $0$ or $1$ and consequently $2b\leq r$. As $s\cdot t=b-n \geq 0$, we find $r\geq 2n$, as annouced in the lemma.

Finally, by contracting $f-F_1,f-F_2,...,f-F_n,F_{n+1},F_{n+2},...,F_{r}$, we obtain a birational morphism $p_0$ of conic bundles which sends $s$ on a section of self-intersection $0$ and thus whose image is $\mathbb{F}_0=\mathbb{P}^1\times\mathbb{P}^1$. Similarly, the morphism $p_1:S\rightarrow \mathbb{F}_1$ is given by the contraction of $f-F_1,f-F_2,...,f-F_{n-1},F_{n},F_{n+1},...,F_{r}$.
\proofend\end{proof}

\bigskip

\section{Finite groups of automorphisms: two representations}
\label{Sec:CBTwoRepr}
\pagetitless{Automorphisms of conic bundles}{Finite groups of automorphisms: two representations}
In this section, we study finite groups of automorphisms of conic bundles, or equivalently, finite subgroups of $\CrP$ that leave invariant a pencil of rational curves.

We can view this kind of group in two ways:
\begin{enumerate}
\item
One way is to consider these groups directly as \cadre{automorphisms of a conic bundle}:

Let $(S,\pi)$ be a conic bundle.
The group of automorphisms of $S$ that leave the conic bundle structure invariant (i.e.\ that send every fibre to another fibre) is denoted by $\AutSp$.
\item
Another way is to consider these groups as \cadre{birational maps of some rational surface}, which leave invariant a pencil of rational curves:

Let $\Lambda$ be a pencil of rational curves of some surface $S$.
We denote by $\Crc{S}{\Lambda}$ the group of birational maps of $S$ that leaves $\Lambda$ invariant. If $\Lambda$ is induced by a morphism $\pi:S \rightarrow \mathbb{P}^1$ (i.e.\ if $\Lambda$ is base-point-free), we write $\Crc{S}{\pi}=\Crc{S}{\Lambda}$.
\end{enumerate}

\bigskip

\begin{Exa}\fauxtitred\upshape
\label{Exa:CrPp}
Let $\pi_1:\mathbb{P}^1\times \mathbb{P}^1 \rightarrow \mathbb{P}^1$ be the projection on the first factor.
Explicitly, $\CrPp$ is the group of birational maps of the form:
$$(x,y)\dasharrow \left(\frac{ax+b}{cx+d},\frac{\alpha(x)y+\beta(x)}{\gamma(x)y+\delta(x)}\right),$$
where $a,b,c,d \in \K$, $\alpha,\beta,\gamma,\delta \in \K(x)$, and $(ad-bc)(\alpha\delta-\beta\gamma)\not=0$.

This group, called the \defn{de Jonqui\`eres group}, is isomorphic to $\PGL(2,\K(x))\rtimes \PGL(2,\K)$.
\end{Exa}

\bigskip

We now show that these two points of view are equivalent. For this we use Proposition \refd{Prp:GS} and a Noether-Enriques theorem (a recent proof of which can be found in \cite{bib:Be1}, Theorem III.4):

\begin{Prp}\fauxtitred
\label{Prp:BeauTsen}
\begin{itemize}
\item
Let $\Lambda$ be a pencil of rational curves of some surface $S$. There exists a birational map $\varphi: S \dasharrow \mathbb{P}^1\times\mathbb{P}^1$ such that $\varphi(\Lambda)$ is the pencil of fibres of $\pi_1:\mathbb{P}^1\times\mathbb{P}^1\rightarrow \mathbb{P}^1$
and $\varphi \Crc{S}{\Lambda} \varphi^{-1}=\CrPp$.
\item
Let $G$ be a finite subgroup of $\CrPp$. There exist a conic bundle $(S,\pi)$ and a birational map of conic bundles $\varphi: S \dasharrow \mathbb{P}^1\times\mathbb{P}^1$ such that $\varphi^{-1} G\varphi \subset \Aut(S,\pi)$.
\item
Let $G \subset \Aut(S,\pi)$ be a finite subgroup of automorphisms of some conic bundle $(S,\pi)$. Then $G$ is conjugate to a subgroup of $\CrPp$ by a birational map of conic bundles.
\end{itemize}
\end{Prp}
\begin{proof}\fauxtitred\upshape
\begin{itemize}
\item
To prove the first assertion, we blow-up the base points of the pencil $\Lambda$ until it becomes base-point-free. We then get a morphism $\pi:\tilde{S}\rightarrow \mathbb{P}^1$, with rational generic fibre, and a birational morphism $\tilde{S}\rightarrow S$ which sends the pencil of fibres of $\pi$ on $\Lambda$ and conjugates $\Crc{S}{\Lambda}$ to $\Crc{\tilde{S}}{\pi}$. By \cite{bib:Be1}, Theorem III.4, there exists a birational map $\varphi:\tilde{S}\dasharrow \mathbb{P}^1\times\mathbb{P}^1$ such that the following diagram commutes:

\begin{center}
\hspace{0 cm}\xymatrix{ \tilde{S}\ar@{-->}^{\varphi}[rr] \ar_{\pi}[rd] & & \mathbb{P}^1\times\mathbb{P}^1 \ar_{\pi_1}[ld] \\
& \mathbb{P}^1. & }
\end{center}
This proves the first assertion.
\item
To prove the second assertion, we observe that by Proposition \refd{Prp:GS}, there exist a rational surface $S$ and a birational map $\varphi_1: S \dasharrow \mathbb{P}^1\times\mathbb{P}^1$ such that $G_1:=\varphi_1^{-1} G\varphi_1 \subset \Aut(S)$. The fibres of $\pi_1$ are sent on a $G_1$-invariant pencil $\Lambda$ of rational curves of $S$. If this pencil has base points, we blow them up (they are $G_1$-invariant). We get
a birational map $\varphi_2:S \dasharrow S'$ which conjugates $G_1$ to some group $G_2 \subset \Aut(S')$ and sends $\Lambda$ on a base-point-free pencil, to which corresponds a morphism $\pi':S'\rightarrow \mathbb{P}^1$. We may suppose that the $G_2$-orbit of any exceptional curve contained in a fibre is not contractible.

This implies that every reducible fibre is a transverse union of two curves of self-intersection $-1$ (see Lemma \refd{Lem:Negativity} in the Appendix). Hence, $(S',\pi')$ is a conic bundle.

\item
The last assertion follows directly from the first, since $G\subset \Aut(S,\pi) \subset \Crc{S}{\pi}$, which is conjugate by a birational map of conic bundles to $\CrPp$.\end{itemize}\proofend
\end{proof}

\bigskip

Note that the first assertion implies that if $(S,\pi)$ is a conic bundle, the group $\Crc{S}{\pi}$ is birationally conjugate to the de Jonqui\`eres group.

We may thus study either finite subgroups of $\AutSp$, for some conic bundle $(S,\pi)$, or finite subgroups of $\CrPp$. The advantage of the first possibility is that we have \emph{biregular morphisms} and can use their action on the set of exceptional curves of the fibration (i.e.\ exceptional divisors contained in a fibre). The advantage of the second possibility is that we can use the structure of $\CrPp \cong \PGL(2,\K(x))\rtimes \PGL(2,\K)$ and the nature of the surface $\mathbb{P}^1\times \mathbb{P}^1$.
We shall use both representations to classify finite abelian subgroups of $\CrPp$.

\section{The exact sequence}
\pagetitless{Automorphisms of conic bundles}{The exact sequence}
\label{Sec:Exactsequence}
Let $G \subset \Crc{S}{\pi}$ be some subgroup acting (biregularly) on a conic bundle $(S,\pi)$. We have a natural homomorphism
$\overline{\pi}: G
\rightarrow \Aut(\mathbb{P}^1)= \PGLn{2}$ that satisfies $\overline{\pi}(g)\pi=\pi g$, for every $g \in G$.
We observe that the group $G' =\ker \overline{\pi}$ of automorphims 
that leaves every fibre invariant embeds in the subgroup
$\PGL(2,\K(x))$ of automorphisms of the generic fibre $\mathbb{P}^1(\K(x))$.

We use the exact sequence
\begin{equation}
\label{eq:ExactSeqCB}
1 \rightarrow G' \rightarrow G \stackrel{\overline{\pi}}{\rightarrow} \overline{\pi}(G) \rightarrow 1\end{equation} to restrict the structure of $G$. 

When $G$ is abelian and finite, so are $G'$ and $\overline{\pi}(G)$. The finite abelian subgroups of $\PGLn{2}$ and $\PGL(2,\K(x))$ are 
either cyclic or isomorphic to $(\Z{2})^2$ (see Lemma 
\refd{Lem:AutP1} in the Appendix). We examine all possibilities in Chapter \refd{Chap:FiniteAbConicBundle}.

\section{Minimality and contraction of curves - the twisting of singular fibres}
\pagetitless{Automorphisms of conic bundles}{Minimality and contraction of curves - the twisting of singular fibres}
\label{Sec:Minimality}
Let us also consider the action of some group on the conic bundle and define the notion of minimality in this context:

\begin{Def} \fauxtitred
\begin{itemize}
\item
Let $G \subset \Aut(S,\pi)$ be a group of automorphisms of a conic bundle. We say that a birational morphism of conic bundles $\varphi:S \rightarrow S'$ is \defn{$G$-equivariant} if the $G$-action on $S'$ induced by $\varphi$ is biregular (it is clear that it preserves the conic bundle structure).
\item
Let $(S,\pi)$ be a conic bundle and $G \subset \Aut(S,\pi)$ a group of automorphisms.
We say that the triple $(G,S,\pi)$ is \defn{minimal} if any $G$-equivariant birational morphism of conic bundles $\varphi:S\rightarrow S'$ is an isomorphism.\end{itemize}\end{Def}

\begin{remark}\fauxtitred
If $G \subset \Aut(S,\pi)$ is such that the pair $(G,S)$ is minimal, so is the triple $(G,S,\pi)$. The converse is not true in general (see Section \refd{SubSec:DelPezzo6Kappa}).
\end{remark}

As we mentioned above, any automorphism of the conic bundle acts on the set of singular fibres and on its irreducible components. Here is a basic simple result:
\begin{Lem}\fauxtitre
\label{Lem:MinTripl}
 Let $G \subset \Aut(S,\pi)$ be a finite group of automorphisms of a conic bundle.
The following two conditions are equivalent:
\index{Conic bundles!minimality}
\begin{itemize}
\item[\upshape 1.]
The triple $(G,S,\pi)$ is minimal.
\item[\upshape 2.]
For any singular fibre $\{F_1,F_2\}$ of $\pi$ (where $F_1$ and $F_2$ are two exceptional curves) there exists $g \in G$ such that $g(F_1)=F_2$ (and consequently $g(F_2)=F_1$).
\end{itemize}
\end{Lem}
\begin{proof}\fauxtitred\upshape
Suppose that the second assertion is true. Let $(S',\pi')$ be another conic bundle and $\varphi:S\rightarrow S'$ be a $G$-equivariant birational morphism of conic bundles. Suppose that $\varphi$ contracts some exceptional curve $F$. This curve must be an irreducible component of some singular fibre, since $\varphi$ leaves invariant the conic bundle structure. And its orbit by $G$ must be contracted, as $\varphi$ is $G$-equivariant. The second condition then implies that the two components of the singular fibre containing $F$ must be contracted, which is impossible. So $\varphi$ is an isomorphism and the triple $(G,S,\pi)$ is minimal.

Conversely, suppose that the second assertion is false. Let $\{F_1,F_2\}$ be a singular fibre of $\pi$ such that no element of $G$ sends $F_1$ on $F_2$. Then, the orbit of $F_1$ by $G$ consists of the disjoint union of exceptional curves lying on distinct singular fibres. So we may contract the orbit by $G$ of $F_1$ to get a $G$-equivariant birational morphism of conic bundles which is not an isomorphism. The triple $(G,S,\pi)$ is then not minimal.\proofend\end{proof}

\bigskip

The permutation of two components of a singular fibre is thus very important. For this reason, we introduce some terminology:
\begin{Def} \fauxtitre
Let $g\in \Aut(S,\pi)$ be an automorphism of the conic bundle $(S,\pi)$. Let $F=\{F_1,F_2\}$ be a singular fibre. We say that $g$ \defn{twists} the singular fibre $F$ if $g(F_1)=F_2$ {\upshape (}and consequently $g(F_2)=F_1${\upshape )}.
\end{Def}
We need only consider groups of automorphisms of conic bundles in which every singular fibre is twisted by some element of the group.
\begin{remark}\fauxtitred
An automorphism of a conic bundle with a non-trivial action on the base of the fibration may twist at most two singular fibres. However, an automorphism with a trivial action on the fibration may twist a large number of fibres. (We will see in Lemma \refd{Lem:DeJI} that this number is even in both cases.)
\end{remark}

\bigskip

The following result is a consequence of Lemma \refd{Prp:IskWithoutG}.
\begin{Lem}\fauxtitre
\label{Lem:GoingToF0F1}
Let $(S,\pi)$ be a conic bundle and $g \in \Aut(S,\pi)$ be an automorphism of finite order which twists at least one singular fibre.
\begin{itemize}
\item[\upshape 1.]
There exist two birational morphisms of conic bundles $p_0: S \rightarrow \mathbb{F}_0$ and $p_1:S \rightarrow \mathbb{F}_1$  (which are not necessarily $g$-equivariant).
\item[\upshape 2.]
Let $-n$ be the minimal self-intersection of sections of $\pi$ and let $r$ be the number of singular fibres of $\pi$. Then, $r\geq 2n\geq 2$
\end{itemize}
\end{Lem}
\begin{proof}\fauxtitred\upshape
Note that any section touches exactly one component of each singular fibre. As $g$ twists some singular fibre, its action on the set of sections of $S$ fixes no element. The number of sections of minimal self-intersection is then greater than $1$ and we apply Lemma \refd{Prp:IskWithoutG} to get the result.\proofend
\end{proof}
\begin{remark}\fauxtitred
A result of the same kind can be found in \cite{bib:Isk1}, Theorem 1.1.\end{remark}
\chapter{Finite abelian groups of automorphisms of surfaces of large degree}
\label{Chap:LargeDegreeSurfaces}
\ChapterBegin{In this chapter, we describe the finite abelian groups of automorphisms of surfaces which have a large degree. We find that these are birationally conjugate to subgroups of $\Aut(\Pn)$ or $\Aut(\mathbb{P}^1\times\mathbb{P}^1)$ (Proposition \refd{Prp:LargeDegree}, proved at the end of the chapter).}

By definition, the degree of a surface is the square of its canonical divisor. We begin with the following obvious result: blowing-down a curve decreases the degree:
\begin{Lem}\fauxtitre
\label{Lem:Degree}
Let $\eta:\tilde{S}\rightarrow S$ be the blow-down of some exceptional curve of $S$. Then, $K_{\tilde{S}}^2=K_S^2-1$.
\end{Lem}
\begin{proof}\fauxtitred\upshape
This follows from the fact that $K_{\tilde{S}}=\eta^{*}(K_S)+E$, where $E$ denotes the exceptional curve contracted. We have $K_{\tilde{S}}^2=\eta^{*}(K_S)^2+E^2+2E\cdot \eta^{*}(K_S)$. As $E^2=-1$ and $E\cdot \eta^{*}(K_S)=0$, we get the stated result.\proofend\end{proof}
\begin{remark}
Note that $K_{\Pn}^2=9$ and $K_{\mathbb{F}_n}^2=8$ for any integer $n$. The maximal degree of rational surfaces is thus $9$, since the minimal surfaces are the Hirzeburch surfaces and\hspace{0.2 cm}$\Pn$.\end{remark}

\bigskip

In this chapter, we first describe the automorphism groups of minimal surfaces and then those of surfaces of degree $\geq 5$. We prove at the end of this chapter (Proposition \refd{Prp:LargeDegree}) that any finite abelian group of automorphisms of a rational surface of degree $\geq 5$ is birationaly conjugate to a subgroup of $\Aut(\mathbb{P}^1\times\mathbb{P}^1)$ or $\Aut(\Pn)$.

\bigskip

Note that if $(S,\pi)$ is a conic bundle and $r$ denotes the number of singular fibres of $\pi$, then the degree of $S$ is $8-r$. (This follows directly from Lemma \refd{Lem:Degree} and the remark that follows it).

\break

\section{Automorphisms of\hspace{0.2 cm}$\Pn$}
\pagetitless{Finite abelian groups of automorphisms of surfaces of large degree}{Automorphisms of\hspace{0.2 cm}$\Pn$}\index{Del Pezzo surfaces!$\mathbb{P}^2$|idb}
\label{Sec:AutP2}
The simplest examples of birational maps of the plane are linear automorphisms, which are maps of the form $(x:y:z) \mapsto (a_{11}x+a_{12}y+a_{13}z:a_{21}x+a_{22}y+a_{23}z:a_{31}x+a_{32}y+a_{33}z)$, for some $(a_{ij})_{i,j=1}^3 \in \PGLn{3}$. 
\begin{Def} \fauxtitre
A \defn{diagonal automorphism of\hspace{0.2 cm}$\Pn$} is an automorphism $(x:y:z) \mapsto (\alpha x: \beta y: \gamma z)$ with $\alpha,\beta,\gamma \in \K^{*}$.
We will denote such an automorphism by $\Diag{\alpha}{\beta}{\gamma}$, and denote by $\mathcal{T}$ the diagonal torus of $\PGLn{3}$, generated by all these automorphisms.
\end{Def}

The group of automorphisms of\hspace{0.2 cm}$\Pn$ is well known to be $\PGLn{3}$ and the conjugacy classes of its finite abelian subgroups are also classical (we use the notation $\zeta_n=e^{2\im \pi/n}$ and $\omega=\zeta_3$):
\begin{Prp}\fauxtitre
\label{Prp:PGL3Aut}
Every finite abelian subgroup of\hspace{0.2cm}$\PGLn{3}$ is conjugate, in this group, to one of the following:
\begin{itemize}
\item[\upshape 1.]
A cyclic group isomorphic to $\Z{n}$, generated by $\Diag{1}{\zeta_n^a}{\zeta_n^b}$, with $\gcd(n,a,b)=1$.
\item[\upshape 2.]
A diagonal group, isomorphic to $\Z{n} \times \Z{m}$, where $n$ divides $m$, generated by \\
$\Diag{1}{1}{\zeta_n}$ and  $\Diag{1}{\zeta_m^a}{\zeta_m^b}$, with $\gcd(m,a,b)=1$.
\item[\upshape 3.]
The special group $V_9$, isomorphic to $\Z{3} \times \Z{3}$, generated by
 $\Diag{1}{\omega}{\omega^2}$ \\ and $(x:y:z) \mapsto (y:z:x)$.
\end{itemize}
\end{Prp}\begin{proof}\fauxtitred\upshape
It is clear that every element of $\PGLn{3}$ of finite order is diagonalizable, so the cyclic groups are conjugate to the first case.

On computing the commutator of an element we see that, except for the conjugacy class of $V_9$, every finite abelian group is diagonalizable, and thus is contained in the diagonal torus $\mathcal{T}$, of rank $2$. This provides the second case. \proofend
\end{proof}

\bigskip

\begin{remarks}\fauxtitred
\begin{itemize}
\item
Note that there are distinct conjugacy classes in $\PGLn{3}$ of cyclic groups of the same order. For example the groups respectively generated by $\Diag{\zeta_n}{1}{1}$ and $\Diag{\zeta_n}{\zeta_n^2}{1}$ are not conjugate, if $n \geq 3$. This can be proved by counting the number of fixed points of the two groups.
\item
The same is true of non-cyclic groups, for example let $G_1$ be the group generated by $\Diag{-\omega}{1}{1}$ and $\Diag{1}{-1}{1}$ and $G_2$ the group generated by $\Diag{-\omega}{\omega^2}{1}$ and $\Diag{1}{-1}{1}$. Both are isomorphic to $\Z{2}\times\Z{6}$, but they are not conjugate, since some elements of order $6$ of $G_1$ fix infinitely many points, but each element of order $6$ of $G_2$ fixes only a finite number of points.
\end{itemize}
\end{remarks}
\bigskip

The conjugation in the Cremona group is different. Some groups are conjugate in the Cremona group, but not in $\PGLn{3}$. We make this explicit in the following proposition:

\begin{Prp}\fauxtitre
\label{Prp:PGL3Cr}
Every finite abelian subgroup of\hspace{0.2cm}$\PGLn{3}$ is conjugate, in the Cremona group $\CrP$, to one and only one of the following:
\begin{itemize}
\item[${\num{0.n}}$]
A cyclic group isomorphic to $\Z{n}$, generated by $\Diag{\zeta_n}{1}{1}$.
\item[${\num{0.mn}}$]
A diagonal group, isomorphic to $\Z{n} \times \Z{m}$, where $n$ divides $m$, generated by\\ $\Diag{1}{\zeta_n}{1}$ and $\Diag{\zeta_m}{1}{1}$.
\item[$\num{0.V9}$]
 The special group $V_9$, isomorphic to $\Z{3} \times \Z{3}$, generated by $\Diag{1}{\omega}{\omega^2}$\\  and $(x:y:z) \mapsto (y:z:x)$.
\end{itemize}
Thus, except for the third case, two isomorphic finite abelian sugroups of $\PGLn{3}$ are conjugate in $\CrP$.
\end{Prp}\begin{proof}\fauxtitred\upshape
We recall that the automorphisms of\hspace{0.2 cm}$\Pn$ of a given finite order are all conjugate in the Cremona group (see \cite{bib:BeB}). From this, we deduce that all cyclic groups of $\PGLn{3}$ of order $n$ are conjugate in $\CrP$ to ${\nump{0.n}}$.

Let us prove that the group $V_9$ is not diagonalizable, even in $\CrP$. Note that this group has no fixed points, so it cannot be conjugate by a birational map to a group of automorphisms having some fixed points (see \cite{bib:KoS}, Proposition A.2). As diagonal subgroups fix the points $(1:0:0)$, $(0:1:0)$ and $(0:0:1)$, $V_9$ cannot be diagonalizable in $\CrP$.

Let $G$ be a finite abelian subgroup of $\PGLn{3}$, neither cyclic nor conjugate to $V_9$. We may assume (see Proposition \refd{Prp:PGL3Aut}) that it is diagonal and is isomorphic to $\Z{n} \times \Z{m}$, where $n$ divides $m$. Note that the group $\GLZ{2} \subset \CrP$ of birational maps of the form $(x,y) \dasharrow (x^ay^b,x^cy^d)$ normalizes the torus $\mathcal{T}$ of diagonal automorphisms.\ We conjugate the group $G$ by an element of $\GLZ{2}$ (see \cite{bib:BeB}) so that it contains $\Diag{\zeta_m}{1}{1}$. As the group remains diagonal and contains a subgroup isomorphic to $\Z{n} \times \Z{n}$, it contains all elements of $\mathcal{T}$ of order $n$. Hence it contains $\Diag{1}{\zeta_n}{1}$, and we get $\nump{0.mn}$.\proofend
\end{proof}

\bigskip

We will need the following proposition to decide whether or not a group of birational maps of\hspace{0.2 cm}$\Pn$ is conjugate to a linear automorphism.
We begin by giving some classical definitions:
\begin{Def} \fauxtitred
\begin{itemize}\item
Using the notation of \cite{bib:SeR}, we denote by $C^n(P_1^{k_1},P_2^{k_2},...,P_r^{k_r})$ the linear system of curves of\hspace{0.2 cm}$\Pn$ of degree $n$ passing through the point $P_i$ with multiplicity $\geq k_i$ for $i=1,...,r$, where $P_1,...,P_r$ are points of\hspace{0.2 cm}$\Pn$ or infinitely near and $k_1,...,k_r$ are positive integers.
\item
Let $\Lambda$ be a linear system of curves of\hspace{0.2 cm}$\Pn$, without fixed components and $\varphi \in \CrP$ be a birational map. Note that $\varphi$ is an isomorphism between two open dense subsets of\hspace{0.2 cm}$\Pn$, say $U$ and $V$. We define \defn{the homoloidal transform of $\Lambda$ by $\varphi$} to be the linear system generated by the adherence in $\Pn$ of general elements of $\Lambda \cap U$ by $\varphi$. We denote it by $\varphi(\Lambda)$.
\item
We say that a birational map $\varphi \in \CrP$ {\upshape (}respectively a subgroup $G \subset \CrP${\upshape )} \emph{leaves a linear system $\Lambda$ invariant} if $\varphi(\Lambda)=\Lambda$ (respectively if $g(\Lambda)=\Lambda$ for every $g\in G$).
\item
A \defn{homoloidal net} is the homoloidal transform of the linear system of lines of\hspace{0.2 cm}$\Pn$ by some birational map of\hspace{0.2 cm}$\Pn$.
\end{itemize}
\end{Def}

\bigskip

\begin{Prp} \PropRef{\cite{bib:SeR}, VII, 7.}
\label{Prp:HomoNet}
The linear system $C^n(P_1^{k_1},P_2^{k_2},...,P_r^{k_r})$ is a homoloidal net if and only if 
\begin{itemize}
\item
$\sum_{i=1}^r k_i=3n-3$.
\item
$\sum_{i=1}^r (k_i)^2=n^2-1$.
\item
The system is generically irreducible.
\end{itemize}
\end{Prp}

\bigskip

These reminders allow us to give the following basic simple result:
\begin{Prp}\fauxtitre
\label{Prp:HomolNetInvariant}
A subgroup $G$ of $\CrP$ is conjugate in $\CrP$ to a subgroup of $\Aut(\Pn)$ if and only if $G$ leaves a homoloidal net invariant.
\end{Prp}\begin{proof}\fauxtitred\upshape
We denote by $H$ the linear system of lines of\hspace{0.2 cm}$\Pn$.

Suppose that $\varphi\in \CrP$ is such that $\varphi G \varphi^{-1} \subset \Aut(\Pn)$.\ Then, for every $g\in G$, we have $\varphi g \varphi^{-1}(H)=H$, whence $g (\varphi^{-1}(H))=\varphi^{-1}(H)$. The group $G$ thus leaves the homoloidal net $\varphi^{-1}(H)$ invariant.

Conversely, let us suppose that $G$ leaves invariant the homoloidal net $\Lambda$, associated to the birational map $\lambda$ (so that $\lambda^{-1}(H)=\Lambda$). Then, for every $g\in G$, we have $\lambda g\lambda^{-1}(H)=H$, which implies that $\lambda g \lambda^{-1} \in \Aut(\Pn)$. The group $\lambda G \lambda^{-1}$ is a subgroup of $\Aut(\Pn)$.\proofend
\end{proof}

\bigskip

\section{Automorphisms of Hirzeburch surfaces}
\label{Sec:SmoothFib}
\pagetitless{Finite abelian groups of automorphisms of surfaces of large degree}{Automorphisms of Hirzeburch surfaces}
Let us first see that every finite abelian group of automorphisms of a Hirzebruch surface is birationally conjugate to a subgroup of $\Aut(\mathbb{P}^1\times\mathbb{P}^1)$:
\begin{Prp}\fauxtitre
\label{Prp:FnAb}
Let $G \subset \Aut(\mathbb{F}_n)$ be a finite abelian subgroup of automorphisms of $\mathbb{F}_n$, for some integer $n>0$.
 Then, $G$ is conjugate in the de Jonqui\`eres group to a finite group of automorphisms of $\mathbb{F}_0=\mathbb{P}^1\times\mathbb{P}^1$.
\end{Prp}
\begin{proof}\fauxtitred\upshape
Let $G \subset \Aut(\Fn)$ be a finite abelian group, with $n\geq 1$.  Note that $G$ preserves the unique ruling of $\mathbb{F}_n$, which is a smooth conic bundle structure on the surface. We have the exact sequence (see Section \refd{Sec:Exactsequence})
\begin{equation}
1 \rightarrow G' \rightarrow G \stackrel{\overline{\pi}}{\rightarrow} \overline{\pi}(G) \rightarrow 1.
\tag{\refd{eq:ExactSeqCB}}\end{equation}
We denote by $E\subset \mathbb{F}_n$ the unique section
of self-intersection $-n$, which is neccessarily invariant by $G$. In this case, the group $G'=\ker\overline{\pi}$ must fix a point in each fibre (the point lying on $E$) and acts cyclically on each of them. This implies that $G'$ is cyclic and if it is not trivial, it fixes two points in each fibre, one on $E$, the other on another section.

\break

As the group $\overline{\pi}(G) \subset \PGLn{2}$ is abelian, it is isomorphic to a cyclic group or to $(\Z{2})^2$.
\begin{enumerate}
\item
\textit{If $\overline{\pi}(G)$ is a cyclic group,} two fibres at least are invariant by $G$. The group $G$ fixes two points in one such fibre. We can blow-up the point that does not lie on $E$
and blow-down the corresponding fibre to get a group of automorphisms of $\mathbb{F}_{n-1}$. 
We do this $n$ times and finally obtain that $G$ is conjugate in the de Jonqui\`eres group to a group of automorphisms of $\mathbb{F}_0=\mathbb{P}^1 \times
\mathbb{P}^1$.
\item
\textit{If $\overline{\pi}(G)$ is isomorphic to $(\Z{2})^2$,} there exist two fibres $F,F'$ of $\pi$ whose union  is invariant by $G$. 

Let $G_F \subset G$ be the subgroup of $G$ of elements that leave $F$ invariant. This group is of index $2$ in $G$ and thus is normal. As $G_F$ fixes the point $F\cap E$ in $F$, it acts cyclically on it. There exists another point $P \in F$, $P\notin E$, which is fixed by $G_F$. The orbit of $P$ by $G$ consists of two points, $P$ and $P'$, such that $P' \in F'$, $P' \notin E$.
We blow-up these two points and blow-down the strict transforms of $F$ and $F'$ to get 
a group of automorphisms of $\mathbb{F}_{n-2}$. We do this $\lfloor n/2\rfloor$ times to obtain $G$ as a group of automorphisms of $\mathbb{F}_0$ or
$\mathbb{F}_1$.

If $n$ is even, we get in this manner a group of automorphisms of $\mathbb{F}_0=\mathbb{P}^1 \times
\mathbb{P}^1$. 

Note that $n$ cannot be odd, if the group $\overline{\pi}(G)$ is not cyclic.\ Otherwise, we could conjugate  $G$ to a group of automorphisms of $\mathbb{F}_1$ and then to a group of automorphisms of\hspace{0.2 cm}$\Pn$ by blowing-down the exceptional section. We get an abelian subgroup of $\PGLn{3}$ that fixes a point, and thus a group with at least three fixed points, whose action
on the pencil of lines passing through one of the three fixed points is cyclic (see Proposition \refd{Prp:PGL3Cr}). \proofend
\end{enumerate}
\end{proof}

\bigskip

\subsection{Automorphisms of $\mathbb{P}^1\times\mathbb{P}^1$}
\label{SubSec:AutP1}\index{Del Pezzo surfaces!$\mathbb{P}^1\times\mathbb{P}^1$|idb}
We now investigate the automorphisms of $\mathbb{P}^1 \times \mathbb{P}^1$.
It is well known that $\Aut(\mathbb{P}^1 \times \mathbb{P}^1) = (\PGLn{2} \times \PGLn{2}) \rtimes \Z{2}$, where $\Z{2}$ acts on $\PGLn{2} \times \PGLn{2}$ by exchanging the two components. Taking some affine coordinates $x$ and $y$, the two copies of $\PGLn{2}$ act respectively on $x$ and $y$ and the group $\Z{2}$ is the permutation of $x$ and $y$. 

\bigskip

\begin{Exa}\fauxtitred\upshape
\label{Exa:TorP1}
Note that $\Aut(\mathbb{P}^1 \times \mathbb{P}^1)$ contains the group $(\K^{*})^2 \rtimes \Z{2}$, where $(\K^{*})^2$ is generated by automorphisms of the form $(x,y)\mapsto (\alpha x,\beta y)$, $\alpha,\beta \in \K^{*}$, and $\Z{2}$ is generated by the automorphism $(x,y) \mapsto (y,x)$.

The birational map $(x,y) \dasharrow (x:y:1)$ from $\mathbb{P}^1\times \mathbb{P}^1$ to $\Pn$ conjugates $(\K^{*})^2 \rtimes \Z{2}$ to the group of linear automorphisms of\hspace{0.2 cm}$\Pn$ generated by $(x:y:z)\mapsto (\alpha x:\beta y:z)$, $\alpha,\beta \in \K^{*}$ and $(x:y:z) \mapsto (y:x:z)$.
\end{Exa}

\bigskip

We can describe the finite abelian subgroups of $\Aut(\mathbb{P}^1 \times \mathbb{P}^1)$:

\break

\begin{Prp}\fauxtitre
\label{Prp:AutP1P1}
Let $G \subset \Aut(\mathbb{P}^1 \times \mathbb{P}^1)$ be a finite abelian group. Up to conjugation in $\Aut(\mathbb{P}^1 \times \mathbb{P}^1)$, $G$ is one of the following:\\
\begin{tabular}{llll}
$[1]$ & \multicolumn{3}{l}{a subgroup of the group $(\K^{*})^2 \rtimes \Z{2}$ of Example $\refd{Exa:TorP1}$}\\ 
$[2]$ & $G \cong (\Z{2})^2\times\Z{n} $& generated by &$(x,y) \mapsto (\pm x^{\pm 1},y)$ and
$(x,y) \mapsto (x,\zeta_n y)$\\
$[3]$ & $G \cong \Z{2}\times\Z{2n} $& generated by &$(x,y) \mapsto (x^{-1},y)$ and
$(x,y) \mapsto (-x,\zeta_{2n} y)$\\
$[4]$ & $G \cong (\Z{2})^2 $& given by &$(x,y) \mapsto (\pm x^{\pm 1},y)$\\
$[5]$ & $G \cong (\Z{2})^2 $& generated by &$(x,y) \mapsto (-x,-y)$ and
$(x,y) \mapsto (x^{-1},y)$\\
$[6]$ & $G \cong (\Z{2})^2 $& generated by &$(x,y) \mapsto (-x,-y)$ and
$(x,y) \mapsto (x^{-1},y^{-1})$\\
$[7]$ & $G \cong (\Z{2})^3 $& generated by &$(x,y) \mapsto (\pm x,\pm y)$ and
$(x,y) \mapsto (x^{-1},y)$\\
$[8]$ & $G \cong (\Z{2})^4$& given by &$(x,y) \mapsto (\pm x^{\pm 1},\pm y^{\pm 1})$\\
$[9]$ & $G \cong \Z{2}\times \Z{4}$& generated by & $(x,y) \mapsto (x^{-1},y^{-1})$, and $(x,y) \mapsto (-y,x)$\\
$[10]$ & $G \cong (\Z{2})^3$& generated by & $(x,y) \mapsto (-x,-y)$, $(x,y) \mapsto (x^{-1},y^{-1})$, 
\\
& & & and $(x,y) \mapsto (y,x)$
\end{tabular}\\
{\upshape (}in $[2]$ and $[3]$, $n$ is an integer $>1$ and $\zeta_n=e^{2\im \pi/n}${\upshape )}.
\end{Prp}
\begin{remark}\fauxtitred
It can be proved that cases $[2],...,[10]$ represent distinct conjugacy classes in $\Aut(\mathbb{P}^1\times\mathbb{P}^1)$. We leave this as an exercise, since we do not need it here.\end{remark}
\begin{proof}\fauxtitred\upshape
We first consider a group $G \subset \PGLn{2} \times \PGLn{2}$. We denote by $\pi_1$ and $\pi_2$ the projections $\pi_i: \PGLn{2} \times \PGLn{2} \rightarrow \PGLn{2}$ on the $i$-th factor. Note that $\pi_1(G)$ and $\pi_2(G)$ are finite abelian subgroups of $\PGLn{2}$ and are therefore conjugate either to a diagonal cyclic group or to the group $x\dasharrow \pm x^{\pm 1}$, isomorphic to $(\Z{2})^2$. Let us enumerate the possible cases.
\begin{itemize}
\item
{\bf Both groups $\pi_1(G)$ and $\pi_2(G)$ are cyclic.}\\
In this case, the group $G$ is conjugate to a subgroup of the diagonal torus $(\K^{*})^2$ of automorphisms of the form $(x,y)\mapsto (\alpha x,\beta y)$, $\alpha,\beta \in \K^{*}$. This is part of the first possibility given in the proposition. Note that this diagonal torus is conjugate in $\CrP$ to that of $\Aut(\Pn)$.
\item
{\bf Exactly one of the two groups $\pi_1(G)$ and $\pi_2(G)$ is cyclic.}\\
Up to conjugation in $\Aut(\mathbb{P}^1 \times \mathbb{P}^1)$, we may suppose that $\pi_2(G)$ is cyclic, generated by $y \mapsto \zeta_m y$ and is isomorphic to $\Z{m}$, for some integer $m\geq 1$, and that $\pi_1(G)$ is the group $x\dasharrow \pm x^{\pm 1}$, isomorphic to $(\Z{2})^2$.
We use the exact sequence
\begin{center}$1 \rightarrow G \cap \ker \pi_2 \rightarrow G \stackrel{\pi_2}{\rightarrow} \pi_2(G) \rightarrow 1$.\end{center}
Note that $G$ is a subgroup of $\pi_1(G) \times \pi_2(G) \cong (\Z{2})^2\times\Z{m}$.
As $\pi_2(G)$ is cyclic, and $G$ is not, the group $G \cap \ker \pi_2$ is not trivial. Hence it is either isomorphic to $(\Z{2})^2$ or to $\Z{2}$. In the first case, the group $G$ is the entire group $\pi_1(G) \times \pi_2(G)$, given in $[2]$ if $n>1$ and in $[4]$ if $n=1$. In the second case, using the exact sequence, we see that $G$ is isomorphic to $\Z{2} \times \Z{m}$ or to $\Z{2m}$. As $G$ is not cyclic, $m$ must be even and $G$ isomorphic to $\Z{2} \times \Z{m}$. Up to conjugation, we may assume that $G \cap \ker \pi_2$ is generated by $x \mapsto x^{-1}$, and we get case $[3]$ if $m>2$ and $[5]$ if $m=2$.
\item
{\bf Both groups $\pi_1(G)$ and $\pi_2(G)$ are isomorphic to $(\Z{2})^2$.}\\
Up to conjugation, the groups $\pi_1(G)$ and $\pi_2(G)$ are respectively $x\dasharrow \pm x^{\pm 1}$ and $y\dasharrow \pm y^{\pm 1}$. 
We use once again that $G$ is contained in $\pi_1(G) \times \pi_2(G) \cong (\Z{2})^4$ and the exact sequence
\begin{center}$1 \rightarrow G \cap \ker \pi_2 \rightarrow G \stackrel{\pi_2}{\rightarrow} \pi_2(G) \rightarrow 1$.\end{center}
The possibilities for $G \cap \ker \pi_2$ are to be isomorphic to $\{1\}$, to $\Z{2}$ or to $(\Z{2})^2$. These three possibilities give groups conjugate respectively to $[6]$, $[7]$ and $[8]$.
\end{itemize}

We now have to investigate the groups not contained in $\PGLn{2} \times \PGLn{2}$. Any element $\varphi\in \Aut(\mathbb{P}^1 \times \mathbb{P}^1)$ not contained in $\PGLn{2} \times \PGLn{2}$ is written $\varphi:(x,y)\mapsto (\alpha(y),\beta(x))$, where $\alpha,\beta \in \Aut(\mathbb{P}^1)$. Conjugating this element by $(x,y)\mapsto (\beta(x),y)$, we may suppose that $\beta(x)=x$. Supposing that $\varphi$ is of finite order, as $\varphi^2(x,y)=(\alpha(x),\alpha(y))$, we see that $\alpha$ is of finite order and thus is diagonalizable. We conjugate $\varphi$ by some element $(x,y)\mapsto (\gamma^{-1}(x),\gamma(y))$ and get $(x,y) \mapsto (\gamma\alpha\gamma^{-1}(y),x)$. Hence, {\it any element of finite order of $\Aut(\mathbb{P}^1\times \mathbb{P}^1) \backslash \PGLn{2} \times \PGLn{2}$ is conjugate to some element $(x,y)\mapsto (\lambda y,x)$, for some $\lambda \in \K^{*}$ of finite order.}

We see then that if $G \subset \Aut(\mathbb{P}^1\times \mathbb{P}^1)$, not contained in $\PGLn{2} \times \PGLn{2}$, is a cyclic group, it is conjugate to the first case of the proposition (and thus conjugate to some group of automorphisms of\hspace{0.2 cm}$\Pn$). If $G$ is not cyclic, it must be generated by the group $H =G\cap (\PGLn{2} \times \PGLn{2})$ and some element conjugate to $(x,y)\mapsto (\lambda y,x)$, for some $\lambda \in \K^{*}$ of finite order. As the group $G$ is abelian, every element of $H$ is of the form $(x,y) \mapsto (\alpha(x),\alpha(y))$, where $\alpha(\lambda x)=\lambda \alpha(x)$. Computing the group of elements of $\PGL(2,\K)$ which commute with $\alpha$, there are three possible cases:
\begin{itemize}
\item
$\lambda=1$\\
In this case, we can conjugate the group by $(x,y) \mapsto (\gamma(x),\gamma(y))$ and suppose that $H$ is either diagonal or equal to the group generated by $(x,y) \mapsto (-x,-y)$ and $(x,y) \mapsto (x^{-1},y^{-1})$. In the first situation, the group is contained in $(\K^{*})^2 \rtimes \Z{2}$ (first case of the proposition) and the second situation gives $[10]$.
\item
$\lambda=-1$\\
Note that the group $H$ contains the square of $(x,y) \mapsto (-y,x)$, which is $(x,y) \mapsto (-x,-y)$. 

The group $H$ is either cyclic or generated by $(x,y) \mapsto (-x,-y)$ and $(x,y) \mapsto (x^{-1},y^{-1})$. 
If $H$ is cyclic, it is diagonal, since it contains $(x,y) \mapsto (-x,-y)$, so $G$ is contained in $(\K^{*})^2 \rtimes \Z{2}$. The second possibility gives $[9]$. 
\item
$\lambda \not=\pm 1$\\
The group $H$ must be diagonal and then $G$ is contained in $(\K^{*})^2 \rtimes \Z{2}$.\proofend
\end{itemize}
\end{proof}

\bigskip

We now study the birational equivalences between the finite abelian groups of automorphisms of $\mathbb{P}^1\times\mathbb{P}^1$, enumerated in Proposition \refd{Prp:AutP1P1}:
\begin{Prp}\fauxtitred
\label{Prp:AutP1Bir}
\begin{itemize}
\item[1.]
Any finite abelian subgroup of $\Aut(\mathbb{P}^1 \times\mathbb{P}^1)$ leaves invariant a pencil of rational curves.
\item[2.]
Up to birational conjugation, the non-trivial finite abelian subgroups of $\Aut(\mathbb{P}^1\times\mathbb{P}^1)$ are conjugate to one \emph{and only one} of the following:\end{itemize}
\begin{flushleft}\begin{tabular}{llll}
${\num{0.n}}$ & $G \cong \Z{n} $& g.b.&$(x,y) \mapsto (\zeta_n x,y)\hts$\\
${\num{0.mn}}$ & $G \cong \Z{n}\times\Z{m}$& g.b.&$(x,y) \mapsto (\zeta_n x,y)$ and
$(x,y) \mapsto (x,\zeta_m y)\hts$\\
${\num{P1.22n}}$ & $G \cong \Z{2}\times\Z{2n} $& g.b.&$(x,y) \mapsto (x^{-1},y)$ and
$(x,y) \mapsto (-x,\zeta_{2n} y)\hts$\\
${\num{P1.222n}}$ & $G \cong (\Z{2})^2\times\Z{2n} $& g.b.&$(x,y) \mapsto (\pm x^{\pm 1},y)\hts$ and
$(x,y) \mapsto (x,\zeta_{2n} y)$\\
${\num{P1.22.1}}$ & $G \cong (\Z{2})^2 $& g.b. &$(x,y) \mapsto (\pm x^{\pm 1},y)\hts$\\
${\num{P1.222}}$ & $G \cong (\Z{2})^3 $& g.b.&$(x,y) \mapsto (\pm x,\pm y)$ and
$(x,y) \mapsto (x^{-1},y)\hts$\\
${\num{P1.2222}}$ & $G \cong (\Z{2})^4$& g.b. &$(x,y) \mapsto (\pm x^{\pm 1},\pm y^{\pm 1})\hts$\\
${\num{P1s.24}}$ & $G \cong \Z{2}\times \Z{4}$& g.b. & $(x,y) \mapsto (x^{-1},y^{-1})$ and $(x,y) \mapsto (-y,x)$\\
&conjugate to the group & g.b. & $(x,y)\dasharrow (-x^{-1},-x^{-1}y)$ and \\
& & & $(x,y)\dasharrow (x^{-1},y^{-1})$\\
${\num{P1s.222}}$ & $G \cong (\Z{2})^3$& g.b. & $(x,y) \mapsto (-x,-y)$, $(x,y) \mapsto (x^{-1},y^{-1})$, 
\\
& & & and $(x,y) \mapsto (y,x)$ \\
&conjugate to the group & g.b. &  $(x,y) \dasharrow (x^{-1},y^{-1})$, $(x,y) \mapsto (x,-y)$, \\ 
& & & and $(x,y) \mapsto (x,x/y)$
\end{tabular}\\
{\upshape (}where $n$ is an integer $>1$, $n$ divides $m$ and $\zeta_n=e^{2\im \pi/n}${\upshape )}.\end{flushleft}
\begin{itemize}
\item[3.]
The groups $\nump{0.n}$ and $\nump{0.mn}$ are conjugate to subgroups of $\Aut(\Pn)$, since the others are not. 
\end{itemize}
\end{Prp}
\begin{remark}\fauxtitred
We write the groups $\nump{0.n}$ and $\nump{0.mn}$ differently from the others, using the notation for subgroups of $\Aut(\Pn)$ (see Proposition \refd{Prp:PGL3Cr}).
\end{remark}
\begin{proof}\fauxtitred\upshape
Note that the first assertion is not immediately obvious from the classification of Proposition \refd{Prp:AutP1P1}. Indeed, the groups $[2],[3],...,[8]$ leave invariant the fibres of the two projections, but this is not true for the other cases. However, the second assertion of this proposition directly implies the first one, since any of the groups given above leaves invariant the fibres of the first projection.

To prove assertion $2$, we use the classification of Proposition \refd{Prp:AutP1P1} and conjugate those groups to the groups of the above list. Incidentally, we see that some groups non-conjugate in $\Aut(\mathbb{P}^1\times\mathbb{P}^1)$ are birationally conjugate. (For example [4], [5] and [6].)
\begin{itemize}
\item[[1\hspace{-0.15cm}]]
Let $G$ be a finite abelian subgroup of the group $(\K^{*})^2 \rtimes \Z{2}$ of Example $\refd{Exa:TorP1}$ (case [1]), conjugate to a subgroup of automorphisms of $\mathbb{P}^2$. By the classification of finite abelian subgroups of $\Aut(\Pn)$ given in Proposition \refd{Prp:PGL3Cr}, as $G$ has fixed points, it is diagonalizable and is then conjugate to $\nump{0.n}$ or $\nump{0.mn}$. Note that these two classes of groups leave the two projections invariant.
\item[[2\hspace{-0.15cm}]]
Let $G \cong (\Z{2})^2\times\Z{m}$, generated by $(x,y) \mapsto (\pm x^{\pm 1},y)$ and
$(x,y) \mapsto (x,\zeta_m y)$. If $m$ is even, we directly get $\nump{P1.222n}$, for $n=\frac{m}{2}$. Otherwise, we conjugate the group by the birational map $\varphi:(x,y) \dasharrow (x, y(x+x^{-1}))$. Note that $\varphi$ commutes with $(x,y) \mapsto (x^{-1},y)$ and $(x,y) \mapsto (x,\zeta_m y)$, and conjugates $(x,y) \mapsto (-x,y)$ to $(x,y) \mapsto (-x,-y)$. We get the group $\nump{P1.22n}$, for $n=m$.
\item[[3\hspace{-0.15cm}]]
This group is equal to $\nump{P1.22n}$.
\item[[4-6\hspace{-0.15cm}]]\hspace{0.05cm}
Note that the groups $[4]$, $[5]$ and $[6]$ have the same action on the first coordinate. As there is a section of the morphism induced by this action, we can see that the three groups are conjugate in $\CrPp$, by using Lemma \refd{Lem:splitsequence} (proved later). We give below the explicit conjugation between the groups, which are then all conjugate to $\nump{P1.22.1}$, which is $[4]$. Explicitly:
\begin{itemize}
\item 
We conjugate $[4]$ by the birational map $\varphi:(x,y) \dasharrow (x,\frac{y+x^2}{y+x^{-2}})$. This map commutes with $(x,y) \mapsto (-x,y)$ and conjugates $(x,y) \mapsto (x^{-1},y)$ to $(x,y)\mapsto (x^{-1},y^{-1})$. The group obtained is conjugate to $[5]$ by an automorphism of $\mathbb{P}^1\times\mathbb{P}^1$.
\item
The birational map $\varphi:(x,y) \dasharrow (x,x\frac{y+x^2}{y+x^{-2}})$ conjugates $[4]$ to $[6]$. Indeed, $\varphi$ conjugates $(x,y) \mapsto (-x,y)$ to $(x,y)\mapsto (-x,-y)$ and conjugates $(x,y) \mapsto (x^{-1},y)$ to $(x,y) \mapsto (x^{-1},y^{-1})$.
\end{itemize}
\item[[7-8\hspace{-0.15cm}]]
The groups [7] and [8] are equal respectively to $\nump{P1.222}$ and $\nump{P1.2222}$.
\item[[9\hspace{-0.15cm}]]
Via the birational map $(x,y)\dasharrow (x:y:1)$ from $\mathbb{P}^1 \times\mathbb{P}^1$ to $\Pn$, the group $[9]$ is conjugate to the group of birational maps of\hspace{0.2 cm}$\Pn$ generated by $(x:y:z) \dasharrow (yz:xz:xy)$ and $(x:y:z) \mapsto (-y:x:z)$, which leaves invariant the pencil of lines passing through $(0:0:1)$ (lines of the form $ax+by=0$, $(a:b) \in \mathbb{P}^1$). We conjugate this group by the birational morphism $(x:y:z) \dasharrow (x:y)\times (x:z)$ and get the group of birational maps of $\mathbb{P}^1\times\mathbb{P}^1$ generated by $(x,y)\dasharrow (-x^{-1},-x^{-1}y)$ and $(x,y)\dasharrow (x^{-1},y^{-1})$. This gives $\nump{P1s.24}$.
\item[[10\hspace{-0.15cm}]]
We conjugate the group [10] by the birational map $(x,y) \dasharrow (x/y,x)$ of $\mathbb{P}^1 \times\mathbb{P}^1$ and get the group generated by $(x,y) \dasharrow (x^{-1},yx^{-1})$, $(x,y) \mapsto (x,-y)$, $(x,y) \mapsto (x^{-1},y^{-1})$, which leaves invariant the first projection. This yields $\nump{P1s.222}$
\end{itemize}

We now prove that distinct groups of the list are not birationally conjugate.

First of all, $\nump{0.n}$ and $\nump{0.mn}$ fix at least one point of $\mathbb{P}^1\times\mathbb{P}^1$. Since the other groups of the list don't fix any point, they are not conjugate to $\nump{0.n}$ or $\nump{0.mn}$.

Consider the other groups. There are only two pairs of isomorphic distinct groups, namely:
\begin{itemize}
\item
$\nump{P1s.222}$ and $\nump{P1.222}$, both isomorphic to $(\Z{2})^3$,
\item
$\nump{P1s.24}$ and $\nump{P1.22n}$, with $n=2$, both isomorphic to $\Z{2}\times\Z{4}$.
\end{itemize}
We use the invariant defined in Section \refd{Sec:PairInvariantPencils}:

Both $\nump{P1.222}$ and $\nump{P1.22n}$ leave invariant two pencils of rational curves (the fibres of the two projections $\mathbb{P}^1\times\mathbb{P}^1\rightarrow \mathbb{P}^1$) which intersect freely in exactly one point. We prove that this is not the case for $\nump{P1s.222}$ and $\nump{P1s.24}$; this shows that distinct groups are not birationally conjugate.

The proof is similar for both cases. Take $G\subset \Aut(\mathbb{P}^1 \times \mathbb{P}^1)$ to be the representation of one of these two groups as a group of automorphisms of $\mathbb{P}^1\times\mathbb{P}^1$ which permute the two fibrations. We have then $\Pic{\mathbb{P}^1 \times \mathbb{P}^1}^G=\mathbb{Z}d$, where $d=-\frac{1}{2}K_{\mathbb{P}^1\times\mathbb{P}^1}$ denotes the diagonal of $\mathbb{P}^1 \times \mathbb{P}^1$.
Suppose that there exist two $G$-invariant pencils $\Lambda_1=n_1 d$ and $\Lambda_2=n_2 d$ of rational curves, for some positive integers $n_1,n_2$. The intersection $\Lambda_1 \cdot \Lambda_2=2 n_1 n_2$ is an even integer. Note that the fixed part of the intersection is also even, as the order of $G$ is $8$ and $G$ acts without fixed point on $\mathbb{P}^1\times\mathbb{P}^1$. The free part of the intersection is then also an even integer and hence is not $1$.

We now prove the last assertion.
The finite abelian groups of automorphisms of\hspace{0.2 cm}$\Pn$ are conjugate either to $\nump{0.n}$, $\nump{0.mn}$, or to the group $\nump{0.V9}$, isomorphic to $(\Z{3})^2$ (see Proposition \refd{Prp:PGL3Cr}). As no group of the list $[3],...,[10]$ given above is isomorphic to $(\Z{3})^2$, we have the third assertion.
\proofend
\end{proof}

\section{Automorphisms of surfaces of degree $5$, $6$ or $7$}
\pagetitless{Finite abelian groups of automorphisms of surfaces of large degree}{Automorphisms of surfaces of degree $5$, $6$ or $7$}
\index{Del Pezzo surfaces!of degree $5$}\index{Del Pezzo surfaces!of degree $6$}\index{Del Pezzo surfaces!of degree $7$}
We have described all the finite abelian groups of automorphisms of surfaces of degree $8$ or $9$. We now treat the next cases.

We begin with a simple result which will often be used in the study of automorphisms of conic bundles and will help us to prove Proposition \refd{Prp:Degree567}.

\begin{Lem}\fauxtitre
\label{Lem:Degree567CB}
\index{Conic bundles!automorphisms|idi}
Let $G \subset \Aut(S,\pi)$ be a finite abelian group of automorphisms of the conic bundle $(S,\pi)$, such that:
\begin{itemize}
\item
 $\pi$ has at most $3$ singular fibres (or equivalently such that $K_S^2\geq 5$);
 \item
the triple $(G,S,\pi)$ is minimal.
 \end{itemize}
 Then, $S$ is a Del Pezzo surface of degree $5$ or $6$.\index{Del Pezzo surfaces!of degree $5$}\index{Del Pezzo surfaces!of degree $6$}
\end{Lem}
\begin{proof}\fauxtitred\upshape
Let $-n$ be the minimal self-intersection of sections of $\pi$ and let $r$ be the number of singular fibres of $\pi$. 
As $(G,S,\pi)$ is minimal, every singular fibre is twisted by some element of $G$. Using Lemma \refd{Lem:GoingToF0F1}, we get $r\geq 2n\geq 2$, whence $r =2$ or $3$ and $n=1$, and the existence of some birational morphism of conic bundles (not necessarily $G$-equivariant) $p_1:S\rightarrow \mathbb{F}_1$. So the surface is obtained by the blow-up of $2$ or $3$ points of $\mathbb{F}_1$, not on the exceptional section, and thus by blowing-up $3$ or  $4$ points of\hspace{0.2 cm}$\Pn$, no $3$ of which are collinear (otherwise we would have a section of self-intersection $-2$). The surface is then a Del Pezzo surface of degree $6$ or $5$.
\proofend
\end{proof}

\bigskip

\begin{Prp}\fauxtitre
\label{Prp:Degree567}
Let $S$ be a rational surface of degree $5$, $6$ or $7$. Let $G \subset \Aut(S)$ be a finite subgroup of automorphisms of $S$ such that the pair $(G,S)$ is minimal.

Then, $S$ is a Del Pezzo surface of degree $5$ or $6$.\index{Del Pezzo surfaces!of degree $5$}\index{Del Pezzo surfaces!of degree $6$}\index{Del Pezzo surfaces!of degree $7$}
\end{Prp}
\begin{proof}\fauxtitred\upshape
We first prove that $S$ is a Del Pezzo surface.
Since the pair $(G,S)$ is minimal, either $\rkPic{S}^G=1$ or $G$ preserves a conic bundle structure (Proposition \refd{Prp:TwoCases}). In the first case, $S$ is a Del Pezzo surface (Lemma \refd{Lem:rkDelPezzo}), and we are done. The second case follows from Lemma \refd{Lem:Degree567CB}.

Finally, note that the group of automorphisms of the Del Pezzo surface of degree $7$ is not minimal. This has already been annouced and is easy to see. Indeed, there are three exceptional divisors on the surface (see Proposition \refd{Prp:DescExcCurvDelPezzo}), and only one intersects the other two. This one is then invariant, so we may contract it and get an equivariant birational morphism to $\mathbb{P}^1\times\mathbb{P}^1$.\proofend
\end{proof}

\subsection{The Del Pezzo surface of degree $6$}
\index{Del Pezzo surfaces!of degree $6$|idb}
\label{Subsec:DelPezzo6}
There is a single isomorphism class of Del Pezzo surfaces of degree $6$, since all sets of three non-collinear points of\hspace{0.2 cm}$\Pn$ are equivalent under the action of linear automorphisms.

Consider the surface $S$ of degree $6$ defined by the blow-up of the points $A_1=(1:0:0)$, $A_2=(0:1:0)$ and $A_3=(0:0:1)$. We can view it in $\mathbb{P}^6$ as the image of the rational map $(x:y:z) \dasharrow (x^2 y:x^2z:xy^2:xyz:xz^2:y^2 z:yz^2)$, given by the linear system of cubics passing through $A_1$, $A_2$ and $A_3$. We may also view it in $\Pn \times \Pn$, defined as $\{ (x:y:z) \times (u:v:w)\ | \ ux=vy=wz\}$, where the blow-down is the projection on one copy of\hspace{0.2 cm}$\Pn$, explicitly $p: (x:y:z) \times (u:v:w) \mapsto (x:y:z)$.

If $A_i A_j$ denotes the line passing through $A_i$ and $A_j$, for $1 \leq i,j \leq 3$, the 6 exceptional divisors of $S$ are (see Proposition \refd{Prp:DescExcCurvDelPezzo}):\index{Del Pezzo surfaces!curves on the surfaces|idi}
\begin{center}
\hspace{-3cm}
\begin{tabular}{ll}
$E_1= \{(1:0:0) \times (0:a:b) \ | \ (a:b) \in \mathbb{P}^1\}= p^{-1}(A_1)$\\
$E_2= \{(0:1:0) \times (a:0:b) \ | \ (a:b) \in \mathbb{P}^1\}= p^{-1}(A_2)$\\
$E_3= \{(0:0:1) \times (a:b:0) \ | \ (a:b) \in \mathbb{P}^1\}= p^{-1}(A_3)$\\
$D_{23}=\{(0:a:b) \times (1:0:0) \ | \ (a:b) \in \mathbb{P}^1\}=  \tilde{p}^{-1}(A_2 A_3)$\\
$D_{13}=\{(a:0:b) \times (0:1:0) \ | \ (a:b) \in \mathbb{P}^1\}= \tilde{p}^{-1}(A_1 A_3)$\\
$D_{12}=\{(a:b:0) \times (0:0:1) \ | \ (a:b) \in \mathbb{P}^1\}=  \tilde{p}^{-1}(A_2 A_2)$.\end{tabular}\\
\begin{flushright} \vspace{-3.5cm}
\includegraphics[width=30.00mm]{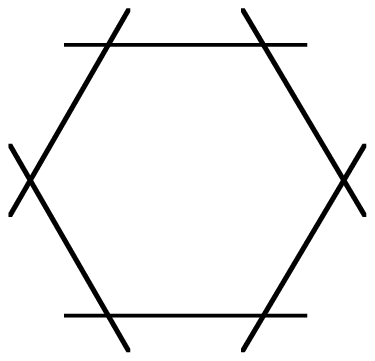}%
\drawat{-18.1mm}{24mm}{$D_{23}$}%
\drawat{-6.64mm}{18mm}{$E_3$}%
\drawat{-6.64mm}{6.6mm}{$D_{13}$}%
\drawat{-17.8mm}{0.1mm}{$E_1$}%
\drawat{-29.0mm}{6.6mm}{$D_{12}$}%
\drawat{-28.0mm}{18mm}{$E_2$}
\end{flushright}\end{center}
The Picard group of $S$ is of rank $4$, generated by $E_1$, $E_2$, $E_3$ and $L=p^{-1}(H)$, where $H$ denotes a general line of\hspace{0.2 cm}$\Pn$. The exceptional divisors $D_{12}$, $D_{13}$ and $D_{23}$ are linearly equivalent respectively to $L-E_1-E_2$, $L-E_1-E_3$ and $L-E_2-E_3$.

\index{Del Pezzo surfaces!automorphisms}
\begin{Lem}\fauxtitre
\label{Lem:AutDP6}
\begin{itemize}
\item
The exceptional divisors form a hexagon: it is connected and each divisor touches two others.
\item
The action of $\Aut(S)$ on the hexagon gives rise to the exact sequence 
\begin{center}$1\rightarrow (\K^{*})^2 \rightarrow \Aut(S) \stackrel{\rho}{\rightarrow} \Sym_3 \times \Z{2} \rightarrow 1$.\end{center}
\item
The exact sequence splits and $\Aut(S)=(\K^{*})^2 \rtimes (\Sym_3 \times \Z{2})$, where:
\begin{itemize}
\item
$(\K^{*})^2$ is generated by automorphisms of the form \\
$(x:y:z)\times (u:v:w) \mapsto (x:\alpha y: \beta z) \times (\alpha\beta u:\beta v:\alpha w)$, $\alpha,\beta \in \K^{*}$.
\item
$(\K^{*})^2 \rtimes \Sym_3$ is the lift on $S$ of the group of automorphisms of\hspace{0.2 cm}$\Pn$ that leave invariant the set $\{A_1,A_2,A_3\}$. 
\item
$\Z{2}$ is generated by the automorphism \begin{center}$(x:y:z) \times (u:v:w) \mapsto (u:v:w) \times (x:y:z)$,\end{center} which corresponds to the standard quadratic transformation\\ $(x:y:z) \dasharrow (yz:xz:xy)$ of\hspace{0.2 cm}$\Pn$. It exchanges $E_i$ and $D_{jk}$, for $\{i,j,k\}=\{1,2,3\}$.
\item
$\Sym_3$ acts on the torus $(\K^{*})^2$ by permuting the coordinates, and the action of $\Z{2}$ is the inversion. 
\end{itemize}
\end{itemize}
\end{Lem}
\begin{proof}\fauxtitred\upshape
The first assertion follows directly from the description of exceptional divisors given above. By rotating the hexagon we find $E_1$, $D_{12}$, $E_2$, $D_{23}$, $E_3$, $D_{13}$ and then $E_1$ again.

As $\Aut(S)$ preserves the exceptional divisors and the intersection form, it must preserve the hexagon. So the action of $\Aut(S)$ on the hexagon gives rise to a homomorphism \begin{center}$\rho:\Aut(S)\rightarrow \Sym_3 \times \Z{2}$.\end{center}

As any element of the kernel leaves invariant every exceptional divisor, it comes from an automorphism of\hspace{0.2 cm}$\Pn$ that fixes the three points $A_1,A_2$ and $A_3$. The kernel of $\rho$ thus consists of automorphisms of the form $(x:y:z)\times (u:v:w) \mapsto (x:\alpha y: \beta z) \times (\alpha\beta u:\beta v:\alpha w)$, with $\alpha,\beta \in \K^{*}$, and is the lift of the torus $\mathcal{T}$ of diagonal automorphisms of\hspace{0.2 cm}$\Pn$.

Note that the group $\Sym_3$ of permutations of the variables $x$, $y$ and $z$ (and the corresponding variables $u,v$ and $w$), generated by the two automorphisms 
\begin{center}\begin{tabular}{l}$(x:y:z)\times (u:v:w) \mapsto (y:x:z) \times (v:u:w),$\\
$(x:y:z)\times (u:v:w) \mapsto (z:y:x) \times (w:v:u),$\end{tabular}
\end{center} is sent by $\rho$ on $\Sym_3$. The group generated by the automorphism 
\begin{center}$(x:y:z) \times (u:v:w) \mapsto (u:v:w) \times (x:y:z)$\end{center} is sent by $\rho$ on $\Z{2}$. This gives the surjectivity of $\rho$ and an obvious section. The other assertions are evident.\proofend
\end{proof}

\bigskip

\begin{Prp}\fauxtitre
\label{Prp:DP6conjLin}
Let $G$ be an abelian subgroup of automorphisms of the Del Pezzo surface $S$ of degree $6$ such that $\rkPic{S}^{G}=1$. Then $G$ is birationally conjugate to a subgroup of $\Aut(\Pn)$ of order $6$.
\end{Prp}\begin{proof}\fauxtitred\upshape
Let $G$ be an abelian subgroup of $\Aut(S)$ such that $\rkPic{S}^{G}=1$. In Lemma \refd{Lem:SizeOrbits}, we proved that the size of the orbits on the exceptional divisors are divisible by the degree of the surface. The action of $G$ on the $6$ exceptional divisors is thus transitive, which implies that the image of $G$ by $\rho$ is isomorphic to $\Z{6}$, since the group is abelian. Let $g \in G$ be such that $\rho(g)$ generates the group $\rho(G)$. We have then
\begin{center} $g:(x:y:z) \times (u:v:w) \mapsto (\alpha v:\beta w:u) \times (\beta y:\alpha z:\alpha\beta x)$,\end{center} where $\alpha,\beta \in \K^{*}$. As the only element of the kernel of $\rho$ which commutes with $g$ is the identity, $G$ must be cyclic, generated by $g$. Conjugating it by \begin{center}$p:(x:y:z) \times (u:v:w) \mapsto (\beta x: y:\alpha z) \times (\alpha u:\alpha\beta v:\beta w)$,\end{center}
we may assume that $\alpha=\beta=1$.

Consider the birational map $g'=p gp^{-1}$ of\hspace{0.2 cm}$\Pn$, explicitly $g':(x:y:z) \dasharrow (xz:xy: yz)$. As $g$ is an automorphism of the surface, it fixes the canonical divisor $K_S$, so the birational map
$g'$ leaves invariant the linear system of cubics of\hspace{0.2 cm}$\Pn$ passing through $A_1,A_2$ and $A_3$, i.e.\ $g'(C^3(A_1,A_2,A_3))=C^3(A_1,A_2,A_3)$ (this can also be verified directly).

Note that $g'$ fixes exactly one point of\hspace{0.2 cm}$\Pn$, namely $P=(1:1:1)$, and that its action on the projective tangent space $\mathbb{P}(T_{P}(\Pn))$ of\hspace{0.2 cm}$\Pn$ at $P$ is of order $3$, with two fixed points, corresponding to the lines $(x-y)+\omega^k(y-z)=0$, where $\omega=e^{2\im\pi/3}$, $k=1,2$. Hence, by choosing $P'$ to be one of the two points infinitely near to $P$ fixed by $g'$, the linear system $C^3(A_1,A_2,A_3,P^2\leftarrow P')$ is invariant by $g'$. This is a homoloidal net, by Proposition \refd{Prp:HomoNet}. The birational map $g'$ is conjugate to the linear automorphisms of order $6$ (see Proposition \refd{Prp:HomolNetInvariant}).\proofend
\end{proof}

\bigskip

\index{Conic bundles!on the Del Pezzo surface of degree $6$|idb}
Consider the three conic bundle structures of $S$ (see Proposition \refd{Prp:NbConicDP}) and the corresponding groups of automorphisms. We know (Proposition \refd{Prp:InterDPCB}) that their finite subgroups are not minimal, but we will use some of these results later, when looking at other groups of automorphisms of conic bundles (see for example Section \refd{SubSec:DelPezzo6Kappa}).

Let $\pi_1:S\rightarrow \mathbb{P}^1$ be the morphism defined by:
\begin{center}
$\pi_1((x:y:z) \times (u:v:w))=\left\{\begin{array}{lll}
(y:z)& \mbox{ if } &(x:y:z) \not= (1:0:0),\\
(w:v)& \mbox{ if } &(u:v:w) \not= (1:0:0).\end{array}\right.$
\end{center}
Note that $p$ sends the fibres of $\pi_1$ on lines of\hspace{0.2 cm}$\Pn$ passing through $A_1$. There are exactly two singular fibres of this fibration, namely
\begin{center}
\begin{tabular}{lll}
$\pi_1^{-1}(1:0)=\{E_2,D_{12}\}$& and &$\pi_1^{-1}(0:1)=\{E_3,D_{13}\}$;
\end{tabular}\end{center}
$E_1$ and $D_{23}$ are sections of $\pi_1$.

As in Lemma \refd{Lem:AutDP6}, we determine the structure of the group of automorphisms of $S$ that leave the conic bundle invariant:

\begin{Lem}\fauxtitred
\label{Lem:AutDP6CB}
\begin{itemize}
\item
The group $\Aut(S,\pi_1)$ of automorphisms of the conic bundle acts on the hexagon and leaves invariant the set $\{E_1,D_{23}\}$.
\item
The action on the hexagon gives rise to the exact sequence 
\begin{center}
$1\rightarrow (\K^{*})^2 \rightarrow \Aut(S,\pi_1) \rightarrow (\Z{2})^2 \rightarrow 1.$
\end{center}
\item
The exact sequence splits and $\Aut(S,\pi_1)= (\K^{*})^2\rtimes (\Z{2})^2$, where 
\begin{itemize}
\item
$(\K^{*})^2$ is generated by automorphisms of the form\\
$(x:y:z)\times (u:v:w) \mapsto (x:\alpha y: \beta z) \times (\alpha\beta u:\beta v:\alpha w)$,
$\alpha,\beta \in \K^{*}$.
\item
The group $ (\Z{2})^2$ is generated by the automorphisms
\begin{center}
$(x:y:z) \times (u:v:w) \mapsto (x:z:y) \times (u:w:v)$,
\end{center}
whose action on the set of exceptional divisors is  $(E_2 \ E_3)(D_{12}\ D_{13})$;
\begin{center}
$(x:y:z) \times (u:v:w) \mapsto  (u:v:w) \times (x:y:z),$
\end{center}
whose action is $(E_1 \ D_{23})(E_2\ D_{13})(E_3\ D_{12})$.
\end{itemize}
\end{itemize}
\end{Lem}
\begin{proof}\fauxtitred\upshape
Since $\Aut(S)$ acts on the hexagon, so does $\Aut(S,\pi_1) \subset \Aut(S)$. The group $\Aut(S,\pi_1)$ sends a fibre on a fibre and hence a section on a section. The set $\{E_1,D_{23}\}$ is therefore invariant.

It is clear that the group $(\K^{*})^2$ leaves the conic bundle invariant, so it is the kernel of the action of $\Aut(S,\pi_1)$ on the hexagon. As the set $\{E_1,D_{23}\}$ is invariant, the image is contained in the group $(\Z{2})^2$ generated by $(E_2 \ E_3)(D_{12}\ D_{13})$ and $(E_1 \ D_{23})(E_2\ D_{13})(E_3\ D_{12})$. The rest of the lemma follows directly.\proofend
\end{proof}

\bigskip

By permuting coordinates, we have two other conic bundle structures on the surface $S$, given by the following morphisms $\pi_2,\pi_3:S \rightarrow \mathbb{P}^1$:
\begin{center}
$\pi_2((x:y:z) \times (u:v:w))=\left\{\begin{array}{lll}
(x:z)& \mbox{ if }& (x:y:z) \not= (0:1:0),\\
(w:u)& \mbox{ if }& (u:v:w) \not= (0:1:0).\end{array}\right.$\vspace{0.1cm}\\
$\pi_3((x:y:z) \times (u:v:w))=\left\{\begin{array}{lll}
(x:y)& \mbox{ if }& (x:y:z) \not= (0:0:1),\\
(v:u)& \mbox{ if }& (u:v:w) \not= (0:0:1).\end{array}\right.$
\end{center}

\begin{Lem}\fauxtitre
\label{Lem:DP6CBPasMin}
The pair $(\Aut(S,\pi_1),S)$ is not minimal. Furthermore, $\Aut(S,\pi_1)$ is conjugate by the birational morphism $q \mapsto (\pi_2(q),\pi_3(q))$ to a subgroup of $\Aut(\mathbb{P}^1\times \mathbb{P}^1)$.

Similarly, the birational morphism $q \mapsto (\pi_1(q),\pi_3(q))$ {\upshape (}respectively $q \mapsto (\pi_1(q),\pi_2(q)${\upshape )} conjugates $\Aut(S,\pi_2)$ {\upshape (}respectively $\Aut(S,\pi_3)${\upshape )} to a subgroup of $\Aut(\mathbb{P}^1\times \mathbb{P}^1)$.
\end{Lem}
\begin{proof}\fauxtitred\upshape
We observe that the union of the sections $E_1$ and $D_{23}$ is invariant by the action of the whole group $\Aut(S,\pi_1)$. As these two exceptional divisors don't intersect, we can contract both of them and get a birational $\Aut(S,\pi_1)$-equivariant morphism from $S$ to $\mathbb{P}^1 \times \mathbb{P}^1$. The pair $(\Aut(S,\pi_1),S)$ is then not minimal. 
Explicitly, the birational morphism is given by $q \mapsto (\pi_2(q),\pi_3(q))$, as stated in the lemma.

We get the other cases by permuting coordinates.\proofend
\end{proof}
\begin{remark}\fauxtitred The subgroup of $\Aut(\mathbb{P}^1 \times \mathbb{P}^1)$ obtained in this manner doesn't leave any of the two fibrations of $\mathbb{P}^1 \times \mathbb{P}^1$ invariant.\end{remark}

\bigskip

\subsection{The Del Pezzo surface of degree $5$}
\index{Del Pezzo surfaces!of degree $5$|idb}
There is a single isomorphism class of Del Pezzo surfaces of degree $5$, since all sets of four non collinear points of\hspace{0.2 cm}$\Pn$ are equivalent under the action of linear automorphisms. Let $S$ be the surface obtained by the blow-up of the four points $A_1=(1:0:0)$, $A_2=(0:1:0)$, $A_3=(0:0:1)$ and $A_4=(1:1:1)$, and let $\pi:S\rightarrow \Pn$ be the blow-down. 

The exceptional divisors on $S$ are (see Proposition \refd{Prp:DescExcCurvDelPezzo}):

\begin{center}
\hspace{-3cm}
\begin{tabular}{ll}\index{Del Pezzo surfaces!curves on the surfaces|idi}
$E_1=\pi^{-1}(A_1),...,E_4=\pi^{-1}(A_4)$, the $4$ pull-backs \\
of the points $A_1$,...,$A_4$;\\
$D_{ij}=\tilde{\pi^{-1}}(A_i A_j)$ for $1\leq i,j\leq 4$, $i\not=j$, the $6$ strict pull-backs\\ of the lines through $2$ of the $A_i$'s.\end{tabular}\\
\begin{flushright} \vspace{-2.5cm}
\includegraphics[width=30.00mm]{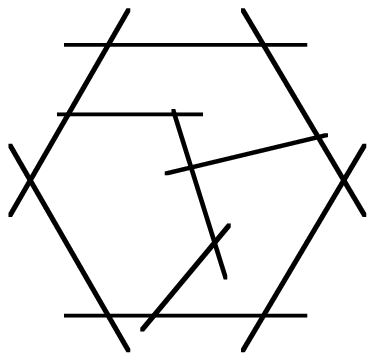}%
\drawat{-18.1mm}{24.5mm}{$D_{23}$}%
\drawat{-6.6mm}{18.5mm}{$E_3$}%
\drawat{-6.64mm}{6.6mm}{$D_{13}$}%
\drawat{-17.8mm}{0.1mm}{$E_1$}%
\drawat{-29.8mm}{6.6mm}{$D_{12}$}%
\drawat{-29.5mm}{18.5mm}{$E_2$}%
\drawat{-22mm}{19.5mm}{$D_{24}$}%
\drawat{-13.5mm}{17mm}{$D_{34}$}%
\drawat{-14mm}{11.5mm}{$E_4$}%
\drawat{-21mm}{7.5mm}{$D_{14}$}
\end{flushright}\end{center}

There are thus $10$ exceptional divisors. Each intersects three others:
\begin{itemize}
\item
Each $E_i$ intersects the strict transforms of the $3$ lines passing through $A_i$.
\item
The divisor $D_{ij}$ intersects $E_i$ and $E_j$ and the strict transform of the line passing through the other two points.
\end{itemize}

\break

\begin{Prp}\fauxtitred
\begin{itemize}
\item
There are $5$ sets of $4$ skew exceptional divisors on $S$, namely
\begin{itemize}
\item[]
$F_1=\{E_1,D_{23},D_{24},D_{34}\},$
\item[]
$F_2=\{E_2,D_{13},D_{14},D_{34}\},$
\item[]
$F_3=\{E_3,D_{12},D_{14},D_{24}\},$
\item[]
$F_4=\{E_4,D_{12},D_{13},D_{23}\},$
\item[]
$F_5=\{E_1,E_2,E_3,E_4\}$.
\end{itemize}
\item
The action on these sets gives rise to an isomorphism $\rho:\Aut(S) \mapsto \Sym_{5}$. 
\end{itemize}
\end{Prp}
\begin{proof}\fauxtitred\upshape
\begin{itemize}
\item
The first assertion follows by direct observation, or from Proposition \refd{Prp:DescExcCurvDelPezzo}. 
\item
Let us prove that $\rho$ is injective. An automorphism of $S$ that fixes the $5$ sets must fix the $10$ exceptional lines; it is therefore conjugate to a linear automorphism of\hspace{0.2 cm}$\Pn$ which fixes the four points, so it is the identity. 
\item
Let us prove that $\rho$ is surjective.
\begin{itemize}
\item
The group $\Aut(S)$ contains the lift on $S$ of the group $\Sym_4$ of automorphisms of\hspace{0.2 cm}$\Pn$ that leave the set $\{A_1,A_2,A_3,A_4\}$ invariant. This group embeds naturally in $\Sym_5$ as the group of permutations of $F_1,...,F_4$, which fix $F_5$.
\item
Let us show that the lift on $S$ of the birational map $\alpha:(x:y:z) \dasharrow (yz:xz:xy)$ gives an automorphism of $S$ whose image by $\rho$ is the transposition $(F_4 \ F_5)$.

Let $S_6$ be the Del Pezzo surface of degree $6$ obtained by blowing-up $A_1,A_2,A_3 \in \Pn$. As we saw in Section \refd{Subsec:DelPezzo6}, the lift of $\alpha$ on the surface $S_6$ is an automorphism whose action on the $6$ exceptional divisors is $(E_1\ D_{23})(E_2\ D_{13})(E_3\ D_{12})$.\ Note that $\alpha$ fixes the point $A_4$, so the blow-up of this point on $S_6$ gives rise to an automorphism of $S$ whose action on the $10$ exceptional divisors is the same: $(E_1\ D_{23})(E_2\ D_{13})(E_3\ D_{12})$ (it leaves invariant $E_4, D_{14},D_{24},D_{34}$). Its image by $\rho$ is therefore the transposition $(F_4\ F_5)$ as stated above.
\end{itemize}
The image of $\rho$ in $\Sym_5$ contains both $\Sym_4$ and the permutation $(F_4\ F_5)$, and is therefore the whole group $\Sym_5$.
\proofend
\end{itemize}
\end{proof}

\bigskip

\begin{Prp}\fauxtitre
\label{Prp:DP5conjLin}\index{Del Pezzo surfaces!automorphisms}
Let $G$ be an abelian subgroup of automorphisms of the Del Pezzo surface $S$ of degree $5$ such that $\rkPic{S}^{G}=1$. Then $G$ is conjugate, in $\CrP$, to a group of automorphisms of\hspace{0.2 cm}$\Pn$ of order $5$.
\end{Prp}\begin{proof}\fauxtitred\upshape
Let $G$ be an abelian subgroup of $\Aut(S)$ such that $\rkPic{S}^{G}=1$. By Lemma \refd{Lem:SizeOrbits}, the order of $G$ is divisible by $5$. This implies that $G$ is a cyclic subgroup of $\Sym_5$ of order $5$. We may assume that it is generated by the element $(F_1 \ F_2 \ F_3 \ F_4 \ F_5)=(F_1 \ F_2 \ F_3 \ F_4)(F_4 \ F_5)$. 

The conjugation by $\pi$ of the generator of $G$ to a birational map of\hspace{0.2 cm}$\Pn$ yields $h:(x:y:z) \dasharrow (xy:y(x-z):x(y-z))$ which has two fixed points, namely $(\zeta+1:\zeta:1)$, where $\zeta^2-\zeta-1=0$. Denoting one of them by $P$, the linear system $C^3(A_1,A_2,A_3,A_4,P^2)$ is invariant by $h$.
As this is a homoloidal net (see \refd{Prp:HomoNet}), the birational map $h$ is conjugate to the linear automorphisms of order $5$, by Proposition \refd{Prp:HomolNetInvariant}.\proofend
\end{proof}

\bigskip

\section{Conjugation to a subgroup of $\Aut(\mathbb{P}^1\times\mathbb{P}^1)$ or $\Aut(\Pn)$}
\pagetitless{Finite abelian groups of automorphisms of surfaces of large degree}{Conjugation to a subgroup of $\Aut(\mathbb{P}^1\times\mathbb{P}^1)$ or $\Aut(\Pn)$}
We are now equipped to prove the following result, announced at the beginning of this chapter:
\begin{Prp}\fauxtitre
\label{Prp:LargeDegree}
Let $G \subset \Aut(S)$ be a finite abelian group of automorphisms of a surface $S$ and suppose that $K_S^2\geq 5$. Then, $G$ is birationally conjugate to a subgroup of $\Aut(\mathbb{P}^1\times\mathbb{P}^1)$ or of $\Aut(\Pn)$.

More precisely:
\begin{itemize}
\item
If $G$ is cyclic, it is conjugate both to a subgroup of $\Aut(\mathbb{P}^1\times\mathbb{P}^1)$ {\emph{and}} to a subgroup of $\Aut(\Pn)$.
\item
If $G\not\cong (\Z{3})^2$, it is conjugate to a subgroup of $\Aut(\mathbb{P}^1\times\mathbb{P}^1)$.
\end{itemize}
\end{Prp}
\begin{proof}\fauxtitred\upshape
We first prove that $G$ is birationally conjugate to a group of automorphisms of $\Aut(\mathbb{P}^1\times\mathbb{P}^1)$ or $\Aut(\Pn)$.
We assume that the pair $(G,S)$ is minimal by blowing-down sets of disjoint $G$-invariant exceptional curves (this may increase $K_S^2$).
\begin{itemize}
\item
If the degree of $S$ is $8$, then $S$ is either $\Pn$ or a Hirzeburch surface $\mathbb{F}_n$, for some integer $n\geq 0$. If $S=\Pn$, we are done. Otherwise, Proposition \refd{Prp:FnAb} tells us that $G$ is birationally conjugate to a subgroup of $\Aut(\mathbb{P}^1\times\mathbb{P}^1)$. 
\item
If the degree of $S$ is $5$, $6$ or $7$, Proposition \refd{Prp:Degree567} tells us that $S$ is a Del Pezzo surface of degree $5$ or $6$.
As $(G,S)$ is minimal, either $\rkPic{S}^G=1$, or $G$ preserves a conic bundle structure on $S$. Proposition \refd{Prp:InterDPCB} tells us that the second case is not possible.

We then use Proposition \refd{Prp:DP6conjLin} if $K_S^2=6$ (or Proposition \refd{Prp:DP5conjLin} if $K_S^2=5$ ) to see that $G$ is conjugate to a subgroup of $\Aut(\Pn)$.\end{itemize}

\noindent We now prove the last two assertions:
\begin{itemize}
\item
By Proposition \refd{Prp:PGL3Cr}, any finite abelian subgroup of $\Aut(\Pn)$ is conjugate either to the group $V_9$, isomorphic to $(\Z{3})^2$, or to a diagonal group.
\item
By Proposition \refd{Prp:AutP1Bir}, any finite abelian subgroup of $\Aut(\mathbb{P}^1\times\mathbb{P}^1)$ is conjugate either to a diagonal group, or to a non-cyclic group. 
\end{itemize}
Since the diagonal tori of $\Aut(\mathbb{P}^1\times\mathbb{P}^1)$ and $\Aut(\Pn)$ are conjugate, we get the assertions.\proofend
\end{proof}

\chapter{Finite abelian groups of automorphisms of conic bundles}
\label{Chap:FiniteAbConicBundle}
\index{Conic bundles!automorphisms|idb}
\section{Two fundamental examples}
\pagetitless{Finite abelian groups of automorphisms of conic bundles}{Two fundamental examples}
In this section, we give three important examples that illustrate almost all possible finite abelian groups of automorphisms of conic bundles (see 
Proposition \refd{Prp:TwocasesKappaTwist}). 

\subsection{Automorphisms of the Del Pezzo surface of degree 6}\index{Del Pezzo surfaces!of degree $6$|idb}\index{Del Pezzo surfaces!automorphisms|idi}
\index{Conic bundles!on the Del Pezzo surface of degree $6$|idb}
\label{SubSec:DelPezzo6Kappa}
We denote by $S_6$ the Del Pezzo surface of degree $6$ viewed in $\Pn \times \Pn$ as
\begin{center}
$S_6=\{ (x:y:z) \times (u:v:w)\ | \ ux=vy=wz\}$,
\end{center}
and choose $(S_6,\pi)$ to be the conic bundle induced by 
\begin{center}
$\pi((x:y:z) \times (u:v:w))=\left\{\begin{array}{ll}
(y:z)& \mbox{ if } (x:y:z) \not= (1:0:0),\\
(w:v)& \mbox{ if } (u:v:w) \not= (1:0:0).\end{array}\right.$
\end{center}
(See Section \refd{Subsec:DelPezzo6} for a more precise description of the surface and its three conic bundle structures.)

\bigskip

For any $\alpha,\beta \in \K^{*}$, we define $\kappa_{\alpha,\beta}$ to be the following automorphism of $S_6$:
\begin{center}
$\kappa_{\alpha,\beta}:(x:y:z) \times (u:v:w) \mapsto (u:\alpha w:\beta v)\times (x:\alpha^{-1} z:\beta^{-1} y)$.\end{center}

Note that $\kappa_{\alpha,\beta}$ twists the two singular fibres of $\pi$ (see Lemma \refd{Lem:KappaDP6} below); its action on the fibration is $(x_1:x_2) \mapsto (\alpha x_1:\beta x_2)$ and 
\begin{center}$\kappa_{\alpha,\beta}^2((x:y:z) \times (u:v:w))=(x:\alpha\beta^{-1} y:\alpha^{-1}\beta z)\times(u:\alpha^{-1}\beta v:\alpha\beta^{-1} w)$.\end{center} So $\kappa_{\alpha,\beta}$ is an involution if and only if its action on the fibration is trivial.

\break

\begin{Lem}\fauxtitre
\label{Lem:KappaDP6}
Let $g \in \Aut(S_6,\pi)$ be an automorphism of the conic bundle. The following conditions are equivalent:
\begin{itemize}
\item
The triple $(<g>,S_6,\pi)$ is minimal.
\item
$g$ twists the two singular fibres of $\pi$.
\item
The action of $g$ on the exceptional divisors of $S_6$ is $(E_1 \ D_{23})(E_2\ D_{12})(E_3\ D_{13})$.
\item
$g=\kappa_{\alpha,\beta}$ for some $\alpha,\beta \in \K^{*}$.
\end{itemize}
\end{Lem}
\begin{proof}\fauxtitred\upshape
By Lemma \refd{Lem:AutDP6CB} the action of $\Aut(S_6,\pi)$ on the exceptional curves is isomorphic to $(\Z{2})^2$ and hence the possible actions of $g\not=1$ are these:
\begin{enumerate}
\item
$(E_2 \ E_3)(D_{12}\ D_{13})$,
\item
$(E_1 \ D_{23})(E_2\ D_{13})(E_3\ D_{12})$,
\item
$(E_1 \ D_{23})(E_2\ D_{12})(E_3\ D_{13})$.
\end{enumerate}
In the first two cases, the triple $(<g>,S_6,\pi)$ is not minimal. Indeed, the blow-down of respectively $\{E_2,E_3\}$ and $\{E_2,D_{13}\}$ gives a $g$-equivariant birational morphism of conic bundles.

Hence, if $(<g>,S_6,\pi)$ is minimal, its action on the exceptional curves is the third one above, as stated in the proposition, and it then twists the two singular fibres of $\pi$. Conversely if $g$ twists the two singular fibres of $\pi$, the triple $(<g>,S_6,\pi)$ is minimal (by Lemma \refd{Lem:MinTripl}).

It remains to see that the last assertion is equivalent to the others. This follows from the semi-direct product structure of $\Aut(S_6,\pi)$ (see Lemma \refd{Lem:AutDP6CB}).\proofend
\end{proof}

\begin{remark}\fauxtitred
The pair $(\Aut(S_6,\pi),S_6)$ is not minimal (this was proved in Lemma 
\refd{Lem:DP6CBPasMin} and also in Proposition \refd{Prp:InterDPCB}). Thus $<\kappa_{\alpha,\beta}>$ is an example of a group whose action on the surface is not minimal, but whose action on a conic bundle is minimal.
\end{remark}

\subsection{Twisting de Jonqui\`eres involutions}
\index{Conic bundles!twisting involutions|idb}
Let us first recall what we mean by a \emph{hyperelliptic curve}:
\begin{Def} \fauxtitre
\index{Curves!hyperelliptic!definition|idb}We say that an irreducible curve is \defn{hyperelliptic} if it is isomorphic to some (ramified) double covering of $\mathbb{P}^1$.
\end{Def}
By the Riemann-Hurwitz formula, this double covering is ramified over an even number $2k$ of points, where $k\geq 1$ and the genus of the curve is $k-1$. In our definition, rational and elliptic curves are then hyperelliptic curves (these are respectively the cases $k=1,2$). \index{Curves!elliptic!viewed as hyperelliptic curves|idb}
Note that the isomorphism class of the hyperelliptic curve depends only on the linear equivalence of the set of ramification points, a fact which follows from the unicity of the double covering.\index{Curves!hyperelliptic!isomorphism class|idb} We describe the automorphisms of hyperelliptic curves in Section \refd{Sec:HyperellipticCurves}, Proposition \refd{Prp:Authyperellic}.

\begin{Lem}\titreProp{de Jonqui\`eres involutions twisting a conic bundle}\\
\label{Lem:DeJI}
Let $g \in \Aut(S,\pi)$ be an automorphism of finite order of the conic bundle such that 
$g$ twists at least one singular fibre of $\pi$.

Then, the following properties are equivalent:
\begin{itemize}
\item[\upshape 1.]
$g$ is an involution;
\item[\upshape 2.]
$\overline{\pi}(g)=1$, i.e.\ $g$ has a trivial action on the basis of the fibration;
\item[\upshape 3.]
the set of points of $S$ fixed by $g$ is a smooth hyperelliptic curve of genus $k-1$, a double covering of $\mathbb{P}^1$ by means of $\pi$, ramified over $2k$ points, plus perhaps a finite number of isolated points, which are the singular points of the singular fibres not twisted by $g$;\index{Curves!hyperelliptic!birational maps that fix a hyperelliptic curve}
\item[\upshape 4.]
$g$ fixes some curve of positive genus.
\end{itemize}
Furthermore, if the three conditions above are satisifed, the number of singular fibres of $\pi$ twisted by $g$ is $2k$.
\end{Lem}
\begin{proof}\fauxtitred\upshape
\begin{itemize}
\item[$1\Rightarrow 2$]
By contracting some exceptional curves, we may assume that the triple $(<g>,S,\pi)$ is minimal. Suppose that $g$ is an involution and $\overline{\pi}(g)\not=1$. Then  $g$ may twist only two singular fibres, which are the fibres of the two points of $\mathbb{P}^1$ fixed by $\overline{\pi}(g)$. Hence, the number of singular fibres is $\leq 2$. Lemma \refd{Lem:Degree567CB} tells us that $S$ is a Del Pezzo surface of degree $6$ and then Lemma \refd{Lem:KappaDP6} shows that
\begin{center}
$\begin{array}{lllll}
g=\kappa_{\alpha,\beta}:&\h{2}(x:y:z) \times (u:v:w) \mapsto (u:\alpha w:\beta v)&\h{2}\times &\h{2}(x:\alpha^{-1} z:\beta^{-1} y), \vspace{0.1 cm}\\
g^2=\kappa_{\alpha,\beta}^2:&\h{2}(x:y:z) \times (u:v:w) \mapsto (x:\alpha\beta^{-1} y:\alpha^{-1}\beta z)&\h{2}\times &\h{2}(u:\alpha^{-1}\beta v:\alpha\beta^{-1} w),\end{array}$\end{center} for some $\alpha,\beta \in \K^{*}$. As $\overline{\pi}(g)=(x_1:x_2) \mapsto (\alpha x_1:\beta x_2)$, $\beta \not=\alpha$, so $g$ is not an involution.
\item[$2\Rightarrow 3$]
Suppose first that $(<g>,S,\pi)$ is minimal.
As $\overline{\pi}(g)=1$, every fibre twisted by a power of $g$ is twisted by $g$. The minimality of the triple $(<g>,S,\pi)$ thus implies (Lemma \refd{Lem:MinTripl}) that $g$ twists every singular fibre of $\pi$. Therefore, on a singular fibre there is one point fixed by $g$ (the singular point of the fibre) and on a general fibre there are two fixed points. The set of points of $S$ fixed by $g$ is therefore a smooth irreducible curve. The projection $\pi$ gives it as a double covering of $\mathbb{P}^1$ ramified over the points whose fibres are singular and twisted by $g$. By the Hurwitz formula, this number is even.

The situation when $(<g>,S,\pi)$ is not minimal is obtained from this one, by blowing-up some fixed points. This adds in each new singular fibre (not twisted by the involution) an isolated point, which is the singular point of the singular fibre.

We then get the third assertion and the final remark. 
\item[$3\Rightarrow 2,1$]
We observe first of all that $\overline{\pi}(g)=1$. Then, $g^2$ leaves invariant every component of every singular fibre of $\pi$.
Blowing-down one irreducible component in each singular fibre we get a $g^2$-equivariant birational morphism $p:S\rightarrow \mathbb{F}_n$, for some integer $n\geq 0$.  If $g^2$ is not the identity, it is (by hypothesis) an automorphism of finite order acting on $\mathbb{F}_n$, trivially on the fibration. Its set of fixed points is then the union of two sections. But this is not possible, since the fixed points of $g^2$ must contain those of $g$.\proofend
\end{itemize}
\end{proof}

\bigskip

\begin{Def} \titreProp{twisting $(S,\pi)$-de Jonqui\`eres involutions}\\
\label{Def:DeJI}
Let $(S,\pi)$ be a conic bundle. We say that an automorphism $g \in \Aut(S,\pi)$ that satisfies the conditions of Lemma \refd{Lem:DeJI} is a \defn{de Jonqui\`eres involution twisting the conic bundle $(S,\pi)$}, or equivalently a \defn{twisting $(S,\pi)$-de Jonqui\`eres involution}.
\end{Def}

\begin{remark}\fauxtitred
Note that this sort of involution is a \defn{de Jonqui\`eres involution} in the classical sense (an involution which leaves invariant a pencil of rational curves, see Definition \refd{Def:DeJonqInv}). See Section \refd{Sec:DeJonqInv} for further relations between the two definitions.\end{remark}

\bigskip

\begin{Exa}\titreProp{twisting de Jonqui\`eres involutions}\upshape\\
\label{Exa:ThreeTwistDeJI}
Let $a_1,...,a_n,b_1,...,b_n,c_1,...,c_n \in \K$ be distinct and let $S \subset \mathbb{P}^1 \times \mathbb{P}^2$ be the smooth surface defined by the equation
\begin{center}$y_1^2\prod_{i=1}^n (x_1-a_i x_2)+y_2^2\prod_{i=1}^n (x_1-b_i x_2)+y_3^2\prod_{i=1}^n (x_1-c_i x_2)=0$.\end{center}

Let $\pi:S \rightarrow \mathbb{P}^1$ denote the projection $(x_1:x_2)\times (y_1:y_2:y_3)\mapsto (x_1:x_2)$ on the first factor. It induces a conic bundle structure $(S,\pi)$ with $3n$ singular fibres, which are the fibres of the points $\{(a_i:1)\}_{i=1,...,n}$, $\{(b_i:1)\}_{i=1,...,n}$ and $\{(c_i:1)\}_{i=1,...,n}$. \\
The group of automorphisms of the form
\begin{center} $(x_1:x_2) \times (y_1:y_2:y_3) \mapsto (x_1:x_2) \times (\pm y_1:\pm y_2:y_3)$\end{center} contains three distinct de Jonqui\`eres involutions twisting the conic bundle $(S,\pi)$:

\begin{itemize}
\item
The involution $(x_1:x_2) \times (y_1:y_2:y_3) \mapsto (x_1:x_2) \times (- y_1:y_2:y_3)$ twists the singular fibre of $\{(b_i:1)\}$ and $\{(c_i:1)\}$ for $i=1,...,n$. It leaves invariant the two irreducible components of the singular fibre of $\{(a_i:1)\}$ for $i=1,...,n$. Its set of fixed points is the hyperelliptic curve of $S$ of equation
$y_2^2\prod_{i=1}^n (x_1-b_i x_2)+y_3^2\prod_{i=1}^n (x_1-c_i x_2)=0 ,y_1=0$.\index{Curves!hyperelliptic!birational maps that fix a hyperelliptic curve}
\item
The situation is similar for the involutions $(x_1:x_2) \times (y_1:y_2:y_3) \mapsto (x_1:x_2) \times (y_1:-y_2:y_3)$ and $(x_1:x_2) \times (y_1:y_2:y_3) \mapsto (x_1:x_2) \times (y_1:y_2:-y_3)$.
\end{itemize}\end{Exa}

\bigskip

We give now some further simple results on twisting de Jonqui\`eres involutions.
\begin{Lem}\fauxtitre
\label{Lem:NoRoot}
Let $(S,\pi)$ be some conic bundle. No twisting $(S,\pi)$-de Jonqui\`eres involution has a root in $\Aut(S,\pi)$ which acts trivially on the fibration.
\end{Lem}
\begin{proof}\fauxtitred\upshape
Such a root must twist a singular fibre and so (Lemma \ref{Lem:DeJI}) is an involution.\proofend
\end{proof}
\begin{remark}Some roots in $\Aut(S,\pi)$ which do not act trivially on the fibration may exist (see Section \ref{Sec:RootDeJ}).\end{remark}
\begin{Lem}\fauxtitre
\label{Lem:CurveFixDeJ}\index{Curves!of positive genus!birational maps that fix a curve of positive genus}
Let $(S,\pi)$ be some conic bundle and let $g \in \Aut(S,\pi)$ an element of finite order. The following conditions are equivalent.
\begin{itemize}
\item
$g$ is a twisting $(S,\pi)$-de Jonqui\`eres involution that twists at least $4$ singular fibres.
\item
 $g$ fixes a curve of positive genus. 
\end{itemize}
\end{Lem}
\begin{proof}\fauxtitred\upshape
Lemma \ref{Lem:DeJI} shows that the first assertion implies the second.

Suppose that $g$ fixes a curve of positive genus. It is clear that $g$ acts trivially on the fibration, and that it fixes $2$ points on a general fibre.
The curve fixed by $g$ must then be a smooth hyperelliptic curve, and we get the first assertion using Lemma \ref{Lem:DeJI}.
\proofend
\end{proof}

\bigskip

\subsection{How the two examples determine almost all cases}
As we mentioned above, the automorphisms that twist some singular fibre are fundamental (see Section \refd{Sec:Minimality}, and in particular Lemma \refd{Lem:MinTripl}). Let us see that these are essentially the two examples given above:

\begin{Prp}\fauxtitre
\label{Prp:TwocasesKappaTwist}
Let $(S,\pi)$ be a conic bundle and $g \in \Aut(S,\pi)$ be an automorphism of finite order which twists at least one singular fibre.
\begin{itemize}
\item[\upshape 1.]
If $\overline{\pi}(g)=1$, i.e.\ if the action on the basis is trivial, then $g$ is a twisting $(S,\pi)$-de Jonqui\`eres involution.
\item[\upshape 2.]
If the cyclic group $G$ generated by $g$ does not contain any twisting $(S,\pi)$-de Jonqui\`eres involution, there exists some $g$-equivariant birational morphism of conic bundles  $\eta:S\rightarrow S_6$, where $S_6$ is the Del Pezzo surface of degree $6$\index{Del Pezzo surfaces!of degree $6$} and $\eta g\eta^{-1}=\kappa_{\alpha,\beta}$ for some $\alpha,\beta \in \K^{*}$.\\
Moreover, we have
\begin{itemize}
\item[\upshape (a)]
$\overline{\pi}(g)\not=1${\upshape;}
\item[\upshape (b)]
$g$ twists exactly two singular fibres;
\item[\upshape (c)]
the order of $g$ is $4n$, for some positive integer $n$;
\item[\upshape (d)]
the exact sequence induced by $\overline{\pi}$ is:
$1 \rightarrow G' \rightarrow G \stackrel{\overline{\pi}}{\rightarrow} \Z{2n} \rightarrow 1$,\\
where $G' \cong \Z{2}$.
\end{itemize}
\end{itemize}
\end{Prp}
\begin{proof}\fauxtitred\upshape
The first part of the lemma follows from Lemma \refd{Lem:DeJI}.

Let us prove the second part. The fact that $\overline{\pi}(g)\not=1$ follows from the first part. If necessary, we blow-down some components of singular fibres and suppose that the triple $(G,S,\pi)$ is minimal. Note that if some element $h\in G$ twists a singular fibre, it must act non-trivially on the fibration, since we assumed that no twisting $(S,\pi)$-de Jonqui\`eres involution belongs to $G$.  Then,  a fibre twisted by $h$ is also twisted by $g$. There are then at most $2$ singular fibres of $\pi$ and $g$ twists each one, by minimality. Lemma \refd{Lem:Degree567CB} tells us that $S$ is a Del Pezzo surface of degree $6$ and Lemma \refd{Lem:KappaDP6} shows that
\begin{center}
$\begin{array}{lllllll}
g=&\kappa_{\alpha,\beta}:&(x:y:z) \times (u:v:w)& \mapsto&\h{1} (u:\alpha w:\beta v)&\times &\h{1}(x:\alpha^{-1} z:\beta^{-1} y), \vspace{0.1cm}\\
g^2=&\kappa_{\alpha,\beta}^2:&(x:y:z) \times (u:v:w)& \mapsto&\h{1} (x:\alpha\beta^{-1} y:\alpha^{-1}\beta z)&\times &\h{1}(u:\alpha^{-1}\beta v:\alpha\beta^{-1} w),\end{array}$\end{center} for some $\alpha,\beta \in \K^{*}$.

This implies that $g$ is of even order. Denoting this order by $2m$, we see that $\alpha\beta^{-1}$ is a primitive $m$-th root of the unity. As the action of $g$ on the fibration is $(x_1:x_2) \mapsto (\alpha x_1:\beta x_2)$, it is cyclic of order $m$. The fact that $m$ is an even integer follows from the fact that $g^m \in G'$ does not twist any of the singular fibres, since $g^m$ is not a twisting $(S,\pi)$-de Jonqui\`eres involution. This gives the annouced result.
\proofend
\end{proof}
\section{The example Cs.24}
\label{SubSec:ExampleCs24} 
\pagetitless{Finite abelian groups of automorphisms of conic bundles}{The example Cs.24}
\index{Conic bundles!the example Cs.24|idb}\ChapterBegin{
We now give  the most important example of this paper. This is the only finite abelian subgroup of the Cremona group which is not conjugate to a group of automorphisms of\hspace{0.2 cm}$\Pn$ or $\mathbb{P}^1\times\mathbb{P}^1$ but whose elements do not fix any curve of positive genus\index{Curves!of positive genus!birational maps that do not fix a curve of positive genus!the example Cs.24|idb} (see  Section \refd{Sec:Results}, Theorem \refThmAbelGenusCurve).}

\bigskip

Let $S_6 \subset \Pn \times \Pn$ be the Del Pezzo surface of degree $6$\index{Del Pezzo surfaces!of degree $6$}\index{Del Pezzo surfaces!conic bundles on the surfaces|idi}\index{Conic bundles!on Del Pezzo surfaces|idi} defined by 
\begin{center}$S_6=\{ (x:y:z) \times (u:v:w)\ | \ ux=yv=zw\}$\end{center} and let us keep the notation of Sections \refd{Subsec:DelPezzo6} and  \refd{SubSec:DelPezzo6Kappa}. We denote by $\eta:S\rightarrow S_6$ the blow-up of $A_4,A_5 \in S_6$ defined by:
\begin{center}$\begin{array}{ccccc}A_4=&(0:1:1) &\times& (1:0:0)& \in D_{23},\\
A_5=&(1:0:0) &\times& (0:1:-1)& \in E_1.\end{array}$\end{center}

By convention we again denote by $E_1,E_2,E_3,D_{12},D_{13},D_{23}$ the total pull-backs by $\eta$ of these divisors of $S_6$. We denote by $\tilde{E_1}$ and $\tilde{D_{23}}$ the strict pull-backs of $E_1$ and $D_{23}$ by $\eta$. (Note that for the other exceptional divisors, the strict and total pull-backs are the same.) Let us illustrate the situation:
\begin{center}
\includegraphics[width=90.00mm]{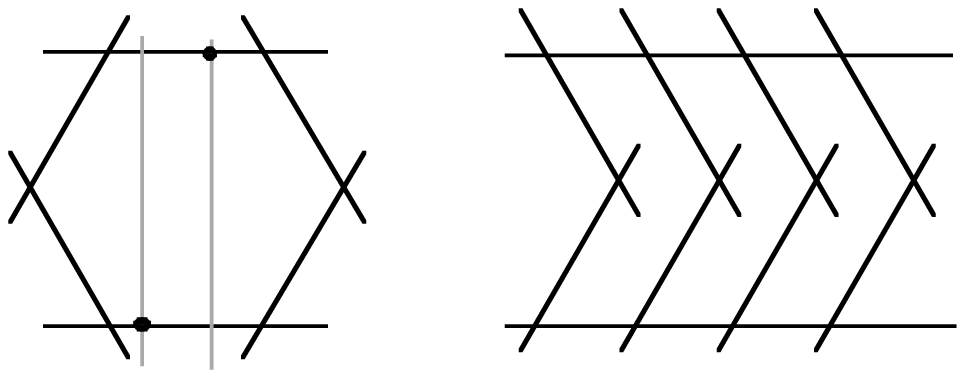}%
\drawat{-63.78mm}{30mm}{$D_{23}$}%
\drawat{-60.8mm}{23.25mm}{$E_3$}%
\drawat{-60.93mm}{8.12mm}{$D_{13}$}%
\drawat{-64.07mm}{1.0mm}{$E_1$}%
\drawat{-88mm}{8.12mm}{$D_{12}$}%
\drawat{-87.0mm}{23.25mm}{$E_2$}%
\drawat{-74.47mm}{6.16mm}{$A_4$}%
\drawat{-68.49mm}{25.03mm}{$A_5$}%
\drawat{-8.94mm}{30mm}{$\tilde{D_{23}}$}%
\drawat{-8mm}{0.5mm}{$\tilde{E_1}$}%
\drawat{-35.66mm}{22.6mm}{$E_2$}%
\drawat{-26.96mm}{22.6mm}{$D_{15}$}%
\drawat{-18.14mm}{22.6mm}{$E_4$}%
\drawat{-9.55mm}{22.6mm}{$E_3$}%
\drawat{-38mm}{6.9mm}{$D_{12}$}%
\drawat{-28.8mm}{6.9mm}{$E_5$}%
\drawat{-21mm}{6.9mm}{$D_{14}$}%
\drawat{-11.9mm}{6.9mm}{$D_{13}$}%
\end{center}

Let $\pi'$ denote the morphism $S_6 \rightarrow \mathbb{P}^1$ defined in Section \refd{SubSec:DelPezzo6Kappa}. The morphism $\pi=\pi' \circ \eta$ gives the surface $S$ a conic bundle structure $(S,\pi)$. It has $4$ singular fibres, which are the fibres of $(-1:1)$, $(0:1)$, $(1:1)$ and $(1:0)$.
We denote by $f$ the divisor of $S$ corresponding to the fibre of $\pi$ and we set $E_4=\eta^{-1}(A_4)$, $E_5=\eta^{-1}(A_5)$. Note that $E_4$ is one of the components of the singular fibre of $(1:1)$; we denote by $D_{14}=f-E_4$ the other component, which is the strict pull-back by $\eta$ of $\pi'^{-1}(1:1)$. Similarly, we denote by $D_{15}$ the divisor $f-E_5$, so that the singular fibre of $(-1:1)$ is $\{E_5,D_{15}\}$.

There are $8$ exceptional curves on $S$, which are $E_2,E_3,E_4,E_5,D_{12},D_{13},D_{14},D_{15}$, and two rational irreducible curves of self-intersection $-2$, which are $\tilde{E_1}=E_1-E_5$ and $\tilde{D_{23}}=D_{23}-E_4$. (It can in fact be proved that there are no other rational irreducible curves of negative self-intersection.)

Note that the group of automorphisms of $S$ that leave invariant every curve of negative self-intersection is isomorphic to $\K^{*}$ and corresponds to automorphisms of\hspace{0.2 cm}$\Pn$ of the form $(x:y:z)\mapsto (\alpha x:y:z)$, for $\alpha \in \K^{*}$. Indeed, such automorphisms are the lifts of automorphisms of $S_6$ leaving invariant every exceptional curve (which are of the form $(x:y:z) \times (u:v:w) \mapsto (x:\alpha y:\beta z) \times (u:\alpha^{-1} v:\beta^{-1} w)$, for $\alpha,\beta \in \K^{*}$) and which fix both points $A_4$ and $A_5$. 

\break

\begin{Lem}\fauxtitre
\label{Lem:Cs24Properties}
Let $G \subset \CrP$ be the following group:
\begin{center}\begin{tabular}{ccrl}
${\num{Cs.24}}$ &\h{1} $G \cong \Z{2}\times\Z{4}$& g.b. &$g_1':(x:y:z)\dasharrow (yz:xy:-xz)$ \\
&\h{1} & and& $g_2':(x:y:z)\dasharrow ( yz(y-z):xz(y+z):xy(y+z))$ \end{tabular}\end{center}
Then:
\begin{itemize}
\item
The lift of $G$ on $S$ is a subgroup of $\Aut(S,\pi)$. We denote  the lift of $g_1'$ (respectively $g_2'$) on $S$ by $g_1$ (respectively $g_2$).
\item
The action of $g_1$ and $g_2$ on the set of irreducible rational curves of negative self-intersection is respectively:
\begin{center}$(\tilde{E_1}\ \tilde{D_{23}})(E_2\ D_{12})(E_3\ D_{13})(E_4\ E_5)(D_{14}\ D_{15})$, \vspace{0.1 cm}\\
$(\tilde{E_1}\ \tilde{D_{23}})(E_2\ D_{13})(E_3\ D_{12})(E_4\ D_{14})(E_5\ D_{15})$.\end{center}
\item
Both $g_1$ and $g_2$ are elements of order $4$ and $(g_1')^2=(g_2')^2=(x:y:z)\mapsto (-x:y:z)$. Thus $(g_1)^2=(g_2)^2 \in \ker \overline{\pi}$ is an automorphism of $S$ which leaves invariant every divisor of negative self-intersection.
\item
$G$ contains no twisting $(S,\pi)$-de Jonqui\`eres involution. In particular, no element of $G$ fixes a curve of positive genus.
\item
The pair $(G,S)$ and the triple $(G,S,\pi)$ are both minimal.
\end{itemize}\end{Lem}
\begin{proof}\fauxtitred\upshape
Observe first that $g_1'$ and $g_2'$ preserve the pencil of lines of\hspace{0.2 cm}$\Pn$ passing through the point $A_1=(1:0:0)$ (of the form $ay+bz=0$, $(a:b) \in \mathbb{P}^1$), so $g_1,g_2 \in \Crc{S}{\pi}$. Then, we compute $(g_1')^2=(g_2')^2=(x:y:z)\mapsto (-x:y:z)$. This implies that both $g_1'$ and $g_2'$ are birational maps of order $4$.
\begin{itemize}
\item
Note that the lift of $g_1'$ on the surface $S_6$ is the automorphism 
\begin{center}
$\kappa_{1,-1}:(x:y:z) \times (u:v:w) \mapsto (u: w: -v)\times (x: z:-y)$\end{center}
(see Section \refd{SubSec:DelPezzo6Kappa}). As this automorphism permutes $A_4$ and $A_5$, its lift on $S$ is biregular. The action on the divisors with negative self-intersection is deduced from that of $\kappa_{1,-1}$ (see Lemma \refd{Lem:KappaDP6}).
\item
Let us compute the involution $g_3'=g_1'g_2'=(x:y:z) \dasharrow (x(y+z):z(y-z):-y(y-z))$. The associated linear system is $\{ax(y+z)+(by+cz)(y-z)=0\ | \ (a:b:c) \in \Pn\}$, which is the linear sytem of conics passing through $(0:1:1)$ and $A_1=(1:0:0)$, with tangent $y+z=0$ at this point. Blowing-up these three points (two on $\Pn$ and one in the blow-up of $A_1$), we get an automorphism of some rational surface $S'$. As the points $A_2=(0:1:0)$ and $A_3=(0:0:1)$ are permuted by $g_3'$, we can also blow them up and again get an automorphism.

The isomorphism class of the surface obtained is independent of the order of the points blown-up. We first blow-up $A_1,A_2,A_3$ and get $S_6$. Then, we blow-up the two other base points of $g_3'$, which are in fact $A_4$ (the point $(0:1:-1)$) and $A_5$ (the point infinitely near $A_1$ corresponding to the tangent $y+z=0$). This shows that $g_3$, and therefore $g_2$, belong to $\Aut(S,\pi)$.

Since $g_3'$ permutes the points $A_2$ and $A_3$, $g_3=g_1g_2$ permutes the divisors $E_2$ and $E_3$. It also permutes $D_{12}$ and $D_{13}$, since $g_3'$  leaves invariant the lines through $A_1$. It therefore leaves $\tilde{E_1}$ and $\tilde{D_{23}}$ invariant, since $E_2$ and $E_3$ touch $\tilde{D_{23}}$ but not $E_1$. The remaining exceptional divisors are $E_4,E_5,D_{14},D_{15}$. Either $g_1g_2$ leaves all four invariant, or it acts as $(E_4\ D_{15})(E_5\ D_{14})$ (using the intersection with $\tilde{E_1}$ and $\tilde{D_{23}}$). As $A_4$ and $A_5$ are base point of $g_1'g_2'$, $E_4$ and $E_5$ are not invariant. Thus, $g_1g_2$ acts on the irreducible rational curves of negative self-intersection as:
\begin{center}
$(E_2\ E_3)(D_{12}\ D_{13})(E_4\ D_{15})(E_5\ D_{14}).$\end{center}
We get the action of $g_2$ by composing this one with that of $g_1$.
\item
We see then that $G$ contains no twisting $(S,\pi)$-de Jonqui\`eres involution, since $\ker \overline{\pi}$ is generated by $g_1^2=g_2^2$ (which is a de Jonqui\`eres involution that does not twist any singular fibre).
\item
Note that the orbits of the action of $G$ on the exceptional divisors of $S$ are \begin{center}$\{E_2,E_3,D_{12},D_{13}\}$ and $\{E_4,E_5,D_{14},D_{15}\}$.\end{center} As these orbits cannot be contracted, the pair $(G,S)$ is minimal, and so is the triple $(G,S,\pi)$.\proofend\end{itemize}
\end{proof}

\begin{remark}\fauxtitred
From now on, $\nump{Cs.24}$ will denote the triple $(G,S,\pi)$ {\upshape (}or only the pair $(G,S)${\upshape)} of  Lemma \refd{Lem:Cs24Properties}, or possibly the corresponding conjugacy class in the Cremona group.
\end{remark}

\bigskip

Let us now prove that this group is not birationally conjugate to a group of automorphisms of a surface of larger degree:
\begin{Lem}\fauxtitre
\label{Lem:Cs24notconjugate}
The pair $\nump{Cs.24}$ is not birationally  conjugate to a group of automorphisms of\hspace{0.2 cm}$\Pn$ or $\mathbb{P}^1\times\mathbb{P}^1$.
\end{Lem}
\begin{proof}\fauxtitred\upshape
Up to birational conjugation, any subgroup of $\Aut(\Pn)$ and $\Aut(\mathbb{P}^1\times\mathbb{P}^1)$ isomorphic to $\Z{2}\times\Z{4}$ is conjugate to one of the following two groups of automorphisms of $\mathbb{P}^1\times\mathbb{P}^1$ (see Propositions \refd{Prp:PGL3Cr} and  \refd{Prp:AutP1Bir}):
\begin{center}\begin{tabular}{lll}
${\nump{0.mn}}$, for $m=4,n=2:$ & $G_1$ generated by &$(x,y) \mapsto (- x,y)$ and
$(x,y) \mapsto (x,\im y)$,\vspace{0.1 cm}\\
${\nump{P1s.24}}$ &  $G_2$  generated by & $(x,y) \mapsto (x^{-1},y^{-1})$ and $(x,y) \mapsto (-y,x)$.\end{tabular}\end{center}
We denote the $2$-torsion of $G_1$  by $H_1$ and that of $G_2$ by $H_2$. The following simple observation will be useful:
\begin{itemize}
\item
Both $H_1$ and $G_1$ fix $4$ points of $\mathbb{P}^1\times\mathbb{P}^1$.
\item
Neither $H_2$ nor $G_2$ fixes any point of $\mathbb{P}^1\times\mathbb{P}^1$.
\end{itemize}

Let $(G,S)$ be the pair $\nump{Cs.24}$ defined above. We use the notation $g_1,g_2$ as above. We denote by $H$ the $2$-torsion of $G$, generated by $(g_1)^2=(g_2)^2$ and $g_3=g_1g_2$.

We will see that $H$ fixes $4$ points of $S$ and that $G$ does not fix any point of $S$. This implies that $G$ is not birationally conjugate to $G_1$ or $G_2$, since the existence of fixed points is an invariant (see Section \refd{Sec:ExistenceFixedPoints}).

Note that $g_1^2$ and $g_3$ are respectively the lifts of the following birational maps of\hspace{0.2 cm}$\Pn$:\begin{center}
$\begin{array}{rl}
(g_1')^2=(g_2')^2=&(x:y:z)\mapsto (-x:y:z),\\
g_3=g_1g_2=&(x:y:z) \dasharrow(x(y+z):z(y-z):-y(y-z)).\end{array}$\end{center}
Note that $g_3$ acts cyclically on the fibration  with order $2$. Its set of fixed points is the disjoint union of the fibres of $(1:\pm \im)$, which we see on $\Pn$ as the lines $y^2+z^2=0$. The element $(g_1)^2$ leaves each fibre invariant and fixes exactly two points on each fibre, which are the intersection of the fibre with $\tilde{E_1}$ and $\tilde{D_{23}}$. We see two of the four fixed points on $\Pn$, namely the points $(0:1:\pm \im)$.

Note that $g_1$ permutes the fibres of $(1:\im)$ and $(1:-\im)$, so $G$ fixes no point of $S$.\proofend\end{proof}

\bigskip

\section{Groups without a twisting de Jonqui\`eres involution}
\pagetitless{Finite abelian groups of automorphisms of conic bundles}{Groups without a twisting de Jonqui\`eres involution}
\index{Conic bundles!twisting involutions|idi}
In Definition \refd{Def:DeJI}, we defined what we mean by de Jonqui\`eres involutions twisting some conic bundle. In this section, we consider groups that do not contain any of these involutions.
\begin{Prp}\fauxtitre
\label{Prp:necDeJi}
Let $(S,\pi)$ be a conic bundle and let $G \subset \Aut(S,\pi)$ be a finite abelian group such that:
\begin{itemize}
\item
no twisting $(S,\pi)$-de Jonqui\`eres involution belongs to $G$;
\item
the triple $(G,S,\pi)$ is minimal.
\end{itemize}
Then, one of the following occurs:
\begin{itemize}
\item
The fibration is smooth, i.e.\ $S$ is a minimal ruled surface $\mathbb{F}_n$ for some integer $n\geq 0$.
\item
$S$ is the Del Pezzo surface of degree $6$.\index{Del Pezzo surfaces!of degree $6$}
\item 
The triple $(G,S,\pi)$ is the triple $\nump{Cs.24}$ of Section \refd{SubSec:ExampleCs24}.
\end{itemize}
\end{Prp}
\begin{remark}\fauxtitred
This proposition is very useful, as we have then only to consider groups of automorphisms of conic bundles which contain a twisting de Jonqui\`eres involution, or the three cases stated above.
Note that the cases are simplified, using birational conjugation, in Proposition \refd{Prp:CBdistCases}.
\end{remark}
\begin{proof}\fauxtitred\upshape
We assume that the fibration is not smooth, and show that either $S$ is the Del Pezzo surface of degree $6$, or the pair $(G,S)$ is isomorphic to the pair $\nump{Cs.24}$ of Section \refd{SubSec:ExampleCs24}. 

As the triple $(G,S,\pi)$ is minimal, there exists for every singular fibre $F=\{F_1,F_2\}$ some element $g \in G$ such that $g(F_1)=F_2$ and $g(F_2)=F_1$ (Lemma \refd{Lem:MinTripl}). 
This means that $\pi(F)$ is fixed by $\overline{\pi}(g)$. Furthermore, any automorphism $g \in G$ twists either no singular fibre or exactly two (see Proposition \refd{Prp:TwocasesKappaTwist}, assertion $2(b)$) which correspond to the fibres of the two fixed points of $\overline{\pi}(g)\in \PGL(2,\K)$. 

If the number of singular fibres is exactly two, $S$ has degree $6$, and by Proposition 
\refd{Prp:TwocasesKappaTwist}, $S$ is a Del Pezzo surface of degree $6$, and we are done.

Suppose now that the number of singular fibres is larger than two. This implies that $\overline{\pi}(G)$ is not a cyclic group (otherwise the non-trivial elements of $\overline{\pi}(G)$ have the same two fixed points: there are then at most two singular fibres). Therefore, $\overline{\pi}(G)$ is isomorphic to $(\Z{2})^2$ (see Lemma \ref{Lem:AutP1} in Appendix). By a judicious choice of coordinates we may suppose that \begin{center}$\overline{\pi}(G)=\left\{\left(\begin{array}{cc}1 & 0\\ 0 & 1\end{array}\right),\left(\begin{array}{cc}-1 & 0\\ 0 & 1\end{array}\right),\left(\begin{array}{cc}0 & 1\\ 1 & 0\end{array}\right),\left(\begin{array}{cc}0 & -1\\ 1 & 0\end{array}\right)\right\}.$\end{center}
As a singular fibre corresponds to a fixed point of one of the three elements of order $2$ of $\overline{\pi}(G)$, only the fibres of $(0:1),(1:0),(1:1),(-1:1),(\im:1),(-\im:1)$ can be singular. Since this group acts transitively on the sets $\{(1:0),(0:1)\}$, $\{(1:\pm 1)\}$ and $\{(1:\pm \im)\}$, there are $4$ or $6$ singular fibres.

We denote by $g_1$ an element of $G$ which twists the two singular fibres of $(1:0)$ and $(0:1)$. Note that $g_1$ is of order $4$  and $(g_1)^2\in G'$  (see Proposition \refd{Prp:TwocasesKappaTwist}, assertion 2). Since by hypothesis $(g_1)^2$ is not a twisting $(S,\pi)$-de Jonqui\`eres involution, it leaves invariant each component of every singular fibre of $\pi$.
Let $\eta:S\rightarrow S_6$ denote the birational $g_1$-equivariant morphism given by Proposition \refd{Prp:TwocasesKappaTwist}, which conjugates $g_1$ to the automorphism
\begin{center}$\eta g_1 \eta^{-1}=\kappa_{\alpha,\beta}:(x:y:z) \times (u:v:w) \mapsto (u:\alpha w:\beta v)\times (x:\alpha^{-1} z:\beta^{-1} y)$\end{center}
of the Del Pezzo surface $S_6$ of degree $6$, for some $\alpha,\beta \in \K^{*}$. In fact, as $\overline{\pi}(g_1)$ has order $2$, we have $\beta=-\alpha$, so $\eta g_1 \eta^{-1}=\kappa_{\alpha,-\alpha}$.
The points blown-up by $\eta$ are fixed by 
\begin{center}$\eta (g_1)^2 \eta^{-1}=(\kappa_{\alpha,-\alpha})^2:(x:y:z) \times (u:v:w)\mapsto (x:-y:- z)\times(u:-v:-w)$, \end{center}
and therefore belong to the curves
\begin{center} $\begin{array}{rlll}
& E_1&=& \{(1:0:0) \times (0:a:b) \ | \ (a:b) \in \mathbb{P}^1\} \\
\mbox{and} &D_{23}&=&\{(0:a:b) \times (1:0:0) \ | \ (a:b) \in \mathbb{P}^1\}.\end{array}$\end{center} As these points consist of orbits of $\eta g_1\eta^{-1}$, half of them lie in $E_1$ and the other half in $D_{23}$. In fact, up to a change of coordinates, $(x,y,z)\times (u,v,w) \leftrightarrow (u,v,w)\times (x,y,z)$, the points that may be blown-up by $\eta$ are
\begin{center}
$\begin{array}{rlllll}
A_4=&(0:1:1) &\times& (1:0:0)& \in D_{23},\\
\kappa_{\alpha,-\alpha}(A_4)=A_5=&(1:0:0) &\times& (0:1:-1)& \in E_1,\\
A_6=&(0:1:\im) &\times&( 1:0:0)& \in D_{23},\\
\kappa_{\alpha,-\alpha}(A_6)=A_7=&(1:0:0) &\times& (0:1:\im)& \in E_1.\end{array}$
\end{center}
The strict pull-backs $\tilde{E_1}$ and $\tilde{D_{23}}$ by $\eta$ of $E_1$ and $D_{23}$ respectively then have self-intersection $-2$ or $-3$ in $S$, depending on the number of points blown-up.
By convention we again denote by $E_1,E_2,E_3,D_{12},D_{13},D_{23}$ the total pull-backs by $\eta$ of these divisors. (Note that for $E_2,E_3,D_{12},D_{13}$, the strict and the total pull-backs are the same.)
We set $E_4=p^{-1}(A_4)$,..., $E_7=p^{-1}(A_7)$ and denote by $f$ the divisor class of the fibre of the conic bundle. 

\begin{enumerate}
\item
Suppose that $\eta$ is the blow-up of $A_4$ and $A_5$. 
The Picard group of $S$ is then generated by $E_1,E_2,...,E_5$ and $f$.

Since we assumed that $(G,S,\pi)$ is minimal, the singular fibres of $(1:1)$ and $(-1:1)$ must be twisted. As we assumed that no de Jonqui\`eres involution twisting the conic bundle belongs to the group, one element $g_2$ twists these two singular fibres and acts cyclically with order $2$ on the basis of the fibration, with action $(x_1:x_2) \mapsto (x_2:x_1)$. As $g_1$ and $g_2$ twist some singular fibre, both must invert the two curves of self-intersection $-2$, namely $\tilde{E_1}$ and $\tilde{D_{23}}$. The action of $g_1$ and $g_2$ on the irreducibles rational curves of negative self-intersection are then respectively
\begin{center}$\begin{array}{l}(\tilde{E_1}\ \tilde{D_{23}})(E_2\ D_{12})(E_3\ D_{13})(E_4\ E_5)(D_{14}\ D_{15}), \vspace{0.1 cm}\\
(\tilde{E_1}\ \tilde{D_{23}})(E_2\ D_{13})(E_3\ D_{12})(E_4\ D_{14})(E_5\ D_{15}).\end{array}$\end{center}

The surface is that of Section \refd{SubSec:ExampleCs24} and the elements $g_1$ and $g_2$ have the same action on $\Pic{S}$ as the two automorphisms given the same name in that example. Note that the group $H$ of automorphisms of $S$ that leave invariant every curve of negative self-intersection is isomorphic to $\K^{*}$ and corresponds to automorphisms of\hspace{0.2 cm}$\Pn$ of the form $(x:y:z)\mapsto (\alpha x:y:z)$, for any $\alpha \in \K^{*}$. Then, $g_1$ and $g_2$ are equal to the following:
\begin{center}$\begin{array}{l}g_1:(x:y:z)\dasharrow (\mu yz:xy:-xz),\\
g_2:(x:y:z)\dasharrow (\nu yz(y-z):xz(y+z):xy(y+z)),\end{array}$\end{center} 
for some $\mu,\nu \in \K^{*}$

As $g_1g_2(x:y:z)=(\mu x(y+z):\nu z(y-z):-\nu y(y-z))$ and $g_2g_1(x:y:z)=(\nu x(y+z):\mu z(y-z):-\mu y(y-z))$ must be the same by hypothesis, we get $\mu^2=\nu^2$.

We observe that $\overline{\pi}(g_1)$ and $\overline{\pi}(g_2)$ generate $\overline{\pi}(G)\cong (\Z{2})^2$; on the other hand, an element of $G'$ does not twist a singular fibre by hypothesis and hence belongs to $H$. As the only elements of $H$ which commute with $g_1$ are $id$ and $(g_1)^2:(x:y:z)\mapsto (-x:y:z)$, we see that $g_1$ and $g_2$ generate the whole group $G$.

Conjugating $g_1$ and $g_2$ by $(x:y:z) \mapsto (\alpha x: y:z)$, where $\alpha \in \K^{*}, \alpha^2=\mu$, we may suppose that $\mu=1$. So $\nu=\pm 1$ and we get in both cases the same group, because $(g_1)^2(x:y:z)=(-x:y:z)$. The pair $(G,S)$ is then isomorphic to that of $\nump{Cs.24}$ (see Section \refd{SubSec:ExampleCs24}).
\item
Suppose that $\eta$ is the blow-up of $A_6$ and $A_7$. 
We get a case isomorphic to the previous one, using the automorphism $(x:y:z) \times (u:v:w) \mapsto (x:y:\im z)\times (u:v:-\im w)$ of $S_6$.
\item
Suppose that $\eta$ is the blow-up of $A_4,A_5,A_6$ and $A_7$. 
The Picard group of $S$ is then generated by $E_1,E_2,...,E_5,E_6,E_7$ and $f$.

Since we assumed that the group is minimal, there must be two elements $g_2,g_3 \in G$ that twist respectively the fibres of $(\pm 1:1)$ and those of $(\pm \im:1)$.
As in the previous example, the three actions of these elements on the basis are of order $2$, and the three elements transpose $\tilde{E_1}$ and $\tilde{D_{23}}$. The actions of $g_1,g_2$ and $g_3$ on the divisors of negative self-intersection are then respectively
\begin{center}
$\begin{array}{l}(\tilde{E_1}\ \tilde{D_{23}})(E_2\ D_{12})(E_3\ D_{13})(E_4\ E_5)(D_{14}\ D_{15})(E_6\ E_7)(D_{16}\ D_{17}), \vspace{0.1 cm}\\
(\tilde{E_1}\ \tilde{D_{23}})(E_2\ D_{13})(E_3\ D_{12})(E_4\ D_{14})(E_5\ D_{15})(E_6\ E_7)(D_{16}\ D_{17}), \vspace{0.1 cm}\\
(\tilde{E_1}\ \tilde{D_{23}})(E_2\ D_{13})(E_3\ D_{12})(E_4\ E_5)(D_{14}\ D_{15})(E_6\ D_{16})(E_7\ D_{17}).\end{array}$
\end{center}
We see that $g_1$ and $g_2$ correspond to elements given in the previous case. The element $g_3$ corresponds to the birational map 
\begin{center}$g_3:(x:y:z)\dasharrow (\zeta yz(y+\im z):xz(y-\im z),-xy(y-\im z))$,\end{center} for some $\zeta \in \K^{*}$ (the conjugation of $g_2$ by the automorphism $(x:y:z) \times (u:v:w) \mapsto (x:y:\im z)\times (u:v:-\im w)$ of $S_6$). Let us calculate some compositions:
\begin{center}
$\begin{array}{rrrrrrr}
g_1g_3(x:y:z)=(\hspace{-0.2 cm}&\h{1} -\mu x(y-\im z) &:& \zeta z(y+\im z) &:& \zeta y(y+\im z)&\hspace{-0.3 cm}),\\
g_3g_1(x:y:z)=(\hspace{-0.2 cm}&\h{1} \zeta x(y-\im z) &:& \mu z(y+\im z) &:& \mu y(y+\im z)&\hspace{-0.3 cm}),\\
g_2g_3(x:y:z)=(\hspace{-0.2 cm}&\h{1} \nu x(y-\im z)(y+z) &:& -\zeta y(y+\im z)(y-z) &:& \zeta z(y+\im z)(y-z)&\hspace{-0.3 cm}),\\
g_3g_2(x:y:z)=(\hspace{-0.2 cm}&\h{1}\zeta x(y-\im z)(y+z) &:& -\nu y(y+\im z)(y-z) &:& \nu z(y+\im z)(y-z)&\hspace{-0.3 cm}).\\
\end{array}$\end{center}
Thus, $g_1$ and $g_3$ commute if and only if $\zeta^2=-\mu^2$; and $g_2$ and $g_3$ commute if and only if $\zeta^2=\nu^2$. In the preceding case, we saw that
$g_1$ and $g_2$ commute if and only if $\mu^2=\nu^2$, so the three elements $g_1$, $g_2$ and $g_3$ cannot commute pairwise.\proofend
\end{enumerate}
\end{proof}

\bigskip

\section{The distinct cases}
\pagetitless{Finite abelian groups of automorphisms of conic bundles}{The distinct cases}
The following proposition summarises the results obtained so far:
\begin{Prp}\fauxtitre
\label{Prp:CBdistCases}\index{Conic bundles!automorphisms}
Let $(S,\pi)$ be a conic bundle and let $G \subset \Aut(S,\pi)$ be a finite abelian group such that the triple $(G,S,\pi)$ is minimal. 

Then, up to birational conjugation, one of the following occurs:
\begin{itemize}
\item
Some twisting $(S,\pi)$-de Jonqui\`eres involution belongs to $G$;
\item
$(S,\pi)=(\mathbb{P}^1\times \mathbb{P}^1,\pi_1)$;
\item
the group $G$ is isomorphic to $\Z{2}\times\Z{4}$ and is either $\nump{Cs.24}$ or $\nump{P1s.24}$.
\end{itemize}
Moreover, in the first case, $G'$ contains the twisting $(S,\pi)$-de Jonqui\`eres involution and is isomorphic either to $\Z{2}$, or to $(\Z{2})^2$.
\end{Prp}
\begin{remark}\fauxtitred
The situations  in which $G'\cong \Z{2}$ are described in Sections \refd{Sec:DeJonqInv}, \refd{Sec:ExtCyclDeJ} and \refd{Sec:ExtC2-C2C2}. The case when $G'\cong (\Z{2})^2$ arises in Section \refd{Sec:GroupsC2C2inPGL}.
\end{remark}
\begin{proof}\fauxtitred\upshape
By Proposition \refd{Prp:necDeJi}, if no twisting $(S,\pi)$-de Jonqui\`eres involution belongs to $G$, one of the following occurs:
\begin{itemize}
\item
The fibration is smooth, i.e.\ $S$ is a minimal ruled surface $\mathbb{F}_n$ for some integer $n\geq 0$.
\item
The surface $S$ is the Del Pezzo surface of degree $6$.\index{Del Pezzo surfaces!of degree $6$}
\item
The pair $(G,S)$ is isomorphic to the pair $\nump{Cs.24}$ of Section \refd{SubSec:ExampleCs24}.
\end{itemize}
Using Proposition \refd{Prp:FnAb}, the first case can be reduced to $n=0$, so $(S,\pi)=(\mathbb{P}^1\times \mathbb{P}^1,\pi_1)$.

In the second case, the group $G$ is birationally conjugate to a subgroup of $\Aut(\mathbb{P}^1 \times \mathbb{P}^1)$ (this was proved in Lemma 
\refd{Lem:DP6CBPasMin} and also in Proposition \refd{Prp:InterDPCB}). We can use Proposition 
\refd{Prp:AutP1Bir} to get either a subgroup of $\Aut(\mathbb{P}^1 \times \mathbb{P}^1,\pi_1)$ or one of the two groups  $\nump{P1s.222}$ and $\nump{P1s.24}$.

Note that the group ${\nump{P1s.222}}$ is birationally conjugate to the group generated by \begin{center}$(x,y) \dasharrow (x^{-1},y^{-1})$, $(x,y) \mapsto (x,-y)$, and $(x,y) \mapsto (x,x/y)$\end{center} (see Proposition 
\refd{Prp:AutP1Bir}) and this group is then conjugate to a group of automorphisms of some conic bundle $(S,\pi)$ (using Proposition \refd{Prp:BeauTsen}). This group contains some twisting $(S,\pi)$-de Jonqui\`eres involution, which corresponds to the involution $(x,y) \mapsto (x,x/y)$ (this will be proved in Proposition \refd{Prp:TwistingOrNot}).

Let us prove the last assertion. Recall that a twisting $(S,\pi)$-de Jonqui\`eres involution $\sigma$ leaves invariant any fibre of the conic bundle and so belongs to $G'=G \cap \ker \overline{\pi}$. Note that no root of this element can act trivially on the fibration (see Lemma 
\refd{Lem:NoRoot}), so $G'$ is isomorphic either to $\Z{2}$ or to $(\Z{2})^2$. 
\proofend
\end{proof}

As the subgroups of $\Aut(\mathbb{P}^1 \times \mathbb{P}^1,\pi_1)$ have been classified up to birational equivalence (see Proposition \refd{Prp:AutP1Bir}), we need only consider groups of automorphisms of some conic bundles $(S,\pi)$ containing some twisting $(S,\pi)$-de Jonqui\`eres involution. The exact sequence (introduced in Section \refd{Sec:Exactsequence})
\begin{equation}
1 \rightarrow G' \rightarrow G \stackrel{\overline{\pi}}{\rightarrow} \overline{\pi}(G) \rightarrow 1
\tag{\refd{eq:ExactSeqCB}}\end{equation}
and the restriction imposed on $G'$ by the proposition reduce the possibilities for $G$.

\section{The de Jonqui\`eres involutions}
\index{Involutions!De Jonqui\`eres|idb}
\pagetitless{Finite abelian groups of automorphisms of conic bundles}{The de Jonqui\`eres involutions}
\label{Sec:DeJonqInv}
There are three kinds of elements of order $2$ of $\CrP$, classically called \defn{Geiser}, \defn{Bertini} and \defn{de Jonqui\`eres} involutions (see for example \cite{bib:SeR}, or more recently \cite{bib:BaB}).
We describe Geiser and Bertini involutions, acting minimally respectively on Del Pezzo surfaces of degree $2$ and $1$, in Sections \refd{Sec:DP2} and \refd{Sec:DP1}.

\begin{Def} \fauxtitre
\label{Def:DeJonqInv}
An element of order $2$ of the Cremona group is a \defn{de Jonqui\`eres involution} if it leaves invariant a pencil of rational curves.
\end{Def}
Equivalently (see Section \refd{Sec:CBTwoRepr}), a de Jonqui\`ere involution is an involution birationally conjugate to an element of the group $\CrPp$.
To describe these elements, we use the algebraic structure of this group, which is $\PGL(2,\K(x))\rtimes \PGL(2,\K)$ (see Example \refd{Exa:CrPp}). We begin by noting the following result:
\begin{Prp}\fauxtitre
\label{Prp:DeJconjPGL2C}
Every de Jonqui\`eres involution is birationally conjugate to an element of \begin{center}$\PGL(2,\K(x)) \subset \PGL(2,\K(x))\rtimes \PGL(2,\K)=\CrPp$.\end{center}
\end{Prp}
\begin{proof}\fauxtitred\upshape
Let $\sigma$ be a de Jonqui\`eres involution. By Proposition \refd{Prp:BeauTsen}, up to birational conjugation, $\sigma$ acts biregularly on some conic bundle $(S,\pi)$.\ Applying Proposition \refd{Prp:CBdistCases}, birationally we get two cases: either $\sigma$ is a twisting $(S,\pi)$-de Jonqui\`eres involution or $(S,\pi)=(\mathbb{P}^1\times \mathbb{P}^1,\pi_1)$.
\begin{itemize}
\item
In the first case, $\sigma$ acts trivially on the fibration, and is then conjugate to an involutive element of $\CrPp$ which acts trivially on the fibration of $\pi_1$ (by Proposition \refd{Prp:BeauTsen}). This element belongs to $\PGL(2,\K(x))$.
\item
In the second case, $\sigma$ is birationally conjugate to the element $(x,y) \mapsto (x,-y)$ (see Proposition \refd{Prp:AutP1Bir}), which belongs to $\PGL(2,\K(x))$.\proofend
\end{itemize}
\end{proof}

\bigskip

To study the de Jonqui\`eres involutions, we classify them up to conjugations, first of all in the group $\PGL(2,\K(x))$, then in the de Jonqui\`eres group and then in the Cremona group. The three situations yield different results.

We proceed to describe the conjugacy classes of elements of order $2$ of $\PGL(2,\K(x))$:
\begin{Prp}\titreProp{classification of involutions in $\PGL(2,\K(x))$}
\index{Involutions!De Jonqui\`eres!in $\PGL(2,\K(x))$|idb}
\label{Prp:InvDeJPgl2} 
\begin{itemize}
\item[\upshape 1.]
Any involution belonging to $\PGL(2,\K(x))$ is conjugate in this group to an element of the form
\begin{center}$\sigma_{g(x)}=\left(\begin{array}{cc}0 & g(x) \\ 1 & 0\end{array}\right),$\end{center} for some $g(x)\in\K(x)\backslash\{0\}$.
\item[\upshape 2.]
The involutions $\sigma_{g(x)}$ and $\sigma_{\lambda(x)^2g(x)}$ are conjugate by the element \begin{center}$\left(\begin{array}{cc}\lambda(x) & 0 \\ 0& 1\end{array}\right),$\end{center} for any $\lambda(x)\in\K(x)\backslash\{0\}$.
\item[\upshape 3.]
The morphism $\delta: \PGL(2,\K(x)) \rightarrow \K(x)^{*}/\K(x)^{*2}$ given by the determinant induces a bijection between the set of conjugacy classes of involutions of $\PGL(2,\K(x))$ and $\K(x)^{*}/\K(x)^{*2}$.
\item[\upshape 4.]
An involution which does not belong to the kernel of $\delta$ has no roots in the group $\PGL(2,\K(x))$.
\end{itemize}
\end{Prp}
\begin{proof}\fauxtitred\upshape
Let $\sigma \in \PGL(2,\K(x))$ be an involution. Choosing some general vector $v \in \K(x)\times \K(x)$, in the basis $(v, \sigma(v))$, the involution is of the form $\sigma_{g(x)}$, for some $g(x)\in\K(x)\backslash\{0\}$, so we have assertion $1$.

Assertion $2$ follows directly by a simple calculation. 

Note that the morphism $\delta$ is invariant by conjugation.\ Assertion $3$ thus follows from assertion $2$.

As any square of $\PGL(2,\K(x))$ belongs to the kernel of $\delta$, an involution which does not belong to the kernel of $\delta$ cannot have a square root in the group $\PGL(2,\K(x))$. 
Let us now prove that this sort of involution cannot have a root of any order. We write the root of the involution $\sigma \in \PGL(2,\K(x))$ as the class of a matrix $A \in GL(2,\K(x))$ and have $A^{2n}=\lambda I$, for some $\lambda \in \K(x)^{*}$. The minimal polynomial $\min(A)$ of $A$ is irreducible of degree $2$ if $\delta(\sigma)\not=1$; moreover, $\min(A)$ divides the polynomial $y^{2n}-\lambda \in \K(x)[y]$. As $\det(A)^{2n}=((det(A))^n)^2=\lambda^2$, there exists $g \in \K(x)^{*}$ such that $\lambda=g^n$; this implies that $y^{2n}-\lambda=\prod_{k=1}^n(y^2-e^{2\im k \pi/n} g)\hts$ , whence $\min(A)=(y^2-e^{2\im k \pi/n} g)$ for a suitable $k$, and finally that $A^2=e^{2\im k \pi/n} g I$.
\proofend
\end{proof}

\bigskip

\break

We now give the conjugacy classes of involutions belonging to $\PGL(2,\K(x))$ in the de Jonqui\`eres group $\PGL(2,\K(x))\rtimes \PGL(2,\K)$. 
\begin{Prp}\titreProp{Conjugacy classes of involutions of $\PGL(2,\K(x))$ in the group $\PGL(2,\K(x))\rtimes \PGL(2,\K)$)}\index{Involutions!De Jonqui\`eres!in $\PGL(2,\K(x))$|idi}
\label{Prp:classDeJIJonqG}
\begin{itemize}
\item[\upshape 1.]
Any involutive element of $\PGL(2,\K(x))$ is conjugate in the de Jonqui\`eres group 
\begin{center}$\CrPp=\PGL(2,\K(x))\rtimes \PGL(2,\K)$\end{center} either to the automorphism
\begin{center}$\sigma_0:(x_1:x_2)\times (y_1:y_2) \dasharrow (x_1:x_2) \times(y_2:y_1)$\end{center}
of $\mathbb{P}^1\times\mathbb{P}^1$, or to the birational map 
\begin{center}$\sigma_{a,b}:(x_1:x_2)\times (y_1:y_2) \dasharrow (x_1:x_2) \times(y_2 \prod_{i=1}^k(x_1-b_i x_2):y_1 \prod_{i=1}^k(x_1-a_i x_2))$\end{center}of $\mathbb{P}^1\times\mathbb{P}^1$, for some integer $k\geq 1$ and some $a_1,...,a_k,b_1,...,b_k \in \K$ all distinct. 
\item[\upshape 2.]
The set of fixed points of $\sigma_{0}$ is the disjoint union of the two rational curves $y_1=y_2$ and $y_1=-y_2$.
\item[\upshape 3.]
The set of fixed points of $\sigma_{a,b}$ is the smooth hyperelliptic curve $\Gamma \subset \mathbb{P}^1\times\mathbb{P}^1$ of equation \index{Curves!hyperelliptic!birational maps that fix a hyperelliptic curve}
\begin{center}$y_1^2\prod_{i=1}^k(x_1-a_i x_2)=y_2^2\prod_{i=1}^k(x_1-b_i x_2)$.\end{center}
The projection $(x_1:x_2)\times (y_1:y_2) \mapsto (x_1:x_2)$ gives $\Gamma$ as a double covering of $\mathbb{P}^1$ ramified over the $2k$ points $(a_1:1),...,(a_k:1),(b_1:1),...,(b_k:1)$.
\item[\upshape 4.]
$\sigma_0$ is not conjugate to any $\sigma_{a,b}$ in $\CrPp$.
\item[\upshape 5.]
The linear equivalence of the set $\{(a_1:1),...,(a_k:1),(b_1:1),...,(b_k:1)\}$, and the isomorphism class of the hyperelliptic curve fixed by $\sigma_{a,b}$, determine the conjugacy class of $\sigma_{a,b}$ in $\CrPp$.\index{Curves!hyperelliptic!isomorphism class|idi}
\item[\upshape 6.]
The morphism $\delta: \PGL(2,\K(x)) \rightarrow \K(x)^{*}/\K(x)^{*2}$ given by the determinant induces a bijection between the set of conjugacy classes of involutions of $\PGL(2,\K(x))$ in the group $\PGL(2,\K(x))\rtimes \PGL(2,\K)$ and  the equivalence classes of $\K(x)^{*}/\K(x)^{*2}$, up to linear action of $\PGL(2,\K)$.
\end{itemize}
\end{Prp}
\begin{proof}\fauxtitred\upshape
By Proposition \refd{Prp:InvDeJPgl2}, any involution belonging to $\PGL(2,\K(x))$ is conjugate in this group to an involution $\sigma_{g(x)}$, for some $g(x)\in \K(x)\backslash\{0\}$. Furthermore, we may suppose that\begin{center}$g(x)=1$\hspace{0.3 cm} or\hspace{0.3 cm} $g(x)=\frac{\prod_{i=1}^l(x-b_i)}{\prod_{i=1}^k(x-a_i)}$,\end{center} for some $a_1,...,a_k,b_1,...,b_l \in \K$ all distinct and $l=k$ or $l=k+1$.

Writing $\sigma_{g(x)}$ as a birational map of $\mathbb{P}^1 \times \mathbb{P}^1$, we get $\sigma_0$ if $g(x)=1$, $\sigma_{a,b}$ if $l=k$ and the birational map
$$(x_1:x_2)\times (y_1:y_2) \dasharrow (x_1:x_2) \times(y_2 \prod_{i=1}^l(x_1-b_i x_2):y_1 x_2\prod_{i=1}^k(x_1-a_i x_2)),$$ if $l=k+1$. It suffices then to make a change of coordinates on $(x_1:x_2)$ to obtain the standard form given in assertion~1.

Assertions $2$ and $3$ follow directly by calculation.

Let us prove assertion $4$. Note that $\sigma_0$ commutes with $\PGL(2,\K)$. The conjugacy class of this element in $\PGL(2,\K(x))\rtimes \PGL(2,\K)$ is the same as its conjugacy class in $\PGL(2,\K(x))$. By Proposition \refd{Prp:InvDeJPgl2}, $\sigma_0$ is not conjugate to any $\sigma_{a,b}$ in this group, since $\delta(\sigma_0)\not=\delta(\sigma_{a,b})$.

Note that the hyperelliptic curve fixed by $\sigma_{a,b}$, for $a=\{a_1,...,a_k\}$ and $b=\{b_1,...b_k\}$, is rational if $k=1$, elliptic if $k=2$ and of genus $>1$ otherwise.\index{Curves!elliptic!viewed as hyperelliptic curves|idi}\index{Curves!elliptic!birational maps that fix an elliptic curve|idi}
Using the invariant defined in Section \refd{Sec:CurveFixedPoints} (the normalized fixed curve), we see that two involutions $\sigma_{a,b}$ that are conjugate in $\CrPp$ (and in the Cremona group) must have isomorphic fixed curves.
We recall that the isomorphism class of a hyperelliptic curve is defined by the linear equivalence of the set of ramification points.

To prove assertion $5$, it remains then to show that two involutions $\sigma_{a,b}$ and $\sigma_{a',b'}$ which have isomorphic curves of fixed points in $\mathbb{P}^1\times \mathbb{P}^1$ are conjugate in $\CrPp$.
After a change of coordinates, we may assume that the sets of ramification are the same, i.e.\ that $\{a_1,...,a_k,b_1,...,b_k\}=\{a_1',...,a_k',b_1',...,b_k'\}$. Note that conjugating by 
\begin{center}$(x_1:x_2)\times (y_1:y_2) \mapsto (x_1:x_2)\times ((x_1-a_i x_2)y_1:(x_1-b_j x_2)y_2)$\end{center} amounts to exchanging the roles of $a_i$ and $b_j$. This shows that $\sigma_{a,b}$ and $\sigma_{a',b'}$ are conjugate in $\CrPp$.

Let us prove the last assertion. The determinant of $\sigma_0$ is equal to $1$ and that of $\sigma_{a,b}$ is equal to $-\prod_{i=1}^k(x_1-a_i x_2)\cdot \prod_{i=1}^k(x_1-b_i x_2)$. Moreover, two elements of $\PGL(2,\K(x))$ that are conjugate in $\PGL(2,\K(x))\rtimes \PGL(2,\K)$ have the same determinant, up to an action of $\PGL(2,\K)$ on $\K(x)$.  Assertion $6$ now follows from assertions $4$ and $5$.
 \proofend
\end{proof}

We now prove that the conjugacy classes of involutions in the de Jonqui\`eres group are the trace of their conjugacy classes in the Cremona group, except for the elements that do not fix a curve of positive genus: all involutions that do not fix a curve of positive genus are conjugate in the Cremona group, but this is not true in the de Jonqui\`eres group: $\sigma_0$, $\sigma_{a_1,b_1}$, $a_1,b_1 \in \K$ and $(x,y)\mapsto (-x,y)$ represent three different conjugacy classes of the de Jonqui\`eres group.\index{Curves!of positive genus!birational maps that do not fix a curve of positive genus!involutions|idb}

 The following result gives in particular the classification of de Jonqui\`eres involutions in the Cremona group, a famous result which was proved in \cite{bib:BaB}, Proposition 2.7.

\begin{Prp}\titreProp{classification of involutions of $\PGL(2,\K(x))$ in the Cremona group}
\label{Prp:classDeJICremona}\index{Involutions!De Jonqui\`eres!in $\PGL(2,\K(x))$|idi}
\begin{itemize}
\item[\upshape 1.]
Any de Jonqui\`eres involution is birationally conjugate to some $\sigma_{a,b}\in \PGL(2,\K(x))$.
\item[\upshape 2.]
The normalized fixed curve (invariant defined in Section \refd{Sec:CurveFixedPoints}) determines the conjugacy class of the de Jonqui\`eres involutions in the Cremona group.
\end{itemize}
\end{Prp}
\begin{proof}\fauxtitred\upshape
By Proposition \refd{Prp:DeJconjPGL2C}, any de Jonqui\`eres involution is conjugate to an element of $\PGL(2,\K(x))$. Note that $\sigma_0$ (whose set of fixed points consists of two distinct rational curves) is conjugate to an involution $\sigma_{a,b}$ whose set of fixed points is an irreducible rational curve.

Indeed, $\sigma_{a,b}$ is conjugate to an automorphism of some conic bundle $(S,\pi)$, which acts trivially on the fibration. If the fixed curve is rational, the number of singular fibres twisted is at most $2$, hence by Proposition \refd{Prp:LargeDegree}, $\sigma_{a,b}$ is conjugate to an involution of $\Aut(\mathbb{P}^1\times\mathbb{P}^1)$, itself conjugate to $\sigma_0$ (see Proposition \refd{Prp:AutP1Bir}).

Assertion $2$ follows then from Proposition \refd{Prp:classDeJIJonqG}, assertion $5$.
\proofend
\end{proof}
\begin{remark}\fauxtitred
The birational conjugation from $\sigma_0$ to $\sigma_{a,b}$ can be determined explicitly.
 Indeed, $\sigma_0$ is conjugate in $\PGL(2,\K(x))$ to the birational map $(x,y)\mapsto (x,-y)$, itself conjugate by the birational map $(x,y) \dasharrow (x^2-y^2,x-y)$ to $(x,y) \dasharrow (x,x/y)$. This last map corresponds to the birational map $(x_1:x_2) \times (y_1:y_2) \dasharrow (x_1:x_2) \times (y_2x_1:y_1x_2)$, which is conjugate to $\sigma_{a_1,b_1}$, for any distinct $a_1,b_1 \in \K$.
\end{remark}

As we remarked earlier (see Proposition \refd{Prp:BeauTsen}), involutive elements of $\PGL(2,\K(x))$ correspond to involutive automorphisms of conic bundles with a trivial action on the fibration. We determine in which cases involutions of this sort are twisting de Jonqui\`eres involutions.
\begin{Prp}\fauxtitre
\label{Prp:TwistingOrNot}
Let $\sigma \in \PGL(2,\K(x))\subset \CrPp$ be an involution and let $\varphi: (\mathbb{P}^1\times\mathbb{P}^1,\pi_1) \dasharrow (S,\pi)$ be a birational map of conic bundles such that $\sigma'=\varphi \sigma\varphi^{-1}\in \Aut(S,\pi)$. Let $\Delta \subset S$ denote the union of the singular points of each singular fibre of $\pi$.
Then, one of the two following situations occurs:

\begin{itemize}
\item[\upshape 1.]

\begin{itemize}
\item
$\sigma'$ is a twisting $(S,\pi)$-de Jonqui\`eres involution.
\item
The set of points of $S$ fixed by $\sigma'$ is the union of $\Delta$ and a smooth irreducible hyperelliptic curve.\index{Curves!hyperelliptic!birational maps that fix a hyperelliptic curve}
\item
$\sigma$ is conjugate in $\CrPp$ to the birational map 
\begin{center}$\sigma_{a,b}:(x_1:x_2)\times (y_1:y_2) \dasharrow (x_1:x_2) \times(y_2 \prod_{i=1}^k(x_1-b_i x_2):y_1 \prod_{i=1}^k(x_1-a_i x_2))$,\end{center}for some integer $k\geq 1$ and some $a_1,...,a_k,b_1,...,b_k \in \K$ all distinct. \end{itemize}
\item[\upshape 2.]
\begin{itemize}
\item
$\sigma'$ is an involution which leaves invariant every irreducible component of each singular fibre of $\pi$.
\item
The set of points of $S$ fixed by $\sigma'$ is the union of $\Delta$ and two disjoint rational curves.
\item
$\sigma$ is conjugate in $\CrPp$ to the automorphism 
\begin{center}$\sigma_0:(x_1:x_2)\times (y_1:y_2) \dasharrow (x_1:x_2) \times(y_2:y_1).$\end{center}
\end{itemize}
\end{itemize}
Therefore, $\sigma'$ is a twisting $(S,\pi)$-de Jonqui\`eres involution if and only if the set of points of $\mathbb{P}^1\times\mathbb{P}^1$ fixed by $\sigma$ is a hyperelliptic curve, i.e.\ if and only if $\delta(\sigma)\not=1$.
\end{Prp}
\begin{proof}\fauxtitred\upshape
Note first of all that $\sigma'$ is either a twisting $(S,\pi)$-de Jonqui\`eres involution or an involution that leaves invariant every irreducible component of each singular fibre of $\pi$. The two cases mentioned above are then distinct for $\sigma'$. We will see further that the distinction between these two cases depends only on $\sigma$; this is the main interest of the proposition.

Lemma \refd{Lem:DeJI} gives information on the set of points of $S$ fixed by $\sigma'$, in the case where this is a twisting $(S,\pi)$-de Jonqui\`eres involution. If $\sigma'$ does not twist any singular fibre, we blow-down one irreducible component in every singular fibre, and conjugate $\sigma'$ to an involution of some Hirzebruch surface, with trivial action on the fibration. The set of points fixed by this involution is the union of two rational sections. This gives the above result for the set of points of $S$ fixed by $\sigma'$.

The set of points of $\mathbb{P}^1\times\mathbb{P}^1$ fixed by any involutive element of $\PGL(2,\K(x))$ intersects a general fibre in two points. The set of fixed points is then either a hyperelliptic curve, or two rational curves, plus some isolated points. This distinction clearly doesn't depend on a birational conjugation that preserves the general fibre.

It was proved in Proposition \refd{Prp:classDeJIJonqG} that the set of fixed points of $\sigma_{a,b}$ is a hyperelliptic curve and that the set of fixed points of $\sigma_0$ is the union of the two rational sections $y_1=y_2$ and $y_1=-y_2$. This shows that the distinction between twisting and non-twisting involutions is the same as the conjugation to $\sigma_0$ or to some $\sigma_{a,b}$ (it was proved in Proposition \refd{Prp:classDeJIJonqG} that any involutive element of $\PGL(2,\K(x))$ is conjugate to some $\sigma_{a,b}$ or to $\sigma_0$, the two cases being distinct).

The last assertion follows directly from the preceding ones and Proposition \refd{Prp:InvDeJPgl2}.
\proofend
\end{proof}

\bigskip 

Since we have seen that $\sigma_0$ is conjugate to some $\sigma_{a_1,b_1}$ (not in the de Jonqui\`eres group but in the Cremona group), we get the following result:
\begin{Prp}\fauxtitre
\label{Prp:DeJconjTwist}
Every de Jonqui\`eres involution is birationally conjugate to a twisting $(S,\pi)$-de Jonqui\`eres involution, for some conic bundle $(S,\pi)$.
\end{Prp}
\begin{proof}\fauxtitred\upshape
By Proposition \refd{Prp:classDeJICremona},
any de Jonqui\`eres involution is birationally conjugate to some involution $\sigma_{a,b}\in \PGL(2,\K(x))$. By Proposition \refd{Prp:BeauTsen}, $\sigma_{a,b}$ is conjugate to an involution of some conic bundle $(S,\pi)$, with a trivial action on the fibration. Proposition \refd{Prp:TwistingOrNot}, tells us that this is a twisting $(S,\pi)$-de Jonqui\`eres involution.\proofend
\end{proof}

\section{Extension of a cyclic group by a de Jonqui\`eres involution}
\pagetitless{Finite abelian groups of automorphisms of conic bundles}{Extension of a cyclic group by a de Jonqui\`eres involution}
\label{Sec:ExtCyclDeJ}
Let us study the finite abelian groups $G$ of automorphisms of a conic bundle such that $G'$ is generated by a twisting de Jonqui\`eres involution and $\overline{\pi}(G)$ is a cyclic group, isomorphic to $\Z{n}$. 
Using the exact sequence (introduced in Section \refd{Sec:Exactsequence})
\begin{equation}
1 \rightarrow G' \rightarrow G \stackrel{\overline{\pi}}{\rightarrow} \overline{\pi}(G) \rightarrow 1,
\tag{\refd{eq:ExactSeqCB}}\end{equation}
we see that $G$ is isomorphic to $\Z{2n}$ or to $\Z{2}\times \Z{n}$ (the two cases are the same if and only if $n$ is odd). This depends on the splitting of the exact sequence. We accordingly separate our investigation into two cases, in Sections \refd{Sec:Exactseqsplit} and \refd{Sec:RootDeJ}.

\subsection{The case when the exact sequence splits}
\label{Sec:Exactseqsplit}
We study in this section the simplest case, when the exact sequence ({\refd{eq:ExactSeqCB}) splits.

We use the following purely algebraic lemma which follows directly from \cite{bib:Be2}, Lemma 2.5, and from Proposition 2.6:
\begin{Lem}
\label{Lem:splitsequence}
Let $G \subset \PGL(2,\K(x))\rtimes \PGL(2,\K)$ be a finite abelian group and suppose that the exact sequence
\begin{center}
$1 \rightarrow G' \rightarrow G \stackrel{\overline{\pi}}{\rightarrow} \overline{\pi}(G) \rightarrow 1$\end{center} splits. 
Then $G$ is conjugate, in the group $\PGL(2,\K(x))\rtimes \PGL(2,\K)$, to the group\begin{center} $\{ (g,h) \ |\ g \in G', h \in \overline{\pi}(G)\} \subset \PGL(2,\K(x))\rtimes \PGL(2,\K)$,\end{center} the direct product of $G'$ and $\overline{\pi}(G)$.
\end{Lem}
\begin{proof}\fauxtitred\upshape
All the ideas can be found in $\cite{bib:Be2}$. As in the proof of $\cite{bib:Be2}$, Proposition 2.6, we use the cohomology exact sequence associated to the exact sequence
\begin{center}$1 \rightarrow \K(x)^{*} \rightarrow GL(2,\K(x)) \rightarrow \PGL(2,\K(x)) \rightarrow 1$,\end{center} to deduce that the cohomology set $H^1(H,\PGL(2,\K(x)))$ is reduced to a single element, for any finite subgroup $H\subset \PGL(2,\K)$. This is then the case for $H^1(\overline{\pi}(G),\PGL(2,\K(x)))$ and we can apply $\cite{bib:Be2}$, Lemma 2.5 directly to obtain the desired result.
\proofend
\end{proof}

\bigskip

We use this powerful lemma to give the following classification:
\begin{Prp}\fauxtitre
\label{Prp:DirectProduct}
Let $G \subset \CrPp $ be a finite abelian subgroup of the de Jonqui\`eres group such that
\begin{itemize}
\item
The exact sequence (\refd{eq:ExactSeqCB}) splits.
\item
$\overline{\pi}(G)\cong \Z{n}$, for some integer $n>1$.
\item
$G'\cong \Z{2}$ is generated by a de Jonqui\`eres involution, whose set of fixed points is a hyperelliptic curve of genus $k-1$.\index{Curves!hyperelliptic!birational maps that fix a hyperelliptic curve}
\end{itemize}
Then, $k\geq 2$ is divisible by $n$, so the genus $k-1$ of the curve is equal to $n-1 \pmod{n}$, and $G\cong \Z{2}\times \Z{n}$ is conjugate in the de Jonqui\`eres group to the group generated by
\begin{center}
$\begin{array}{rrcrl}
{\num{C.2,\h{0.5}n}} \h{2}&\sigma_{a,b}:(x_1:x_2)\times (y_1:y_2)\h{3}& \dasharrow& \h{4}(x_1:x_2) \times\h{3} & (y_2 \prod_{i=1}^k(x_1-b_i x_2):y_1 \prod_{i=1}^k(x_1-a_i x_2)),\\
\h{2}&\rho:(x_1:x_2)\times (y_1:y_2)\h{3}& \mapsto& \h{4}(x_1:\zeta_n x_2)\times\h{3} & (y_1:y_2),\end{array}$\end{center}
where $a_1,...,a_k,b_1,...,b_k \in \K^{*}$ are all distinct, $\zeta_n=e^{2\im \pi/n}$ and the sets $\{a_1,...,a_k\}$ and $\{b_1,...,b_k\}$ are both invariant by multiplication by $\zeta_n$.

Furthermore if $k>2$, two such groups are conjugate in $\CrPp$ {\upshape (}and hence in $\CrP${\upshape )} if and only if they \emph{have the same action on their sets of fixed points} 
(see Definition \refd{Def:actionFixedPoints}). 
Explicitly, the curve $y_1^2 \prod_{i=1}^k(x_1-a_i x_2)= y_2^2 \prod_{i=1}^k(x_1-b_i x_2)$ and the action $(x_1:x_2)\ \mapsto (x_1:\zeta_n x_2)$ on it determine the conjugacy class of the group.
\end{Prp}
\begin{proof}\fauxtitred\upshape
Conjugating $G$ in $\CrPp$ we may assume that $\overline{\pi}(G)$ is diagonal, generated by 
\begin{center}$(x_1:x_2) \mapsto (x_1:\zeta_n x_2)$.\end{center} We conjugate $G$ by some element of $\PGL(2,\K(x))$ and suppose that $G'$ is generated by a de Jonqui\`eres involution of the standard form given in Proposition \refd{Prp:classDeJIJonqG}, assertion $1$. We apply Lemma \refd{Lem:splitsequence} and see that $G$ is generated by $\sigma_{a,b}$ and $\rho$, as affirmed in the statement of the proposition. The fact that the sets $\{a_1,...,a_k\}$ and $\{b_1,...,b_k\}$ are both invariant by multiplication by $\zeta_n$ follows from the hypothesis that $G$ is abelian. Note that the orbits of the multiplication by $\zeta_n$ on $\K$ have size exactly $n$, except for $0$, which is fixed. As all the $a_i$ and $b_j$ are distinct, the number $0$ can appear only once; by counting the size of the orbits we see that $a_1,...,a_k,b_1,...,b_k \in \K^{*}$. This implies that $k$ is divisible by $n$ and so $k>1$.

We now prove the last assertion. Let $G$ (respectively $G'$) be the group generated by $\sigma_{a,b}$ (respectively $\sigma_{a',b'}$) and $\rho$. We denote by $\Gamma$ (respectively $\Gamma'$) the curve of fixed points of $\sigma_{a,b}$ (respectively $\sigma_{a',b'}$). We suppose that $G$ and $G'$ \emph{have the same action on their sets of fixed points}. This means that there exists an isomorphism $\varphi:\Gamma' \rightarrow \Gamma$ and an integer $m$ such that $\rho^m\varphi=\varphi \rho$.

As $k>2$, the genus of $\Gamma$ and $\Gamma'$ is $>1$, so $\varphi$ induces an automorphism $\overline{\varphi}$ of $\mathbb{P}^1$. In this case, we have $\overline{\rho}^m\circ\overline{\varphi}=\overline{\varphi}\circ \overline{\rho}$, where $\overline{\rho}(x_1:x_2)=(x_1:\zeta_n x_2)$. This implies that $\overline{\varphi}$ is of the form $(x_1:x_2) \mapsto (\lambda x_2:x_1)$ or $(x_1:x_2) \mapsto (\lambda x_1:x_2)$, for some $\lambda \in \K^{*}$, and that $m=\pm 1$. Note that the automorphism $\tilde{\varphi}:(x_1:x_2) \times (y_1:y_2) \mapsto \overline{\varphi}(x_1:x_2)\times (y_1:y_2)$ of $\mathbb{P}^1\times\mathbb{P}^1$ sends the ramification points of $\Gamma'$ on those of $\Gamma$ and does not change the standard form: if we conjugate $G'$ by $\tilde{\varphi}$, we get a group $G''$ generated by $\rho$ and some $\sigma_{a'',b''}$, with $\{a'',b''\}=\{a,b\}$.

Let $\psi$ denote the birational map
\begin{center}
$\psi:(x_1:x_2)\times (y_1:y_2) \dasharrow (x_1:x_2)\times (y_1 \prod_{i=1}^k(x_1-a_i x_2) :y_2 \prod_{i=1}^k(x_1-a_i'' x_2))$\end{center}
of $\mathbb{P}^1\times\mathbb{P}^1$, which commutes with $\rho$. The two birational maps
\begin{center}
$\begin{array}{rcl}
\psi \circ \sigma_{a,b}(x,y)& = &(x_1:x_2)\times (y_2 \prod_{i=1}^k(x_1-b_i x_2):y_1\prod_{i=1}^k(x_1-a''_i x_2) )\vspace{0.1 cm}\\
\sigma_{a'',b''}\circ \psi(x,y)& = &(x_1:x_2)\times (y_2 \prod_{i=1}^k(x_1-b''_i x_2):y_1\prod_{i=1}^k(x_1-a_i x_2) )\end{array}$\end{center}
are equal, since $\{a'',b''\}=\{a,b\}$. The groups $G$ and $G'$ are thus conjugate in $\CrPp$ by the birational map $\psi$.
\proofend
\end{proof}
\begin{remark}\fauxtitred
If $k=2$, the curve is elliptic\index{Curves!elliptic!viewed as hyperelliptic curves}\index{Curves!elliptic!birational maps that fix an elliptic curve} and an automorphism of it does not necessarily preserve the fibration. We do not know whether the action on the sets of fixed points determines the conjugacy classes in this case.
\end{remark}

\bigskip

\subsection{When the exact sequence does not split - roots of de Jonqui\`eres involutions}
\label{Sec:RootDeJ}
\index{Involutions!De Jonqui\`eres!roots|idb}
We continue the decription of the groups $G$ of automorphisms of a conic bundle such that $G'$ is generated by a twisting de Jonqui\`eres involution and $\overline{\pi}(G)$ is a cyclic group, isomorphic to $\Z{n}$. In this section, we consider the case in which the exact sequence
\begin{equation}
1 \rightarrow G' \rightarrow G \stackrel{\overline{\pi}}{\rightarrow} \overline{\pi}(G) \rightarrow 1
\tag{\refd{eq:ExactSeqCB}}\end{equation}
does not split, i.e.\ in which $G$ is isomorphic to $\Z{2n}$ but not to $\Z{2}\times \Z{n}$, so that $n$ is even.
We then have to classify roots of de Jonqui\`eres involutions of even order, that belong to the de Jonqui\`eres group. Note that the roots of odd order were classified in Proposition \refd{Prp:DirectProduct} because then $G \cong \Z{2}\times \Z{n}$. 
\begin{remark}\fauxtitred
Any root of a de Jonqui\`eres involution {\upshape (}not necessarily in $\CrPp${\upshape )} generates a finite group which acts biregularly on a rational surface (see Proposition \refd{Prp:TwoCases}). If the root preserve some conic bundle structure, it is conjugate to an element of the de Jonqui\`eres group, described for odd roots in Section \refd{Sec:Exactseqsplit} and for even roots in this section.

If the group does not leave invariant a pencil of rational curves, then it is conjugate to an automorphism of a Del Pezzo surface\index{Del Pezzo surfaces!automorphisms|idi}. We will see further (Theorem \refThmCyclicGroup\ in Chapter 10)  that the de Jonqui\`eres involutions acting on these surfaces necessarily have a rational or elliptic curve of fixed points.\index{Curves!elliptic!viewed as hyperelliptic curves}\index{Curves!elliptic!birational maps that fix an elliptic curve}

\end{remark}
The following result describes the geometry of the situation:
\begin{Prp}\titreProp{roots of even order of de Jonqui\`eres involutions}\\
\label{Prp:rootsDeJ}
\index{Involutions!De Jonqui\`eres!roots|idb}
Let $(S,\pi)$ be a conic bundle, let $\alpha \in \Aut(S,\pi)$ and let $G$ be the group generated by $\alpha$. We assume that $\alpha^n$ is a twisting $(S,\pi)$-de Jonqui\`eres involution, for some even integer $n\geq 2$ and that the triple $(G,S,\pi)$ is minimal. Let $\Gamma$ denote the smooth hyperelliptic curve fixed by $\alpha^n$, a double covering of $\mathbb{P}^1$ ramified over $2k$ points.\index{Curves!hyperelliptic!birational maps that fix a hyperelliptic curve}

\begin{itemize}
\item[\upshape 1.]
The action of $\alpha$ on the basis has order $n$.
\item[\upshape 2.]
The action of $\alpha$ on the set of ramification has no fixed points. In particular, $n$ divides $2k$.
\item[\upshape 3.]
Let $r$ denote the number of orbits of $\alpha$ on the ramification points (so that $rn=2k$) and $l$ denote the number of singular fibres of $\pi$. Then, one of the following holds:
\begin{itemize}
\item[{\upshape (a)}]
$l=2k$, $r$ is even, and $\alpha$ acts on $\Gamma$ with $4$ fixed points;
\item[{\upshape (b)}]
$l=2k+1$, $r$ is odd, and $\alpha$ acts on $\Gamma$ with $2$ fixed points;
\item[{\upshape (c)}]
$l=2k+2$, $r$ is even, and $\alpha$ fixes no point of $\Gamma$.
\end{itemize}
\item[\upshape 4.]
The involution $\alpha^n$ twists $2k$ of the $l$ singular fibres of $\pi$, and $\alpha$ twists the $l-2k$ others.
\end{itemize}

\end{Prp}
\begin{proof}\fauxtitred\upshape
As an $(S,\pi)$-twisting de Jonqui\`eres involution has no root in $\PGL(2,\K(x))$ (Proposition \refd{Prp:InvDeJPgl2}, assertion $4$), the kernel $G'$ of $\overline{\pi}$ is generated by the de Jonqui\`eres involution $\alpha^n$; the action of $\alpha$ on the basis is then exactly of order $n$. This proves assertion $1$.

If $\alpha$ fixes a ramification point, the singular fibre of this point is invariant by $\alpha$. But $\alpha^n$ must twist this singular fibre, and this is not possible since $n$ is an even integer.
The action of $\alpha$ on the ramification points thus has only orbits of size $n$. This implies that $n$ divides the number of ramification points, which is $2k$. This yields assertion $2$.

Let $p_0:S\rightarrow \mathbb{P}^1\times\mathbb{P}^1$ be the birational morphism given by Lemma \refd{Lem:GoingToF0F1}, denote by $F_1,...,F_l$ the divisor classes of the exceptional curves contracted by $p_0$ (each belonging to a singular fibre of $\pi$), and let $f$ denote the divisor class of the fibre of $\pi$.
Up to a change in the numbering of the $F_i$'s, we may assume that $\pi(F_1),...,\pi(F_{2k}) \in \mathbb{P}^1$ are the ramification points of the double covering $\pi:\Gamma\rightarrow\mathbb{P}^1$ and that the action of $\alpha$ on these points is as follows:
\begin{center}
$\pi(\alpha(F_i))=\left\{\begin{array}{ll}
\pi(F_{i+1}) & \mbox{if }i \mbox{ is not divisible by }n,\\
\pi(F_{i+1-n})& \mbox{otherwise.}\end{array}\right.$
\end{center}

Since the triple $(G,S,\pi)$ is minimal by hypothesis, every singular fibre of $\pi$ must be twisted by some element of $G$. If $l>2k$, the singular fibres $\{F_i,f-F_i\}$ for $2k<i\leq l$ are not twisted by $\alpha^n$, and hence must be twisted by some element that does not belong to $G'$. This implies that the singular fibres are invariant by $\alpha$ and twisted by it. As $\alpha$ leaves invariant exactly two fibres of $\pi$, the number $l$ of singular fibres of $\pi$ is equal to $2k,2k+1$ or $2k+2$ and $\alpha(F_i)=f-F_i$ if $i>2k$, as stated above. We get then assertion $4$.

Let $E$ denote the divisor class of the strict pull-back of the general section of self-intersection $0$ of $\mathbb{P}^1\times\mathbb{P}^1$. We have $E\cdot F_i=0$ for $1\leq i \leq l$ and $K_S=-2E-2f+\sum_{i=1}^l F_i$. The image of $E$ by $\alpha$ is then \begin{center}$\alpha(E)=\frac{1}{2}\alpha(-K_S-2f+\sum_{i=1}^l F_i)=E+\frac{1}{2}\sum_{i=1}^l(\alpha(F_i)-F_i)$.\end{center} 

Let us compute the number of times that the divisor $f$ appears in $\sum_{i=1}^l(\alpha(F_i)-F_i)$. This number must be even by the above relation. We will say that $\alpha$ makes \emph{a change at $F_i$} if $\alpha(F_i)=f-F_{j}$ for some $j=i+1 \pmod{n}$. Note that the number of times that $f$ appears in $\sum_{i=1}^n(\alpha(F_i)-F_i)$ is exactly the number of changes that $\alpha$ makes at $F_1,...,F_n$. Since $\alpha^n(F_1)=f-F_1$, this number is odd. The situation is similar for the other orbits of size $n$ of $\alpha$. As $\alpha(F_i)-F_i=f-2F_i$ for $i>2k$, the number of times that $f$ appears in $\sum_{i=1}^l(\alpha(F_i)-F_i)$ is equal to $r+(l-2k) \pmod{2}$. So $l-2k=r \pmod{2}$, which gives the three cases stated in assertion $3$.
It remains to prove the relation on the fixed points of the action of $\alpha$ on the curve $\Gamma$.

As the action of $\alpha$ on the basis has order $n$, it leaves exactly two fibres invariant (which do not contain a point of ramification, see above). Let $F$ be a fibre of $\pi$ invariant by $\alpha$; we describe the situation, which depends on the singularity of the fibre:
\begin{itemize}
\item
If $F$ is non-singular, then $F \cong \mathbb{P}^1$. The action of $\alpha$ on $F$ has finite order and is then diagonalizable. It therefore fixes two points of $F$ which are also fixed by $\alpha^n$ and are then $\Gamma\cap F$.
\item
If $F$ is singular, then $\alpha$ twists it (assertion $4$ proved above). So $\alpha$ fixes exactly one point of $F$, which is the singular point. But this point does not belong to $\Gamma$ (see Lemma \refd{Lem:DeJI}), so the two points of $\Gamma \cap F$ are exchanged by $\alpha$.
\end{itemize}
From this observation, we see that the number of points of $\Gamma$ fixed by $\alpha$ is $4$ if $\alpha$ does not twist any singular fibre (the case $l=2k$), is $2$ if $\alpha$ twists exactly one singular fibre (the case $l=2k+1$) and is $0$ is $\alpha$ twists two singular fibres (the case $l=2k+2$). We thus get assertion~3.\proofend
\end{proof}
\begin{remarks}\fauxtitred
\begin{itemize}
\item
The cases $l=2k,2k+1,2k+2$ represent distinct conjugacy classes of groups, since the action on the fixed curve has respectively $0$, $2$, and $4$ fixed points.
\item
In the case of a root of a de Jonqui\`eres involution whose set of fixed points does not contain a curve of positive genus (the case $k=1$), this proposition proves that any root belonging to the de Jonqui\`eres group acts on a conic bundle with $l=3$ singular fibres. By Proposition \refd{Prp:LargeDegree}, this root is then conjugate to a element of $\Aut(\Pn)$ and also to an element of $\Aut(\mathbb{P}^1\times\mathbb{P}^1)$. This will be generalised to any root of some element conjugate to a linear automorphism, in Theorem \refThmRootsLinearAutomorphisms\ (Chapter 10).
\end{itemize}
\end{remarks}

\bigskip

Proposition \refd{Prp:rootsDeJ} describes the geometry on the conic bundle and implies that the action of $\alpha$ on the curve $\Gamma$ fixes no ramification point. Conversely, given some $\Gamma$ and such an action on it, we prove in the following proposition that there exists a root $\alpha$ of the de Jonqui\`eres involution corresponding to $\Gamma$.
\begin{Prp}\titreProp{existence of the roots}\\
\index{Involutions!De Jonqui\`eres!roots|idi}
\label{Prp:ExistRoot}
Let $\beta \in \PGL(2,\K)$ be an automorphism of $\mathbb{P}^1$ of even order $n\geq 2$ and let $\{p_1,...,p_{rn}\}$ be a set of $rn$ distinct points in $\mathbb{P}^1$, invariant by $\beta$.
Assume that no point $p_i$ is fixed by $\beta$.

Then, for $l=rn+1$ if $r$ is odd (respectively for $l=rn$ and $l=rn+2$ if $r$ is even), there exists a conic bundle $(S,\pi)$, which is the blow-up of $l$ points of $\mathbb{P}^1\times\mathbb{P}^1$, and $\alpha \in \Aut(S,\pi)$ such that
\begin{itemize}
\item
the triple $(<\alpha>,S,\pi)$ is minimal;
\item\index{Curves!hyperelliptic!birational maps that fix a hyperelliptic curve}
$\alpha^n$ is a twisting $(S,\pi)$-de Jonqui\`eres involution, whose set of fixed points is the disjoint union of $l-rn$ points and a smooth hyperelliptic curve $\Gamma$, such that $\pi:\Gamma \rightarrow \mathbb{P}^1$ is a double covering ramified over $p_1,...,p_{rn}$.
\item
The action of $\overline{\alpha}$ on $\mathbb{P}^1$ corresponds to $\beta$.
\end{itemize}
\end{Prp}
\begin{proof}\fauxtitred\upshape
Our aim is to find $l$ points $P_1,...,P_l \in \mathbb{P}^1\times\mathbb{P}^1$ such that:
\begin{itemize}
\item
The blow-up $p:S \rightarrow \mathbb{P}^1\times\mathbb{P}^1$ of these points gives a conic bundle $(S,\pi)$, where $\pi=\pi_1 \circ p$.
\item
There exists an $\alpha \in \Aut(S,\pi)$ which satisfies the conditions stated in the proposition.
\end{itemize}

Firstly, we describe conditions on the action of $\alpha$ on $\Pic{S}$, necessary for the existence of this automorphism. Secondly, we view these conditions as conditions on the birational map $\alpha' =p \alpha p^{-1}$ of $\mathbb{P}^1\times\mathbb{P}^1$. Thirdly we prove that such a map $\alpha'$ exists, and finally we verify that the map $\alpha'$ obtained gives an automorphism $\alpha$ satisfying the conclusions of the proposition.
\bigskip

Up a permutation of the $p_i$'s, we may assume that for $i=1,...,rn$,
\begin{center}
$\beta(p_i)=\left\{\begin{array}{ll}
p_{i+1} & \mbox{if }i \mbox{ is not divisible by }n,\\
p_{i+1-n} & \mbox{otherwise.}\end{array}\right.$
\end{center} If $l>rn$, we denote by $p_{rn+1}$ one of the points of $\mathbb{P}^1$ fixed by $\beta$. If $l=rn+2$, we denote the other one by $p_{rn+2}$. We impose $\pi_1(P_i)=p_i$ for $i=1,...,l$. We set $F_i=p^{-1}(P_i)$ for $i=1,...,l$ and $f$ denotes the divisor class of the fibre of $\pi$.

The conditions which we want to impose on the action of $\alpha$ on $\Pic{S}$ are these: 
\begin{equation}
\label{conditionsroot}
\alpha(F_i)=\left\{\begin{array}{ll}
f-F_{i}& \mbox{if }i>rn,\\
f-F_{i+1-n} & \mbox{if }i=n,2n,...,rn,\\
F_{i+1} & \mbox{otherwise.}
\end{array}\right.
\end{equation}
Note that we need only find an automorphism $\alpha \in \Aut(S,\pi)$ with the same action as $\beta$ on $\mathbb{P}^1$ and satisfying the conditions (\refd{conditionsroot}). Indeed, in this case $\alpha^n$ acts trivially on the fibration and twists the fibres of $p_1,...,p_{rn}$, so it fixes a hyperelliptic curve given by $\pi$ as the double covering of $\mathbb{P}^1$ ramified over $p_1,...,p_{rn}$. The triple $(<\alpha>,S,\pi)$ is minimal since $\alpha$ twists the fibres of $p_i$ for $rn<i\leq l$.

Let $E$ be the divisor class of the strict pull-back by $p$ of the general section of self-intersection $0$ of $\mathbb{P}^1\times\mathbb{P}^1$. As proved in Proposition \refd{Prp:rootsDeJ}, we have the following relations:
\begin{center}$\begin{array}{lllll}
\alpha(E)&=&E+\frac{1}{2}\sum_{i=1}^l(\alpha(F_i)-F_i)&=&E+df-\sum_{i=0}^{r-1}F_{in+1}-\sum_{i=rn+1}^l F_i,\vspace{0.1 cm}\\
\alpha^{-1}(E)&=&E+\frac{1}{2}\sum_{i=1}^l(\alpha^{-1}(F_i)-F_i)&=&E+df-\sum_{i=1}^rF_{in}-\sum_{i=rn+1}^l F_i,\end{array}$\end{center}
where $d=\frac{r+(l-rn)}{2}$.
The birational map $\alpha'=p \alpha p^{-1}$ of $\mathbb{P}^1\times\mathbb{P}^1$ is then
\begin{center}
$\alpha': (x_1:x_2)\times(y_1:y_2) \dasharrow \beta(x_1:x_2) \times (\lambda_1(x_1,x_2,y_1,y_2):\lambda_2(x_1,x_2,y_1,y_2))$,\end{center} where $\lambda_1,\lambda_2$ generate a pencil $\Lambda$ of curves of bidegree $(d,1)$. We set $P_i=(p_i,q_i) \in \mathbb{P}^1\times\mathbb{P}^1$, for $i=1,...,l$. The elements $q_1,...,q_l \in \mathbb{P}^1$, $(\lambda_1,\lambda_2) \in \mathbb{P}^{4d+3}$ belong to the projective variety $(\mathbb{P}^1)^l\times\mathbb{P}^{4d+3}$. Note that $\alpha'$ is birational if and only if $\lambda_1$ and $\lambda_2$ are not collinear. Indeed, in such a case the intersection of the pencils $\beta(x_1,x_2)$ and $\Lambda$ is $1$. We impose these conditions:
\begin{itemize}
\item
For $k=n,2n,...,rn$, $P_k$ is a base point of $\Lambda$ and $\alpha'(p_k \times \mathbb{P}^1)=P_{k+1-n}$.\item
For $k<rn$, and $k$ not divisible by $n$, $\alpha'(P_k)=P_{k+1}$.
\item
For $rn<k\leq l$, $P_k$ is a base point of $\Lambda$ and $\alpha'(p_k \times \mathbb{P}^1)=P_k$.
\end{itemize}
For any $k \in \{1,...,l\}$, imposing that $P_k$ be a base point implies two conditions, namely $\lambda_1(p_k,q_k)=\lambda_2(p_k,q_k)=0$. Then, imposing that $\lambda(p_k \times \mathbb{P}^1)$ be equal to some point gives only one more condition.

The number of conditions is thus 
\begin{center}$3r+(n-1)r+ 3(l-rn)=rn+2r+3(l-rn)=2(r+(l-rn))+l=4d+l.$\end{center}
The dimension of the solutions is therefore at least $3$, so we can choose one general solution, which gives a birational map $\alpha'$ whose $2d$ base points (the number of these follows from the degree of the map) are exactly $\{P_k \ |\ k=n,...,rn\} \cup \{P_k \ | \ rn<k\leq l\}$. 
As stated above, the blow-up $p:S \rightarrow \mathbb{P}^1\times\mathbb{P}^1$ of $P_1,...,P_l$ gives a conic bundle $(S,\pi)$, where $\pi=\pi_1 \circ p$. As $\alpha'$ leaves invariant the pencil of lines $\{ax_1+bx_2=0\ | \ (a:b)\in\mathbb{P}^1\}$, the element  $\alpha=p^{-1}\alpha'p$ belongs to $\Crc{S}{\pi}$. It remains to show that $\alpha$ satisfies the condition above and is biregular.

As $\alpha$ keeps invariant the fibration induced by $\pi$, we need only look at the image of the fibre of $x\in \mathbb{P}^1$, for any $x$, and show that $\alpha$ is an isomorphism from $\pi^{-1}(x)$ to $\pi^{-1}(\beta(x))$ :
\begin{itemize}
\item
If the fibre of $x$ by $\pi_1$  does not contain a base point of $\alpha'$, the map $\alpha'$ is an isomorphism from the fibre of $x$ to the fibre  of $\beta(x)$.
\begin{itemize}
\item
If the fibre does not contain any point $P_k$, for $k=1,...,l$, then $\alpha$ is also an isomorphism from the smooth fibre $\pi^{-1}(x)$ to the fibre $\pi^{-1}(\beta(x))$.
\item
If the fibre contains one point $P_k$ (with $k<rn$, $k\not=0\pmod{n}$, since the fibre does not contain a base point), then the fibre of $\pi(x)$ contains $P_{k+1}$, which is the image by $\alpha'$ of $P_k$. As $p$ blows-up $P_k$ and $P_{k+1}$, $\alpha$ is an isomorphism from the singular fibre $\{F_k,f-F_k\}$ to the singular fibre $\{F_{k+1},f-F_{k+1}\}$, and sends $F_k$ on $F_{k+1}$ and $f-F_k$ on $f-F_{k+1}$.
\end{itemize}
\item
Supose that the fibre of $x$ by $\pi_1$ contains some basis point $P_k$. The fibre $\pi_1^{-1}(x)$ is then contracted on some point $P_j$ (with $j=k+n-1$ if $k\leq rn$ and $j=k$ otherwise). 

Since $p$ blows-up all the base points of $\alpha'$ (all lying on $\mathbb{P}^1\times\mathbb{P}^1$), the birational map $\alpha'\circ p$ is a morphism from $S$ to $\mathbb{P}^1\times\mathbb{P}^1$, which is then surjective, as it is birational.

Thus, the map $\alpha'\circ p$ sends $F_k$ on the fibre of $\pi^{-1}(\beta(x))$ and its restriction to $F_k$ is an isomorphism. As $\alpha'$ blows-up $P_j$, which is the image of the curve $f-F_k$, $\alpha$ sends $f-F_k$ isomorphically on $F_j$.
\end{itemize}
We now have all the conditions on $\alpha$ that induce that $\alpha \in \Aut(S,\pi)$ and $\alpha^n$ is a $(S,\pi)$-twisting de Jonqui\`eres involution, that twists the fibres $\pi^{-1}(p_k)$, for $k=1,...,rn$.\proofend
\end{proof}

\bigskip

\begin{Exa}\fauxtitred\upshape
\label{Exa:Root2oddorder}
Let $n$ be an odd integer and let $g \in \K(x)$ be a rational function. Then, the map
\begin{center}$\alpha:(x,y)\dasharrow (\zeta_{2n}\cdot x, -g(x^n)\cdot \frac{y+g(-x^n)}{y+g(x^n)})$\end{center} (where $\zeta_{2n}=e^{2\im\pi/2n}$) is a birational map which is a $2n$-th root of the de Jonqui\`eres involution \begin{center}$\alpha^{2n}:(x,y)\dasharrow (x, \frac{g(x^n)\cdot g(-x^n)}{y})$\end{center} whose set of fixed points is $y^2=g(x^n)\cdot g(-x^n)$. 

The form of $\alpha^{2n}$ follows from the computation of $\alpha^2:(x,y)\dasharrow (\zeta_{n}\cdot x, \frac{g(x^n)\cdot g(-x^n)}{y})$, which is the composition of $(x,y)\dasharrow (x, \frac{g(x^n)\cdot g(-x^n)}{y})$ and $(x,y)\mapsto (\zeta_{n} x,y)$, two elements which commute.
\end{Exa}

\bigskip

\begin{Exa}\fauxtitred\upshape
\label{Exa:4throotDeJonq}
A $4$-th root of a de Jonqui\`eres involution:
\begin{center}
$\begin{array}{llcccl}
\alpha&:(x,y) \dasharrow (\h{3}&\im x&\h{3},\h{3}&\frac{(x+1)((\sqrt{2}-1)-x)y+(x^4-1)}{y+(x+1)((\sqrt{2}-1)-x))}&\h{4}),\vspace{0.1 cm}\\
\alpha^2&:(x,y) \dasharrow (\h{3}&-x&\h{3},\h{3}&\frac{-\im(x+1)(x-\im)y+x^4-1}{y-\im(x+1)(x-\im)}&\h{4}),\vspace{0.1 cm}\\
\alpha^4&:(x,y) \dasharrow (\h{3}&x&\h{3},\h{3}&\frac{x^4-1}{y}&\h{4}),\vspace{0.1 cm}\\
\alpha^8&:(x,y) \dasharrow (\h{3}&x&\h{3},\h{3}&y&\h{4}).
\end{array}$
\end{center}
\end{Exa}

\bigskip

\bigskip

Let us prove that square roots are conjugate to Example \refd{Exa:Root2oddorder}:
\begin{Prp}\fauxtitre
\label{Prp:SquareRoots}
Let $a \in \PGL(2,\K(x))\rtimes \PGL(2,\K)$ be such that $\overline{\pi}(a)\not=1$ and that $a^2 \in \PGL(2,\K(x))$ is an involution.

Then, up to conjugation in the group $\PGL(2,\K(x))\rtimes \PGL(2,\K)$, we have
\begin{center}$a=\Bigg(\left(\begin{array}{cc}\nu(x)& -\nu(x)\nu(-x) \\ 1 &-\nu(x)\end{array}\right),\left(\begin{array}{cc}-1 & 0 \\ 0& 1\end{array}\right)\Bigg),$\end{center}
for some $\nu(x) \in \K(x)\backslash\{0\}$.
\end{Prp}
\begin{proof}\fauxtitred\upshape
To simplify the notation, we will denote in this proof  the function $\varphi(-x)$ by $\tilde{\varphi}$, for any $\varphi=\varphi(x) \in \K(x)$. 

Up to conjugation, we may suppose that:

\break

\begin{center}$\begin{array}{ccccccc}a&=&\Bigg(&\h{3}\left(\begin{array}{cc}\alpha & \beta \\ \gamma& \delta\end{array}\right)&\h{3},\h{3}&\left(\begin{array}{cc}-1 & 0 \\ 0& 1\end{array}\right)&\h{3}\Bigg),\vspace{0.1cm}\\
a^2&=&\Bigg(&\h{3}\left(\begin{array}{cc}0 & \rho \\ 1& 0\end{array}\right)&\h{3},\h{3}&\left(\begin{array}{cc}1 & 0 \\ 0& 1\end{array}\right)&\h{3}\Bigg),\end{array}$\end{center} for some elements $\alpha,\beta,\gamma,\delta,\rho \in \K(x)\backslash\{0\}$.

From the relation \begin{center}$\left(\begin{array}{cc}0 & \rho \\ 1& 0\end{array}\right)=\left(\begin{array}{cc}\alpha & \beta \\ \gamma& \delta\end{array}\right)\left(\begin{array}{cc}\tilde{\alpha} & \tilde{\beta} \\ \tilde{\gamma}& \tilde{\delta}\end{array}\right)$\end{center}
 we get 
\begin{center}$\alpha\tilde{\alpha}+\beta\tilde{\gamma}=\delta\tilde{\delta}+\tilde{\beta}\gamma=0$,\end{center}
and we see that $\alpha\tilde{\alpha}=\delta\tilde{\delta}$. Writing this as $(-\frac{\delta(x)}{\alpha(x)})(-\frac{\delta(-x)}{\alpha(-x)})=1$ and using Hilbert $90$ (see \cite{bib:Serre}, page 158), there exists an element $\lambda(x) \in \K(x)\backslash\{0\}$ such that $\frac{\lambda(x)}{\lambda(-x)}=-\frac{\delta(x)}{\alpha(x)}$.
Conjugating then $a$ by the element $\Bigg(\h{1}\left(\begin{array}{cc}\lambda & 0 \\ 0& 1\end{array}\right),\left(\begin{array}{cc}1 & 0 \\ 0& 1\end{array}\right)\h{1}\Bigg)$ we may suppose that 
\begin{center}$a=\Bigg(\h{1}\left(\begin{array}{cc}\alpha\lambda& \beta\lambda\tilde{\lambda} \\ \gamma &\delta\tilde{\lambda}\end{array}\right),\left(\begin{array}{cc}-1 & 0 \\ 0& 1\end{array}\right)\h{1}\Bigg)=\Bigg(\h{1}\left(\begin{array}{cc}\nu& -\nu\tilde{\nu} \\ 1 &-\nu\end{array}\right),\left(\begin{array}{cc}-1 & 0 \\ 0& 1\end{array}\right)\h{1}\Bigg)$,\end{center}
where $\nu=\frac{\alpha\lambda}{\gamma}$.
\proofend
\end{proof}
\begin{remark}\fauxtitred
For $n$ odd, we do not know whether any $2n$-th root of a de Jonqui\`eres involution is conjugate to the form given Example \refd{Exa:Root2oddorder}
\end{remark}

\bigskip

\section{Extension of a group $(\Z{2})^2$ by a de Jonqui\`eres involution}
\pagetitless{Finite abelian groups of automorphisms of conic bundles}{Extension of a group $(\Z{2})^2$ by a de Jonqui\`eres involution}
\label{Sec:ExtC2-C2C2}
We now study the finite abelian groups $G$ of automorphisms of a conic bundle such that $G'$ is generated by a twisting de Jonqui\`eres involution and $\overline{\pi}(G)$ is isomorphic to $(\Z{2})^2$. 
Using the exact sequence (introduced in Section \refd{Sec:Exactsequence})
\begin{equation}
1 \rightarrow G' \rightarrow G \stackrel{\overline{\pi}}{\rightarrow} \overline{\pi}(G) \rightarrow 1
\tag{\refd{eq:ExactSeqCB}}\end{equation}
we see that $G$ is either isomorphic to $(\Z{2})^3$ or to $\Z{2}\times\Z{4}$. As in section \refd{Sec:ExtCyclDeJ}, this depends on the splitting of the exact sequence. 

\break

\begin{Prp}\fauxtitre
\label{Prp:DirectProductC2C2C2}
Let $G \subset \CrPp $ be a finite abelian subgroup of the de Jonqui\`eres group such that
\begin{itemize}
\item
The exact sequence (\refd{eq:ExactSeqCB}) splits.
\item
$\overline{\pi}(G)\cong (\Z{2})^2$.
\item
$G'\cong \Z{2}$ is generated by a de Jonqui\`eres involution, whose set of fixed points is a hyperelliptic curve of genus $k-1$.\index{Curves!hyperelliptic!birational maps that fix a hyperelliptic curve}
\end{itemize}
Then, $k\geq 4$ is divisible by $4$, and $G\cong (\Z{2})^3$ is conjugate in $\CrPp$ to the group generated by
\begin{center}
$\begin{array}{rrcccl}
{\num{C.2,\h{0.5}22}} &\sigma:(x_1:x_2)\times (y_1:y_2)\h{1}& \dasharrow \h{1}&(x_1:x_2) \h{2}&\times\h{2}&(y_2 \prod_{i=1}^{k/4}P(b_i):y_1 \prod_{i=1}^{k/4}P(a_i)),\\
&\rho_1:(x_1:x_2)\times (y_1:y_2)\h{1}& \mapsto \h{1}&(x_1:- x_2)\h{2}&\times\h{2}& (y_1:y_2),\\
&\rho_2:(x_1:x_2)\times (y_1:y_2)\h{1}& \mapsto \h{1}&(x_2: x_1)\h{2}&\times\h{2}& (y_1:y_2),\end{array}$\end{center}
where $P(a)=(x_1-ax_2)(x_1+ax_2)(x_1-a^{-1}x_2)(x_1+a^{-1}x_2) \in \K[x_1,x_2]$ and $a_1$,..., $a_{k/4}$, $b_1$,..., $b_{k/4} \in \K\backslash\{0,\pm 1\}$ are all distinct.

Furthermore, two such groups are conjugate in $\CrPp$ {\upshape(}and hence in $\CrP${\upshape )} if and only if they \emph{have the same action on their sets of fixed points} 
{\upshape (}see Definition \refd{Def:actionFixedPoints}\ {\upshape )}. 
Explicitly, the curve $y_1^2 \prod_{i=1}^{k/4}P(a_i)= y_2^2 \prod_{i=1}^{k/4}P(b_i)$ and the action of $(x_1:x_2)\ \mapsto (x_1: -x_2)$ and $(x_1:x_2)\ \mapsto (x_2: x_1)$ on it determine the conjugacy class of the group.
\end{Prp}
\begin{proof}\fauxtitred\upshape
Conjugating $G$ in $\CrPp$, we may suppose that $\overline{\pi}(G)$ is generated by $(x_1:x_2) \mapsto (x_1:- x_2)$ and $(x_1:x_2) \mapsto (x_2: x_1)$. We then conjugate $G$ by some element of $\PGL(2,\K(x))$ and suppose that $G'$ is generated by a de Jonqui\`eres involution of the standard form given in Proposition \refd{Prp:classDeJIJonqG}, assertion $1$. We then apply Lemma \refd{Lem:splitsequence}, which implies that, up to conjugation in the de Jonqui\`eres group, $G$ is generated by the automorphisms $\rho_1,\rho_2$ as written above, and the birational map
\begin{center}$\sigma_{\alpha,\beta}:(x_1:x_2)\times (y_1:y_2) \dasharrow (x_1:x_2) \times(y_2 \prod_{i=1}^k(x_1-\beta_i x_2):y_1 \prod_{i=1}^k(x_1-\alpha_i x_2))$\end{center}of $\mathbb{P}^1\times\mathbb{P}^1$, for some $\alpha_1,...,\alpha_k,\beta_1,...,\beta_k \in \K$ all distinct. 

As $\rho_1$ commutes with $\sigma_{\alpha,\beta}$, we find that the sets $\{\alpha_1,...,\alpha_k\}$ and $\{\beta_1,...,\beta_k\}$ are invariant by multiplication by $-1$. With a little more calculation, we find that both sets are also invariant by the action $x \mapsto x^{-1}$ and that neither $0,1$ nor $-1$ belongs to the sets.
We thus get the announced result.

The proof of the last assertion is the same as in the proof of Proposition \refd{Prp:DirectProduct}. 
Let $G$ (respectively $G'$) be the group generated by $\sigma_{a,b}$ (respectively $\sigma_{a',b'}$) and $\rho_1,\rho_2$. Let $\Gamma$ (respectively $\Gamma'$) denote the curve of fixed points of $\sigma_{a,b}$ (respectively $\sigma_{a',b'}$). We suppose that $G$ and $G'$ \emph{have the same action on their sets of fixed points}. This means that there exists an isomorphism $\varphi:\Gamma' \rightarrow \Gamma$ such that $<\rho_1,\rho_2>\varphi=\varphi <\rho_1,\rho_2>$.

As $k>1$, $\varphi$ induces an automorphism $\overline{\varphi}$ of $\mathbb{P}^1$. The element $\overline{\varphi}\in \PGL(2,\K)$ normalizes the group isomorphic to $(\Z{2})^2$ generated by $(x_1:x_2) \mapsto (-x_1:x_2)$ and $(x_1:x_2) \mapsto (x_2:x_1)$. Note that the automorphism $\tilde{\varphi}:(x_1:x_2) \times (y_1:y_2) \mapsto \overline{\varphi}(x_1:x_2)\times (y_1:y_2)$ of $\mathbb{P}^1\times\mathbb{P}^1$ sends the ramification points of $\Gamma'$ on those of $\Gamma$ and does not change the standard form: we conjugate $G'$ by $\tilde{\varphi}$ and get a group $G''$ generated by $\rho_1,\rho_2$ and some $\sigma_{a'',b''}$, with $\{a'',b''\}=\{a,b\}$.

Let $\psi$ denote the birational map
\begin{center}
$\psi:(x_1:x_2)\times (y_1:y_2) \dasharrow (x_1:x_2)\times (y_1 \prod_{i=1}^k(x_1-a_i x_2) :y_2 \prod_{i=1}^k(x_1-a_i'' x_2))$\end{center}
of $\mathbb{P}^1\times\mathbb{P}^1$, which commutes with $\rho_1$ and $\rho_2$. Note that the two birational maps
\begin{center}
$\begin{array}{rcl}
\psi \circ \sigma_{a,b}(x,y)& = &(x_1:x_2)\times (y_2 \prod_{i=1}^k(x_1-b_i x_2):y_1\prod_{i=1}^k(x_1-a''_i x_2) ), \vspace{0.1 cm}\\
\sigma_{a'',b''}\circ \psi(x,y)& = &(x_1:x_2)\times (y_2 \prod_{i=1}^k(x_1-b''_i x_2):y_1\prod_{i=1}^k(x_1-a_i x_2) ),\end{array}$\end{center}
are equal, since $\{a'',b''\}=\{a,b\}$. The groups $G$ and $G'$ are then conjugate in $\CrPp$ by the birational map $\psi$.\proofend
\end{proof}

\bigskip

We now consider the case when the exact sequence (\refd{eq:ExactSeqCB}) does not split.
\begin{Prp}\fauxtitre
\label{Prp:DoesnotsplitC2-C2C2}
Let $G \subset \CrPp $ be a finite abelian subgroup of the de Jonqui\`eres group such that
\begin{itemize}
\item
The exact sequence (\refd{eq:ExactSeqCB}) does not split.
\item
$\overline{\pi}(G)\cong (\Z{2})^2$.
\item
$G'\cong \Z{2}$ is generated by a de Jonqui\`eres involution, whose set of fixed points is a hyperelliptic curve of genus $k-1$.\index{Curves!hyperelliptic!birational maps that fix a hyperelliptic curve}
\end{itemize}
Then, $k\geq 2$ is divisible by $2$, $G\cong (\Z{4}\times\Z{2})$ and up to conjugation in the de Jonqui\`eres group,  $G$ is generated by the following commuting elements $\alpha,\rho \in \CrPp$ of order respectively $4$ and $2$: 
\begin{center}
$\begin{array}{rrclcccc}
{\num{C.24}} & \alpha:(x,y)& \dasharrow &(\h{3}&-x&\h{3},\h{3}&\nu(x)\cdot \frac{y-\nu(-x)}{y-\nu(x)}&\h{4}),\vspace{0.1cm}\\
& \alpha^2:(x,y)& \dasharrow &(\h{3}&x&\h{3},\h{3}&\frac{\nu(x)\nu(-x)}{y}&\h{4}),\vspace{0.1cm}\\
&\rho:(x,y)&\mapsto &(\h{3}&x^{-1}&\h{3},\h{3}&y&\h{4}),\end{array}$\end{center}
where $\nu(x) \in \K(x)^{*}$ is invariant by the action of $x \mapsto x^{-1}$, $\nu(x)\not=\nu(-x)$.
\end{Prp}
\begin{proof}\fauxtitred\upshape
As $G$ is abelian and the exact sequence does not split, we have $G \cong \Z{4}\times\Z{2}$, generated by some elements $\alpha,\rho$ of order respectively $4$ and $2$ such that:
\begin{itemize}
\item $\overline{\pi}(\alpha),\overline{\pi}(\rho)$ are two commuting involutions of $\PGL(2,\K)$;
\item $\alpha^2$ is the generator of $G'$.
\end{itemize}
Up to conjugation in the de Jonqui\`eres group, we may suppose that $\alpha^2(x,y)=(x,g(x)/y)$ (see Section \refd{Sec:DeJonqInv}) and that the action of $\rho$ on the fibration is $x \mapsto x^{-1}$.
 
We now denote by $H$ the group isomorphic to $(\Z{2})^2$ generated by $\alpha^2$ and $\rho$. Note that $\overline{\pi}(H)=<\overline{\pi}(\rho)>\cong \Z{2}$, $H'=<\alpha^2> \cong \Z{2}$ and the exact sequence for the group $H$ splits. Using Proposition \refd{Prp:DirectProduct}, we see that $k\geq 2$ is divisible by $2$. We then use Lemma \refd{Lem:splitsequence} and may suppose, up to conjugation in the de Jonqui\`eres group, that $\rho(x,y)=(x^{-1},y)$. As $\rho$ commutes with $\alpha^2$, we see that $g(x^{-1})=g(x)$. 

Up to conjugation by some element of $\PGLn{2}$ acting on $\K(x)$ (which changes $g$ by a change of variable), we may suppose that the action of $\alpha$ on the fibration is $x\mapsto -x$.

We proceed exactly as in the proof of Proposition \refd{Prp:SquareRoots}. We first write $$a=\Bigg(\h{1}\left(\begin{array}{cc}\alpha & \beta \\ \gamma& \delta\end{array}\right),\left(\begin{array}{cc}-1 & 0 \\ 0& 1\end{array}\right)\h{1}\Bigg),$$ for some elements $\alpha,\beta,\gamma,\delta \in \K(x)\backslash\{0\}$.
Note that here, since $a$ commutes with $\rho$, we may suppose that $\alpha,\beta,\gamma,\delta \in L$, where $L=\K(x+x^{-1})$ is the subfield of $\K(x)$ invariant by $x\mapsto x^{-1}$.

After some calculations (see the proof of Proposition \refd{Prp:SquareRoots}), we obtain \begin{center}$\alpha(x)\alpha(-x)+\beta(x)\gamma(-x)=\delta(x)\delta(-x)+\beta(-x)\gamma(x)=0$,\end{center}
and see that $\alpha(x)\alpha(-x)=\delta(x)\delta(-x)$. Writing this as $(-\frac{\delta(x)}{\alpha(x)})(-\frac{\delta(-x)}{\alpha(-x)})=1$ and using Hilbert $90$ (\cite{bib:Serre}, page 158), there exists an element $\lambda(x) \in L\backslash\{0\}$ such that $\frac{\lambda(x)}{\lambda(-x)}=-\frac{\delta(x)}{\alpha(x)}$.

We conjugate the group $G$ by the element $\Bigg(\h{1}\left(\begin{array}{cc}\lambda & 0 \\ 0& 1\end{array}\right),\left(\begin{array}{cc}1 & 0 \\ 0& 1\end{array}\right)\h{1}\Bigg)$ and may assume that 
\begin{center}$a=\Bigg(\h{1}\left(\begin{array}{cc}\nu(x)& -\nu(x)\nu(-x) \\ 1 &-\nu(x)\end{array}\right),\left(\begin{array}{cc}-1 & 0 \\ 0& 1\end{array}\right)\h{1}\Bigg)$,\end{center}
where $\nu=\frac{\alpha\lambda}{\gamma}$ (see Proposition \refd{Prp:SquareRoots} again). We obtain the announced result. The fact that $\nu(x)\not=\nu(-x)$ follows from the hypothesis that the set of fixed points fixed by $\alpha^2$ is a hyperelliptic curve.
\proofend
\end{proof}

\section{Groups $(\Z{2})^2 \subset \PGL(2,\K(x))$}
\pagetitless{Finite abelian groups of automorphisms of conic bundles}{Groups $(\Z{2})^2 \subset \PGL(2,\K(x))$}
\label{Sec:GroupsC2C2inPGL}
In this section, we study the groups generated by two commuting de Jonqui\`eres involutions, having a trivial action on the fibration. This amounts to studying the subgroups of $ \PGL(2,\K(x))$ isomorphic to $(\Z{2})^2$.

The following result is proved in \cite{bib:Be2}, Lemma 2.1:
\begin{Prp}\fauxtitred
\label{Prp:C2C2PGL2Ct}
\begin{itemize}
\item
If $G \subset \PGL(2,\K(x))$ is a group isomorphic to $(\Z{2})^2$, then up to conjugation in this group, $G$ is generated by 
\begin{center}$\left(\begin{array}{cc}0 & g(x) \\ 1 & 0\end{array}\right)$\hspace{1 cm} and \hspace{1 cm}$\left(\begin{array}{cc}h(x) & -g(x) \\ 1 & -h(x)\end{array}\right)$\end{center} for some $g(x),h(x)\in\K(x)\backslash\{0\}$.
\item
The morphism $\delta: \PGL(2,\K(x)) \rightarrow \K(x)^{*}/\K(x)^{*2}$ given by the determinant induces a bijection between the set of conjugacy classes of subgroups of $\PGL(2,\K(x))$ isomorphic to $(\Z{2})^2$ and the set of subgroups of order $\leq 4$ of $\K(x)^{*}/\K(x)^{*2}$.\proofend
\end{itemize}
\end{Prp}

\bigskip

We now give some information on the conjugacy classes in the de Jonqui\`eres group:
\begin{Prp}\fauxtitre
\label{Prp:C2C2PGL2CtinJonq}
The morphism $\delta: \PGL(2,\K(x)) \rightarrow \K(x)^{*}/\K(x)^{*2}$ given by the determinant induces a bijection between the set of conjugacy classes of subgroups of $\PGL(2,\K(x))$  isomorphic to $(\Z{2})^2$ in the group $\PGL(2,\K(x))\rtimes \PGL(2,\K)$ and the equivalence classes of subgroups of order $\leq 4$ of $\K(x)^{*}/\K(x)^{*2}$, up to the action of $\PGL(2,\K)$ on $\K(x)$.
\end{Prp}
\begin{proof}\fauxtitred\upshape
This follows from the fact that if we take $\alpha \in \PGLn{2}$, $g \in \mathrm{\PGL}(2,\K(t))$, then $\delta(\alpha g \alpha^{-1})=\alpha (\delta(g))$.\proofend

\end{proof}

\bigskip

Note that, by Proposition \refd{Prp:BeauTsen}, a group $G$ of this form is conjugate to a group $\tilde{G}$ of automorphisms of a conic bundle $(S,\pi)$, generated by two involutions whose actions on the basis are trivial. Proposition \refd{Prp:TwistingOrNot}, tells us that each involution $\sigma \in \PGL(2,\K(x))$ corresponds to a twisting $(S,\pi)$-de Jonqui\`eres involution if and only if $\delta(\sigma)\not=1$. Then, one (and only one) of the following occurs:
\begin{center}
\begin{tabular}{ll}
$\delta(G)=1$ & $\tilde{G}$ contains no twisting $(S,\pi)$-de Jonqui\`eres involution;\\
$\delta(G)\cong \Z{2}$ & $\tilde{G}$ contains exactly two twisting $(S,\pi)$-de Jonqui\`eres involutions;\\
$\delta(G)\cong (\Z{2})^2$& $\tilde{G}$ contains exactly three twisting $(S,\pi)$-de Jonqui\`eres involutions.
\end{tabular}
\end{center}

Let us now look at the conjugacy classes of these groups in the Cremona group:

\begin{Prp}\fauxtitre
\label{Prp:C2C2PGL2CtinCremona}
Let $G\subset \PGL(2,\K(x))$ be a group isomorphic to $(\Z{2})^2$.
\begin{itemize}
\item
\index{Curves!of positive genus!birational maps that do not fix a curve of positive genus}
If no element of $G$ fixes a curve of positive genus, then $G$ is birationally conjugate to a subgroup of $\Aut(\mathbb{P}^1\times\mathbb{P}^1,\pi_1)$. More precisely, $G$ is birationally conjugate either to $(x,y)\mapsto (\pm x^{\pm 1},y)$ or to $(x,y) \mapsto (\pm x,\pm y)$.
\item
\index{Curves!of positive genus!birational maps that fix a curve of positive genus}
If exactly one element of $G$ fixes a curve of positive genus, this curve is elliptic and the group is \index{Curves!elliptic!birational maps that fix an elliptic curve}
birationally conjugate to a subgroup $H$ of the de Jonqui\`eres group that satisfies
\begin{center} $\overline{\pi}(H)\cong H'\cong \Z{2}$.\end{center} In particular, $H$ does not act trivially on the fibration.
\item
If at least two elements of $G$ each fix a curve of positive genus, then any subgroup of $\PGL(2,\K(t))\rtimes \PGLn{2}$ which is birationally conjugate to $G$ belongs to $\PGL(2,\K(t))$ (i.e.\ acts trivially on the fibration).
\end{itemize}
\end{Prp}
\begin{proof}\fauxtitred\upshape
Let us prove the first assertion. As mentioned above, the group $G$ is conjugate to a group $\tilde{G}\subset \Aut(S,\pi)$ of automorphisms of some conic bundle $(S,\pi)$, generated by two de Jonqui\`eres involutions acting trivially on the fibration. The hypothesis in the first assertion imposes that each involution twists either $0$ or $2$ fibres of the fibration. Suppose that the triple $(\tilde{G},S,\pi)$ is minimal; then one of the following situations occurs:
\begin{itemize}
\item
The fibration is smooth, no element of $\tilde{G}$ is twisting, $\delta(G)=1$.
\item
There are two singular fibres, twisted by two of the involutions of $\tilde{G}$, $\delta(G)\cong \Z{2}$.
\item
There are three singular fibres, each involution twists two of these, $\delta(G)\cong(\Z{2})^2$.
\end{itemize}
In any case, the number of singular fibres is at most $3$, so the group is birationally conjugate to a subgroup of automorphisms of $\mathbb{P}^1\times\mathbb{P}^1$ (see Proposition \refd{Prp:LargeDegree}). We use the classification of groups of this kind (see Proposition \refd{Prp:AutP1Bir}) to see that there are only two possibilities, up to conjugation, namely $(x,y)\mapsto (\pm x^{\pm 1},y)$ or $(x,y) \mapsto (\pm x,\pm y)$.

To prove the second assertion, we need properties of Del Pezzo surfaces of degree $4$\index{Del Pezzo surfaces!of degree $4$|idi}. This uses Lemma \refd{Lem:EllRatRat} and will be proved in Chapter \refd{Chap:ConjBetweenTheCases}. 

Let us prove the last assertion. Note that an element of the de Jonqui\`eres group that does not act trivially on the fibration does not fix a curve of positive genus. Then, if $H\subset \PGL(2,\K(x))\rtimes \PGL(2,\K)$ is birationally conjugate to $G$, it contains two involutions that each fix a curve of positive genus, and hence that both act trivially on the fibration. The group generated by these two involutions, which is $H$, belongs then to $\PGL(2,\K(x))$.
\proofend
\end{proof}

\subsection{Extension of some group of $\PGL(2,\K)$ by two de Jonqui\`eres involutions}
\label{Sec:ExtensionC2C2}
Let us study the finite abelian groups $G$ of automorphisms of a conic bundle such that $G'$ is generated by two twisting de Jonqui\`eres involution and $\overline{\pi}(G)$ is not trivial (i.e.\ is either cyclic of order $n>1$ or isomorphic to $(\Z{2})^2$). 
Using the exact sequence (introduced in Section~\refd{Sec:Exactsequence})
\begin{equation}
1 \rightarrow G' \rightarrow G \stackrel{\overline{\pi}}{\rightarrow} \overline{\pi}(G) \rightarrow 1
\tag{\refd{eq:ExactSeqCB}}\end{equation}
we see that one of the following occurs:
\begin{center}\begin{tabular}{ll}
$\overline{\pi}(G)\cong \Z{n}$: &  \h{3}$G$ is isomorphic either to $(\Z{2})^2\h{0.5}\times\h{0.5}\Z{n}$ or to $\Z{2}\h{0.5}\times\h{0.5}\Z{2n}$,\\
$\overline{\pi}(G)\cong(\Z{2})^2$: & \h{3}$G$ is isomorphic either to $(\Z{2})^4$, to $(\Z{2})^2\h{0.5}\times\h{0.5} \Z{4}$ or to $\Z{4}\h{0.5}\times\h{0.5}\Z{4}$.\end{tabular}\end{center}
 As in the previous sections, this depends on whether or not the exact sequence splits. We describe the cases in the following two propositions:

 \begin{Prp}\fauxtitre
 \label{Prp:SplitC2C2}
 Let $G \subset \CrPp$ be a finite abelian subgroup of the de Jonqui\`eres group such that $G'\cong (\Z{2})^2$, $\overline{\pi}(G)\not=1$ and the exact sequence (\refd{eq:ExactSeqCB}) splits.
 
Then, up to conjugation in the de Jonqui\`eres group, $G$ is one of the following:
 \begin{itemize}
\item
 $\overline{\pi}(G)\cong \Z{n}$, $n>1$:  
 ${G \cong (\Z{2})^2\times \Z{n}}$ is generated by
\begin{center}
$\begin{array}{rrclcccc}
\num{C.22,\h{0.5}n}&(x,y)& \h{3}\dasharrow\h{3} &(&\h{3}x&\h{3},\h{3}&g(x)/y&\h{4}),\\
&(x,y)& \h{3}\dasharrow\h{3} &(&\h{3}x&\h{3},\h{3}&(h(x)y-g(x))/(y-h(x))&\h{4}),\\
&(x,y)& \h{3}\dasharrow\h{3} &(&\h{3}\zeta_n x&\h{3},\h{3}&y&\h{4}),\end{array}$\end{center}
for some $g(x),h(x) \in \K(x) \backslash \{0\}$ invariant by $\zeta_n$.  The genus of the hyperelliptic curve fixed by any involution of $G'$ is equal to $-1 \pmod{n}$.\index{Curves!hyperelliptic!birational maps that fix a hyperelliptic curve}
\item
 $\overline{\pi}(G)\cong (\Z{2})^2$:  
 ${G \cong (\Z{2})^4}$ is generated by
\begin{center}
$\begin{array}{rrclcccc}
\num{C.22,\h{0.5}22}&(x,y)& \h{3}\dasharrow\h{3} &(&\h{3}x&\h{3},\h{3}&g(x)/y&\h{4}),\\
&(x,y)& \h{3}\dasharrow\h{3} &(&\h{3}x&\h{3},\h{3}&(h(x)y-g(x))/(y-h(x))&\h{4}),\\
&(x,y)& \h{3}\dasharrow\h{3} &(&\h{3}\pm x^{\pm 1}&\h{3},\h{3}&y&\h{4}),\end{array}$\end{center}
for some $g(x),h(x) \in \K(x) \backslash \{0\}$ invariant by $x \mapsto \pm x^{\pm 1}$. The genus of the hyperelliptic curve fixed by any involution of $G'$ is equal to $3 \pmod{4}$.\index{Curves!hyperelliptic!birational maps that fix a hyperelliptic curve}
\end{itemize}
\end{Prp}
 \begin{remark}\fauxtitred
 The groups called here $\nump{C.22,\h{0.5}22}$ have already been described in \cite{bib:Be2}. In fact we use the method of this article here also for the other groups.
 \end{remark}
 \begin{proof}\fauxtitred\upshape
Conjugating $G$ in $\CrPp$, we may suppose that $\overline{\pi}(G)$ is generated either by $x\mapsto \zeta_n x$ or by $x\mapsto \pm x^{\pm 1}$.

We then conjugate $G$ by some element of $\PGL(2,\K(x))$ and suppose that $G'$ is of the standard form given in Proposition \refd{Prp:C2C2PGL2Ct}. We apply Lemma \refd{Lem:splitsequence}, which implies that up to conjugation in the de Jonqui\`eres group, $G$ is of one of the two forms given above. The invariance by the actions of $\overline{\pi}(G)$ follows from the fact that $G$ is abelian.

The relations on the genus of the fixed curves follow from Propositions \refd{Prp:DirectProduct} and \refd{Prp:DirectProductC2C2C2}.
\proofend
 \end{proof}

\bigskip

 \begin{Prp}\fauxtitre
 \label{Prp:NotSplitC2C2}
 Let $G \subset \CrPp$ be a finite abelian subgroup of the de Jonqui\`eres group such that $G'\cong (\Z{2})^2$ and the exact sequence (\refd{eq:ExactSeqCB}) does not split.
 
Then, up to conjugation in the de Jonqui\`eres group, $G$ is one of the following:
 \begin{itemize}
\item
 $\overline{\pi}(G)\cong \Z{2n}$, $n>1${\upshape :} 
 ${G \cong \Z{4n}\times \Z{2}}$ is generated by elements $\alpha$ and $\beta$ of order respectively $4n$ and $2${\upshape :}
\begin{center}
$\begin{array}{rrrclcccc}
$\num{C.2n2}$&\alpha:&(x,y)& \h{3}\dasharrow\h{3} &(&\h{3}\zeta_n x&\h{3},\h{3}&(a(x)y+b(x))/(c(x)y+d(x))&\h{4}),\\
&\alpha^n:&(x,y)& \h{3}\dasharrow\h{3} &(&\h{3}x&\h{3},\h{3}&g(x)/y&\h{4}),\\
&\beta:& (x,y)& \h{3}\dasharrow\h{3} &(&\h{3}x&\h{3},\h{3}&(h(x)y-g(x))/(y-h(x))&\h{4}),\end{array}$\end{center}
for some $a,b,c,d,g,h \in \K(x) \backslash \{0\}$ such that $\alpha$ and $\beta$ commute. 
\item
 $\overline{\pi}(G)\cong (\Z{2})^2${\upshape :}  
 ${G \cong (\Z{4})\times (\Z{2})^2}$ is generated by elements $\alpha$, $\beta$ and $\gamma$ of order respectively $4$, $2$ and $2${\upshape :}
\begin{center}
$\begin{array}{rrrclcccc}
$\num{C.422}$&\alpha:&(x,y)& \h{3}\dasharrow\h{3} &(&\h{3}x^{-1}&\h{3},\h{3}&(a(x)y+b(x))/(c(x)y+d(x))&\h{4}),\\
&\alpha^2:&(x,y)& \h{3}\dasharrow\h{3} &(&\h{3}x&\h{3},\h{3}&g(x)/y&\h{4}),\\
&\beta:& (x,y)& \h{3}\dasharrow\h{3} &(&\h{3}x&\h{3},\h{3}&(h(x)y-g(x))/(y-h(x))&\h{4}),\\
&\gamma:& (x,y)& \h{3}\dasharrow\h{3} &(&\h{3}-x&\h{3},\h{3}&y&\h{4}),\end{array}$\end{center}
for some $a,b,c,d,g,h \in \K(x) \backslash \{0\}$ such that $\alpha$ and $\beta$ commute and $g,h$ invariant by $x \mapsto -x$.
\item
 $\overline{\pi}(G)\cong (\Z{2})^2${\upshape :}  
 ${G \cong (\Z{4})^2}$ is generated by elements $\alpha$, $\beta$ of order $4${\upshape :}
\begin{center}
$\begin{array}{rrrclcccc}
$\num{C.44}$&\alpha:&(x,y)& \h{3}\dasharrow\h{3} &(&\h{3}-x&\h{3},\h{3}&(a(x)y+b(x))/(c(x)y+d(x))&\h{4}),\\
&\alpha^2:&(x,y)& \h{3}\dasharrow\h{3} &(&\h{3}x&\h{3},\h{3}&g(x)/y&\h{4}),\\
&\beta:& (x,y)& \h{3}\dasharrow\h{3} &(&\h{3}x^{-1}&\h{3},\h{3}&(p(x)y+q(x))/(r(x)y+s(x))&\h{4}),\\
&\beta^2:& (x,y)& \h{3}\dasharrow\h{3} &(&\h{3}x&\h{3},\h{3}&(h(x)y-g(x))/(y-h(x))&\h{4}),\end{array}$\end{center}
for some $a,b,c,d,g,h,p,q,r,s \in \K(x) \backslash \{0\}$ such that $\alpha$ and $\beta$ commute.
\end{itemize}
\end{Prp}
 \begin{proof}\fauxtitred\upshape
 Since the exact sequence does not split, $\overline{\pi}(G)\not=1$ and at least one root of even order of a de Jonqui\`eres involution of $G'$ belongs to $G$.
 
 If $\overline{\pi}(G)$ is cyclic, then $\overline{\pi}(G)\cong \Z{2n}$, for some integer $n>1$, and one such root $\alpha$ (and its inverse) belongs to $G$. Moreover, $\overline{\pi}(\alpha)$ generates $\overline{\pi}(G)$. Up to conjugation in the de Jonqui\`eres group, we suppose that $\overline{\pi}(\alpha)$ is $x\mapsto \zeta_{2n} x$. We then conjugate $G$ by some element of $\PGL(2,\K(x))$ and may assume that $G'$ is of the standard form given in Proposition \refd{Prp:C2C2PGL2Ct}.  We obtain $\nump{C.2n2}$. \vspace{0.1 cm}

 If $\overline{\pi}(G)\cong (\Z{2})^2$, we first conjugate the group $G$ in the de Jonqui\`eres group, so that $\overline{\pi}(G)$ is the group $x\mapsto \pm x^{\pm 1}$. Then, we conjugate $G$ by some element of $\PGL(2,\K(x))$ and may assume that $G'$ is of the standard form given in Proposition \refd{Prp:C2C2PGL2Ct}. Two cases are possible:
 \begin{itemize}
 \item
\emph{Every involution of $G$ acts trivially on the fibration.}\\
 This implies that the pull-back by $\overline{\pi}$ of every group of order $2$ of $(\Z{2})^2$ is a group isomorphic to $\Z{4} \times \Z{2}$ of the form above (for $n=2$).  The group $G$ is then isomorphic to $(\Z{4})^2$, and generated by two elements $\alpha,\beta$, such that $\alpha^2,\beta^2$ are commuting involutions that generate $G'$. We obtain $\nump{C.44}$.
 \item
\emph{One of the involutions of $G$ acts non-trivially on the fibration.}\\
We denote this involution by $\gamma$ and the group $\overline{\pi}^{-1}(\overline{\pi}(\gamma))$ by $H$. Then, the exact sequence associated to $H$ splits. We apply Proposition \refd{Prp:SplitC2C2}, which implies that up to conjugation in the de Jonqui\`eres group, $H \cong (\Z{2})^3$ is generated by these three involutions:
\begin{center}
$\begin{array}{rrclcccc}
\rho:&(x,y)& \h{3}\dasharrow\h{3} &(&\h{3}x&\h{3},\h{3}&g(x)/y&\h{4}),\\
\beta:&(x,y)& \h{3}\dasharrow\h{3} &(&\h{3}x&\h{3},\h{3}&(h(x)y-g(x))/(y-h(x))&\h{4}),\\
\gamma:& (x,y)& \h{3}\dasharrow\h{3} &(&\h{3}-x&\h{3},\h{3}&y&\h{4}),\end{array}$\end{center}

for some $g(x),h(x) \in \K(x) \backslash \{0\}$ invariant by $x\mapsto -x$. 

Furthermore, $G$ contains some element $\alpha$ of order $4$, such that $\alpha^2 \in G'$ and $\overline{\pi}(\alpha)$ is $x \mapsto x^{-1}$. Without loss of generality, we may suppose that $\alpha^2=\rho$ and obtain $\nump{C.422}$.
 \end{itemize}\proofend
 \end{proof}
 
 \bigskip

\begin{Exa}\fauxtitred\upshape
\label{Exa:C4C2ellratdoesnotsplit}
Let $G\cong \Z{4}\times\Z{2}$ be the group generated by $\alpha$, $\beta$, elements of order respectively $4$ and $2$:
\begin{center}
$\begin{array}{llcccl}
\alpha&:(x,y) \dasharrow (\h{3}&-x&\h{3},\h{3}&(1+\im x)\frac{y-(1-\im x)}{y-(1+\im x)}&\h{4})\\
\alpha^2&:(x,y) \dasharrow (\h{3}&x&\h{3},\h{3}&\frac{x^2+1}{y}&\h{4})\\
\beta&:(x,y) \dasharrow (\h{3}&x&\h{3},\h{3}&\frac{y(1+x)-(1+x^2)}{y-(1+x)}&\h{4})\\
\alpha^4=\beta^2&:(x,y) \dasharrow (\h{3}&x&\h{3},\h{3}&y&\h{4})
\end{array}$
\end{center}
The elements $\alpha^2$ and $\beta$ are de Jonqui\`eres involutions and the action of $\alpha$ on the fibration is cyclic of order $2$.

Both $\alpha^2$ and $\beta$ fix a rational curve (ramified over respectively $\{\pm \im\}$ and $\{0,\infty\}$), and $\alpha^2\beta$ fixes an elliptic curve (ramified over $\{0,\pm \im, \infty\}$).\index{Curves!elliptic!birational maps that fix an elliptic curve}
\end{Exa}
\begin{remark}\fauxtitred
We do not know whether a group of the form $\nump{C.44}$ exists.\end{remark}

\bigskip

Let us conclude  this long and difficult chapter by mentioning its most surprising result. This is clearly the existence of roots of de Jonqui\`eres involutions of any order (see Proposition \refd{Prp:ExistRoot}). This result (and most of the others) do not appear anywhere in the literature.

 One question now arises. Taking two roots of the same de Jonqui\`eres involution, that have the same action on the hyperelliptic curve fixed by it, are the roots (or at least the groups generated by these roots) birationally conjugate?\index{Curves!hyperelliptic!birational maps that fix a hyperelliptic curve}
\chapter{Finite abelian groups of automorphisms of Del Pezzo surfaces of degree $\leq 4$}\index{Del Pezzo surfaces!automorphisms|idb}
\label{Chap:RankOne}
\ChapterBegin{We describe in this chapter the pairs $(G,S)$, where $S$ is a rational surface and $G \subset \Aut(S)$ is a finite abelian group such that $\rkPic{S}^G=1$ (case $2$ of Proposition \refd{Prp:TwoCases}). As was proved in Lemma \refd{Lem:rkDelPezzo}, the surface $S$ is of Del Pezzo type. Since the finite abelian subgroups of automorphisms of surfaces of degree $\geq 5$ were described in Chapter \refd{Chap:LargeDegreeSurfaces} (see in particular Proposition 
\refd{Prp:LargeDegree}), the remaining cases are Del Pezzo surfaces of degree $1$,$2$,$3$ or $4$.} 

\section{Del Pezzo surfaces of degree $4$}\index{Del Pezzo surfaces!of degree $4$|idb}
\pagetitless{Finite abelian groups of automorphisms of Del Pezzo surfaces of degree $\leq 4$}{Del Pezzo surfaces of degree $4$}
\label{Sec:DP4}
\subsection{The structure of the group of automorphisms}
There are many different isomorphism classes of Del Pezzo surfaces of degree $4$. These surfaces canonically embed in $\mathbb{P}^4$ as the intersection of two quadrics (see for example \cite{bib:Be1}, Proposition IV.16). 

The literature on these surfaces and their group of automorphisms is dense. A. Wiman \cite{bib:Wim} gave the groups of automorphisms and the equations (that may be found also before in the works of S. Kantor and B. Segre).
It was known that the automorphism group of every quartic Del Pezzo surface contains a normal subgroup isomorphic to $(\Z{2})^4$ (a result proved again in \cite{bib:Koi} and \cite{bib:Ho1}). 
Recently, A. Beauville  \cite{bib:Be2} proved that the conjugacy class of this group in $\CrP$ determines the surface up to isomorphism.

We will go a little farther, following another approach, and use extensively the structure of the Picard group, to make explicit the group of automorphisms of every such surface in a natural way, using the conic bundle structures of the surface. Furthermore, we will consider not only the entire group of automorphisms but also in its abelian subgroups.

Let $S$ be a Del Pezzo surface of degree $4$, obtained by the blow-up of five points of\hspace{0.2 cm}$\Pn$ in general position and let $\pi: S \rightarrow \Pn$ be the blow-down. With a suitable choice of coordinates of\hspace{0.2 cm}$\Pn$, we may assume that the five points are $A_1=(1:0:0)$, $A_2=(0:1:0)$, $A_3=(0:0:1)$, $A_4=(1:1:1)$ and $A_5=(a:b:c)$, a point not aligned with any two of the others.
Note that the linear equivalence of the set $\{A_1,...,A_5\} \subset \Pn$ determines the isomorphism class of $S$. This follows from Proposition \refd{Prp:IsoDPlinear}.

The exceptional divisors on $S$ are (see Proposition 
\refd{Prp:DescExcCurvDelPezzo}):\index{Del Pezzo surfaces!curves on the surfaces|idi}
\begin{itemize}
\item
$E_1=\pi^{-1}(A_1),...,E_5=\pi^{-1}(A_5)$, the $5$ pull-backs of the points $A_1$,...,$A_5$.
\item
$D_{ij}=\tilde{\pi^{-1}}(A_iA_j)$ for $1\leq i,j\leq 5$, $i\not=j$, the $10$ strict pull-backs of the lines through $2$ of the $A_i$'s.
\item
$C=\tilde{\pi^{-1}}(\mathcal{C})$, the strict pull-back of the conic $\mathcal{C}$ of\hspace{0.2 cm}$\Pn$ through $A_1,A_2,...,A_5$.
\end{itemize}
There are thus $16$ exceptional divisors, and each intersects five others:
\begin{itemize}
\item
Each $E_i$ intersects $C$ and the $4$ strict transforms of lines through $A_i$.
\item
The divisor $C$ intersects $E_1,E_2,...,E_5$.
\item
Each divisor $D_{ij}$ intersects $E_i$ and $E_j$ and the transforms of the three lines through two of the other three points.
\end{itemize}
Let $L$ be the pull-back by $\pi$ of a line of\hspace{0.2 cm}$\Pn$. The Picard group of $S$ is generated by $E_1,E_2,E_3,E_4,E_5$ and $L$.
\begin{Def} \fauxtitre
A pair $\{A,B\}$ of divisors of $S$ is an \defn{exceptional pair} if $A$ and $B$ are divisors, $A^2=B^2=0$ and $(A+B)=-K_S$. This implies that $AB=2$ and $AK_S=BK_S=-2$.\end{Def} 
\begin{remark}\fauxtitred
In fact, the divisors of an exceptional pair are effective (see below). The conditions $A^2=0$ and $A(A+K_S)=-2$ are equivalent to saying that $A$ is a base-point-free pencil of rational curves, and so the same holds for $-K_S-A$. Each exceptional pair thus corresponds to two conic bundle structures, as no curve of self-intersection $<-1$ belongs to the surface (see Proposition 
\index{Conic bundles!on Del Pezzo surfaces}
\refd{Prp:NbConicDP}).
\end{remark}
\begin{Lem}\fauxtitred
\begin{itemize}
\item
There are $5$ exceptional pairs on the surface $S$, namely
$\{L-E_i,2L-\sum_{j\not= i} E_j\}=\{L-E_i,-K_S-L+E_i\}$, for $i=1,...,5$.
\item
Geometrically, these correspond to the lines through one of the five blown-up points, and the conics through the other four points.
\end{itemize}
\end{Lem}\begin{proof}\fauxtitred\upshape
Let $\{A,B\}$ denote an exceptional pair, where $A=mL-\sum_{i=1}^5 a_i E_i$ and $B=nL-\sum_{i=1}^5 b_i E_i$, for some $m,n,a_1,...,a_5,b_1,...,b_5 \in \mathbb{Z}$. From the relations $A^2=0$ and $AK_S=-2$ we get:
\begin{equation}\label{eqexcpair}\begin{array}{lll}
\sum_{i=1}^5 {a_i}^2&=&m^2, \vspace{0.1 cm}\\
\sum_{i=1}^5 a_i&=&3m-2.\end{array}
\end{equation} 
Note that $(\sum_{i=1}^5 a_i)^2\leq 5\sum_{i=1}^5 {a_i}^2$ (see Lemma \refd{Lem:Tecaicarr} in the Appendix), which implies here that $(3m-2)^2\leq 5m^2$, that is $4(m^2-3m+1)\leq 0$. As $m$ is an integer, we get $1\leq m\leq 2$.

From $A+B=-K_S$, we get $m+n=3$. Up to permutation of $A$ and $B$, we may suppose that $m=1$ and $n=2$. Replacing $m=1$ in (\refd{eqexcpair}), we see that there exists $i \in \{1,...,5\}$ such that $A=L-E_i$. The pair $\{A,B\}$ must then be equal to $\{L-E_i,2L-\sum_{j\not= i} E_j\}$.

The second assertion follows directly from the first.\proofend
\end{proof}

\bigskip
We use the notion of exceptional pairs to obtain the structure of the group $\Aut(S)$ in a very natural manner. The action of $\Aut(S)$ on the exceptional pairs gives rise to a homomorphism $\rho:\Aut(S)\rightarrow \Sym_5$. Denoting its image and kernel respectively by $H_{S}$ and $G_{S}$, we have the exact sequence
\begin{equation}\label{exGH}
1 \rightarrow G_{S} \rightarrow \Aut(S) \stackrel{\rho}{\rightarrow} H_{S} \rightarrow 1.
\end{equation} 

As every element of $G_{S}$ leaves the set $\{L-E_i,-K_S-L+E_i\}$ invariant for $i=1,...,5$, we can define a homomorphism $\gamma:G_{S}\rightarrow (\mathbb{F}_2)^5$: put $\gamma(g)=(a_1,...,a_5)$, where $a_i=1$ if $g$ permutes $L-E_i$ and $-K_S-L+E_i$, and $a_i=0$ if $g$ fixes both divisors. Note that the homomorphism $\gamma$ is injective: any element of $\Aut(S)$ which fixes all the divisors $L-E_i$, for $i=1,...,5$, also fixes $K_S+\sum_{i=1}^5(L-E_i)=2L$ and thus all the divisors of $\Pic{S}$.

\begin{remark}\fauxtitred
 We use here the notation of \cite{bib:Be2} for the group $G_{S}$. In this article, A. Beauville uses the canonical embedding of $S$ in $\mathbb{P}^4$ as the intersection of two quadrics. By a linear change of coordinates, the two quadrics can be chosen to be $\sum {x_i}^2=0$ and $\sum \lambda_i {x_i}^2$. The group $G_{S}$ is then the $2$-torsion of the maximal diagonal torus of $\PGLn{5}$. For further details on this point of view, see Proposition \refd{Prp:ClassDelPezzoQSpec}.\end{remark}
 
 \bigskip

\begin{Prp}\titreProp{Structure of $\Aut(S)$}\label{Prp:AutSF24S5}\index{Del Pezzo surfaces!automorphisms|idi}
\begin{itemize}
\item[\upshape 1.] The exact sequence (\refd{exGH}) splits.
\item[\upshape 2.] $\gamma$ induces an isomorphism $G_{S}\rightarrow \{ (r_1,r_2,r_3,r_4,r_5) \in (\mathbb{F}_2)^5 \ | \ \sum_{i=1}^5 r_i=0\}\cong (\mathbb{F}_2)^4$.
\item[\upshape 3.] $H_{S}\cong \{ h\in \PGLn{3}\ |\ h(\{A_1,A_2,A_3,A_4,A_5\})=\{A_1,A_2,A_3,A_4,A_5\} \}\subset \Sym_5$.
\item[\upshape 4.] $\Aut(S) \cong G_{S} \rtimes H_{S}$, where the group $H_{S} \subset \Sym_5$ acts on $G_{S}$ by permuting the $5$ coordinates of its elements.
\end{itemize}
\end{Prp}\begin{proof}\fauxtitred\upshape
We denote by $\Aut(\Pic{S})^{*}$ the group of automorphisms of $\Pic{S}$ that preserve the intersection, and fix $K_S$ (this is classically called the \defn{Weyl group} of $\Pic{S}$).  
The action of these automorphisms on the exceptional pairs gives rise to a homomorphism $\overline{\rho}:\Aut(\Pic{S})^{*} \rightarrow \Sym_5$. Note that this homomorphism is surjective: for $\sigma \in \Sym_5$
, the element of $\Aut(\Pic{S})^{*} $ which sends $E_i$ on $E_{\sigma(i)}$ and which fixes $L$, is sent by $\overline{\rho}$ on $\sigma$.

Denoting by $\overline{G_{S}}$ the kernel of $\overline{\rho}$, we have the exact sequences
\begin{center}$\begin{array}{cccccccc}
1 &\rightarrow& \overline{G_{S}}& \rightarrow& \Aut(\Pic{S})^{*}& \stackrel{\overline{\rho}}{\rightarrow} &\Sym_5& \rightarrow 1 \\
&& \cup & & \cup & & \cup & \\
1 &\rightarrow& G_{S}& \rightarrow& \Aut(S)& \stackrel{\rho}{\rightarrow} &H_{S}& \rightarrow 1\end{array}$.
\end{center}
As every element of $\overline{G_{S}}$ leaves the set $\{L-E_i,-K_S-L+E_i\}$ invariant for $i=1,...,5$, we can extend the homomorphism $\gamma:G_{S}\rightarrow (\mathbb{F}_2)^5$, arising from the action on the set of exceptional pairs, to an injective homomorphism $\gamma:\overline{G_{S}}\rightarrow (\mathbb{F}_2)^5$.

Consider the quadratic involution $\lambda:(x:y:z) \dasharrow (ayz:bxz:cxy)$ of\hspace{0.2 cm}$\Pn$ whose base points are $A_1,A_2,A_3$ and which exchanges $A_4$ and $A_5$. Note that the blow-up of $A_1,A_2,A_3$ gives the Del Pezzo surface of degree $6$, where the lift of $\lambda$ acts biregularly and whose action on the exceptional divisors is $(E_1\ D_{23})(E_2\ D_{13})(E_3\ D_{12})$ (see Section \refd{Subsec:DelPezzo6}).

 As $\lambda$ permutes $A_4$ and $A_5$, the map $\pi^{-1}\lambda\pi$ is an automorphism of $S$, and induces the permutation
\begin{center} $(E_1\ D_{23})(E_2\ D_{13})(E_3\ D_{12})(E_4\ E_5)(D_{14}\ D_{15})(D_{24}\ D_{25})(D_{34}\ D_{35})(D_{45}\ C)$\end{center} on the exceptional divisors.
Its action on the Picard group of $S$ is then given by
\begin{center}
$\left(\begin{array}{cccccc}
0 & -1 & -1& 0 & 0 & -1\\
-1& 0 & -1& 0 & 0 & -1\\
-1 & -1& 0& 0 & 0 & -1\\
0 & 0& 0& 0 & 1 & 0\\
0 & 0& 0& 1& 0 & 0\\
1 & 1& 1& 0 & 0 & 2\end{array}\right)$,\end{center}
with respect to the basis $(E_1,E_2,E_3,E_4,E_5,L)$. This automorphism is in the kernel of $\rho$ and its image by $\gamma$ is the element $(0,0,0,1,1) \in (\mathbb{F}_2)^5$.

By changing the role of the $A_i$'s and proceeding in the same manner, we get $10$ automorphisms of $S$ contained in $G_{S}$ whose image by $\gamma$ are elements with two $1$'s and three $0$'s. Then, $G_{S}$ contains the group generated by these elements, whose image by $\gamma$ is \begin{center}$\{ (r_1,r_2,r_3,r_4,r_5) \in (\mathbb{F}_2)^5 \ | \ \sum_{i=1}^5 r_i=0\}\cong (\mathbb{F}_2)^4.$\end{center} We then have $(\mathbb{F}_2)^4\subset G_S\subset \overline{G_S} \subset (\mathbb{F}_2)^5$. 

Let us prove that $\overline{G_S} \not= (\mathbb{F}_2)^5$, which implies that $G_S=\overline{G_S}\cong (\mathbb{F}_2)^4$. No element $h \in \overline{G_S}$ is such that $\gamma(h)=(1,1,1,1,1)$. Indeed, $h$ would send $L=\frac{1}{2}(K_S+\sum_{i=1}^5 (L-E_i))$ on the divisor $\frac{1}{2}(K_S+\sum_{i=1}^5 (-K_S-L+E_i))=\frac{1}{2}(-2L-3 K_S)$, which doesn't belong to $\Pic{S}$. This gives part $2$ of the proposition.

Let $\sigma \in H_{S} \subset \Sym_5$ be a permutation of the set $\{1,...,5\}$ in the image of $\rho$ and $g$ be an automorphism of $S$ such that $\rho(g)=\sigma$. Let $\alpha$ be the element of $\Aut(\Pic{S})^{*}$ that sends $E_{i}$ on $E_{\sigma(i)}$ and fixes $L$. Since $g\alpha^{-1} \in \Aut(\Pic{S})^{*}$ fixes the five exceptional pairs, it belongs to $\overline{G_{S}}$ and hence to $G_{S}$. This implies that $\alpha$ is in fact an automorphism of $S$ and that $r=\pi \alpha \pi^{-1}$ is an automorphism of\hspace{0.2 cm}$\Pn$ such that $r(A_i)=A_{\sigma(i)}$ for $i=1,...,5$.
Conversely, it is clear that every automorphism $r$ of\hspace{0.2 cm}$\Pn$ which leaves the set $\{A_1,A_2,A_3,A_4,A_5\}$ invariant lifts to the automorphism $\pi^{-1}r\pi$ of $S$ whose action on the exceptional pairs is the same as that of $r$ on the set $\{A_1,A_2,A_3,A_4,A_5\}$.
This gives part $3$ and then part $1$ of the proposition.

Part $4$ obviously follows from the others.
\proofend\end{proof}

\bigskip

\subsection{The geometry of elements of $G_{S}=(\mathbb{F}_2)^4$}
In Proposition \refd{Prp:AutSF24S5}, we proved that the homomorphism $\gamma:G_{S}\rightarrow (\mathbb{F}_2)^5$, arising from the action on the exceptional pairs, is an isomorphism from $G_{S}\subset \Aut(S)$ to $\{ (r_1,r_2,r_3,r_4,r_5) \in (\mathbb{F}_2)^5 \ | \ \sum_{i=1}^5 r_i=0\}$.

There are then $10$ elements with two $1$'s and three $0$'s, and $5$ elements with four $1$'s and one $0$. By convention, we will call these two kinds of automorphisms respectively \defn{quadratic involutions} and \defn{cubic involutions}, associating them with their corresponding birational maps of\hspace{0.2 cm}$\Pn$. We will describe the geometry of these elements in the following two propositions. By a linear change of coordinates, we may restrict our attention to the elements $(0,0,0,1,1)$ and $(0,1,1,1,1)$.

We again suppose that $S$ is the blow-up of the points $A_1=(1:0:0)$, $A_2=(0:1:0)$, $A_3=(0:0:1)$, $A_4=(1:1:1)$ and $A_5=(a:b:c)$ on $\Pn$ and denote by $\pi: S \rightarrow \Pn$ the corresponding blow-down.
\begin{Prp}\titreProp{quadratic involutions}\\
\index{Involutions!De Jonqui\`eres!quadratic|idb}
\label{Prp:GeoQuad}
Let $g \in G_{S}$ be the automorphism of $S$ which is the lift on $S$ of the birational map $g'=\pi g \pi^{-1}$ of\hspace{0.2 cm}$\Pn$. If $\gamma(g)=(0,0,0,1,1)$, then:
\begin{itemize}
\item[\upshape 1.]
Both $g$ and $g'$ are de Jonqui\`eres involutions.
\item[\upshape 2.]
The birational map $g'$ is the quadratic involution of\hspace{0.2 cm}$\Pn$ given by $(x:y:z) \dasharrow (ayz:bxz:cxy)$.
\item[\upshape 3.]
Both $g'$ and $g$ have four fixed points, that are respectively $(\pm \sqrt{a}:\pm \sqrt{b}:\pm \sqrt{c}) \in \Pn$ and their images by $\pi^{-1}$.
\item[\upshape 4.]
Any two of the four points of\hspace{0.2 cm}$\Pn$ fixed by $g'$ are aligned with $A_1$, $A_2$ or $A_3$.
\item[\upshape 5.]
The automorphism $g$ induces the permutation
\begin{center} $(E_1\ D_{23})(E_2\ D_{13})(E_3\ D_{12})(E_4\ E_5)(D_{14}\ D_{15})(D_{24}\ D_{25})(D_{34}\ D_{35})(D_{45}\ C)$\end{center} on the exceptional divisors
and its action on the Picard group of $S$ is given by
\begin{center}
$\left(\begin{array}{cccccc}
0 & -1 & -1& 0 & 0 & -1\\
-1& 0 & -1& 0 & 0 & -1\\
-1 & -1& 0& 0 & 0 & -1\\
0 & 0& 0& 0 & 1 & 0\\
0 & 0& 0& 1& 0 & 0\\
1 & 1& 1& 0 & 0 & 2\end{array}\right)$,\end{center}
with respect to the basis $(E_1,E_2,E_3,E_4,E_5,L)$.
\item[\upshape 6.]
The pair $(<g>,S)$ is not minimal.
\item[\upshape 7.]
The map $g'$ (and hence $g$) is conjugate in $\CrP$ to the linear involutions of\hspace{0.2 cm}$\Pn$.
\end{itemize} 
\end{Prp}\begin{proof}\fauxtitred\upshape
Since $g$ leaves invariant the divisor $L-E_1$, which is a pencil of rational curves, $g$ (and hence $g'$) is a de Jonqui\`eres involution.

Assertions $2$ and $5$ were established in the proof of Proposition \refd{Prp:AutSF24S5}. Assertions $3$ and $4$ follow by writing $g'$ as $(x:y:z) \dasharrow (ayz:bxz:cxy)$.

Assertion $5$ follows from assertion $3$. Indeed, the automorphism $g$ leaves invariant the set of disjoint exceptional curves $\{E_1,D_{23},E_4,E_5\}$. We blow-down these $4$ exceptional curves and get a $g$-equivariant birational morphism $p$ from $S$ to $\mathbb{P}^1\times \mathbb{P}^1$; explicitly, \begin{center}$p\pi^{-1}:(x:y:z) \dasharrow (x:y)\times(x:z)$.\end{center}

As $g$ fixes the divisor $2L-E_1-E_2$, the birational map $g'$ leaves the linear system $C^2(A_1,A_2)$ invariant and hence also the homoloidal net $C^2(A_1,A_2,P)$, where $P$ is one of the four fixed points of\hspace{0.2 cm}$\Pn$. This implies (Proposition \refd{Prp:HomolNetInvariant}) that $g'$ is conjugate in $\CrP$ to the linear involutions. Note that this result can also be proved by the above conjugation of $g$ to an involutive element of $\Aut(\mathbb{P}^1\times\mathbb{P}^1)$.\proofend
\end{proof}

\bigskip

\begin{Prp}\titreProp{cubic involutions}\\
\label{Prp:GeoCub}\index{Involutions!De Jonqui\`eres!cubic|idb}
Let $g \in G_{S}$ be the automorphism of $S$ such that $\gamma(g)=(0,1,1,1,1)$, the lift on $S$ of the birational map $g'=\pi g \pi^{-1}$ of\hspace{0.2 cm}$\Pn$. Then:
\begin{itemize}
\item[\upshape 1.]
The birational map $g'$ is the cubic de Jonquières involution of\hspace{0.2 cm}$\Pn$ given by 
\begin{center}$\begin{array}{lrccccccl}
(x:y:z) \dasharrow \big(\h{3}&-ayz(&\h{3}(c-b)x&\h{2}+\h{2}&(a-c)y&\h{2}+\h{2}&(b-a)z&\h{3})&\h{2}:\\ 
\h{3}& y(&\h{3}a(c-b)yz&\h{2}+\h{2}&b(a-c)xz&\h{2}+\h{2}&c(b-a)xy&\h{3})&\h{2}:\\
\h{3}& z(&\h{3}a(c-b)yz&\h{2}+\h{2}&b(a-c)xz&\h{2}+\h{2}&c(b-a)xy&\h{3})&\h{2}\big)\end{array}$\end{center}
\item[\upshape 2.]
The map $g'$ fixes the pencil of lines of\hspace{0.2 cm}$\Pn$ passing through $A_1$. It also fixes the pencil of conics through $A_2,A_3,A_4$ and $A_5$. It thus preserves the conic bundle structures induced by both pencils and is a twisting de Jonqui\`eres involution in both cases.
\item[\upshape 3.]
The automorphism $g$ induces the permutation
\begin{center} $(E_1\ C)(E_2\ D_{12})(E_3\ D_{13})(E_4\ D_{14})(E_5\ D_{15})(D_{23}\ D_{45})(D_{24}\ D_{35})(D_{25}\ D_{34})$\end{center} on the exceptional divisors of $S$
and its action on the Picard group is given by
\begin{center}
$\left(\begin{array}{cccccc}
-1 & -1 & -1& -1 & -1 & -2\\
-1& -1 & 0& 0 & 0 & -1\\
-1 & 0& -1& 0 & 0 & -1\\
-1 & 0& 0& -1 & 0 & -1\\
-1 & 0& 0& 0& -1 & -1\\
2 & 1& 1& 1 & 1 & 3\end{array}\right)$,\end{center} with respect to the basis $(E_1,E_2,E_3,E_4,E_5,L)$.
\item[\upshape 4.]
The curve of\hspace{0.2 cm}$\Pn$ fixed by $g'$ is a smooth cubic of\hspace{0.2 cm}$\Pn$ passing through $A_1,...,A_5$, with equation
\begin{center} $b(a-c)x^2z+c(b-a)x^2y+a(a-c)y^2z+a(b-a)yz^2+2a(c-b)xyz=0,$\end{center}
and is isomorphic to the elliptic curve $y^2z=x(x-z)(bx-cz)$, with invariant $(b:c) \in \mathbb{P}^1_{/{S_3}}$.\index{Curves!elliptic!birational maps that fix an elliptic curve}
\item[\upshape 5.]
The pair $(<g>,S)$ is minimal; $\rkPic{S}^g=2$, generated by $K_S$ and $L-E_1$.
\item[\upshape 6.]
The map $g'$ (and hence $g$) is not birationally conjugate to a linear involution of\hspace{0.2 cm}$\Pn$.
\item[\upshape 7.]
Two involutions of this kind are conjugate if and only if their curves of fixed points are isomorphic.
\end{itemize} 
\end{Prp}\begin{proof}\fauxtitred\upshape
Assertion $1$ can be proved by calculating the composition of two quadratic involutions, using Proposition \refd{Prp:GeoQuad}. The calculation of the fixed points of $g'$ using this explicit form provides the equation of assertion $4$. We also see that $g'$ fixes the pencil of lines of\hspace{0.2 cm}$\Pn$ through $A_1=(1:0:0)$ (which are of the form $ay+bz=0,\ \ (a:b)\in\mathbb{P}^1$).

The trivial action of $g$ on the exceptional pair $\{L-E_1,-K_S-L+E_1\}$ implies that $g$ is an automorphism of the two conic bundles induced respectively by $L-E_1$ and $-K_S-L+E_1$. 
Since $g'$ fixes an irreducible elliptic curve it is a twisting de Jonqui\`eres involution, in both cases, twisting the $4$ singular fibres of both conic bundles. These are: 
\begin{center}
\begin{tabular}{p{14 cm}}
{\it $\{D_{1i},E_i\}$ for $i=2,...,5$: the singular fibres of the conic bundle induced by $L-E_1$;}\\
{\it $\{D_{23},D_{45}\}$, $\{D_{24},D_{35}\}$, $\{D_{25},D_{34}\}$ $\{C,E_1\}$, those  of the conic bundle induced by $2L-E_2-E_3-E_4-E_5$.}\end{tabular}\end{center}

Assertion $3$ follows directly from this observation. (Note that it can also be verified by using the composition of two quadratic involutions.)
The morphism $S\rightarrow \mathbb{P}^1$ induced by $L-E_1$ (which corresponds to the projection $(x:y:z) \dasharrow (y:z)$) gives the curve of fixed points as a double covering of $\mathbb{P}^1$ over the points $(1:0)$, $(0:1)$, $(1:1)$ and $(b:c)$ (corresponding to the singular fibres containing respectively $E_2$, $E_3$, $E_4$ and $E_5$); we then obtain the last part of assertion $4$.

As every exceptional curve is sent on another one which intersects it, the pair $(<g>,S)$ is minimal. Hence, the rank of $\Pic{S}^g$ is at most $2$. As $g$ fixes the divisors $K_S$ and $L-E_1$, we get assertion $5$.

Since no linear automorphism of\hspace{0.2 cm}$\Pn$ fixes a non-rational curve, the maps $g$ and $g'$ cannot be conjugate to a linear involution.

The last point is proven in Proposition \refd{Prp:classDeJICremona} and in \cite{bib:BaB}, Proposition 2.7. \proofend
\end{proof}

\begin{remark}\fauxtitred
The last assertion implies that the involutions $(0,1,1,1,1)$ of the two surfaces given by the blow-up of respectively $A_1,...,A_4,(a:b:c)$ and $A_1,...,A_4,(a':b':c')$ are then conjugate in $\CrP$ if and only if $(b:c)$ and $(b':c')$ are equivalent under the action of $S_3$ on $\mathbb{P}^1$, i.e.\ if and only if $(b':c') \in \{(b:c), (c:b), (b:b-c), (b-c:b), (c:c-b), (c-b:c)\}$. 
\end{remark}
\bigskip

\subsection{The embedding of $S$ in $\mathbb{P}^4$}
Proposition \refd{Prp:GeoCub} describes the five cubic involutions of $S$. We now use these to give explicitly the classical embedding of $S$ as the intersection of two quadrics in $\mathbb{P}^4$ (see \cite{bib:Ho1} and \cite{bib:Be2}).
\begin{Prp}\fauxtitre
\label{Prp:embedtwoquadrics}
Let $g_1,...,g_5 \in G_S$ denote the cubic involutions corresponding via $\gamma$ to \begin{center}$(0,1,1,1,1),(1,0,1,1,1),...,(1,1,1,1,0)$\end{center} and let $F_i$ denote the equation of the elliptic curve of\hspace{0.2 cm}$\Pn$ fixed by $\pi g_i\pi^{-1}$.\index{Curves!elliptic!birational maps that fix an elliptic curve} Then:
\begin{itemize}
\item
The equations $F_1,...,F_5$ are:
\begin{itemize}
\item[]
$F_1=b(a-c)x^2z+c(b-a)x^2y+a(a-c)y^2z+a(b-a)yz^2+2a(c-b)xyz,$
\item[]
$F_2=a(b-c)y^2z+c(a-b)y^2x+b(b-c)x^2z+b(a-b)xz^2+2b(c-a)xyz,$
\item[]
$F_3=b(c-a)z^2x+a(b-c)z^2y+c(c-a)y^2x+c(b-c)yx^2+2c(a-b)xyz,$
\item[]
$F_4=bcx^2(z-y)+abz^2(y-x)+acy^2(x-z),$
\item[]
$F_5=ayz(y-z)+bxz(z-x)+cxy(x-y).$
\end{itemize}
\item
The map $\varphi:(x:y:z) \dasharrow (F_1:F_2:F_3:F_4:F_5)$ is a birational map from $\Pn$ to the surface $S'\subset \mathbb{P}^4$ defined by the equations:
\begin{itemize}
\item[]
$cx_1^2-ax_3^2+(a-c)x_4^2-ac(a-c)x_5^2=0$,
\item[]
$cx_2^2-bx_3^2+(c-b)x_4^2-bc(c-b)x_5^2=0$,
\end{itemize}
and the map $\varphi\pi:S\rightarrow S'$ is an isomorphism whose linear system is $|-K_S|$.
\item
The group $G_S$ is sent by the isomorphism $\varphi\pi$ on the $2$-torsion of the maximal diagonal torus of $\PGLn{5}$. An element $g \in G_S$ such that $\gamma(g)=(a_1,a_2,a_3,a_4,a_5) \in (\mathbb{F}_2)^5$ corresponds to the automorphism $(x_1:...:x_5) \mapsto ((-1)^{a_1}x_1:...:(-1)^{a_5}x_5)$.
\end{itemize}
\end{Prp}
\begin{proof}\fauxtitred\upshape
The equation of $F_1$ was given in Proposition \refd{Prp:GeoCub}. The equation of $F_2$ is obtained by exchanging the roles of $A_1$ and $A_2$, which consists in permuting the coordinates $x$ and $y$ and the coordinates $a$ and $b$. We obtain $F_3$, $F_4$ and $F_5$ in the same manner.

The linear system generated by $F_1,F_2,...,F_5$ is the linear system of cubics passing through $A_1,...,A_5$. The map $\varphi\pi$ is then defined by $|-K_S|$ and is an isomorphism with a non-singular surface of $\mathbb{P}^4$ (see \cite{bib:Kol}, Theorem III.3.5), since $-K_S$ is very ample. The equations follow by direct calculation.

We observe directly that the $2$-torsion of the maximal diagonal torus of $\PGLn{5}$ acts on the surface $S'$ and that the automorphism $(x_1:...:x_5) \mapsto (-x_1:x_2:...:x_5)$ fixes the trace in $S$ of the hyperplane $x_1=0$ of $\mathbb{P}^4$. This set of fixed points corresponds to the curve fixed by the cubic involution of $S$ corresponding to $(0,1,1,1,1)$. We do this for all cubic involutions and obtain a correspondence between the group generated by the cubic involutions and the $2$-torsion of the maximal diagonal torus of $\PGLn{5}$.
\proofend
\end{proof}

\bigskip

\subsection{The subgroups of $G_{S} \subset \Aut(S)$}
As $G_{S}$ is isomorphic to $(\mathbb{F}_2)^4$, its subgroups are isomorphic to $(\mathbb{F}_2)^r$, where $1\leq r\leq 4$. The condition on the fixed part of the Picard group gives:
\begin{Lem}\fauxtitre
\label{Lem:TwoCubicInv}
Let $H$ be a subgroup of $G_{S}$. Then $H$ contains two distinct cubic involutions if and only if $\rkPic{S}^{H}=1$. 
\end{Lem}\begin{proof}\fauxtitred\upshape
Set $V_i=K_S+2(L-E_i)$ for $i=1,...,5$. Note that $\mathcal{B}=(K_S,V_1,V_2,V_3,V_4,V_5)$ is a basis of $\Pic{S} \otimes \mathbb{Q}$ and that every element of $G_{S}$ is diagonal with respect to this basis. In fact, an element $g \in G_{S}$ such that $\gamma(g)=(r_1,r_2,r_3,r_4,r_5) \in (\mathbb{F}_2)^5$ corresponds to the diagonal matrix whose diagonal is $(1,(-1)^{r_1},(-1)^{r_2},(-1)^{r_3},(-1)^{r_4},(-1)^{r_5})$.

From this, we deduce that the subgroup of $\Pic{S}$ fixed by two distinct cubic involutions has rank $1$.

Conversely, if a subgroup $H$ of $G_{S}$ contains no cubic involution, it fixes at least two of the $V_i$ and if it contains exactly one cubic involution, it fixes one of the $V_i$, and therefore $\rkPic{S}^{H}>1$.\proofend
\end{proof}

\bigskip

We have then several subgroups $H \subset G_{S}$ such that $\rkPic{S}^{H}=1$, isomorphic to $(\Z{2})^2,(\Z{2})^3$ and $(\Z{2})^4$. To decide whether these groups are conjugate, we will use the following Proposition.

\bigskip

\begin{Prp}\fauxtitre
\label{Prp:conjisoA1A2}
Let $S$ be the blow-up of $A_1,...,A_4,A_5 \in \Pn$ and $S'$ the blow-up of $A_1,...,A_4,A_5' \in \Pn$; these are two Del Pezzo surfaces of degree $4$.

Let $G\subset \Aut(S)$, $G'\subset \Aut(S')$ be two subgroups isomorphic to $(\Z{2})^2$, generated by two cubic involutions.
 Suppose that there exists a birational map $\varphi:S\dasharrow S'$  that conjugates $G$ to $G'$.
Then $\varphi$ is an isomorphism.
\end{Prp}
\begin{remark}\fauxtitred
\begin{itemize}
\item
This proposition is a generalisation of \cite{bib:Be2}, Proposition 4.3, which is a particular case of Theorem 3.3 of \cite{bib:Isk2}. This theorem cannot be applied here since the group $G$ fixes $4$ points on each surface. Despite this, our approach is inspired by \cite{bib:Be2}. 
\item
The classification of elementary links (see \cite{bib:Isk5}, Theorem 2.6) shows that $\varphi$ is the composition of elementary links, which may start by a Bertini or a Geiser involution (from $S$ to itself), or the blow-up of one point, and then some links of conic bundles. Using these might provide another proof, but we believe that it would be longer and more difficult.
\item
Here, we have examples of groups that \emph{have the same action on their sets of fixed points} (see Section \refd{Sec:ActionFixedCurves}) but are not conjugate:

Take for example $A_5=(a:b:c)$ and $A_5'=(c^2:ab:ac)$. The curves fixed by $g_1$ and $g_1'$ are isomorphic {\upshape (}the invariant is $(b:c)${\upshape )}, the curves fixed by $g_2$ and $g_2'$ are isomorphic {\upshape(}$(a:c) \equiv (c:a)${\upshape)} and both $g_1g_2$ and $g_1'g_2'$ fix $4$ points.  Furthermore, every invariant of Chapter \refd{Chap:ConjugacyInvariants} is the same for both groups.

But the sets $\{A_1,A_2,A_3,A_4,A_5\}$ and $\{A_1,A_2,A_3,A_4,A_5'\}$ are not linearly equivalent for a general point $(a:b:c)\in \Pn$, so the surface obtained by blowing-up are not isomorphic (see Proposition \refd{Prp:IsoDPlinear}), hence the  groups $G$ and $G'$ are not birationally conjugate.
\end{itemize}
\end{remark}
\begin{proof}\fauxtitred\upshape
Let $g_1$ and $g_2$ be the two cubic involutions of $G$. We denote by $g_1'$ and $g_2'$ the cubic involutions of $G'$ such that $\varphi$ conjugates $g_1$ to $g_1'$ and $g_2$ to $g_2'$, i.e.\  $g_i'=\varphi g_i \varphi^{-1}$ for $i=1,2$.

Without loss of generality, we may assume that  $g_1$ and $g_1'$ both correspond to the cubic involution $(0,1,1,1,1)$ and that $g_2$ and $g_2'$ both correspond to the cubic involution $(1,0,1,1,1)$. 

We solve the indeterminacy and get a diagram
\begin{center}
\hspace{0 cm}\xymatrix{ & \tilde{S} \ar_{\eta}[rd] \ar_{\epsilon}[ld] & \\
S\ar@{-->}^{\varphi}[rr] \ar_{\pi}[d] & & S' \ar_{\pi'}[d] \\
\Pn& & \Pn}
\end{center}
where $\eta$ and $\epsilon$ are birational morphisms, and $\eta^{-1}$ and $\epsilon^{-1}$ consist of sequences of blow-ups of points. We have some effective divisors $F_1,...,F_r$ and $F_1',...,F_r'$ such that 
\begin{center} $\Pic{\tilde{S}}=\epsilon^{*} \Pic{S} \oplus \sum_{i=1}^{r} \mathbb{Z}F_i$, where $F_i^2=-1$, $F_i\cdot F_j=0$ for $i\not=j$, $1\leq i,j \leq r$,\vspace{0.1 cm}\\
$\Pic{\tilde{S}}=\eta^{*} \Pic{S'} \oplus \sum_{i=1}^{r} \mathbb{Z}F_i'$, where $F_i'^2=-1$, $F_i'\cdot F_j'=0$ for $i\not=j$, $1\leq i \leq r$,\vspace{0.1 cm}\\
$K_{\tilde{S}}=\epsilon^{*}(K_S)+\sum_{i=1}^{r}F_i=\eta^{*}(K_{S'})+\sum_{i=1}^{r}F_i'$.
\end{center}
We may suppose that there exist $\tilde{g_1},\tilde{g_2}\in \Aut(\tilde{S})$, such that $\epsilon$ conjugates $\tilde{g_i}$ to $g_i$ and $\eta$ conjugates $\tilde{g_i}$ to $g_i'$, for $i=1,2$. Furthermore, $\tilde{g_1}$ and $\tilde{g_2}$ act on the set $\{F_1,...,F_r\}$ (respectively on $\{F_1',...,F_r'\}$).

Note that $\Pic{S}^{g_1}=\mathbb{Z} (L-E_1) \oplus \mathbb{Z} (-K_S-L+E_1)$, which implies that
\begin{equation}
\label{etaLE1}\begin{array}{rcc}
\eta^{*}(L-E_1)&=&k \epsilon^{*}(L-E_1)+ l \epsilon^{*}(-K_S-L+E_1)-\sum_{i=1}^r s_i F_i
\end{array}\end{equation}
for some integers $k,l,(s_i)_{i=1}^r$. On computing the intersection of these effective divisors with the effective divisors $\epsilon^{*}(L-E_1), \epsilon^{*}(-K_S-L+E_1), F_1,...,F_r$, we see that all these integers are non-negative.

\bigskip

Let us first prove that
\begin{equation}
\label{assert1Prpg1g2}\begin{array}{rcl}
s_i \not=0& \mbox{implies that } F_i \mbox{ is invariant by }  \tilde{g_1}.
\end{array}\end{equation}

Computing the intersection of $\eta^{*}(L-E_1)$ with \begin{center}$-K_{\tilde{S}}=\epsilon^{*}(-K_S)-\sum_{i=1}^r F_i=\eta^{*}(-K_{S'})-\sum_{i=1}^r F_i'$\end{center} yields the equation
\begin{equation}\label{eqsumab1} \begin{array}{ccc}\sum_{i=1}^r s_i&=&2(k+l-1).\end{array}\end{equation}
The set of points of $\tilde{S}$ fixed by the automorphism $\tilde{g_1}$ contains a curve numerically equivalent to $\Gamma=\epsilon^{*}(-K_S)-\sum v_i F_i$, each $v_i$ being $1$ if $F_i$ is invariant, and $0$ if not. This follows from the fact that the projection of this curve on $\Pn$ is a smooth cubic curve passing through all the points $A_1,...,A_5$ (by Proposition \refd{Prp:GeoCub}).

The curves on $S'$ equivalent to $L-E_1$ are the strict pull-backs by $\pi'$ of lines passing through $A_1$. This is a pencil of rational curves, each invariant by $g_1'$, but not fixed (by Proposition \refd{Prp:GeoCub}). The intersection of $\eta^{*}(L-E_1)$ with $\Gamma$, which is $2(k+l)-\sum_{i=1}^r v_i s_i$, must then be exactly two. Using (\refd{eqsumab1}) we see that $v_i=1$ if $s_i>0$, so $F_i$ is invariant by $\tilde{g_1}$. We have thus proved assertion (\refd{assert1Prpg1g2}).

\bigskip

Let us calculate $\eta^{*}(L-E_1)+\tilde{g_2}(\eta^{*}(L-E_1))=\eta^{*}(-K_{S'})$, using (\refd{etaLE1}). We find \begin{equation}\label{etaKS}
\eta^{*}(-K_{S'})=(k+l) \epsilon^{*}(-K_S)-\sum_{i=1}^r s_i (F_i+\tilde{g_2}(F_i)).\end{equation}

Remember that $g_1$ and $g_2$ commute. Then, for every $i$ such that $s_i>0$, the divisor $F_i$ is invariant by $\tilde{g_1}$, and so must be the divisor $g_2(F_i)\in \{F_1,...,F_r\}$. So the decomposition of $\eta^{*}(-K_{S'})$ in $\epsilon^{*} \Pic{S} \oplus \sum_{i=1}^{r} \mathbb{Z}F_i$ is a linear combination involving only $\epsilon^{*}(-K_S)$ and divisors $F_i$ invariant by $\tilde{g_1}$. We can repeat this reasoning after exchanging $g_1$ and $g_2$, so $\eta^{*}(-K_{S'})$ is finally a linear combination of $\epsilon^{*}(-K_S)$ and $F_i$'s that are invariant both by $\tilde{g_1}$ and $\tilde{g_2}$.

Note that there are four points of $S$ that are fixed both by $g_1$ and $g_2$. Their images by $\pi$ form, together with $A_1$,...,$A_4,A_5$, the intersection in $\Pn$ of the two curves of fixed points of $g_1$ and $g_2$; they are also the four fixed points of the quadratic involution $g_1g_2$ (Proposition \refd{Prp:GeoQuad}). We may assume that these four points are blown-up by $\epsilon^{-1}$ and also that they are $\epsilon(F_1),\epsilon(F_2),\epsilon(F_3)$ and $\epsilon(F_4)$, numbered so that $\sum_{i=1}^r s_i =\sum_{i=1}^4 s_i $.

We set $m=k+l$; equation (\refd{etaKS}) can then be written \begin{equation}\label{etaKS2}\eta^{*}(-K_{S'})=m \epsilon^{*}(-K_S)-2\sum_{i=1}^4 s_i F_i.\end{equation}

We saw in Proposition \refd{Prp:GeoQuad} that two of the fixed points of the birational map $\pi(g_1 g_2)\pi^{-1}$ of\hspace{0.2 cm}$\Pn$ are aligned with $A_3$, $A_4$ or $A_5$. Then for every $i\not=j \in \{1,...,4\}$, there exists $3 \leq v \leq 5$ such that the divisor $\epsilon^{*}(L-E_v)-F_i-F_j$ is effective. By computing its intersection with $\eta^{*}(-K_{S'})$, and using (\refd{etaKS2}), we see that $s_i+s_j\leq k+l=m$.

Using (\refd{etaKS2}) again, we compute the intersection of $\eta^{*}(-K_{S'})$ with itself and $-K_{\tilde{S}}$, and get two more conditions, that give, with the above inequality:
\begin{equation}
\begin{array}{rcll}
\sum_{i=1}^4 s_i^2&=&m^2-1, \vspace{0.1 cm}\\
 \sum_{i=1}^4 s_i&=&2(m-1),\vspace{0.1 cm}\\
 s_i+s_j &\leq& m, & 1 \leq i,j\leq 4,\ i \not=j.
\end{array}
\end{equation}
These relations imply  that $s_i=0$ for $i=1,...,4$, and that $m=1$ (see Lemma \refd{Lem:ArCondSum} in the Appendix). Then, using (\refd{etaKS2}), we have $\eta^{*}(-K_{S'})=\epsilon^{*}(-K_S)$, from which it follows that $\varphi$ is an isomorphism (Lemma \refd{Lem:epsetaiso} below).\proofend
\end{proof}

\bigskip

\begin{Cor}
\label{Cor:DelP4contrEx}
Let $S$ and $S'$ be two Del Pezzo surfaces of degree $4$ and let $H$ and $H'$ be respectively subgroups of $G_{S}$ and $G_{S'}$ such that $\rkPic{S}^H=\rkPic{S'}^{H'}=1$.

The groups $H$ and $H'$ are conjugate in $\CrP$ if and only if they are conjugate by an automorphism $\varphi:S\rightarrow S'$.
\proofend
\end{Cor}

\bigskip
\begin{Lem}\fauxtitre\label{Lem:epsetaiso}
Let $\varphi:S\dasharrow S'$ a birational map of surfaces.

Let $\tilde{S}$ be a surface, $\epsilon:\tilde{S} \rightarrow S$ and $\eta:\tilde{S} \rightarrow S'$ be two birational morphisms such that the diagram
\begin{center}
\hspace{0 cm}\xymatrix{ & \tilde{S} \ar_{\eta}[rd] \ar_{\epsilon}[ld] & \\
S\ar@{-->}^{\varphi}[rr] & & S' }
\end{center}
commutes.
Then, $\epsilon^{*}(K_{S})=\eta^{*}(K_{S'})$ if and only if $\varphi$ is an isomorphism. 
\end{Lem}
\begin{proof}\fauxtitred\upshape
First of all, if $\varphi$ is an isomorphism, we have $\eta^{*}(K_{S'})=(\epsilon\varphi)^{*}(K_S')=\epsilon^{*}(K_S)$. Let us prove that the converse is also true.

We have some effective divisors $F_1,...,F_r$ and $F_1',...,F_r'$ such that 
\begin{center} $\Pic{\tilde{S}}=\epsilon^{*} \Pic{S} \oplus \sum_{i=1}^{r} \mathbb{Z}F_i$, where $F_i^2=-1$, $F_i\cdot F_j=0$ for $i\not=j$, $1\leq i,j \leq r$,\vspace{0.1 cm}\\
$\Pic{\tilde{S}}=\eta^{*} \Pic{S'} \oplus \sum_{i=1}^{r} \mathbb{Z}F_i'$, where $F_i'^2=-1$, $F_i'\cdot F_j'=0$ for $i\not=j$, $1\leq i \leq r$,\vspace{0.1 cm}\\
$K_{\tilde{S}}=\epsilon^{*}(K_S)+\sum_{i=1}^{r}F_i=\eta^{*}(K_{S'})+\sum_{i=1}^{r}F_i'$.
\end{center}

We write the effective divisor $F_i' =\sum_{j=1}^r a_{ij}F_j +\epsilon^{*}(D_i)$, for some integers $a_{ij}$ and $D_i \in \Pic{S}$. Two situations are possible:
\begin{itemize}
\item
$F_i'$ is equal to some $F_j$, so $D_i=0$;
\item
$D_i$ is an effective divisor of $S$.
\end{itemize}
Since $\epsilon^{*}(K_{S})=\eta^{*}(K_{S'})$, we see that $\sum_{i=1}^r F_i=\sum_{i=1}^r F_i'$. This implies that $\sum_{i=1}^r D_i=0$, so $D_i=0$ for $i=1,..r$, and $F_i'$ is equal to some $F_j$. The sets $\{F_1,...,F_n\}$ and $\{F_1',...,F_n'\}$ are then equal, so $\varphi$ is an isomorphism.\proofend
\end{proof}

\subsection{Abelian subgroups of $\Aut(S)$ not contained in $G_{S}$}
We have seen (Proposition \refd{Prp:AutSF24S5}) that the group $\Aut(S)$ of automorphisms of $S$ is equal to $(\mathbb{F}_2)^4 \rtimes H_{S}$, where $H_{S}$ is the lift on $S$ of automorphisms of\hspace{0.2 cm}$\Pn$ that leave the set $\{A_1,A_2,A_3,A_4$, $A_5\}$ invariant. In general, this group is reduced to the identity, so $\Aut(S)=(\mathbb{F}_2)^4$. The following proposition enumerates the other cases:
\begin{Prp}\fauxtitre
\label{Prp:ClassDelPezzoQSpec} 
Let $S$ be a Del Pezzo surface of degree $4$ with a group of automorphisms not reduced to $(\mathbb{F}_2)^4$. Then:
\begin{itemize}
\item
The surface $S$ is isomorphic to the blow-up of $A_1=(1:0:0)$, $A_2=(0:1:0)$, $A_3=(0:0:1)$, $A_4=(1:1:1)$ and a point $A_5$ that lies on the line $z=x+y$, $A_5=(1:\xi:1+\xi)$, where $\xi \in \K \backslash \{0,1,-1\}$.
\item
$S$ can be embedded in $\mathbb{P}^4$ as the surface defined by the equations:
\begin{center}\begin{itemize}
\item
$(1+\xi)x_1^2-x_3^2-\xi x_4^2+\xi(1+\xi)x_5^2=0$
\item
$(1+\xi)x_2^2-\xi x_3^2-x_4^2+\xi(1+\xi)x_5^2=0$
\end{itemize}\end{center}
\item
The group $H_S$ contains the automorphism $(A_1\ A_2)(A_3\ A_4)$, which corresponds in $\mathbb{P}^4$ to the automorphism \begin{center}$(x_1:x_2:x_3:x_4:x_5) \mapsto (x_2:x_1:x_4:x_3:-x_5)$.\end{center}
\item
The groups $H_{S}$ and $\Aut(S)$ are then:

\begin{tabular}{llll}
${\num{S4.2}}$ & If $\xi^2+1\not=0$, $\xi^2\pm \xi \pm 1\not=0$, then \\ 
& $H_{S}=<(A_1\ A_2)(A_3\ A_4)>$, $\Aut(S) \cong (\mathbb{F}_2)^4 \rtimes \Z{2}$.\\ \\
${\num{S4.3}}$ & If $\xi^2-\xi+1=0$, then \\
& $H_{S}=<(A_1\ A_2\ A_5), (A_1\ A_2)(A_3\ A_4)>$, $\Aut(S)\cong(\mathbb{F}_2)^4 \rtimes \Sym_3$.\\ 
& $(A_1\ A_2\ A_5)$ corresponds in $\mathbb{P}^4$ to \\
& $(x_1:x_2:x_3:x_4:x_5) \mapsto ((\xi-1)x_5:x_1:(\xi-1)x_3:-\xi x_4:-\xi x_2)$.\\ \\
${\num{S4.4}}$ & If $\xi^2+1=0$, then \\ 
& $H_{S}=<(A_1\ A_3\ A_2\ A_4)>$, $\Aut(S)\cong(\mathbb{F}_2)^4 \rtimes \Z{4}$.\\ 
& $(A_1\ A_3\ A_2\ A_4)$ corresponds in $\mathbb{P}^4$ to \\
& $(x_1:x_2:x_3:x_4:x_5) \mapsto (-x_4:-x_3:x_1:x_2:\xi x_5)$.\\ \\
${\num{S4.5}}$ & If $\xi^2+\xi-1=0$, then \\
& $H_{S}=<(A_1\ A_4\ A_3\ A_2\ A_5), (A_1\ A_2)(A_3\ A_4)>$, $\Aut(S)\cong(\mathbb{F}_2)^4 \rtimes D_5$,\\
& where $D_5$ denotes the dihedral group of $10$ elements.\\
& $(A_1\ A_4\ A_3\ A_2\ A_5)$ corresponds in $\mathbb{P}^4$ to \\
& $(x_1:x_2:x_3:x_4:x_5) \mapsto (x_5:x_3:(1+\xi)x_4:(2+\xi)x_1:(2+\xi)x_2).$\\
\\
& If $\xi^2+\xi+1$, the surface is isomorphic to $\nump{S4.3}$.\\
&  If $\xi^2-\xi-1=0$, the surface is isomorphic to $\nump{S4.5}$.
\end{tabular}
\item
Furthermore, the surfaces obtained with $A_5=(1:\xi:1+\xi)$ and $A_5'=(1:\xi':1+\xi')$ are isomorphic if and only if $\xi'\in \{\xi,-\xi,\xi^{-1},-\xi^{-1}\}$. Then, except for the first case, there is only one surface, up to isomorphism, in each case.
\end{itemize}
\end{Prp}\begin{proof}\fauxtitred\upshape
Note first of all that the automorphism $(x:y:z) \mapsto (z-y:z-x:z)$ induces the permutation $(A_1\ A_2)(A_3\ A_4)$ and fixes the points $(1:\xi:1+\xi)$ of the line $z=x+y$. Therefore, if  $A_5$ lies on this line, the group $H_{S}$ of all examples above contains a $2\times 2$-cycle. Write $A_5=(1:\xi:1+\xi)$; the equations of Proposition 
\refd{Prp:embedtwoquadrics} then become the two given above. Using the birational morphism of Proposition \refd{Prp:embedtwoquadrics}, we can transform any automorphism of\hspace{0.2 cm}$\Pn$ leaving the set of the five points invariant into an automorphism of $\mathbb{P}^4$ leaving the equations of $S$ invariant.
\begin{itemize}
\item
Suppose that $H_{S}$ contains a transposition, which we may take to be $(A_4 \ A_5)$. The automorphism must then be $(x:y:z) \mapsto (a x:b y:c z)$, for some $a,b,c \in \K^{*}$. As it sends $A_5=(a:b:c)$ on $(a^2:b^2:c^2)$, the point $A_5$ must be $(1:\pm 1: \pm 1)$, which is not possible as this is aligned with two of the four other points.
\item
Suppose that $H_{S}$ contains a $2\times 2$-cycle, namely $(A_1\ A_2)(A_3\ A_4)$. The automorphism is then $(x:y:z) \mapsto (z-y:z-x:z)$. The point $A_5$ must be one of the fixed points. These are points of the line $z=x+y$ and the point $(1:0:1)$. And $(1:0:1)$ is not possible, since this point is aligned with $A_1$ and $A_3$. We get then one of the cases in the list above.
\item
Suppose that $H_{S}$ contains a $3$-cycle, namely $(A_1\ A_2\ A_5)$. Denote $A_5$ by $(1:\xi:\mu)$. The automorphism must be of the form $(x:y:z) \mapsto (y:(1-\xi)x+\xi y:\mu y+(1-\mu) z)$, since it sends $A_1$ on $A_2$, $A_2$ on $A_5$ and fixes $A_3$ and $A_4$. As $A_5$ must be sent on $A_1$, we have $\xi^2-\xi+1=0$ and $\mu(\xi+1-\mu)=0$. And $\mu\not=0$, so $\mu=\xi+1$ and $\xi^2-\xi+1=0$. The two solutions are primitive $6$-th roots of unity, $e^{2\im \pi/6}$ and its inverse. We get then the two isomorphic surfaces (see below), given in $\nump{S4.3}$.
\item
If $H_{S}$ contains a $4$-cycle, namely $(A_1\ A_3\ A_2\ A_4)$, it must be the automorphism $(x:y:z) \mapsto (y:y-z:y-x)$. So $A_5$ must be one of its three fixed points. These are $(1:1:0)$, and $(1:\xi:1+\xi)$, where $\xi^2+1=0$. As the point $(1:1:0)$ is aligned with $A_1$ and $A_2$, there are only two possible surfaces, which are given in case $\nump{S4.4}$.
\item
Suppose that $H_{S}$ contains a $5$-cycle, namely $(A_1\ A_4\ A_3\ A_2\ A_5)$. Denote $A_5$ by $(1:\xi:\mu)$. The automorphism must be of the form $(x:y:z) \mapsto (x-y:x-\xi y+(\xi-1)z:x - \mu y)$. Since the image of $A_5$ is $(1-\xi:1-\xi^2+\mu\xi-\mu:1-\mu\xi)$ and must be $A_1$, this implies that $\mu=\xi^{-1}$ and $\xi^3-2\xi+1=(\xi-1)(\xi^2+\xi-1)=0$. As $\xi=1$ would give $A_5=A_4$, there are two possibilities for $A_5$: the points $(1:\xi:1+\xi)$, where $\xi^2+\xi-1=0$, given in $\nump{S4.5}$.
\end{itemize}
Note that the automorphism $(x:y:z)\mapsto (y:x:z)$ sends the points $A_1,A_2,A_3,A_4,(1:\xi:1+\xi)$ respectively on the points $A_2,A_1,A_3,A_4,(1:\xi^{-1}:1+\xi^{-1})$. Similarly, the automorphism $(x:y:z) \mapsto (y-z:y:y-x)$ sends the points $A_1,A_2,A_3,A_4,(1:\xi:1+\xi)$ respectively on the points $A_3,A_4,A_1,A_2,(1:-\xi:1-\xi)$.

If we look for other automorphisms (calculation), we see that the sets $\{A_1,A_2,A_3,A_4$, $(1:\xi:1+\xi)\}$ and $\{A_1,A_2,A_3,A_4,(1:\xi':1+\xi')\}$ are linearly equivalent if and only if $\xi'\in \{\xi,-\xi,\xi^{-1},-\xi^{-1}\}$.
\proofend
\end{proof}

\bigskip

\begin{remark}\fauxtitred As five points of the plane in general position determine one conic, an automorphism that leaves the five points invariant is an automorphism of the conic. As all the conics are equivalent, we can also prove Proposition \refd{Prp:ClassDelPezzoQSpec} by looking at subgroups of $\PGLn{2}$ that leave invariant a set of five points in $\mathbb{P}^1$.\end{remark}

\bigskip

\bigskip

We now know the structure of the group of automorphisms of all Del Pezzo surfaces of degree $4$. 

This allows us to prove an amazing result:
\begin{Lem}\fauxtitre
\label{Lem:EllRatRat}
Let $S$ be the Del Pezzo surface of degree $4$ viewed in $\mathbb{P}^4$ as the surface defined by the equations:
\begin{center}\begin{itemize}
\item
$(1+\xi)x_1^2-x_3^2-\xi x_4^2+\xi(1+\xi)x_5^2=(1+\xi)x_2^2-\xi x_3^2-x_4^2+\xi(1+\xi)x_5^2=0$
\end{itemize}\end{center}
for some $\xi \in \K\backslash\{0,\pm 1\}$.

Let $G\cong (\Z{2})^2$ be the group generated by 
\begin{center}
$\alpha:(x_1:x_2:x_3:x_4:x_5) \mapsto (x_1:x_2:x_3:x_4:-x_5)$,\\
$\beta:(x_1:x_2:x_3:x_4:x_5) \mapsto (x_2:x_1:x_4:x_3:-x_5)$.\end{center}

\index{Conic bundles!on Del Pezzo surfaces|idi}
Then, there exist two conic bundle structures on $S$ induced by some morphisms $\pi_1,\pi_2:S\rightarrow \mathbb{P}^1$, such that 
\begin{itemize}
\item
$G \subset \Aut(S,\pi_1)$, $G\subset \Aut(S,\pi_2)$.
\item
$\alpha$ acts trivially on both fibrations  {\upshape (}i.e.\ $\overline{\pi_1}(\alpha)=\overline{\pi_2}(\alpha)=1${\upshape )} and in each case twists the $4$ singular fibres.
\item
$\beta$ does not act trivially on the first one {\upshape (}i.e.\ $\overline{\pi_2}(\beta)\not=1${\upshape )} and acts trivially on the second one {\upshape (}$\overline{\pi_1}(\beta)=1${\upshape )}, twisting two singular fibres.
\end{itemize}
\end{Lem}
\begin{proof}\fauxtitred\upshape
Note that $S$ is the blow-up of $A_1,...,A_4,A_5=(1:\xi:1+\xi) \in \Pn$, and  $\Aut(S) =G_S \rtimes H_S$, where $G_S \cong (\Z{2})^4$ and $H_S$ contains the element $(A_1\ A_2)(A_3\ A_4)$ (see Proposition \refd{Prp:ClassDelPezzoQSpec}).

The element $\alpha$ is a cubic involution (see Proposition \refd{Prp:GeoCub}), and $\beta$ is $(A_1\ A_2)(A_3\ A_4)$ (see Proposition \refd{Prp:ClassDelPezzoQSpec}). Both elements preserve the conic bundles induced by the divisors $L-E_5$ and $2L-E_1-E_2-E_3-E_4$. Consider the action of $\alpha$ and $\beta$ on the singular fibres:

Firstly, the singular fibres of the conic bundle associated to $L-E_5$ are $\{E_1,D_{15}\}$,..., $\{E_4,D_{45}\}$. The automorphism $\alpha$ twists the four singular fibres and $\beta$ permutes the fibres $\{E_1,D_{15}\}$ and $\{E_2,D_{25}\}$, and also the fibres $\{E_3,D_{35}\}$ and $\{E_4,D_{45}\}$, hence does not act trivially on the fibration.

Secondly, the singular fibres of the conic bundle associated to $2L-E_1-E_2-E_3-E_4$ are $\{E_5,C\}$, $\{D_{12},D_{34}\}$, $\{D_{13},D_{24}\}$ and $\{D_{14},D_{23}\}$. Both $\alpha$ and $\beta$ leave these $4$ fibres invariant, and twist respectively $4$ and $2$ fibres.\proofend
\end{proof}
\begin{remark}\fauxtitred
This phenomenon is generalised in Proposition \refd{Prp:HardChgExacrSeq}.
\end{remark}
We shall now determine the abelian groups of automorphisms not contained in $G_{S}$ which satisfy the minimality condition.
\begin{Prp}\fauxtitre
Let $S$ be a Del Pezzo surface of degree $4$ and $H \subset \Aut(S)$ an abelian group not contained in $G_{S}$ such that $\rkPic{S}^{H}=1$. Then, up to a birational conjugation, $S$ is the surface in $\mathbb{P}^4$ defined by equations
\begin{center}\begin{itemize}
\item
$(1+\xi)x_1^2-x_3^2-\xi x_4^2+\xi(1+\xi)x_5^2=0$,
\item
$(1+\xi)x_2^2-\xi x_3^2-x_4^2+\xi(1+\xi)x_5^2=0$,
\end{itemize}\end{center}
where $\xi \in \K \backslash \{0,1,-1\}$,
and one of the following situations arises:\vspace{0.1 cm}

\begin{tabular}{ll}
${\num{4.42}}$ & $H\cong \Z{4} \times \Z{2}$, for any $\xi\in \K\backslash \{0,\pm 1\}$ \\ & generated by $((1,0,0,0,1),(1 \ 2)(3 \ 4))$ and $(1,1,1,1,0) \in \Aut(S)$,\\ & which correspond to $(x_1:x_2:x_3:x_4:x_5) \mapsto (-x_2:x_1:x_4:x_3: - x_5)$\\
& and $(x_1:x_2:x_3:x_4:x_5) \mapsto (x_1:x_2:x_3:x_4:-x_5)$.\\
\\
${\num{4.223}}$ & $H\cong \Z{2}\times\Z{2}\times \Z{3}$, \\
 & $S$ is the surface $\nump{S4.3}$ of Proposition $\refd{Prp:ClassDelPezzoQSpec}$ {\upshape(}$\xi=e^{2\im\pi/6}${\upshape)}.\\ & $H$ is generated by $(1,1,1,0,1),(1,1,0,1,1)$ and $(1\ 2\ 5) \in \Aut(S)$,\\ & which correspond to $(x_1:x_2:x_3:x_4:x_5) \mapsto (x_1:x_2:\pm x_3:\pm x_4:x_5)$\\
& and $(x_1:x_2:x_3:x_4:x_5) \mapsto ((\xi-1)x_5:x_1:(\xi-1)x_3:-\xi x_4:-\xi x_2)$.\\ 
\\
${\num{4.8}}$ & $H\cong \Z{8}$, $S$ is the surface $\nump{S4.4}$ of Proposition $\refd{Prp:ClassDelPezzoQSpec}$ {\upshape(}$\xi=\im${\upshape)}.\\ & $H$ is generated by $((0,0,0,1,1),(1 \ 3 \ 2 \ 4)) \in \Aut(S)$,\\ & which correspond to $(x_1:x_2:x_3:x_4:x_5) \mapsto (x_4:x_3:-x_1:x_2:\xi x_5)$.
\end{tabular}

\end{Prp}\begin{proof}\fauxtitred\upshape
Let $\rho:\Aut(S)\rightarrow \Sym_5$ be the homomorphism given by the action on the exceptional pairs and suppose that $H \subset \Aut(S)$ is an abelian group with $\rkPic{S}^{H}=1$ and that  $\rho(H)$ is not the trivial group.
From Proposition \refd{Prp:ClassDelPezzoQSpec}, the group $\rho(H)$ is cyclic, generated by a $2\times 2$-cycle or an $n$-cycle, where $3\leq n\leq 5$. As $H$ is abelian, the group $H \cap G_{S}$ must be fixed by the action of $\rho(H)$. We will enumerate the cases and use the exact sequence
\begin{center}
$1 \rightarrow H \cap G_{S} \rightarrow H \stackrel{\rho}{\rightarrow} \rho(H) \rightarrow 1$
\end{center}
to specify $H$.

Note that the product of two elements $(\alpha,\beta), (\alpha',\beta') \in \Aut(S)\subset (\mathbb{F}_2)^4 \rtimes \Sym_5$ is $(\alpha+\beta(\alpha'), \beta\beta')$ and the conjugate of an element $(\alpha,\beta)$ by $(\mu,id)$ is $(\mu + \beta(\mu) + \alpha,\beta)$.
\begin{itemize}
\item
If $\rho(H)$ is generated by a $2\times 2$-cycle, namely $(1 \ 2)(3\ 4)$, then $H \cap G_{S} \subset V=\{ (a,a,b,b,0) \ |\ a,b \in \mathbb{F}_2\}\cong (\Z{2})^2$. Let $g=((a,b,c,d,e),(1\ 2)(3\ 4)) \in H$ be an element such that $\rho(g)=(1\ 2)(3\ 4)$. As $\rkPic{S}^{H}=1$ and $V$ fixes the divisor $L-E_5$, we may suppose that $e=1$. (Otherwise, the group $H$ would fix $L-E_5$.) We conjugate the group by the element $((0,b,0,d,b+c),id)$ and may assume that $g=((a+b,0,c+d,0,1),(1\ 2)(3\ 4))$. As $a+b+c+d+e=0$, we have $g=((\alpha,0,1+\alpha,0,1),(1\ 2)(3\ 4))$, where $\alpha=a+b=c+d+1 \in \mathbb{F}_2$.

If $\alpha=1$, then $g=((1,0,0,0,1),(1\ 2)(3\ 4))$ and $g^2=((1,1,0,0,0),id) \in V$. As $g$ fixes the divisor $2L-E_3-E_4$, $H$ cannot be equal to $<g>$. It must thus contain all of $V$ and in particular the element $(1,1,1,1,0)$. This gives case $\nump{4.42}$.

If $\alpha=0$, then $g=((0,0,1,0,1),(1\ 2)(3\ 4))$. We conjugate the group by the automorphism $(x:y:z) \mapsto (y-z:y:y-x)$ of\hspace{0.2 cm}$\Pn$ that corresponds to the permutation $(A_1\ A_3)(A_2\ A_4)$, to get the preceding case, with an another point $A_5$.
\item
If $\rho(H)$ is generated by a $3$-cycle, namely $(1\ 2 \ 5)$, then $H \cap G_{S}$ must be a subgroup of $V= \{ (a,a,b,a+b,a) \ | \ a,b \in \mathbb{F}_2\}$, which is the centraliser of $\Aut(S)$. The order of $H$ must be a multiple of $4$, by Lemma \refd{Lem:SizeOrbits}, hence $H \cap G_{S} =V$. As $H$ contains $V$ and an element $(\alpha, (1\ 2 \ 5))$, it also contains an element $((a,b,0,0,c),(1 \ 2 \ 5))$. We conjugate the group by $((c,a,0,0,b),id)$ to obtain the group generated by $V$ and $(id,(1 \ 2 \ 5))$, and get case $\nump{4.223}$.
\item
If $\rho(H)$ is generated by a $4$-cycle, say $(1\ 3\ 2\ 4)$, then $H\cap G_{S}$ must be a subgroup of the centraliser of $\Aut(S)$, which we denote by $V=<(1,1,1,1,0)>$. Let \begin{center}$g=((a,b,c,d,e),(1\ 3\ 2\ 4)) \in H$\end{center} be such that $\rho(g)=(1\ 3\ 2\ 4)$. We conjugate the group by $((0,b+c,c,a,a+c),id)$ and may suppose that $g=((0,0,0,e,e),(1\ 3\ 2\ 4))$. 

If $e=1$, then $g^4=(1,1,1,1,0)$, so $V \subset <g>$. The group $H$ is then cyclic, generated by $g$, of order $8$. Note that $\Pic{S}^{g^4}=\mathbb{Z}K_S \oplus \mathbb{Z}L-E_5$. Since $g$ sends $L-E_5$ on $-K_S-L+E_5$,
we have $\Pic{S}^{H}=\mathbb{Z}K_S$. This gives case $\nump{4.8}$

If $e=0$, the element $g$ belongs to $H_{S}$, so it fixes the divisors $L$ and $E_5$. As the group $V$ fixes $L-E_5$, the rank of $\Pic{S}^{H}$ cannot be $1$.
\item
If $\rho(H)$ is generated by a $5$-cycle, $H \cap G_{S}= id$, so the group $H$ has order $5$ and the rank of $\Pic{S}^{H}$ cannot be $1$, by Lemma \refd{Lem:SizeOrbits}.\proofend
\end{itemize}
\end{proof}

\bigskip

Let us now summarise our results:
\begin{Prp}\fauxtitre
\label{Prp:DelP4}
Let $G \subset \Aut(S)$ be an finite abelian group of automorphisms of a Del Pezzo surface $S$ of degree $4$, such that $\rkPic{S}^G=1$. 
\index{Del Pezzo surfaces!automorphisms}
\begin{itemize}
\item[\upshape 1.]
Up to isomorphism, 
$S\subset \mathbb{P}^4$ is the blow-up of $(1:0:0)$, $(0:1:0)$, $(0:0:1)$, $(1:1:1)$, $(a:b:c) \in \Pn$, of equations
\begin{center}$\begin{array}{cccccc}cx_1^2&&-ax_3^2&+(a-c)x_4^2&-ac(a-c)x_5^2&=0,\vspace{0.1cm}\\
&cx_2^2&-bx_3^2&+(c-b)x_4^2&-bc(c-b)x_5^2&=0.\end{array}$\end{center}
and one of the following occurs:
\begin{itemize}
\item[$\nump{4.22}$]
$G\cong (\Z{2})^2$, given by $G=\{(x_1:x_2:x_3:x_4:x_5) \mapsto (\pm x_1:\pm x_2:x_3:x_4:x_5)\}$.
\item[$\nump{4.222}$]
$G\cong (\Z{2})^3$, given by $G=\{(x_1:x_2:x_3:x_4:x_5) \mapsto (\pm x_1:\pm x_2:\pm x_3:x_4:x_5)\}$.
\item[$\nump{4.2222}$]
$G\cong (\Z{2})^4$, given by $G=\{(x_1:x_2:x_3:x_4:x_5) \mapsto (\pm x_1:\pm x_2:\pm x_3:\pm x_4:x_5)\}$.
\item[$\nump{4.42}$]
$(a:b:c)=(1:\xi:1+\xi)$, for any $\xi\in \K\backslash \{0,\pm 1\}$\\
$G\cong \Z{4} \times \Z{2}$,  generated by $(x_1:x_2:x_3:x_4:x_5) \mapsto (-x_2:x_1:x_4:x_3: - x_5)$
 and $(x_1:x_2:x_3:x_4:x_5) \mapsto (x_1:x_2:x_3:x_4:-x_5)$.
\item[$\nump{4.223}$]
$(a:b:c)=(1:\zeta_6:1+\zeta_6)$, where $\zeta_6=e^{2\im\pi/6}$ \\
$G\cong (\Z{2})^2\times\Z{3}$, generated  by
$(x_1:x_2:x_3:x_4:x_5) \mapsto (x_1:x_2:\pm x_3:\pm x_4:x_5)$
and $(x_1:x_2:x_3:x_4:x_5) \mapsto ((\zeta_6-1)x_5:x_1:(\zeta_6-1)x_3:-\zeta_6 x_4:-\zeta_6 x_2).$
\item[$\nump{4.8}$]
$(a:b:c)=(1:\im:1+\im)$ \\
$G\cong \Z{8}$ generated by $(x_1:x_2:x_3:x_4:x_5) \mapsto (x_4:x_3:-x_1:x_2:\im x_5)$
\end{itemize}
\item[\upshape 2.]
$G$ leaves invariant a pencil of rational curves of $S$ if and only if the pair $(G,S)$ is conjugate to $\nump{4.22}$, $\nump{4.223}$ or $\nump{4.8}$.
\end{itemize}
\end{Prp}
\begin{proof}\fauxtitred\upshape
The first assertion follows directly from the work done in this chapter.
Let us prove the second one. 
\begin{itemize}
\item
Firstly, suppose that $G$ fixes one point $P$ of $S$. Let us see that $P$ does not lie on an exceptional curve. Remember that the orbits of $G$ on the set of exceptional curves have size at least $4$ (see Lemma \refd{Lem:SizeOrbits}). If $P$ lies on some exceptional curve, it then belongs to at least $4$ of them. But no $4$ exceptional divisors intersect in the same point.

Then, the blow-down of five skew exceptional divisors of $S$ gives a birational morphism (not $G$-equivariant) $\pi:S\rightarrow \Pn$ which is the blow-up of $5$ points in general position (no $3$ collinear) and is an isomorphism of some neighbourhood of $P$ on some neighbourhood of $\pi(P)\in \Pn$.   The linear system of cubics of\hspace{0.2 cm}$\Pn$ passing through the $5$ blown-up points and with multiplicity two at $\pi(P)$ is then a pencil of rational curves invariant by $G$. (Recall that the linear system associated to the anti-canonical divisor is the system of cubics passing through the $5$ blown-up points.)

Observe now that the groups of automorphisms of pairs $\nump{4.22}$ and $\nump{4.223}$ and $\nump{4.8}$ fix some point of $S$.
 \begin{itemize}
 \item[$\nump{4.22}$]
Suppose that $(G,S)$ is the pair $\nump{4.42}$. The set $\Delta$ of points of the surface fixed by the group $\{(x_1:x_2:x_3:x_4:x_5) \mapsto (\pm x_1:\pm x_2:x_3:x_4:x_5)\}$ is the trace of the plane $x_1=x_2=0$ on the surface; this is the intersection of two conics. Note that $\Delta$ is also the set of points of $S$ fixed by the involution $(x_1:x_2:x_3:x_4:x_5) \mapsto (-x_1:-x_2:x_3:x_4:x_5)$ and therefore consists of four points (see Proposition \refd{Prp:GeoQuad}).
\item[$\nump{4.223}$]
Suppose that $(G,S)$ is the pair $\nump{4.223}$. As above, the subgroup $H=\{(x_1:x_2:x_3:x_4:x_5) \mapsto (x_1:x_2:\pm x_3:\pm x_4:x_5)\}$ of $G$ fixes exactly four points of the surface. The element of $G$ which is  $(x_1:x_2:x_3:x_4:x_5) \mapsto ((\zeta_6-1)x_5:x_1:(\zeta_6-1)x_3:-\zeta_6 x_4:-\zeta_6 x_2)$, and has order $3$  thus fixes at least one of the four points fixed by $H$. This implies that the entire group $G$ fixes a point of $S$.
\item[$\nump{4.8}$]
Suppose that $(G,S)$ is the pair $\nump{4.8}$. In this case, $G$ is generated by the element $g:(x_1:x_2:x_3:x_4:x_5) \mapsto (x_4:x_3:-x_1:x_2:\im x_5)$, that fixes the point $p=(\alpha^3:\alpha:1:\alpha^2:0)=(\alpha^2:1:-\alpha^3:\alpha:0)=g(p)$ with $\alpha=e^{2\im\pi/8}$.
\end{itemize}
\item
Secondly, if $G$ leaves invariant a pencil of rational curves, it is birationally conjugate to some group $\tilde{G}\subset \Aut(\tilde{S},\pi)$ of automorphisms of a conic bundle. We suppose that this is the case and get a contradiction, if $(G,S)$ is one of the pairs $\nump{4.222}$, $\nump{4.2222}$, and $\nump{4.42}$. Note that any cubic involution acts trivially on the fibration, since its set of fixed points is an elliptic curve (see Lemma \ref{Lem:CurveFixDeJ}). \index{Curves!elliptic!birational maps that fix an elliptic curve|idi}
\begin{itemize}
\item
If $(G,S)$ is one of the pairs $\nump{4.222}$ or $\nump{4.2222}$, the group  contains at least three cubic involutions, so a subgroup $\tilde{H}\subset \tilde{G}$ isomorphic to $(\Z{2})^3$ acts trivially on the fibration. This is not possible (see Section \refd{Sec:Exactsequence}).
\item
Suppose $(G,S)$ is the pair $\nump{4.42}$. The group $G$ is generated by a cubic involution  $\alpha$ and an element $\beta$ whose square is a quadratic involution:
\begin{center}
$\begin{array}{llll}
\alpha:&(x_1:x_2:x_3:x_4:x_5)& \mapsto &(x_1:x_2:x_3:x_4:-x_5),\\
\beta:&(x_1:x_2:x_3:x_4:x_5)& \mapsto& (-x_2:x_1:x_4:x_3: - x_5),\\
\beta^2:&(x_1:x_2:x_3:x_4:x_5)& \mapsto &(-x_1:-x_2:x_3:x_4:x_5).\end{array}$
\end{center}
Note that $\beta^2$ fixes $4$ points of the surface, whose union is the trace of the plane $x_1=x_2=0$ on $S$ (see above). But none of the four points is fixed by $\alpha$ (whose set of fixed points is the elliptic curve which is the trace of the hyperplane $x_5=0$). The group $H\cong (\Z{2})^2$ generated by $\alpha$ and $\beta^2$ therefore fixes no point of $S$. 

Denoting by $\tilde{\alpha},\tilde{\beta}$ the elements of $\tilde{G}$ corresponding respectively to $\alpha$ and $\beta$, the group $\tilde{H}=<\tilde{\alpha},\tilde{\beta}^2>$ also fixes no point of $\tilde{S}$ (see Section \refd{Sec:ExistenceFixedPoints}). 

Note that any element of $\Aut(\tilde{S},\pi)$ leaves invariant the set of singular points of singular fibres of $\pi$. As $\tilde{H}$ fixes no point of $\tilde{S}$, $\tilde{\beta}^2$ acts non-trivially on the fibration and the two fibres invariant by it are smooth. On each of the two fibres, $\alpha$ fixes exactly two points. As $\alpha$ and $\beta$ commute, $\beta$ acts on the set of these two points, so $\beta^2$ fixes both. This contradicts the fact that $\tilde{H}$ fixes no point of $\tilde{S}$. 
\proofend 
\end{itemize}
\end{itemize}
\end{proof}

\bigskip

\bigskip

\section{Del Pezzo surfaces of degree $3$: cubic surfaces}
\index{Del Pezzo surfaces!of degree $3$|idb}
\pagetitless{Finite abelian groups of automorphisms of Del Pezzo surfaces of degree $\leq 4$}{Del Pezzo surfaces of degree $3$: cubic surfaces}
\label{Sec:Cubic}
Recall that a Del Pezzo surface of degree~3 is obtained by the blow-up of $6$ points $A_1,...,A_6$ in general position in $\Pn$ (no $3$ of them collinear and not all $6$ on a conic). Let $\pi: S\rightarrow \Pn$ be the blow-up of the $A_i$'s. The map $\varphi$ associated to the linear system of cubics passing through $A_1,...,A_6$ gives rise to an embedding of $S$ as a non-singular cubic surface in $\mathbb{P}^3$ (see \cite{bib:Kol}, Theorem III.3.5); let $S$ denote also the image of $\varphi$.
\begin{center}
\hspace{0 cm}\xymatrix{ S\ar_{\pi}[d] \ar[drrr] & & \\
\Pn\ar@{-->}^{\varphi}[rrr]& & & S \subset \mathbb{P}^3 }
\end{center}

The exceptional divisors on $S$ are (see Proposition 
\refd{Prp:DescExcCurvDelPezzo}):
\index{Del Pezzo surfaces!curves on the surfaces|idi}
\begin{itemize}
\item
$E_1=\pi^{-1}(A_1),...,E_6=\pi^{-1}(A_6)$, the $6$ pull-backs of the points $A_1$,...,$A_6$.
\item
$D_{ij}=\tilde{\pi^{-1}}(A_iA_j)$ for $1\leq i,j\leq 6,\ i\not= j$, the $15$ strict pull-backs of 
the lines passing through $2$ of the $A_i$'s.
\item
$C_{i}$ for $1\leq i\leq 6$, the strict pull-backs of the $6$
conics of\hspace{0.2 cm}$\Pn$ passing through $5$ of the $A_i$'s. (Note that $C_{i}$ will denote the conic \textit{not} passing  through $A_i$).
\end{itemize}

There are thus $27$ exceptional divisors, and each of them intersects $10$ others with multiplicity $1$. 
Note that any exceptional divisor of $S\subset \mathbb{P}^3$ is a line of $\mathbb{P}^3$ contained in the surface. We obtain again the famous result on the number of lines on a smooth cubic surface of $\mathbb{P}^3$ (see for example \cite{bib:Seg}).

The groups of automorphisms of non-singular cubic surfaces have been studied by B. Segre \cite{bib:Seg} and later by Hosoh \cite{bib:Ho2} and the list is now exhaustive.
We will restrict ourselves to the abelian case; moreover for our purpose we may assume that the fixed part of the Picard group is of rank one. This will considerably simplify the calculations. We will then give the conjugacy classes of these groups in the Cremona group $\CrP$.

\bigskip

\begin{remark}\fauxtitred
We could use directly the lists of groups of automorphisms of cubic surfaces given in \cite{bib:Ho2}, describe the subgroups, and find minimal actions on the Picard group. In fact, the gain in time is not really significant. We thus prefer not to use the classification.
\end{remark}

\bigskip

Let us recall that every automorphism of a cubic surfaces comes from an automorphism of $\mathbb{P}^3$. We will denote by $\DiaG{\alpha}{\beta}{\gamma}{\delta}$ the automorphism $(w:x:y:z) \mapsto (\alpha w:\beta x:\gamma y: \delta z)$.

\bigskip
\index{Del Pezzo surfaces!automorphisms}
\begin{Prp}\titreProp{automorphisms of cubic surfaces}\\
\label{Prp:Cubics}
Let $G$ be an abelian group of automorphisms of a non-singular cubic surface $S=V(F)$, such that $\rkPic{S}^{G}=1$.
Then, up to isomorphism, we are in one of the following cases:

\begin{tabular}{|llll|}
\hline
{\it name} & {\it structure} & {\it generators} & {\it equation of}\\
{\it of $\mathit (G,S)$} & {\it of $\mathit G$} & {\it of $\mathit G$} & $\mathit F$\\
\hline 
${\num{3.3}}$ & $\Z{3}$& $\DiaG{\omega}{1}{1}{1}$& $w^3+L_3(x,y,z) \hts$ \\
${\num{3.6.1}}$ & $\Z{6} $& $\DiaG{\omega}{1}{1}{-1}$& $w^3+x^3+y^3+xz^2+\lambda yz^2 \hts$ \\
${\num{3.33.1}}$ & $(\Z{3})^2$& $\DiaG{\omega}{1}{1}{1}$, $\DiaG{1}{1}{1}{\omega}$& $w^3+x^3+y^3+z^3\hts$ \\
${\num{3.9}}$ & $\Z{9}$& $\DiaG{\zeta_9}{1}{\omega}{\omega^2}$& $w^3+xz^2+x^2y+y^2z\hts$  \\
${\num{3.33.2}}$ & $(\Z{3})^2$& $\DiaG{\omega}{1}{1}{1}$, $\DiaG{1}{1}{\omega}{\omega^2}$& $w^3+x^3+y^3+z^3+\lambda xyz\hts$  \\
${\num{3.12}}$ & $\Z{12}$& $\DiaG{\omega}{1}{-1}{\im}$& $w^3+x^3+yz^2+y^2x \hts$ \\
${\num{3.36}}$ & $\Z{3}\times\Z{6}$& $\DiaG{\omega}{1}{1}{1}$, $\DiaG{1}{1}{-1}{\omega}$& $w^3+x^3+xy^2+z^3 \hts$ \\
${\num{3.333}}$ & $(\Z{3}) ^3$&$\DiaG{\omega}{1}{1}{1}$, $\DiaG{1}{\omega}{1}{1}$& $w^3+x^3+y^3+z^3\hts$ \\
& &\mbox{and } $\DiaG{1}{1}{\omega}{1}$ &  \\
${\num{3.6.2}}$ & $\Z{6}$ &$\DiaG{1}{-1}{\omega}{\omega^2}$ & $wx^2+w^3+y^3+z^3+\lambda wyz\hts$\\
\hline
\end{tabular} \vspace{0.1 cm}\\
where $L_3$ denotes a non-singular form of degree~3, $\lambda \in \K$ is a parameter such that the surface is non-singular and $\omega=e^{2\im\pi/3},\zeta_9=e^{2\im\pi/9}$. \\
Furthermore, all the cases above are minimal pairs $(G,S)$ with $\rkPic{S}^G=1$.
\end{Prp}
\begin{remark}
Some of the surfaces given above seems to be different but are in fact the same. Indeed, we choosed to put the automorphisms on diagonal form, the equations of the surface may then not be standard.\end{remark}
\begin{proof}\fauxtitred\upshape
Let $G$ be an abelian group of automorphisms of a non-singular cubic surface $S=V(F)$, such that $\rkPic{S}^{G}=1$. 
As $3$ divides the size of the orbits of $G$ on the $27$ lines (see Lemma \refd{Lem:SizeOrbits}), it also divides the order of $G$ which must thus contain an element of order $3$.

\break

There are three kinds of elements of order $3$ in $\PGLn{4}$, depending on the nature of their eigenvalues. Setting as usual $\omega=e^{2\im\pi/3}$, there are elements with one eigenvalue of multiplicity $3$ (conjugate to $\DiaG{1}{1}{1}{\omega}$, or its inverse), elements with two eigenvalues of multiplicity $2$ (conjugate to $\DiaG{1}{1}{\omega}{\omega}$) and elements with three distinct eigenvalues (conjugate to $\DiaG{1}{1}{\omega}{\omega^2}$). We will examine the possibilities for $G$, depending on its elements of order $3$.

\textit{Case a: $G$ contains an element $g$ of order $3$ with an eigenvalue of multiplicity $3$.}\\
We may suppose that $g=\DiaG{\omega}{1}{1}{1}$ and use Lemma \refd{Lem:Cubgexseq}. Let $\Gamma$ denote the elliptic curve fixed by $g$ and 
\begin{center}$1 \rightarrow <g> \rightarrow G \rightarrow H \rightarrow 1$\end{center}be the exact sequence given by this lemma, where $H$ is an abelian subgroup of $\PGLn{3}$ stabilizing the curve $\Gamma$.

If $H$ is the identity, we get $\nump{3.3}$. Otherwise, we use Lemma \refd{Lem:CubP2} to identify the possibilities for $H$ and pull-back all these possibilities to get all the cases given in the proposition, except for the last one. 
If $H$ comes from $\nump{C2}$, we get naturally that $G$ belongs to case $\nump{3.6.1}$. Similarly, we get $\nump{3.33.1}$, $\nump{3.9}$, $\nump{3.33.2}$, $\nump{3.12}$, $\nump{3.36}$ and $\nump{3.333}$ respectively from $\nump{C3.1}$, $\nump{C3.2}$, $\nump{C3.3}$, $\nump{C4}$, $\nump{C6}$ and $\nump{C33.1}$.

\bigskip

\textit{Case b: $G$ contains an element $g$ of order $3$ with two eigenvalues of multiplicity $2$.}\\
With a suitable choice of coordinates, we may assume that $g=\DiaG{1}{1}{\omega}{\omega}$, and that $S$ is the cubic Fermat surface, $F=w^3+x^3+y^3+z^3$. 

The group of automorphisms of $S$ is $(\Z{3})^3 \rtimes \Sym_4$ and the centraliser of $g$ in it is $(\Z{3})^3 \rtimes V$, where $V\cong (\Z{2})^2$ is the subgroup of $\Sym_4$ generated by the two transpositions $(w,\ x)$ and $(y,\ z)$. The structure of this group gives rise to an exact sequence
\begin{center}$\begin{array}{ccccccc}
1\rightarrow &(\Z{3})^3 &\rightarrow& (\Z{3})^3 \rtimes V& \stackrel{\gamma}{\rightarrow}& V &\rightarrow 1 \\
& \cup & & \cup & & \cup \\
1\rightarrow & G \cap (\Z{3})^3& \rightarrow& G& \rightarrow &\gamma(G)& \rightarrow 1.
\end{array}$\end{center}
We may suppose that $G$ contains no element of order $3$ with an eigenvalue of multiplicity $3$, since this case has been examined above. There are then three possibilities for $G \cap (\Z{3})^3$, namely $<g>$, $<g,\DiaG{1}{\omega}{1}{\omega}>$ and $<g,\DiaG{1}{\omega}{\omega}{1}>$. The last is conjugate to the second by the automorphism $(y,\ z)$.

Note that $g$ leaves invariant exactly $9$ of the $27$ lines on the surface; these are $\{w+\omega^i x= y+ \omega^j z=0\}$, for $0\leq i,j \leq 2$. If $G \cap (\Z{3})^3$ is equal to $<g>$, the group $G$ leaves invariant at least one of the $9$ lines, since the order of $\gamma(G)$ is divisible by $4$, so $\rkPic{S}^G>1$.

If $G \cap (\Z{3})^3$ is the group $H=<g,\DiaG{1}{\omega}{1}{\omega}>$, we have $G=H$, since the centraliser of $H$ in $(\Z{3})^3 \rtimes V$ is the group $(\Z{3})^3$. As the set of three skew lines $\{w+\omega^i x=\omega^i y+ z=0\}$ for $0\leq i \leq 2$ is an orbit of $H$, the rank of $\Pic{S}^{G}$ is strictly larger than $1$.

\break

\textit{Case c: $G$ contains an element $g$ of order $3$ with three distinct eigenvalues.}\\
We may suppose that \begin{center}$g=\DiaG{1}{1}{\omega}{\omega^2}$. \end{center}
We then have $F(w,x,\omega y,\omega^2 z)=\omega^i F(w,x,y,z)$, where $i=0,1$ or $2$. The different cases give respectively the equations $F_0,F_1$ and $F_2$ of $S$ listed below:
\begin{center}\begin{itemize}
\item
$F_0=L_3(w,x)+L_1(w,x)yz+y^3+z^3$,
\item
$F_1=y^2L_1(w,x)+L_2(w,x) z +yz^2$,
\item
$F_2=z^2L_1(w,x)+L_2(w,x) y +y^2z$,
\end{itemize}\end{center}
where $L_i$ denotes a form of degree $i$.
Note that every subgroup of $\PGLn{4}$ isomorphic to $(\Z{3})^2$ contains an element with only two distinct eigenvalues, so we assume that $G$ contains only two elements of order $3$, which are $g$ and $g^2$.

The group cannot be reduced to $<g>$. By the Lefschetz formula, the trace of $g$ on $H^{*}(S,\mathbb{Z})$ is the Euler characteristic of its set of fixed points, and hence is strictly larger than $0$, since the set of points of $S$ fixed by $g$ is finite and contains the intersection of $S$ with the line $y=z=0$. The trace of $g$ on $H^{2}(S,\mathbb{Z})=\Pic{S}$ is then strictly larger than $-2$, hence the eigenvalue $1$ must have multiplicity strictly larger than $1$, whence $\Pic{S}^{g}\geq 1$.

The group $G$ thus contains another element, say $h$, of order $\not=1,3$. As it must commute with $g$, we may diagonalise it and suppose that\begin{center} $h=\DiaG{\mu}{\nu}{\xi}{1}$.\end{center}

If $S=V(F_0)$, we deduce first of all that $\xi^3=1$. As $S$ must be non-singular, the equation must be of degree $\geq 2$ in each variable, so either $x^3$ or $wx^2$ must appear in the equation and the same is true of $w^3$ or $w^2x$.
\begin{itemize}
\item
If both $w^3$ and $x^3$ appear, then $h^3=id$, which is excluded.
\item
If $w^2x$ and $wx^2$ appear, we have $\mu\nu^2=\mu^2\nu=1$, so $h$ is once again of order $3$. \item
By permuting $w$ and $x$ if necessary, we suppose that $w^3$ and $wx^2$ appear, which implies that $\mu^3=1$ and $\nu= \pm \mu$. To avoid $h^3=id$, we must take $\nu=-\mu$, so $h$ is of order $2$ or $6$. If $h$ is of order $2$, it is $\DiaG{1}{-1}{1}{1}$. If it of order $6$, as $h^2$ must be $g$ or $g^2$, we get the same group of order $6$, generated by $\DiaG{1}{-1}{\omega}{\omega^2}$. The equation of the surface is then $F_0=wx^2+w^3+y^3+z^3+\lambda wyz$, for some $\lambda \in \K$. This gives $\nump{3.6.2}$.\end{itemize}

If $S=V(F_1)$, the non-singularity of $S$ implies that $w^2z$, $x^2z$ and $yz^2$ appear in the equation, so 
we have $\mu^2=\nu^2=\xi$. The non-singularity of $S$ implies that either $wy^2$ or $xy^2$ appears in $F_1$. By a permutation $(w,\ x)$, we may suppose that $wy^2$ appears, so $\mu \xi^2=\mu^2=\xi$, whence $\xi=\mu^{-1}$ and $\xi^3=\mu^3=1$. To avoid $h^3=1$ we must once again choose $\nu=-\mu$, so $h$ is of order $2$ or $6$. In both cases we get a group of order $6$, generated by $\DiaG{1}{-1}{\omega}{\omega^2}$ and the equation of the surface is then $F_2=wy^2+w^2z+x^2z+yz^2$. As the line $y=z=0$ is fixed by the group, $\rkPic{S}^{G}>1$.

If $S=V(F_2)$, a change of coordinates $(y,\ z)$ takes us back to the case of $F_1$. 
\proofend
\end{proof}

\bigskip

\break

\begin{Lem}\fauxtitre
\label{Lem:Cubgexseq}
Let $G$ be an abelian group of automorphisms of a non-singular cubic surface $S=V(F)$, which contains the element $g(w:x:y:z)=(\omega w:x:y:z)$, where $\omega=e^{2\im\pi/3}$. We then have: \index{Curves!elliptic!birational maps that fix an elliptic curve}
\begin{itemize}
\item[\upshape 1.]
The points of $S$ fixed by the action of $g$ form a non-singular plane cubic curve 
\begin{center} $\Gamma=\{(0:x:y:z) \in \mathbb{P}^3 \ | \ F(0,x,y,z)=0\}$.\end{center}

\item[\upshape 2.]
The group $G$ leaves the plane $w=0$ invariant and its action on this plane gives rise to an exact sequence 
\begin{center}$1 \rightarrow <g> \rightarrow G \rightarrow H \rightarrow 1$,\end{center} where $H$ is a diagonalizable subgroup of $\PGLn{3}$ which leaves the curve $\Gamma$ invariant.
\item[\upshape 3.]
$\rkPic{S}^{G}=1$.
\end{itemize}
\end{Lem}\begin{proof}\fauxtitred\upshape
The set of points of $S$ fixed by $g$ is the intersection of $S$ with the points of $\mathbb{P}^3$ fixed by $g$, that is the plane $w=0$. The first part of the lemma is now clear.

As the group $G$ is abelian, every element of $G$ must leave invariant  the curve $\Gamma$ and hence the plane $w=0$. Elements of $\PGLn{4}$ that fix every point of $\Gamma$ must fix the plane $w=0$ and thus belong to $<g>$. Let us see that $H$ is diagonalizable. Indeed, the commutator of $g$ in $\PGLn{4}$ is the group of elements of the form $\left(\begin{array}{l|l}\lambda & \\ \hline & M\end{array}\right)$, where $\lambda \in \K^{*}$ and $M \in GL(3,\K)$. The abelian subgroups of this groups are  diagonalizable, so $H$ is also.

The third part of the lemma is given by the Lefschetz formula. As the trace of $g$ on $H^{*}(S,\mathbb{Z})$ is the Euler characteristic of its set of fixed points, it must be $0$. The trace of $g$ on $H^{2}(S,\mathbb{Z})=\Pic{S}$ is then $-2$ and its eigenvalues are $1$, three times $\omega$ and $\omega^2$. This gives $\Pic{S}^{g}=1$ and consequently $\Pic{S}^{G}=1$.\proofend
\end{proof}

\bigskip

\begin{Lem}\fauxtitre
\label{Lem:CubP2}
Let $H$ be a non-trivial diagonalizable subgroup of $\PGLn{3}$ and $\Gamma=\{ (x:y:z) \in \Pn \ |\ F(x,y,z)=0\}$ be a smooth cubic curve of\hspace{0.2 cm}$\Pn$ invariant by $H$. Then, up to conjugation in $\PGLn{3}$, $H$ and $\Gamma$ are in the following list:
\begin{center}\begin{tabular}{|llll|}
\hline
{\it name} & {\it structure} & {\it generator} & {\it equation $\mathit F$}\\
{\it of $\mathit H$} & {\it of $\mathit H$} & {\it of $\mathit H$} & {\it of $\mathit \Gamma$}\\
\hline 
${\num{C2}}$ & $\Z{2}$& $\Diag{1}{1}{-1}$& $x^3+y^3+xz^2+\lambda yz^2\hts$ \\
${\num{C3.1}}$ & $\Z{3}$& $\Diag{1}{1}{\omega}$& $x^3+y^3+z^3\hts$ \\
${\num{C3.2}}$ & $\Z{3}$& $\Diag{1}{\omega}{\omega^2}$& $xz^2+x^2y+y^2z$ \\
${\num{C3.3}}$ & $\Z{3}$& $\Diag{1}{\omega}{\omega^2}$& $x^3+y^3+z^3+\lambda xyz\hts$ \\
${\num{C4}}$ & $\Z{4}$& $\Diag{1}{-1}{\im}$& $x^3+yz^2+y^2x\hts$ \\
${\num{C6}}$ & $\Z{6}$& $\Diag{1}{-1}{\omega}$& $x^3+xy^2+z^3$ \\
${\num{C33.1}}$ & $(\Z{3})^2$& $\Diag{1}{1}{\omega}$, $\Diag{1}{\omega}{1}$& $x^3+y^3+z^3\hts$ \\
\hline
\end{tabular}\end{center}
where $\lambda \in \K$ is a parameter such that the curve is smooth, and $\omega=e^{2\im\pi/3}$.
\end{Lem}\begin{proof}\fauxtitred\upshape
We will use frequently that $F$ is of degree $\geq 2$ in each variable, as $\Gamma$ is non-singular and we will also feel free to simplify the equation by scaling the variables. 

Let us begin by assuming that the group $H$ is cyclic, generated by an element $h$ of order $n>1$.
Choosing a coordinate system, we may suppose that $h$ is diagonal.

\bigskip

\textit{Suppose first that two of the eigenvalues of $h$ are equal.}\\
We may suppose that $h=\Diag{1}{1}{\zeta_n}$, where $\zeta_n=e^{2\im\pi/n}$, for some integer $n>1$. We then write $F=\lambda z^3+z^2 L_1(x,y) +z L_2(x,y) + L_3(x,y)$, where $\lambda \in \K$ and $L_i$ is a form of degree $i$. As $\Gamma$ is non-singular, we have $L_3 \not= 0$, and $F(x,y,\alpha z)=F(x,y,z)$. It follows that either $n=2$ or $n=3$ and, up to a change of coordinates, we get respectively $\nump{C2}$ and $\nump{C4}$.

\bigskip

\textit{Suppose now that the three eigenvalues of $h$ are different}.\\
If $h$ is of order $3$, we may assume that $h=\Diag{1}{\omega}{\omega^2}$, then $F(x,\omega y,\omega^2 z)=\omega^i F(x,y,z)$. If $i=0$ or $i=1$, we get respectively cases $\nump{C3.3}$ and $\nump{C3.2}$. If $i=2$, then $F=xy^2+x^2z+yz^2$, so we get $\nump{C3.2}$ by a permutation of the coordinates $y$ and $z$.

If none of the points $(1:0:0),(0:1:0),(0:0:1)$ lies on $\Gamma$, the equation of $F$ contains $x^3,y^3$ and $z^3$, which implies that $h$ is of order $3$, a case examined just above. Let us then suppose that $(1:0:0)\in \Gamma$, so $F=x^2 \cdot L_1(y,z)+ x \cdot L_2(y,z)+L_3(y,z)$, where $L_i$ denotes a form of degree $i$. As the point $(1:0:0)$ is non-singular, either $x^2y$ or $x^2z$ appears in $F$. Up to a permutation of coordinates, we suppose that $x^2y$ appears.
Write $h=\Diag{1}{\mu}{\nu}$, where $1,\mu,\nu$ are all distinct;
we find that $F\circ h (x,y,z)=\mu F(x,y,z)$. 
Using once again the non-singularity of $\Gamma$, we see that either $y^3$ or $y^2z$ appears in the equation (but not $y^2x$, since $\mu^2\not=\mu$). Similarly, either $z^3$ or $z^2y$ or $z^2x$ appears in the equation.
\begin{itemize}
\item
If $y^3$ appears, then $\mu^2=1$, so $\mu=-1$. This implies that $z^2x$ does not appear, since $\nu^2\not=1$.
\begin{itemize}
\item
If $z^3$ appears, then $\nu^3=-1$, so $\nu=-\omega$ or $-\omega^2$ and $F=x^2y+y^3+z^3$. We get $\nump{C6}$ by exchanging the coordinates $x$ and $y$.
\item
If $z^2y$ appears, then $\nu^2=-1$, so $\nu=\pm \im$ and $F=xz^2+x^2y+y^3$; we get $\nump{C4}$ by exchanging $x$ and $y$.
\end{itemize}
\item
If $y^2z$ appears then $\mu^2\nu=\mu$, so $\nu=\mu^{-1}$. Then $z^2y$ does not appear in the equation, since $\nu^2\not=1$.
\begin{itemize}
\item
If $z^3$ appears, then $\nu^3=\mu$, so $\{\mu,\nu\}=\{\im,-\im\}$ and $F=x^2y+y^2z+z^3$. We get $\nump{C4}$ by exchanging $x$ and $z$.
\item
If $z^2x$ appears, then $\nu^2=\mu$, so both $\nu^3=\mu^3$ and $h$ is of order $3$, a case already considered.
\end{itemize}
\end{itemize}
We have given all the cases where $H$ is cyclic; let us give those where $H$ is not cyclic.

The diagonal groups isomorphic to $\Z{n}\times \Z{n}$ contain the element $\Diag{1}{1}{\zeta_n}$, where $\zeta_n=e^{2\im\pi/n}$, so $n\leq 3$. The case $n=2$ is not possible (each variable must be only of degree $2$), and the case $n=3$ gives $\nump{C33.1}$.

If $H$ is diagonal and non-cyclic, it must contain a group $\Z{n}\times\Z{n}$ and as the Fermat cubic contains no other diagonal automorphism than $\Z{3} \times \Z{3}$, we have seen all possibilities.\proofend
\end{proof}

\break

\section{Automorphisms of Del Pezzo surfaces of degree $2$}
\pagetitless{Finite abelian groups of automorphisms of Del Pezzo surfaces of degree $\leq 4$}{Del Pezzo surfaces of degree $2$}\index{Del Pezzo surfaces!of degree $2$|idb}
\label{Sec:DP2}
A Del Pezzo surface of degree $2$\index{Del Pezzo surfaces!of degree $2$|idb} is given by the blow-up of $7$ points of\hspace{0.2 cm}$\Pn$ in general position
(no six of them are on a conic and no three are collinear). Let $A_1,...,A_7 \in \Pn$ be the $7$ points 
and $\pi: S\rightarrow \Pn$ the blow-down. The system of cubics passing through the $7$ points gives 
a morphism $\eta:S\rightarrow\Pn$ whose linear system is $|-K_S|$ and which is a double covering ramified over a smooth quartic curve $\Gamma$ of\hspace{0.2 cm}$\Pn$ (\cite{bib:BaB}, \cite{bib:deF}, \cite{bib:Dol}) :
\begin{center}
\hspace{0 cm}\xymatrix{ S\ar_{\pi}[d] \ar^{\eta}[drrr] & & \\
\Pn\ar@{-->}^{C^3(A_1,...,A_7)}[rrr]& & & \mathbb{P}^2. }
\end{center}
The surface $S$ can then be viewed in the weighted projective space $\mathbb{P}(2,1,1,1)$, as given by the equation $w^2=F(x,y,z)$, where $F$ is the equation of the quartic $\Gamma$ (see \cite{bib:Kol}, Theorem III.3.5 and \cite{bib:KoSC}, Corollary 3.54).

As every automorphism of $S$ must leave the quartic $w=F(x,y,z)=0$ invariant, we have a homomorphism $\rho:\Aut(S) \rightarrow \Aut(\Gamma)$. Let $\sigma$ denote the Geiser involution $\sigma(w,x,y,z)=(-w,x,y,z)$.\index{Involutions!Geiser|idb}

\begin{Prp}\fauxtitre
\label{Prp:exSDP2}
The sequence
\begin{equation}\label{exSDP2}
1 \rightarrow <\sigma> \rightarrow \Aut(S) \stackrel{\rho}{\rightarrow} \Aut(\Gamma) \rightarrow 1
\end{equation} 
is exact. 
\end{Prp}\begin{proof}\fauxtitred\upshape
Clearly, the kernel of $\rho$ is $<\sigma>$. Let us prove that $\rho$ is surjective. 
Note that any automorphism of $\Gamma$ comes from an automorphism of\hspace{0.2 cm}$\Pn$ since $\Gamma$ is canonical. 
Let $\alpha \in \Aut(\Gamma) \subset \PGLn{3}$. As 
$F \circ \alpha=\lambda F$, for some $\lambda \in \K^{*}$, the maps $(w:x:y:z) \mapsto (\pm \sqrt{\lambda}:\alpha(x,y,z))$ are two automorphisms of $S$ that $\rho$ sends on $\alpha$.
\proofend
\end{proof}

\bigskip

\begin{remark}\fauxtitred
In fact, the exact sequence of Proposition \refd{Prp:exSDP2} splits, which implies that  \begin{center}$\Aut(S)\cong <\sigma> \times \Aut(\Gamma)$.\end{center} We received a proof of this beautiful result from Arnaud Beauville in a private communication. We won't include it here, since it is not needed for our classification. 
\end{remark}

\bigskip
\index{Del Pezzo surfaces!curves on the surfaces|idi}
The exceptional divisors on $S$ are (see Proposition 
\refd{Prp:DescExcCurvDelPezzo}):
\begin{itemize}
\item
$E_1=\pi^{-1}(A_1),...,E_7=\pi^{-1}(A_7)$, the $7$ pull-backs of the points $A_1$,...,$A_7$.
\item
$D_{ij}=\tilde{\pi^{-1}}(A_iA_j)$ for $1\leq i,j\leq 7,\ i\not= j$, the $21$ strict pull-backs of 
the lines passing through $2$ of the $A_i$'s.
\item
$C_{ij}$ for $1\leq i,j\leq 7,\ i\not= j$, the strict pull-backs of the $21$
conics of\hspace{0.2 cm}$\Pn$ passing through $5$ of the $A_i$'s. (Note that  $C_{ij}$ will denote the conic \textit{not} passing through $A_i$ or $A_j$.)
\item
$F_1,...,F_7$, the strict pull-backs of the $7$ cubics of\hspace{0.2 cm}$\Pn$ passing through $A_1,...,A_7$, 
with multiplicity two at one of the $A_i$'s.
\end{itemize}
There are thus $56$ exceptional divisors, and each of them intersects $27$ others with multiplicity $1$ and one with
multiplicity $2$ (which is $-K_S$ minus itself and is sent by $\eta$ on the same line of\hspace{0.2 cm}$\Pn$).
Note that the image of every exceptional divisor by
the double covering $\eta$ is a line bitangent to the quartic $\Gamma$. 
The $56$ exceptional divisors of $S$ are then sent by the projection $\eta$ on the 
$28$ bitangents to the curve $\Gamma$.

\bigskip

Let us give a description of the Geiser involution:\index{Involutions!Geiser!matrix|idb}
\begin{Lem}\fauxtitred
\label{Lem:Geiser}
\begin{itemize}
\item[\upshape 1.]
The Geiser involution $\sigma$ switches the pull-backs 
of any point of\hspace{0.2 cm}$\Pn$ by the $2$-covering $\eta$.
\item[\upshape 2.]
It exchanges the two exceptional divisors $D$ and $-K_S-D$, which have the same image by $\eta$,
for any exceptional divisor $D$. 

Explicitly, its orbits are :
\begin{itemize}
\item
$E_i$ and $F_i$, for $i=1,...,7$,
\item
$C_{ij}$ and $D_{ij}$, for $i,j=1,...,7$,\hspace{0.15cm} $i \not= j$.
\end{itemize}
\item[\upshape 3.]
Its action on $\Pic{S}$ is given by the matrix
\begin{center}$\left(\begin{array}{rrrrrrrr}
-2 & -1 & -1 & -1 & -1 & -1 & -1 & -3 \\
-1 & -2 & -1 & -1 & -1 & -1 & -1 & -3 \\
-1 & -1 & -2 & -1 & -1 & -1 & -1 & -3 \\
-1 & -1 & -1 & -2 & -1 & -1 & -1 & -3 \\
-1 & -1 & -1 & -1 & -2 & -1 & -1 & -3 \\
-1 & -1 & -1 & -1 & -1 & -2 & -1 & -3 \\
-1 & -1 & -1 & -1 & -1 & -1 & -2 & -3 \\
3 & 3 & 3 & 3 & 3 & 3 & 3 & 8\end{array}\right)$\end{center}
relative to the basis $(E_1,E_2,...,E_7,L)$, where $L$ is the pull-back by $\pi$ of a line of\hspace{0.2 cm}$\Pn$.
\item[\upshape 4.]
Its eigenvalues on $\Pic{S}$ are $1$ with multiplicity $1$ and $-1$ with multiplicity $7$.
\item[\upshape 5.]
$\rkPic{S}^{\sigma}=1$, so $\sigma$ acts minimally on the surface $S$.
\end{itemize}
\end{Lem}
\begin{proof}\fauxtitred\upshape
The first part follows from the definition of $\sigma$ as $\sigma(w,x,y,z)=(-w,x,y,z)$. We then get 
the second part directly, which implies the third. 

Note that $\sigma(K_S)=K_S$ and $\sigma(K_S+2E_i)=K_S+2(-K_S-E_i)=-(K_S+2E_i)$, for $i=1,...,7$. Hence, the action of $\sigma$ on the basis $(K_S,K_S+2E_1,K_S+2E_2,...,K_S+2E_7)$ of the vector space $\Pic{S}\otimes_{\mathbb{Z}} \mathbb{Q}$ is diagonal with eigenvalues $<1,(-1)^7>$. We thus get  assertion $4$, which implies the last one.
\proofend\end{proof}

\bigskip

Let us denote by $\DiaG{\alpha}{\beta}{\gamma}{\delta}$ the automorphism 
$(w:x:y:z) \mapsto (\alpha w:\beta x:\gamma y:\delta z)$ of $\mathbb{P}(2,1,1,1)$. We see that $\DiaG{\alpha}{\beta}{\gamma}{\delta}$ is sent by
$\eta$ on the automorphism $\Diag{\beta}{\gamma}{\delta}$ of\hspace{0.2 cm}$\Pn$.

\begin{Prp}\titreProp{automorphisms of Del Pezzo surfaces of degree $2$}\\
\label{Prp:DelP2}
\index{Del Pezzo surfaces!automorphisms}
Let $S$ be a Del Pezzo surface of degree $2$ given by the equation $w^2=F(x,y,z)$ and let $G$ 
be an abelian subgroup of $\Aut(S)$ such that $\rkPic{S}^{G}=1$. Then, up to isomorphism we 
are in one of the following cases:
\begin{itemize}
\item
groups containing the Geiser involution $\sigma$:\\
$G=<\sigma>\times H$, where $H \subset
\PGLn{3}$ is an abelian group of automorphisms of the quartic curve $\Gamma \subset \Pn$ defined by the equation $F(x,y,z)=0$:\\
\index{Involutions!Geiser}
\begin{tabular}{|llll|}
\hline
{\it name} & {\it structure} & {\it generators} & {\it equation of}\\
{\it of $\mathit (G,S)$} & {\it of $\mathit G$} & {\it of $\mathit G$} & $\mathit F$\\
\hline 
${\num{2.G}}$ & $\Z{2}$& $\sigma$& $L_4(x,y,z)\hts$ \\
${\num{2.G2}}$ & $(\Z{2})^2$&$\DiaG{\pm 1}{1}{1}{-1}$& $L_4(x,y)+L_2(x,y)z^2+z^4\hts$ \\
${\num{2.G3.1}}$ & $\Z{6}$&$\DiaG{-1}{1}{1}{\omega}$& $L_4(x,y)+z^3L_1(x,y)\hts$ \\
${\num{2.G3.2}}$ & $\Z{6}$&$\DiaG{-1}{1}{\omega}{\omega^2}$& $x(x^3\h{0.4}+\h{0.4}y^3\h{0.4}+\h{0.4}z^3)\h{0.5}+\h{0.5}yzL_1(x^2,yz)\hts$ \\
${\num{2.G4.1}}$ & $\Z{2}\times\Z{4}$&$\DiaG{\pm 1}{1}{1}{\im}$& $L_4(x,y)+z^4\hts$ \\
${\num{2.G4.2}}$ & $\Z{2}\times\Z{4}$&$\DiaG{\pm 1}{1}{-1}{\im}$& $x^3y+y^4+z^4+xyL_1(xy,z^2)\hts$ \\
${\num{2.G6}}$ & $\Z{2}\times\Z{6}$&$\DiaG{\pm 1}{\omega}{1}{-1}$& $x^3y+y^4+z^4+\lambda y^2z^2\hts$ \\
${\num{2.G7}}$ & $\Z{14}$&$\DiaG{-1}{\h{0.5}\zeta_7\h{0.5}}{\h{0.5}(\zeta_7)^4\h{0.5}}{\h{0.5}(\zeta_7)^2}$& $x^3y+y^3z+xz^3\hts$ \\
${\num{2.G8}}$ & $\Z{2}\times\Z{8}$&$\DiaG{\pm 1}{\zeta_8}{-\zeta_8}{1}$& $x^3y+xy^3+z^4\hts$ \\
${\num{2.G9}}$ & $\Z{18}$&$\DiaG{-1}{(\zeta_9)^6}{1}{\zeta_9}$& $x^3y+y^4+xz^3\hts$\\
${\num{2.G12}}$ & $\Z{2}\times\Z{12}$&$\DiaG{\pm 1}{\omega}{1}{\im}$& $x^3y+y^4+z^4\hts$ \\
${\num{2.G22}}$ & $(\Z{2})^3$&$\DiaG{\pm 1}{1}{\pm 1}{\pm 1}$ & $L_2(x^2,y^2,z^2)\hts$ \\
${\num{2.G24}}$ & $(\Z{2})^2\times \Z{4}$&$\DiaG{\pm 1}{1}{\pm 1}{\im}$ & $x^4+y^4+z^4+\lambda x^2y^2\hts$ \\
${\num{2.G44}}$ & $\Z{2}\times (\Z{4})^2$&$\DiaG{\pm 1}{1}{1}{\im}$& $x^4+y^4+z^4$\\
& & and $\DiaG{1}{1}{\im}{1}$& \\
\hline
\end{tabular}
\item
groups not containing the Geiser involution:

\begin{tabular}{|llll|}
\hline
{\it name} & {\it structure} & {\it generators} & {\it equation of}\\
{\it of $\mathit (G,S)$} & {\it of $\mathit G$} & {\it of $\mathit G$} & $\mathit F$\\
\hline 
${\num{2.4}}$ & $\Z{4}$& $\DiaG{1}{1}{1}{\im}$& $L_4(x,y)+z^4\hts$ \\
${\num{2.6}}$ & $\Z{6}$& $\DiaG{-1}{\omega}{1}{-1}$& $x^3y+y^4+z^4+\lambda y^2z^2\hts$ \\
${\num{2.12}}$ & $\Z{12}$& $\DiaG{1}{\omega}{1}{\im}$& $x^3y+y^4+z^4\hts$ \\
${\num{2.24.1}}$ & $\Z{2}\times \Z{4}$& $\DiaG{1}{1}{1}{\im}$, $\DiaG{1}{1}{-1}{1}$& 
$x^4+y^4+z^4+\lambda x^2y^2\hts$ \\
${\num{2.24.2}}$ & $\Z{2}\times \Z{4}$& $\DiaG{1}{1}{1}{\im}$, $\DiaG{-1}{1}{-1}{1}$& 
$x^4+y^4+z^4+\lambda x^2y^2\hts$ \\
${\num{2.44.1}}$ & $(\Z{4})^2$& $\DiaG{1}{1}{1}{\im}$, $\DiaG{1}{1}{\im}{1}$& 
$x^4+y^4+z^4\hts$ \\
${\num{2.44.2}}$ & $(\Z{4})^2$& $\DiaG{1}{1}{1}{\im}$, $\DiaG{-1}{1}{\im}{1}$& 
$x^4+y^4+z^4\hts$\\
\hline
\end{tabular}
\end{itemize}
where $\lambda,\mu \in \K$ are parameters and the $L_i$ are forms of degree $i$, such that the surfaces are smooth. Furthermore, all the cases above are minimal pairs $(G,S)$ with $\rkPic{S}^G=1$.
\end{Prp}
\begin{remark}
Some of the surfaces given above seems to be different but are in fact the same. Indeed, we choosed to put the automorphisms on diagonal form, the equations of the surface may then not be standard.\end{remark}
\begin{proof}\fauxtitred\upshape
The image of $G$ by the homomorphism $\rho:\Aut(S) \rightarrow \Aut(\Gamma) \subset \PGLn{3}$ is an abelian subgroup $H\subset\PGLn{3}$ that leaves $\Gamma$ invariant. We use Lemma \refd{Lem:QuarP2} to identify the possibilities for $H$ and to calculate $\rho^{-1}(H)$. We see in each case that this group is isomorphic to $<\sigma>\times H$ (see the remark that follows Proposition \ref{Prp:exSDP2}). Then, either $G=\rho^{-1}(H)$ or $G$ does not contain $\sigma$.

A simple calculation of the pull-backs of the groups of Lemma \refd{Lem:QuarP2} gives the first list of this proposition, which illustrates the cases where $G=\rho^{-1}(H)$, or equivalently the groups which contains the Geiser involution.

\bigskip

We now assume that $\sigma \notin G$. By Proposition
\refd{Prp:exSDP2}, the group $G$ is isomorphic to $\rho(G)$. The fact that
$\rkPic{S}^G=1$ will restrict the possibilities. Firstly, the order of
$G$ must be divisible by $2$, by Lemma \refd{Lem:SizeOrbits}. Secondly,
the description of automorphisms given in Lemma \refd{Lem:actionNSGeis} will 
help us to decide whether or not $\rkPic{S}^G=1$.

We first suppose that $G$ is cyclic, generated by an
automorphism $\alpha$. Then, $\rho(\alpha)$ is given by cases
$\nump{Q4.1}, \nump{Q4.2}, \nump{Q6},\nump{Q8}$ or $\nump{Q12}$ of Lemma $\refd{Lem:QuarP2}$ and the order of $\alpha$ can be $2$, $4$, $6$, $8$ or $12$. We compute the possibilities for $\alpha$ and use Lemma \refd{Lem:actionNSGeis} to determine its eigenvalues on $\Pic{S}$. Indeed, the rank of $\Pic{S}^\alpha$ is equal to $1$ if and only if the multiplicity of the eigenvalue $1$ is $1$. 

\begin{itemize}
\item
\textit{If the order of $\alpha$ is $2$} and $\rkPic{S}^{\alpha}=1$, 
then for any exceptional divisor $R$,
$R+\alpha(R)$ is a
multiple of $K_S$, which implies that $\alpha(R)$ must be equal to
$-K_S-R$ (using the description of the $56$ exceptional divisors).
As this occurs only when $\alpha=\sigma$, 
the rank of $\Pic{S}^{\alpha}$ is 
strictly larger than $1$.
(We can also verifiy this using Lemma \refd{Lem:actionNSGeis}).
\item
\textit{If the order of $\alpha$ is $4$}, the possibilities for $\alpha \in \Aut(S)$ are: 
\begin{itemize}
\item
$\nump{Q4.1}$: $\alpha=\DiaG{\pm 1}{1}{1}{\im}$, the equation of the surface is $w^2=L_4(x,y)+z^4$,
\item
$\nump{Q4.2}$: $\alpha=\DiaG{\pm 1}{1}{-1}{\im}$, the equation of the surface is $w^2=x^4+y^4+z^4+xyL_1(xy,z^2)$.
\end{itemize}
By Lemma \refd{Lem:actionNSGeis}, $\Pic{S}^{\alpha}=1$
if and only if $\alpha=\DiaG{1}{1}{1}{\im}$. We therefore get $\nump{2.4}$.
\item
\textit{If the order of $\alpha$ is $\geq 6$}, the possibilities for $\alpha \in \Aut(S)$ are: 
\begin{itemize}
\item
$\nump{Q6}$: $\alpha=\DiaG{\pm 1}{\omega}{1}{-1}$, the surface's equation is $w^2=x^3y+y^4+z^4+\lambda y^2z^2$.\\
Lemma \refd{Lem:actionNSGeis} shows that $\rkPic{S}^{\alpha}=1$ if and only if $\alpha=\DiaG{-1}{\omega}{1}{-1}$. We therefore get $\nump{2.6}$
\item
$\nump{Q8}$: $\alpha=\DiaG{\pm 1}{\zeta_8}{-\zeta_8}{1}$, the equation of the surface is $w^2=x^3y+xy^3+z^4$.\\
Lemma \refd{Lem:actionNSGeis} shows that $\rkPic{S}^{\alpha}\geq 2$ in both cases.
\item
$\nump{Q12}$: $\alpha=\DiaG{\pm 1}{\omega}{1}{\im}$, the equation of the surface is $w^2=x^3y+y^4+z^4$. \\
The automorphism $\DiaG{-1}{\omega}{1}{\im}$ does not act minimally on $S$, since it leaves
invariant the two exceptional divisors $y=0,w=\pm z^2$. However, the automorphism
$\alpha=\DiaG{1}{\omega}{1}{\im}$ acts minimally on $S$, since
$\alpha^9=\DiaG{1}{1}{1}{\im}$ does (case $\nump{2.4}$ analysed below). We get then $\nump{2.12}$.
\end{itemize}
\end{itemize}

\break

We have seen all the cases when the group $G$ is cyclic. Let us now investigate a non-cyclic abelian group $G$. We still assume that $G$ does not contain the Geiser involution $\sigma$, 
that the order of the group is
divisible by $2$ (Lemma \refd{Lem:SizeOrbits}), and that $G$ is isomorphic to $\rho(G)$ (Proposition
\refd{Prp:exSDP2}).
\begin{itemize}
\item
\textit{If the group is isomorphic to $\Z{2}\times\Z{2}$}, the group $\rho(G)$ is given by case $\nump{Q22}$, and consists of elements of the form
$\Diag{1}{\pm 1}{\pm 1}$; the surface has equation $w^2=L_2(x^2,y^2,z^2)$, where $L_2$ is a form of degree $2$. There are two kinds of involutions in the group: 
\begin{itemize}
\item
Involutions of the type
$\DiaG{1}{\pm 1}{\pm 1}{\pm 1}$, with eigenvalues $<1^4,(-1)^4>$ (see Lemma \refd{Lem:actionNSGeis}).
\item
Involutions of the type
$\DiaG{-1}{\pm 1}{\pm 1}{\pm 1}$, with eigenvalues $<1^5,(-1)^3>$ (see Lemma \refd{Lem:actionNSGeis}).
\end{itemize}
Note that the action of $G$ on $\Pic{S}$ can be diagonalised, since the group is abelian.
As we assumed that the Geiser involution does not belong to our group, we see that $\rkPic{S}^G>1$, by computing the product of two such elements and its eigenvalues.
\item
\textit{If the group is isomorphic to $\Z{2}\times\Z{4}$}, the group $\rho(G)$ is given by case $\nump{Q24}$,
generated by $\Diag{1}{1}{\im}$ and $\Diag{1}{-1}{1}$, and the equation of the surface is $w^2=x^4+y^4+z^4+\lambda x^2y^2$, for
$\lambda \in \K$. 
\begin{itemize}
\item
If $G$ contains the automorphism $\alpha=\DiaG{1}{1}{1}{\im}$, then $\rkPic{S}^G=1$, since $\rkPic{S}^\alpha=1$ (see above). This yields cases $\nump{2.24.1}$ and $\nump{2.24.2}$.
\item
If $G$ does not contain $\DiaG{1}{1}{1}{\im}$, it is generated by $\alpha=\DiaG{-1}{1}{1}{\im}$ and $\beta=\DiaG{\pm 1}{1}{-1}{1}$. Note that $\alpha$
acts on $\Pic{S}$ with eigenvalues $<1^4,\im^2,(-\im)^2>$. The action of $G$ on $\Pic{S}$ is diagonalizable, and every eigenvalue of $\beta$ is $\pm 1$. If $\rkPic{S}^G=1$, then $\beta$ must acts with eigenvalue $-1$ on a linear space of dimension $3$ or $\Pic{S}^{\alpha}$. The eigenvalues of $\alpha\beta$ must therefore be $<1,(-1)^3,\im^2,-\im^2>$. But these eigenvalues are those of an element conjugate to $\Diag{1}{1}{1}{\im}$ (see Lemma \refd{Lem:actionNSGeis}), so we get the previous case.
\end{itemize}
\item
\textit{If the group is isomorphic to $\Z{4}\times\Z{4}$}, the group $\rho(G)$ is given by case $\nump{Q44}$,
generated by $\Diag{1}{1}{\im}$ and $\Diag{1}{\im}{1}$ and the equation is $w^2=x^4+y^4+z^4$.
Up to a change of coordinates, there are exactly two possibilities for $G$, both with $\rkPic{S}^G=1$, given by $\nump{2.44.1}$ and $\nump{2.44.2}$.\proofend\end{itemize}
\end{proof}

\break

\begin{Lem}\titreProp{action on $\Pic{S}$ of automorphisms of $S$}\\
\label{Lem:actionNSGeis}
The following table gives the eigenvalues of the action of some automorphisms of $S$ on $\Pic{S}$:

\begin{center}$\begin{array}{|lll|}
\hline
{\it automorphism} & \mbox{\it equation }\mathit{F}\mbox{\it\ of the} & \mbox{\it eigenvalues of }\mathit{\alpha}\\
\mathit{\alpha} & \mbox{\it surface }\mathit{w^2=F(x,y,z)} & \mbox{\it on }\mathit{\Pic{S}}\\

\hline 
\DiaG{1}{1}{1}{-1} & L_4(x,y)+L_2(x,y)z^2+z^4 & <1^4,{-1}^4>\hts\\
\DiaG{-1}{1}{1}{-1} & L_4(x,y)+L_2(x,y)z^2+z^4 & <1^5,{-1}^3>\hts\\
\hline
\DiaG{1}{1}{1}{\omega} & L_4(x,y)+L_1(x,y)z^3 & <1^2,(\omega)^3,(\omega^2)^3>\hts\\
\DiaG{1}{1}{\omega}{\omega^2} & x(x^3+y^3+z^3)+yzL_1(x^2,yz) & <1^4,(\omega)^2,(\omega^2)^2>\hts\\
\hline
\DiaG{1}{1}{1}{\im} & L_4(x,y)+z^4 & <1,{-1}^3,\im^2,{-\im}^2>\hts\\
\DiaG{-1}{1}{1}{\im} & L_4(x,y)+z^4 & <1^4,\im^2,{-\im}^2>\hts\\
\DiaG{1}{1}{-1}{\im} & x^4+y^4+z^4+xyL_1(xy,z^2) & <1^3,{-1},\im^2,{-\im}^2>\hts\\
\DiaG{-1}{1}{-1}{\im} & x^4+y^4+z^4+xyL_1(xy,z^2) & <1^2,{-1}^2,\im^2,{-\im}^2>\hts\\
\hline
\DiaG{1}{\omega}{1}{-1} & x^3y+y^4+z^4+\lambda y^2z^2 &
<1^2,(\omega),(\omega^2),(-\omega)^2,(-\omega^2)^2>\hts\\
\DiaG{-1}{\omega}{1}{-1} & x^3y+y^4+z^4+\lambda y^2z^2 & <1,-1,(\omega)^2,(\omega^2)^2,(-\omega),(-\omega^2)>\hts\\
\hline
\DiaG{1}{\zeta_8}{-\zeta_8}{1} & x^3y+xy^3+z^4 &
<1^3,-1,\zeta_8,\zeta_8^3,\zeta_8^5,\zeta_8^7>\hts\\
\DiaG{-1}{\zeta_8}{-\zeta_8}{1} & x^3y+xy^3+z^4 & <1^2,(-1)^2,\zeta_8,\zeta_8^3,\zeta_8^5,\zeta_8^7>\hts\\
\hline
\end{array}$\end{center}
where $\lambda,\mu \in \K$ are parameters, the $L_i$ are forms of degree $i$ such that the surface is smooth, and $\zeta_n=e^{2\im\pi/n},\omega=\zeta_3$.
The sequence  $<\alpha_1^{k_1},...,\alpha_r^{k_r}>$ signifies that the
automorphism has eigenvalue $\alpha_i$ with multiplicity $k_i$ for $i=1,...,r$. 
\end{Lem}

\begin{proof}\fauxtitred\upshape
We will use the Lefschetz formula: the trace of an automorphism on $H^*(S,\mathbb{Z})$
is the Euler characteristic of its set of fixed points. Then, the trace on $H^2(S,\mathbb{Z})=\Pic{S}$ is equal to this
number minus $2$. Note that if we have the eigenvalues of the automorphism $g=\DiaG{1}{\alpha}{\beta}{\gamma}$,
for $\alpha,
\beta, \gamma \in \K$, we get the eigenvalues of $\sigma g=\DiaG{-1}{\alpha}{\beta}{\gamma}$ by multiplying every eigenvalue by $-1$ (except one equal to $1$, corresponding to the canonical divisor). We don't need then to calculate the fixed points of each automorphism (we did this anyway, as a verification).

Note that any point $(w:x:y:z) \in S$ fixed by the automorphism $\DiaG{\lambda}{\alpha}{\beta}{\gamma}$ is sent by the projection $\eta$ on a point $(x:y:z) \in \Pn$ fixed by the automorphism $\Diag{\alpha}{\beta}{\gamma}$. To get the fixed points of $\DiaG{\lambda}{\alpha}{\beta}{\gamma}$, we compute the fixed points of $\Diag{\alpha}{\beta}{\gamma}$ (which are $(1:0:0)$, $(0:1:0)$ and $(0:0:1)$, and the line $x=0$ if $\beta=\gamma$, the line $y=0$ if $\alpha=\gamma$ and the line $z=0$ if $\alpha=\beta$) and see whether the pull-backs of these points are fixed by $\DiaG{\lambda}{\alpha}{\beta}{\gamma}$.

\begin{itemize}
\item
\textit{automorphisms of order $2$}\\
The trace on $\Pic{S}$ of elements of order $2$ gives its eigenvalues directly. The set of points of the surface $w^2=L_4(x,y)+L_2(x,y)z^2+z^4$ fixed by the automorphism $\DiaG{1}{1}{1}{-1}$ is the disjoint union of the elliptic curve $z=0, w^2=L_4(x,y)$\index{Curves!elliptic!birational maps that fix an elliptic curve|idi} and
the two points $(\pm 1:0:0:1)$, so the trace on $\Pic{S}$ is zero and hence the multiplicity of both $1$ and $-1$ is $4$.
\item
\textit{automorphisms of order $3$}\\
As $\omega$ and $\omega^2$ have the same multiplicity, the trace on $\Pic{S}$ directly gives the multiplicities of $1,\omega,\omega^2$.
\begin{itemize}
\item
The set of points of the surface $w^2=L_4(x,y)+L_1(x,y)z^3$ fixed by the automorphism $\DiaG{1}{1}{1}{\omega}$ is the disjoint union of the elliptic curve $z=0,w^2=L_4(x,y)$ \index{Curves!elliptic!birational maps that fix an elliptic curve|idi} and the point $(0:0:0:1)$. Its trace on $\Pic{S}$ is then $-1$ and therefore the eigenvalues are $<1^2,(\omega)^3,(\omega^2)^3>$.
\item
The set of points of the surface $w^2=x(x^3+y^3+z^3)+yzL_1(x^2,yz)$ fixed by the automorphism $\DiaG{1}{1}{\omega}{\omega^2}$ consists of 
the $4$ points $(\pm 1:1:0:0)$, $(0:0:1:0)$ and $(0:0:0:1)$, hence the trace on $\Pic{S}$ is $2$ and the eigenvalues are $<1^4,(\omega)^2,(\omega^2)^2>$.
\end{itemize}
\item
\textit{automorphisms of order $4$}\\
The square of an element of order $4$ is of the same kind as that of $\DiaG{1}{1}{1}{-1}$, whose eigenvalues are $<1^4,(-1)^4>$. Therefore, every automorphism of order $4$ has eigenvalues $\im$ and $-\im$ of multiplicity $2$.
Computing the trace gives the multiplicities of $1$ and $-1$.

\begin{itemize}
\item
The set of points of the surface $w^2=L_4(x,y)+z^4$ fixed by the automorphism $\DiaG{1}{1}{1}{\im}$ is the elliptic curve $z=0,w^2=L_4(x,y)$\index{Curves!elliptic!birational maps that fix an elliptic curve|idi} and thus its trace on $\Pic{S}$ is $-2$. This means that the multiplicities of $1$ and $-1$ are respectively $1$ and $3$.
\item
The set of points of the surface $w^2=x^4+y^4+z^4+xyL_1(xy,z^2)$ fixed by the automorphism $\DiaG{1}{1}{-1}{\im}$ consists of $4$ points, which are $(\pm 1:1:0:0)$ and $(\pm 1:0:1:0)$. The trace on $\Pic{S}$ is $2$ and the multiplicities of $1$ and
$-1$ are respectively $3$ and $1$.
\end{itemize}
\item
\textit{automorphisms of order $6$}\\
Let $S$ be the surface with equation $w^2=x^3y+y^4+z^4+\lambda y^2z^2$ and let us write the eigenvalues of $g=\DiaG{1}{\omega}{1}{-1}$ on $\Pic{S}$ as $<1^a,(-1)^b,(\omega)^c,(\omega^2)^c,(-\omega)^d,(-\omega^2)^d>$, for some non-negative integers $a,b,c,d$. The eigenvalues of $g^2=\DiaG{1}{\omega^2}{1}{1}$ are $<1^2,(\omega)^3,(\omega^2)^3>$ and those of $g^3=\DiaG{1}{1}{1}{-1}$ are $<1^4,(-1)^4>$ (see above). The first relation gives $a+b=2$ and the second $a+2c=4$. As $a\geq 1$ we get $a=2$, $b=0$, $c=1$. The eigenvalues of $g$ are then $<1^2,(\omega),(\omega^2),(-\omega)^2,(-\omega^2)^2>$. We can verify the Lefschetz formula: the set of points of $S$
fixed by the automorphism $g=\DiaG{1}{\omega}{1}{-1}$ consists of the $5$ points $(0:1:0:0)$, $(\pm 1: 0:1:0)$ and $(\pm 1:0:0:1)$ and the trace of $\alpha$ on $\Pic{S}$ is $3$. 
\item
\textit{automorphisms of order $8$}\\
The set of points of the surface $w^2=x^3y+xy^3+z^4$ fixed by the automorphism $g=\DiaG{-1}{\zeta_8}{-\zeta_8}{1}$ is the union of the $2$ points
$(0:1:0:0)$, $(0:0:1:0)$. Its trace on $\Pic{S}$ is then $0$. As
$g^2=\DiaG{1}{\im}{\im}{1}=\DiaG{-1}{1}{1}{-\im}$ has eigenvalues $<1^4,\im^2,(-\im)^2>$ (see above), we deduce that the multiplicity of each of the $4$ primitive $8$-th roots of unity
($\zeta_8,\zeta_8^3,\zeta_8^5$ and $\zeta_8^7$)
as an eigenvalue of $g$ is $1$ and that the sum of multiplicities of $1$ and $-1$ is $4$. As the trace is equal to
zero, $1$ and $-1$ each have multiplicity $2$.
\proofend
\end{itemize}
\end{proof}

\bigskip

\begin{Lem}\fauxtitre
\label{Lem:QuarP2}
Let $H$ be a non-trivial abelian subgroup of $\PGLn{3}$ and $\Gamma=\{ (x:y:z) \in \Pn \ |\ F(x,y,z)=0\}$ a smooth quartic curve of\hspace{0.2 cm}$\Pn$ invariant by $H$. Then, up to conjugation in $\PGLn{3}$, $H$ and $\Gamma$ are in the following list:

\begin{itemize}
\item
$H$ is cyclic:\\
\begin{tabular}{|llll|}
\hline
{\it name} & {\it structure} & {\it generator} & {\it equation $\mathit F$}\\
{\it of $\mathit H$} & {\it of $\mathit H$} & {\it of $\mathit H$} & {\it of $\mathit \Gamma$}\\
\hline 
${\num{Q2}}$ & $\Z{2}$&$\Diag{1}{1}{-1}$& $L_4(x,y)+L_2(x,y)z^2+z^4\hts$ \\
${\num{Q3.1}}$ & $\Z{3}$&$\Diag{1}{1}{\omega}$& $L_4(x,y)+z^3L_1(x,y)\hts$ \\
${\num{Q3.2}}$ & $\Z{3}$&$\Diag{1}{\omega}{\omega^2}$& $x(x^3+y^3+z^3)+yzL_1(x^2,yz)\hts$ \\
${\num{Q4.1}}$ & $\Z{4}$&$\Diag{1}{1}{\im}$& $L_4(x,y)+z^4\hts$ \\
${\num{Q4.2}}$ & $\Z{4}$&$\Diag{1}{-1}{\im}$& $x^4+y^4+z^4+xyL_1(xy,z^2)\hts$ \\
${\num{Q6}}$ & $\Z{6}$&$\Diag{\omega}{1}{-1}$& $x^3y+y^4+z^4+\lambda y^2z^2$ \\
${\num{Q7}}$ & $\Z{7}$&$\Diag{\zeta_7}{(\zeta_7)^4}{(\zeta_7)^2}$& $x^3y+y^3z+xz^3\hts$ \\
${\num{Q8}}$ & $\Z{8}$&$\Diag{\zeta_8}{-\zeta_8}{1}$& $x^3y+xy^3+z^4\hts$ \\
${\num{Q9}}$ & $\Z{9}$&$\Diag{(\zeta_9)^6}{1}{\zeta_9}$& $x^3y+y^4+xz^3$ \\
${\num{Q12}}$ & $\Z{12}$&$\Diag{\omega}{1}{\im}$& $x^3y+y^4+z^4\hts$ \\
\hline
\end{tabular}
\item
$H$ is non-cylic: \\
\begin{tabular}{|llll|}
\hline
{\it name} & {\it structure} & {\it generators} & {\it equation $\mathit F$}\\
{\it of $\mathit H$} & {\it of $\mathit H$} & {\it of $\mathit H$} & {\it of $\mathit \Gamma$}\\
\hline 
${\num{Q22}}$ & $\Z{2}\times \Z{2}$&$\Diag{1}{1}{-1}$, $\Diag{1}{-1}{1}$& $L_2(x^2,y^2,z^2)\hts$ \\
${\num{Q24}}$ & $\Z{2}\times \Z{4}$&$\Diag{1}{-1}{1}$, $\Diag{1}{1}{\im}$& $x^4+y^4+z^4+\lambda x^2y^2\hts$ \\
${\num{Q44}}$ & $(\Z{4})^2$&$\Diag{1}{1}{\im}$, $\Diag{1}{\im}{1}$& $x^4+y^4+z^4$\\
\hline
\end{tabular}
\end{itemize}
where $\lambda,\mu \in \K$ are parameters, the $L_i$ are forms of degree $i$ such that the quartic is smooth, and $\zeta_n=e^{2\im\pi/n},\omega=\zeta_3$.
\end{Lem}
\begin{remark}
This classical result (see \cite{bib:Cia}, \cite{bib:Wim}) has been proved many times for more a hundred years.
\end{remark}
\begin{proof}\fauxtitred\upshape
We suppose first that $H$ is cyclic, generated by a diagonal element $h$ of order $n>1$.

If two eigenvalues of $h$ are equal, we may suppose that $h=\Diag{1}{1}{\alpha}$, where $\alpha$ is a primitive $n$-th root of unity. We write $F=\lambda z^4+z^3 L_1(x,y) +z^2 L_2(x,y) + z L_3(x,y)+L_4(x,y)$, where $\lambda \in \K$ and $L_i$ is a form of degree $i$. As $\Gamma$ is non-singular, we have $L_4 \not= 0$, so either $\alpha^2,\alpha^3$ or $\alpha^4$ is equal to $1$. Up to a change of coordinates, we get respectively $\nump{Q2}$, $\nump{Q3.1}$ or $\nump{Q4.1}$.

We may now suppose that the three eigenvalues of $h$ are different.

If $h$ is of order $3$, say $\Diag{1}{\omega}{\omega^2}$, then $F(x,\omega y,\omega^2 z)=\omega^i F(x,y,z)$. If $i=0$ we get case $\nump{Q3.2}$. If $i=1$ or $i=2$, we return to the case $i=0$ by a permutation of the coordinates $x$, $y$ and $z$.

If $h$ is of order $4$, say $\Diag{1}{-1}{\im}$, then $F(x,-y,\im z)=\im^i F(x,y,z)$. If $i=1$ or $i=3$, the equation $F$ is divisible by $z$ and is singular. If $i=2$, as the variable $z$ is only of degree $2$, the curve is once again singular. Then $i=0$ and we get case $\nump{Q4.2}$.

We may now suppose that the order of $h$ is strictly larger than $4$. Note that the points of\hspace{0.2 cm}$\Pn$ fixed by $h$ are $(1:0:0)$, $(0:1:0)$ and $(0:0:1)$. If none of these $3$ points belongs to $\Gamma$, the equation of $F$ contains $x^4,y^4$ and $z^4$, so $h$ is of order $\leq 4$. Let us then suppose that $A_1=(1:0:0)\in \Gamma$, so $F=x^3 \cdot (\alpha y + \beta z)+ x^2 \cdot L_2(y,z)+x\cdot L_3(y,z)+L_4(y,z)$, where the $L_i$ are forms of degree $i$. As the point $A_1$ is non-singular, either $\alpha$ or $\beta$ is different to zero. We suppose that $\alpha \not=0$ and write $h=\Diag{1}{a}{b}$, so that $F(x,ay,bz)=a F(x,y,z)$. This implies that $\beta=0$ and that $L_2(y,z)$ is a multiple of $z^2$. Using once again the non-singularity of $\Gamma$, the equation $F$ must be of order at least $3$ in each variable, so it must contain either $xy^3$, $y^4$ or $y^3z$. It must also contain either $xz^3$, $yz^3$ or $z^4$. We investigate the $9$ cases:
\begin{itemize}
\item
{\it $\mathit F$ contains $\mathit{ xy^3}$:}\\
In this case, $a^2=1$, so $a=-1$.
\begin{itemize}
\item
If $F$ contains $xz^3$, then $b^3=-1$ and $F$ contains only $x^3y$, $xy^3$ and $xz^3$ and is therefore reducible.
\item
If $F$ contains $yz^3$, then $b^3=1$ and $F=x^3y+xy^3+yz^3$, a reducible curve.
\item
If $F$ contains $z^4$, then $b^4=-1$ and $F=x^3y+xy^3+z^4$ so the group is generated by
$\Diag{1}{-1}{\zeta_8}$ hence also by $\Diag{\zeta_8}{-\zeta_8}{1}$, which is the case $\nump{Q8}$.
\end{itemize}
\item
{\it $\mathit F$ contains $\mathit{y^4}$:}\\
In this case, $a^3=1$, so $a$ is a primitive $3$-rd root of unity, since $a\not=1$.
\begin{itemize}
\item
If $F$ contains $xz^3$, then $b^3=a$, and $F=x^3y+y^4+xz^3$ so we get $\nump{Q9}$.
\item
If $F$ contains $yz^3$, then $b^3=1$ and $h$ is of order $3$.
\item
If $F$ contains $z^4$, then $b^4=a$ so $b^{12}=1$. 
If $b$ is a primitive $12$-th root of unity, $F=xy^3+y^4+z^4$ and we get $\nump{Q12}$,
as the groups generated respectively by $\Diag{1}{(\zeta_{12})^{4}}{\zeta_{12}}$ and
$\Diag{\omega}{1}{\im}$ are the same. If not, it must be a primitive $6$-th root of unity, so $F=xy^3+y^4+z^4+\lambda y^2z^2$ and we get $\nump{Q6}$.
\end{itemize}
\item
{\it $\mathit F$ contains $\mathit{y^3z}$:}\\
In this case, $b^{-1}=a^2$.
\begin{itemize}
\item
If $F$ contains $xz^3$, then $b^3=a$, so $a^7=1$ and $F=x^3y+y^3z+xz^3$; we get $\nump{Q7}$.
\item
If $F$ contains $yz^3$, then $b^3=1$ and $a^2=b^2$, so $a=-b$
and $F$ contains only $x^3y$, $y^3z$ and $yz^3$, so $F$ is divisible by $y$.
\item
If $F$ contains $z^4$, then $b^4=a$, $a^9=1$ and $F=x^3y+y^3z+z^4$; we get $\nump{Q9}$ by a permutation $(x \ z \ y)$.
\end{itemize}
\end{itemize}
We have now enumerated all possible cyclic groups. Let us investigate non-cyclic groups.

The diagonal groups isomorphic to $\Z{n}\times \Z{n}$ contains elements of the form $\Diag{1}{1}{\gamma}$, where $\gamma$ is a primitive $n$-th root of unity, so $n\leq 4$. The case $n=3$ is clearly impossible, and the cases $n=2$ and $n=4$ give respectively $\nump{Q22}$ and $\nump{Q44}$.

If $H$ is diagonal and non-cyclic, it must contain a group $\Z{n}\times\Z{n}$ and hence one of the groups listed in case $\nump{Q22}$. As $F$ must be of degree $\geq 3$ in all variables, it must contain $x^4$, $y^4$ and $z^4$. All diagonal elements must thus be of order $\leq 4$. There is then only one remaining possibility, that is $\nump{Q24}$.

By Proposition \refd{Prp:PGL3Cr}, if $H$ is non-diagonalizable, it must be isomorphic to $\Z{3} \times \Z{3}$, generated 
by $\Diag{1}{\omega}{\omega^2}$ and the automorphism $(x:y:z) \mapsto (y:z:x)$. But this group doesn't leave any quartic 
curve invariant.\proofend
\end{proof}

\bigskip

\begin{Exa}\fauxtitred\upshape
\label{Exa:DP2ConicBundle}
Let $S$ be the following Del Pezzo surface of degree $2$:
\begin{center}$S =\{(w:x:y:z) \in \mathbb{P}(2,1,1,1)\ |\ w^2=L_2(x^2,y^2,z^2)\}$,\end{center}where $L_2\in \K(x,y,z)$ is some form of degree $2$, and $L_2(x^2,y^2,z^2)$ is non-singular.

Let $G \cong (\Z{2})^2$ be the group of automorphisms of the form $(w:x:y:z) \mapsto (w:x:\pm y:\pm z)$.
Then, $G$ acts minimally on $S$ and preserves a conic bundle structure.
\index{Conic bundles!on Del Pezzo surfaces|idi}

Indeed, the eigenvalues of the three involutions of the group on $\Pic{S}$ are $<1^4,(-1)^4>$. This implies that $\rkPic{S}^{G}=2$. As any of the three involutions fixes an elliptic curve\index{Curves!elliptic!birational maps that fix an elliptic curve}, the group is minimal (no such group acts on a Del Pezzo surface of degree $\geq 3$).
\end{Exa}

\bigskip
\bigskip

\section{Automorphisms of Del Pezzo surfaces of degree $1$}
\index{Del Pezzo surfaces!of degree $1$|idb}
\pagetitless{Finite abelian groups of automorphisms of Del Pezzo surfaces of degree $\leq 4$}{Del Pezzo surfaces of degree $1$}
\label{Sec:DP1}
A Del Pezzo surface $S$ of degree $1$ is given by the blow-up of $8$ points of\hspace{0.2 cm}$\Pn$, in general position. The linear system $|-2K_S|$ induces a degree $2$ morphism onto a quadric cone in $Q\subset\mathbb{P}^3$, ramified over the vertex $v$ of $Q$ and a smooth curve $C$ of genus $4$. Moreover $C$ is the intersection of $Q$ with a cubic surface. (See \cite{bib:BaB}, \cite{bib:deF}, \cite{bib:Dol}.)

Note that a quadric cone is isomorphic to the weighted projective plane $\mathbb{P}(1,1,2)$ and the ramification curve $C$ has equation of degree $6$ there. Up to a change of coordinates, we may assume that the surface $S$ has the equation
\begin{center}
$w^2=z^3+F_4(x,y)z+F_6(x,y)$
\end{center}
in the weighted projective space $\mathbb{P}(3,1,1,2)$, where $F_4$ and $F_6$ are forms of respective degree $4$ and $6$ (see \cite{bib:Kol}, Theorem III.3.5 and and \cite{bib:KoSC}, Corollary 3.54). Note that multiple roots of $F_6$ are not roots of $F_4$, since $S$ is non-singular. The point $v=(1:0:0:1)=(-1:0:0:1)$ is the vertex of the quadric.

We denote by $\sigma(w:x:y:z)=(-w:x:y:z)$ the involution associated to the $2$-covering. This is classically called the \defn{Bertini} involution of the surface.\index{Involutions!Bertini|idb}

If $F_4=0$, the surface is a triple covering of $\mathbb{P}(3,1,1)$, ramified over $v$ and the hyperelliptic curve of genus $2$ of equation $z=0$ \index{Curves!hyperelliptic!birational maps that fix a hyperelliptic curve}, $w^2=F_6(x,y)$. In this case we denote the automorphism of order $3$ corresponding to this covering by $\rho(w:x:y:z)=(w:x:y:\omega z)$ (where $\omega=e^{2\im\pi/3}$).

Note that in any case, $|K_S|$ is a pencil of elliptic curves with one base point, $v=(1:0:0:1)$, parametrised by the $(x,y)$-coordinates. Any automorphism of $S$ acts then on the elliptic bundle and fixes the vertex $v$ of $Q$. \index{Curves!elliptic!elliptic bundle|idb}

\break

\begin{Lem}\fauxtitre
\label{Lem:DP1}
Using the above notation for $S$, $\sigma$, $\rho$, $F_4$ and $F_6$, we have:
\begin{itemize}
\item[\upshape 1.]
The action of $\Aut(S)$ on $|K_S|$ gives an exact sequence
\begin{center}
$1\rightarrow G_S \rightarrow \Aut(S) \stackrel{\pi}{\rightarrow} \Aut(\mathbb{P}^1)$,
\end{center}
where \begin{center}$G_S=\left\{\begin{array}{lll}<\sigma,\rho>&\cong \Z{6} &\mbox{ if }F_4 =0, \\
<\sigma>&\cong\Z{2} & \mbox{ otherwise.}\end{array}\right.$\end{center}
\item[\upshape 2.]
The Bertini involution acts minimally on $S$. The eigenvalues of its action on $\Pic{S}$ are $1$ with multiplicity $1$ and $-1$ with multiplicity $8$.\index{Involutions!Bertini}
\item[\upshape 3.]
If $F_4=0$, the automorphism $\rho$ of order $3$ acts minimally on $S$. The eigenvalues of its action on $\Pic{S}$ are $1$ with multiplicity $1$ and $\omega$, $\omega^2$ (where $\omega=e^{2\im\pi/3}$) both of multiplicity $4$.
\item[\upshape 4.]
If $G \subset \Aut(S)$ is abelian, then $\pi(G)$ is cyclic.
\end{itemize}
\end{Lem}
\begin{proof}\fauxtitred\upshape
Note first of all that the double covering of the quadric $Q \cong \mathbb{P}(1,1,2)$ gives an exact sequence
\begin{center}
$1\rightarrow <\sigma> \rightarrow \Aut(S) \rightarrow \Aut(Q)_{C}$,
\end{center}
where $\Aut(Q)_{C}$ denote the automorphisms of $Q$ that leaves invariant the ramification curve $C=\{(x:y:z)\ |\ z^3+zF_4(x,y)+F_6(x,y)=0\}$. (In fact we can prove that the right homomorphism is surjective, but we will not need it here). A quick calculation shows that any element of $\Aut(Q)_{C}$ belongs to $P(GL(2,\K) \times GL(1,\K))$. This implies that $\Aut(S)\subset P(GL(1,\K)\times GL(2,\K) \times GL(1,\K))$.

The first assertion follows by a direct calculation. Let us prove the second one. The set of points of $S$ fixed by the Bertini involution $\sigma$ is the disjoint union of the vertex of the quadric, $(1:0:0:1)$ (which is equal to $(-1:0:0:1)$) and the curve $C$ of genus $4$. By the Lefschetz formula, the trace of the action of $\sigma$ on $H^{*}(S,\mathbb{Z})$ is then equal to $-5$ (the Euler characteristic of the fixed points), and then the trace on $\Pic{S}$ is $-7$. This gives the second assertion.

We once again use the Lefschetz formula to prove the third assertion. As $F_4=0$, the form $F_6$ must be non-singular. The set of points of $S$ fixed by $\rho$ is the disjoint union of the vertex of the quadric, $(1:0:0:1)$ (which is equal to $(1:0:0:\omega)$), and the hyperelliptic curve of genus $2$\index{Curves!hyperelliptic!birational maps that fix a hyperelliptic curve} with equation $z=0$, $w^2=F_6(x,y)$. The trace of $\rho$ on $H^{*}(S,\mathbb{Z})$ is thus equal to $-1$ and hence the trace on $\Pic{S}$ is equal to $-3$. As the eigenvalues of $\rho$ are $\omega$ and $\omega^2$, both with the same multiplicity, and $1$, we get the third assertion.

\break

It remains to prove the last assertion. If $G$ is abelian, then so is $\pi(G)$, which can be either cyclic or conjugate to the group generated by $(x:y) \mapsto (y:x)$ and $(x:y) \mapsto (-x:y)$. But elements of $\Aut(S)$ which are sent on these are of the form $(w:x:y:z) \mapsto (\lambda_1 w: y:x:\mu_1 z)$ and $(w:x:y:z) \mapsto (\lambda_2 w: -x:y:\mu_2 z)$, for some $\lambda_1,\lambda_2,\mu_1,\mu_2 \in \K^{*}$, and don't commute. This proves the last assertion.\proofend
\end{proof}

\bigskip

The $240$ exceptional divisors on $S$ are (see Proposition 
\refd{Prp:DescExcCurvDelPezzo}):
\index{Del Pezzo surfaces!curves on the surfaces|idi}
\begin{itemize}
\item
$E_1=\pi^{-1}(A_1),...,E_8=\pi^{-1}(A_8)$, the $8$ pull-backs of the points $A_1$,...,$A_8$.
\item
The strict pull-backs of the curves of degree $d$ passing through the $A_i$'s with some given multiplicities:
\begin{center}
$\begin{array}{|l|l|l|}
\hline
\mbox{degree } d& \mbox{multiplicities at the }A_i\mbox{'s} & \mbox{number of such curves} \\
\hline
1 & (1,1) & 28\\
2 & (1,1,1,1,1) & 56\\
3 & (2,1,1,1,1,1,1) & 56\\
4 & (2,2,2,1,1,1,1,1) & 56\\
5 & (2,2,2,2,2,2,1,1) & 28\\
6 & (3,2,2,2,2,2,2,2) & 8\\
\hline
\end{array}$\end{center} \end{itemize}

With this information, we can describe the action of the Bertini involution:\index{Involutions!Bertini!matrix|idb}
\begin{Lem}\fauxtitred
\label{Lem:Bertini}
\begin{enumerate}
\item[\upshape 1.]
The Bertini involution $\sigma$ switches the pull-backs 
of any point of the quadric cone by the $2$-covering.
\item[\upshape 2.]
It exchanges
the two exceptional divisors which have the same image, $D$ and $-2K_S-D$, 
for any exceptional divisor $D$. 
\item[\upshape 3.]
Its action on $\Pic{S}$ is given by the matrix
\begin{center}$\left(\begin{array}{rrrrrrrrr}
-3 & -2 & -2 & -2 & -2 & -2 & -2 & -2 & -6 \\
-2 & -3 & -2 & -2 & -2 & -2 & -2 & -2 & -6 \\
-2 & -2 & -3 & -2 & -2 & -2 & -2 & -2 & -6 \\
-2 & -2 & -2 & -3 & -2 & -2 & -2 & -2 & -6 \\
-2 & -2 & -2 & -2 & -3 & -2 & -2 & -2 & -6 \\
-2 & -2 & -2 & -2 & -2 & -3 & -2 & -2 & -6 \\
-2 & -2 & -2 & -2 & -2 & -2 & -3 & -2 & -6 \\
-2 & -2 & -2 & -2 & -2 & -2 & -2 & -3 & -6 \\
6 & 6 & 6 & 6 & 6 & 6 & 6 & 6 & 17\end{array}\right)$\end{center}
relative to the basis $(E_1,E_2,...,E_8,L)$, where $L$ is the pull-back by $\pi$ of a line of\hspace{0.2 cm}$\Pn$.
\item[\upshape 4.]
Its eigenvalues on $\Pic{S}$ are $1$ with multiplicity $1$ and $-1$ with multiplicity $8$.
\item[\upshape 5.]
$\rkPic{S}^{\sigma}=1$, so $\sigma$ acts minimally on the surface $S$.
\end{enumerate}
\end{Lem}
\begin{proof}\fauxtitred\upshape
The first part follows from the definition of $\sigma$ as $\sigma(w:x:y:z)=(-w:x:y:z)$. We then get 
the second part directly, which implies the third. 

Note that $\sigma(K_S)=K_S$ and $\sigma(K_S+2E_i)=K_S+2(-K_S-E_i)=-(K_S+2E_i)$, for $i=1,...,8$. Then, the action of $\Pic{S}$ on the basis $(K_S,K_S+2E_1,K_S+2E_2,...,K_S+2E_8)$ of the vector space $\Pic{S}\otimes_{\mathbb{Z}} \mathbb{Q}$ is diagonal with eigenvalues $<1,(-1)^8>$. We thus get assertion $4$, which implies the last one.
\proofend\end{proof}

\bigskip

\begin{Lem}\fauxtitre
\label{Lem:EllipticRationalDP1}
Let $S\subset \mathbb{P}(3,1,1,2)$ be a Del Pezzo surface of degree $1$,  with equation $w^2=z^3+zF_4(x,y)+F_6(x,y)$ and let $p:S\rightarrow \mathbb{P}^1$ denote the rational map $(w:x:y:z) \dasharrow (x:y)$. 
Let $\Delta$ denote the form $27F_6(x,y)^2+4F_4(x,y)^3$.

Let $(a:b) \in \mathbb{P}^1$. Then, one and only one of the following possibilities occurs:
\begin{itemize}
\item
$(a:b)$ is a root of $\Delta$, and the curve $p^{-1}(a:b)$ is a singular rational curve. 
\item
$(a:b)$ is not a root of $\Delta$, and the curve $p^{-1}(a:b)$ is a smooth elliptic curve. \index{Curves!elliptic!birational maps that fix an elliptic curve}
\end{itemize}
\end{Lem}
\begin{proof}\fauxtitred\upshape
Note that $p^{-1}(a:b)$ represents  an element of the linear system $|-K_S|$ and thus corresponds to a cubic of\hspace{0.2 cm}$\Pn$ passing through the $8$ blown-up points. If the curve of\hspace{0.2 cm}$\Pn$ is smooth, it is an elliptic curve.

Otherwise, it is a rational curve. As the points are in general position, the singularity is not at one of the blown-up points (see Proposition \refd{Prp:ConditionsDelPezzo}), so the curve remains singular in $S$. 

Up to a change of coordinates , we may suppose that $(a:b)=(0:1)$, so the fibre is $\{(w:y:z) \in \mathbb{P}(3,1,2)\ | w^2=z^3 + \alpha zy^4+\beta y^6\}$, where $\alpha=F_4(0,1)$ and $\beta=F_6(0,1)$. Note that the map
\begin{center}$(u,v)\mapsto (u,1,v)$\end{center}
is a isomorphism from $\K^2$ to an affine open subset of $\mathbb{P}(3,1,2)$. The equation of the curve in $\K^2$ becomes $u^2=v^3+\alpha v+\beta$ and it is singular if and only if $v^3+\beta v+\beta$ has a double root in $\K[v]$. This is equivalent to having $27\beta^2+4\alpha^3=0$. Finally, note that the point $(1:0:1)$ is not a singular point of the curve. This gives the result.\proofend
\end{proof}

\bigskip
\bigskip

Let us denote by $\DiaG{\alpha}{\beta}{\gamma}{\delta}$ the automorphism $(w:x:y:z) \mapsto (\alpha w:\beta x:\gamma y:\delta z)$.

\break

\begin{Prp}\titreProp{automorphisms of Del Pezzo surfaces of degree $1$}\\
\label{Prp:DelP1}
\index{Del Pezzo surfaces!automorphisms}
Let $S$ be a Del Pezzo surface of degree $1$ given by the equation 
\begin{center}$w^2=z^3+F_4(x,y)z+F_6(x,y)$ \end{center}
and let $G$ be an abelian subgroup of $\Aut(S)$ such that $\rkPic{S}^{G}=1$. Then, up to isomorphism we 
are in one of the following cases:\\
\begin{tabular}{|lllll|}
\hline
{\it name} & {\it structure} & {\it generators} & {\it equation of}& {\it equation of}\\
{\it of $\mathit (G,S)$} & {\it of $\mathit G$} & {\it of $\mathit G$} & $\mathit F_4$& $\mathit F_6$\\
\hline \index{Involutions!Bertini}
${\num{1.B}}$ & $\Z{2}$& $\sigma$& $L_4(x,y)$ & $L_6(x,y)\hts$ \\
${\num{1.\rho}}$ & $\Z{3}$& $\rho$& $0$ &$L_6(x,y)\hts$ \\
${\num{1.\sigma\rho}}$ & $\Z{6}$& $\sigma\rho$& $0$ & $L_6(x,y)\hts$ \\
${\num{1.B2.1}}$ & $(\Z{2})^2$&$\sigma$, $\DiaG{1}{1}{-1}{1}$& $L_2(x^2,y^2)$ & $L_3(x^2,y^2)\hts$ \\
${\num{1.\sigma\rho 2.1}}$ & $\Z{6}\times \Z{2}$& $\sigma\rho$, $\DiaG{1}{1}{-1}{1}$& $0$ & $L_3(x^2,y^2)\hts$ \\
${\num{1.\rho 2}}$ & $\Z{6}$& $\DiaG{1}{1}{-1}{\omega}$& $0$ & $L_3(x^2,y^2)\hts$ \\
${\num{1.B2.2}}$ & $\Z{4}$& $\DiaG{\im}{1}{-1}{-1}$& $L_2(x^2,y^2)$ & $xyL_2'(x^2,y^2)\hts$ \\
${\num{1.\sigma\rho 2.2}}$ & $\Z{12}$& $\DiaG{\im}{1}{-1}{-\omega}$& $0$ & $xyL_2(x^2,y^2)\hts$ \\
${\num{1.B3.1}}$ & $\Z{6}$&$\DiaG{-1}{1}{\omega}{1}$& $xL_1(x^3,y^3)$ & $L_2(x^3,y^3)\hts$ \\
${\num{1.\sigma\rho 3}}$ & $\Z{6}\times\Z{3}$& $\sigma\rho$, $\DiaG{1}{1}{\omega}{1}$& $0$ &  $L_2(x^3,y^3)\hts$ \\
${\num{1.\rho 3}}$ & $ (\Z{3})^2$& $\rho$, $\DiaG{1}{1}{\omega}{1}$& $0$ & $L_2(x^3,y^3)\hts$ \\
${\num{1.B3.2}}$ & $\Z{6}$& $\DiaG{-1}{1}{\omega}{\omega}$& $\lambda x^2 y^2$ & $L_2(x^3,y^3)\hts$ \\
${\num{1.B4.1}}$ & $\Z{2}\times\Z{4}$&$\sigma$, $\DiaG{1}{1}{\im}{1}$& $L_1(x^4,y^4)$ &  $x^2L_1'(x^4+y^4)\hts$ \\
${\num{1.B4.2}}$ & $\Z{8}$&$\DiaG{\zeta_8}{1}{\im}{-\im}$& $\lambda x^2y^2$ & $xy(x^4+y^4)\hts$ \\
${\num{1.\sigma\rho 4}}$ & $\Z{24}$& $\DiaG{\zeta_{8}}{1}{\im}{-\im\omega}$& $0$& $xy(x^4+y^4)\hts$ \\
${\num{1.5}}$ & $\Z{5}$&$\DiaG{1}{1}{\zeta_5}{1}$& $\lambda x^4$ & $x(\mu x^5+y^5)\hts$ \\
${\num{1.B5}}$ & $\Z{10}$&$\DiaG{-1}{1}{\zeta_5}{1}$& $\lambda x^4$ & $x(\mu x^5+y^5)\hts$ \\
${\num{1.\sigma\rho 5}}$ & $\Z{30}$& $\DiaG{-1}{1}{\zeta_5}{\omega}$ & $0$ & $x(x^5+y^5)\hts$ \\
${\num{1.\rho 5}}$ & $\Z{15}$& $\DiaG{1}{1}{\zeta_5}{\omega}$& $0$ & $x(x^5+y^5)\hts$ \\
${\num{1.6}}$ & $\Z{6}$&$\DiaG{1}{1}{-\omega}{1}$& $\lambda x^4$ & $\mu x^6+y^6\hts$ \\
${\num{1.B6.1}}$ & $\Z{2}\times\Z{6}$&$\sigma$, $\DiaG{1}{1}{-\omega}{1}$& $\lambda x^4$ & $\mu x^6+y^6\hts$ \\
${\num{1.\sigma\rho 6}}$ & $ (\Z{6})^2$& $\sigma\rho$, $\DiaG{1}{1}{-\omega}{1}$&$0$& $x^6+y^6\hts$ \\
${\num{1.\rho 6}}$ & $\Z{3}\times \Z{6}$& $\rho$, $\DiaG{1}{1}{-\omega}{1}$&$0$& $x^6+y^6\hts$ \\
${\num{1.B6.2}}$ & $\Z{2}\times\Z{6}$&$\sigma$, $\DiaG{1}{1}{-\omega}{\omega}$& $\lambda x^2y^2$ & $x^6+y^6\hts$ \\
${\num{1.B10}}$ & $\Z{20}$&$\DiaG{\im}{1}{\zeta_{10}}{-1}$& $x^4$& $xy^5\hts$ \\
${\num{1.B12}}$ & $\Z{2}\times\Z{12}$&$\sigma, \DiaG{\im}{1}{\zeta_{12}}{-1}$& $x^4$ & $y^6\hts$ \\
\hline
\end{tabular}\\
where $L_i$ denotes a form of degree $i$, $\lambda,\mu \in \K$ are parameters, and the equations are such that $S$ is non-singular.

Furthermore, all the cases above are minimal pairs $(G,S)$ with $\rkPic{S}^G=1$.
\end{Prp}
\begin{remark}
Some of the surfaces given above seems to be different but are in fact the same. Indeed, we choosed to put the automorphisms on diagonal form, the equations of the surface may then not be standard.\end{remark}\begin{proof}\fauxtitred\upshape
We use the notation of Lemma \refd{Lem:DP1} for $\pi$ and $G_S$ and recall that
the group $\pi(G)\subset \Aut(\mathbb{P}^1)$ is cyclic. So either $\pi(G)=1$, or $G$ is generated by $\ker \pi \cap G=G_S \cap G$ and one element whose action on $\mathbb{P}^1$ is not trivial. If $G \subset \ker \pi$, we get $\nump{1.B}$, $\nump{1.\sigma\rho}$ or $\nump{1.\rho}$, which act minimally on $S$ (see Lemma \refd{Lem:DP1}).

We may thus suppose that $G$ is generated by $G_S \cap G$ and \begin{center}$(w:x:y:z) \mapsto (aw:x:\zeta_n y:bz)$,\end{center} where $\zeta_n=e^{2\im\pi/n}$, for some integer $n>1$. 
We examine this case in the following manner: 
\begin{itemize}
\item
Firstly, we find the possibilities for $n$ and $F_6,F_4$, which determine the surface $S$.
\item
Secondly, we pull-back the action of $\gamma_n:(x:y) \mapsto (x:\zeta_n y)$ on the elliptic bundle and get an abelian group $\pi^{-1}(<\gamma_n>) \subset \Aut(S)$ which contains $G_S$ and $G$.\index{Curves!elliptic!elliptic bundle}
\item
Thirdly, we look for subgroups $G \subset \pi^{-1}(<\gamma_n>)$ such that $\pi(G)=<\gamma_n>$ and $\rkPic{S}^G=1$, using the Lefschetz formula to compute the eigenvalues of the action of automorphisms of $\Pic{S}$. 
Note that points of $S$ fixed by the automorphism $\DiaG{a}{1}{\zeta_n}{b}$ are contained in the union of the two fibres $\pi^{-1}(0:1)$ and $\pi^{-1}(1:0)$. As the set of fixed points is smooth, the two curves are elliptic (see Lemma \refd{Lem:EllipticRationalDP1}).\index{Curves!elliptic!birational maps that fix an elliptic curve}
\item
Note that the groups which contain $\rho$ or $\sigma$ always occur, since $\rkPic{S}^{\sigma}=\rkPic{S}^{\rho}=1$.
\end{itemize}

We set $F\h{0.2}=\h{0.2}w^2-z^3-zF_4(x,y)-F_6(x,y)$ and remark that $F(aw,x,\zeta_n y,b z)\h{0.6}=\h{0.4}b^3F(w,x,y,z)$, which implies that $F_4(x,\zeta_n y)=b^2F_4(x,y)$ and $F_6(x,\zeta_n y)=b^3F_6(x,y)$. 
If $F_4\not=0$, then $b$ is a power of $\zeta_n$, and then $b^n=1$. Otherwise, the possibilities for $b$ are obtained by multiplication by $\rho$.

We enumerate the possible values for $n$, and in each case the possibilities for $b$. The value of $b$ determines the eigenvalues of $F_4$ and $F_6$, and then the elements $F_4$ and $F_6$ themselves. 

{\it Note that the double roots of $F_6$ are not roots of $F_4$, since $S$ is non-singular.}
\begin{itemize}
\item
\cadre{$n=2$, $b=1$}: {$F_4=L_2(x^2,y^2)$ \it and $F_6=L_3(x^2,y^2)$}.\\
The group $\pi^{-1}(<\gamma_2>)$ is generated by $\sigma$ and $g=\DiaG{1}{1}{-1}{1}$ (and $\rho$ if $L_2=0$), and is isomorphic to $(\Z{2})^2$ (respectively to $\Z{6}\times\Z{2}$). 
The set of points of $S$ fixed by $\DiaG{1}{1}{-1}{1}$ is the disjoint union of the elliptic curve $y=0$ and $2$ or $3$ points ($w=x=0$).\index{Curves!elliptic!birational maps that fix an elliptic curve|idi} The trace on $\Pic{S}$ is then $0$ or $1$. As the trace is an odd integer, the number of points of the set $w=x=0$ is exactly $3$, and then the eigenvalues are $<1^5,(-1)^4>$.

As the eigenvalues of $\sigma$ are $<1,(-1)^8>$, we see that the eigenvalues of both $g$ and $\sigma g$ are $<1^5,(-1)^4>$. 

The group $G$ must then contain either $\sigma$ and $\rho$ (we get $\nump{1.\sigma\rho 2.1}$), or $\sigma$ (we get $\nump{1.B2.1}$) or $\rho$ (we get $\nump{1.\rho 2}$). The group generated by $\DiaG{-1}{1}{-1}{\omega}=\DiaG{1}{-1}{1}{\omega}$ is the same as $\nump{1.\rho 2}$, after a permutation of the coordinates $x$ and $y$.
\item
\cadre{$n=2$, $b=-1$}: {$F_4=L_2(x^2,y^2)$ \it and $F_6=xyL_2'(x^2,y^2)$}.\\
If $F_4\not=0$, the group $\pi^{-1}(<\gamma_2>)$ is cyclic of order $4$, generated by $g=\DiaG{\im}{1}{-1}{-1}$. Otherwise, it is generated by $g$ and $\rho$, and is cyclic of order $12$. 

Then, $G$ contains in any case the Bertini involution. The possible cases are $G=<g>$  (case $\nump{1.B2.2}$) or $G=<g,\rho>$ (case $\nump{1.\sigma\rho 2.2}$).
\item
\cadre{$n=3$, $b=1$}: {$F_4=x L_1(x^3,y^3)$ \it and $F_6=L_2(x^3,y^3)$}.\\
 The group  $\pi^{-1}(<\gamma_3>)$ is generated by $\sigma$ and $g=\DiaG{1}{1}{\omega}{1}$ (and $\rho$ if $F_4=0$), and is isomorphic to $\Z{6}$ (respectively to $\Z{6}\times\Z{3}$). 
 
 The set of points of $S$ fixed by $g$ is the disjoint union of an elliptic curve\index{Curves!elliptic!birational maps that fix an elliptic curve|idi} ($y=0$) and two points (which are $(\pm\sqrt{L_2(0,1)}:0:1:0)$), the trace on $\Pic{S}$ is $0$ and the eigenvalues are $<1^3,(\omega)^3,(\omega^2)^3>$.

If $F_4=0$, the eigenvalues of $g\rho^2=\DiaG{1}{1}{\omega}{\omega^2}=\DiaG{1}{\omega^2}{1}{1}$ are the same as those of $g$ (using a permutation of the coordinates $x$ and $y$).
As the automorphism $g\rho=\DiaG{1}{1}{\omega}{\omega}$ fixes exactly $5$ points (which are $v$, $(\pm\sqrt{L_2(0,1)}:0:1:0)$ and $(\pm\sqrt{L_2(1,0)}:1:0:0)$), its trace on $\Pic{S}$ is $3$ and its eigenvalues are $<1^5,(\omega)^2,(\omega^2)^2>$.
 
In any case, $G$ therefore contains either $\sigma$ or $\rho$. 
\begin{itemize}
\item
If $G$ contains only $\sigma$,  the possibilities for the group are $G=<\sigma,g>=<\sigma g>$ (case $\nump{1.B3.1}$) or $G=<\sigma, \rho g>=<\sigma\rho g>$ (case $\nump{1.B3.2}$) or $G=<\sigma, \rho^2 g>=<\sigma\rho^2 g>$ (we get case $\nump{1.B3.1}$ by a permutation of the coordinates $x$ and $y$).
\item
If $G$ contains $\sigma$ and $\rho$, then $G=<\rho,\sigma,g>=<\rho,\sigma g>$ (case $\nump{1.\sigma\rho 3}$).
\item
The last possibility is $G=<\rho,g>$ (case $\nump{1.\rho 3}$).\end{itemize}
\item
\cadre{$n=3$, $b=\omega$}: {\it  $F_4=\lambda x^2y^2$ \it and $F_6=L_2(x^3,y^3)$}.\\
The case where $F_4=0$ has already been treated above. We may thus suppose that $F_4\not=0$, so the group $\pi^{-1}(<\gamma_3>)\cong \Z{2}\times\Z{3}$ is generated by $\sigma$ and $g=\DiaG{1}{1}{\omega}{\omega}$. The set of points of $S$ fixed by $g$ is the disjoint union of $5$ points (which are $v$, $(\pm\sqrt{L_2(0,1)}:0:1:0)$ and $(\pm\sqrt{L_2(1,0)}:1:0:0)$), the trace on $\Pic{S}$ is $3$ and its eigenvalues are $<1^5,(\omega)^2,(\omega^2)^2>$. So $G$ must be the entire group $\pi^{-1}(<\gamma_3>)$ and we get $\nump{1.B3.2}$.
\item
\cadre{$n=3$, $b=\omega^2$}: {\it  $F_4=y L_1(x^3,y^3)$ \it and $F_6=L_2(x^3,y^3)$}.\\
We get the case $b=1$ by exchanging the coordinates $x$ and $y$.
\item
\cadre{$n=4$, $b=1$}: {$F_4=L_1(x^4,y^4)$ \it and $F_6=x^2L_1'(x^4,y^4)$}.\\
Note that $F_4\not=0$, as $x$ is a double root of $F_6$. The group $\pi^{-1}(<\gamma_4>)\cong \Z{2}\times\Z{4}$ is generated by $\sigma$ and $g=\DiaG{1}{1}{\im}{1}$. We saw above that the eigenvalues of $g^2=\DiaG{1}{1}{-1}{1}$ are $<(1)^5,(-1)^4>$, so $\im$ and $-\im$ are both eigenvalues of multiplicity $2$ of $g$ on $\Pic{S}$. As the set of points of $S$ fixed by $g$ is the disjoint union of an elliptic curve ($y=0$)\index{Curves!elliptic!birational maps that fix an elliptic curve|idi} and the point $(0:0:1:0)$, the trace on $\Pic{S}$ is $-1$ and the eigenvalues of $g$ are $<(1)^2,(-1)^3,(\im)^2,(-\im)^2>$. This implies that $\rkPic{S}^{g}=2$ and $\Pic{S}^{\sigma g}=4$. The only possibility for $G$ is then to be the whole group $\pi^{-1}(<\gamma_4>)$, so we get $\nump{1.B4.1}$. 
\item
\cadre{$n=4$, $b=-1$}: {$F_4=L_1(x^4,y^4)$ \it and $F_6=y^2L_1'(x^4,y^4)$}.\\
We get the previous case by exchanging the coordinates $x$ and $y$.
\item
\cadre{$n=4$, $b=\im$}: {$F_4=\lambda x^2y^2$ \it and $F_6=\mu x^3y^3$,}\\
for some $\lambda,\mu \in \K$. This is impossible, since these polynomials give a singular surface.
\item
\cadre{$n=4$, $b=-\im$}: {$F_4=\lambda x^2y^2$ \it and $F_6=xy(\alpha x^4+\beta y^4)$}.\\
for some $\alpha,\beta,\lambda \in \K$. We have $\alpha\beta \not=0$, since the surface is non-singular, so up to scaling, we may suppose that $\alpha=\beta=1$. 

If $F_4\not=0$, the group $\pi^{-1}(<\gamma_4>)$ is cyclic of order $8$, generated by $g=\DiaG{\zeta_8}{1}{\im}{-\im}$. Otherwise, it is generated by $g$ and $\rho$, and is cyclic of order $24$. 

Then, $G$ contains in any case the Bertini involution. The possible cases are $G=<g>$  (case $\nump{1.B4.2}$) or $G=<g,\rho>$ (case $\nump{1.\sigma\rho 4}$).
\item
\cadre{$n=5$, $b=1$}: {$F_4=\lambda x^4$ \it and $F_6=x(\mu x^5+\alpha y^5)$,}\\
for some $\lambda,\mu,\alpha \in \K$.
As $\alpha\not=0$, since the surface is non-singular we may suppose, up to a linear change of coordinates, that $\alpha=1$. 
The group $\pi^{-1}(<\gamma_5>)$ is generated by $\sigma$ and $g=\DiaG{1}{1}{\zeta_5}{1}$ (and $\rho$ if $F_4=0$) and is isomorphic to $\Z{10}$ (respectively to $\Z{30}$).

The set of points of $S$ fixed by $g$ is the disjoint union of an elliptic curve ($y=0$)\index{Curves!elliptic!birational maps that fix an elliptic curve|idi} and a point (the point $(0:0:1:0)$), so the trace on $\Pic{S}$ is $-1$ and the eigenvalues of $g$ are $1$ with multipicity $1$ and $\zeta_5,(\zeta_5)^2,(\zeta_5)^3,(\zeta_5)^4$ each with multiplicity $2$. Hence $\rkPic{S}^{g}=1$ and then $G=<g>$ is one possibility (case $\nump{1.5}$).  The group $G$ may also be equal to $<\sigma,g>$ (case $\nump{1.B5}$).

If $F_4=0$, we may suppose, up to a linear change of coordinates, that $\mu=1$. Two more possibilities arise for the group $G$, which are $<\sigma,\rho,g>=<\sigma\rho g>$ (case $\nump{1.\sigma\rho 5}$) or  $<\rho,g>=<\rho g>$ (case $\nump{1.\rho 5}$).
\item
\cadre{$n=5$, $b=\zeta_5$}: {$F_4=\lambda x^2y^2$ \it and $F_6=\mu x^3y^3$,}\\
for some $\lambda,\mu \in \K$. This is impossible, since these polynomials give a singular surface.
\item 
\cadre{$n=5$, $b=(\zeta_5)^2$}: {$F_4=\lambda y^4$ \it and $F_6=y(\alpha x^5 +\mu y^5)$,}\\
for some $\lambda,\mu,\alpha \in \K$. The group $\pi^{-1}(<\gamma_5>)$ is generated by $G_S$ and $\DiaG{(\zeta_5)^3\h{0.5}}{\h{0.5}1\h{0.5}}{\h{0.5}\zeta_5\h{0.5}}{\h{0.5}(\zeta_5)^2}\h{0.5}$ $=\DiaG{1}{(\zeta^5)^{-1}}{1}{1}$. Exchanging the coordinates $x$ and $y$, we get the case $b=1$ treated above.
\item 
\cadre{$n=5$, $b=(\zeta_5)^3$}: {$F_4=\lambda x^3y$ \it and $F_6=\mu x^2y^4$,}\\
for some $\lambda,\mu \in \K$. This is also a singular surface.
\item 
\cadre{$n=5$, $b=(\zeta_5)^4$}: {$F_4=\lambda xy^3$ \it and $F_6=\mu x^4y^2$,}\\
for some $\lambda,\mu \in \K$. This is also a singular surface.
\item
\cadre{$n=6$, $b=1$}: {$F_4=\lambda x^4$ \it and $F_6=\mu x^6+\alpha y^6$,}\\
for some $\lambda,\mu,\alpha \in \K$.
As $\alpha\not=0$, since the surface is non-singular we may suppose, up to a linear change of coordinates,  that $\alpha=1$. 
The group $\pi^{-1}(<\gamma_6>)$ is generated by $\sigma$ and $g=\DiaG{1}{1}{-\omega}{1}$ (and $\rho$ if $F_4=0$) and is isomorphic to $\Z{2}\times\Z{6}$ (respectively $(\Z{6})^2$). 

We described above ($n=2$) the eigenvalues of $g^3$, which are $<(1)^5,(-1)^4>$, and those of $g^2$ ($n=3$), which are $<(1)^3,(\omega)^3,(\omega^2)^3>$. This gives two possibilities for the eigenvalues of $g$ on $\Pic{S}$, which are $<1^1,(-1)^2,(\omega)^2,(\omega^2)^2,(-\omega)^1,(-\omega^2)^1>$ or $<1^3,(\omega)^1,(\omega^2)^1,(-\omega)^2,(-\omega^2)^2>$. As the set of points of $S$ fixed by $g$ is an elliptic curve ($y=0$),\index{Curves!elliptic!birational maps that fix an elliptic curve|idi} the trace of $g$ on $\Pic{S}$ is then $-2$, so its eigenvalues on $\Pic{S}$ are then $<1^1,(-1)^2,(\omega)^2,(\omega^2)^2,(-\omega)^1,(-\omega^2)^1>$ and we get $\rkPic{S}^g=1$. The possibility $G=<g>$ is given in $\nump{1.6}$.

Note that the eigenvalues of $g\sigma$ are $<1^3,(\omega)^1,(\omega^2)^1,(-\omega)^2,(-\omega^2)^2>$, (using those of $\sigma$), so $\rkPic{S}^{g\sigma}=3$ and then $G\not=<g\sigma>$. If $F_4\not=0$, the only remaining possibility is $G=<\sigma,g>$ (case $\nump{1.B6.1}$).

We may now suppose that $\lambda=0$ and $F_6=x^6+y^6$.
Let us study the automorphisms of the surface which have the same action as $g$ on $\mathbb{P}^1$ and describe their eigenvalues on $\Pic{S}$. These are the 6 elements $g, g\sigma,g\rho,g\sigma\rho,g\rho^2$ and $g\sigma\rho^2$. We saw above that $\rkPic{S}^{g\sigma}=3$. Note also that $g\sigma\rho^2=\DiaG{-1}{1}{-\omega}{\omega^2}=\DiaG{1}{-\omega^2}{1}{1}$ is conjugate to $g^{-1}$ by the automorphism of $S$ that exchanges the coordinates $x$ and $y$. By the same automorphism (which commutes with $\sigma$), $g\rho^2$ is conjugate to $g^{-1}\sigma=(g\sigma)^{-1}$. It remains to consider $g\rho=\DiaG{1}{1}{-\omega}{\omega}$ and $g\sigma\rho$. In the case $n=3$, we computed the eigenvalues on $\Pic{S}$ of $(g\rho)^2=\DiaG{1}{1}{\omega^2}{\omega^2}$, which are $<(1)^5,(\omega)^2,(\omega^2)^2>$. This implies that $\omega,\omega^2,-\omega,-\omega^2$ are all eigenvalues of $g\rho$ of multiplicity one. As $g\rho$ fixes exactly $3$ points on $S$ ($v$ and $(\pm 1:1:0:0)$), its eigenvalues are thus $<1^3,(-1)^2,(\omega)^1,(\omega^2)^1,(-\omega)^1,(-\omega^2)^1>$, whence $\rkPic{S}^{g\rho}=3$. From this we deduce directly that $g\sigma\rho$ has the same eigenvalues, so $\rkPic{S}^{g\sigma\rho}=3$. Up to conjugation, there is then only one possibility for $G$ with $\ker \pi \cap G=\{1\}$, which is explicitly $<g>$ (case $\nump{1.6}$).

Let us describe the other possibilities:
\begin{itemize}
\item
If $G\cong \Z{2}\times\Z{6}$, it is generated by $\sigma$ and one of the six elements $g, g\sigma,g\rho,g\sigma\rho$, $g\rho^2,g\sigma\rho^2$. Taking $g$ or $g\sigma$, we get $G= <\sigma,g>$ (case $\nump{1.B6.1}$). Taking $g\rho$ or $g\sigma\rho$, we get $G=<\sigma,g\rho>$  (case $\nump{1.B6.1}$). The last case is $G=<\sigma,g\rho^2>$. This group is conjugate to $<\sigma,g>$ by the automorphism of $S$ that exchanges the coordinates $x$ and $y$.
\item
If $G\cong (\Z{6})^2$, it is the entire group $\pi^{-1}(<\gamma_6>)$ (case $\nump{1.\sigma\rho 6}$).
\item
If $G\cong \Z{3}\times\Z{6}$, it is generated by $\rho$ and one of the six elements \begin{center}$g, g\sigma,g\rho,g\sigma\rho,g\rho^2,g\sigma\rho^2$\end{center}. Taking $g$ or $g\rho$ or $g\rho^2$, we obtain the group $G=<g,\rho>$ (case $\nump{1.\rho 6}$). The other possibility is that $G=<\rho,\sigma \rho^2 g>$, which is conjugate to $<\rho,g^{-1}>=<\rho,g>$ by the automorphism of $S$ that exchanges the coordinates $x$ and $y$.
\end{itemize}
\item
\cadre{$n=6$, $b=-1$}: {$F_4=\lambda x^4$ \it and $F_6=\mu x^3y^3$,}\\
for some $\lambda,\mu \in \K$, but this gives a singular surface.

\item
\cadre{$n=6$, $b=\omega=(\zeta_6)^2$}: {$F_4=\lambda y^4$ \it and $F_6=\alpha x^6+\mu y^6$,}\\
for some $\lambda,\mu,\alpha \in \K$. The group $\pi^{-1}(<\gamma_6>)$ is generated by $G_S$ and $\DiaG{1}{1}{\zeta_6}{\omega}=\DiaG{1}{1}{-\omega^2}{\omega}=\DiaG{-1}{-\omega}{1}{1}$. Exchanging the coordinates $x$ and $y$, we get the case $b=1$ treated above.
\item
\cadre{$n=6$, $b=\omega^2=(\zeta_6)^4$}: {$F_4=\lambda x^2y^2$ \it and $F_6=\alpha x^6+\beta y^6$,}\\
for some $\lambda,\alpha,\beta \in \K$.
As $\alpha\beta\not=0$, since the surface is non-singular we may suppose, up to a linear change of coordinates, that $\alpha=\beta=1$. 

The group $\pi^{-1}(<\gamma_6>)$ is generated by $\sigma$ and $g=\DiaG{1}{1}{-\omega}{\omega}$ (and $\rho$ if $F_4=0$) and is isomorphic to $\Z{2}\times\Z{6}$ (respectively $(\Z{6})^2$). 

In the case $n=3$, we computed the eigenvalues on $\Pic{S}$ of $g^2=\DiaG{1}{1}{\omega^2}{\omega^2}$, which are $<(1)^5,(\omega)^2,(\omega^2)^2>$. This implies then that $\omega,\omega^2,-\omega,-\omega^2$ are all eigenvalues of $g$ of multiplicity one. As $g$ fixes exactly $3$ points on $S$ ($v$ and $(\pm 1:1:0:0)$), its eigenvalues are then $<1^3,(-1)^2,(\omega)^1,(\omega^2)^1,(-\omega)^1,(-\omega^2)^1>$ and then $\rkPic{S}^{g}=3$. From this we deduce directly that the eigenvalues of $g\sigma$ are exactly the same, so $\rkPic{S}^{g\sigma}=3$.

If $F_4\not=0$, the only possibility is $G=<\sigma,g>$  (case $\nump{1.B6.2}$). The others cases when $F_4=0$ were treated above ($n=6$, $b=1$).
\item
\cadre{$n=6$, $b=\zeta_6$}: {$F_4=\lambda x^2y^2$ \it and $F_6=\mu x^3y^3$,}\\
for some $\lambda,\mu \in \K$. We also get a singular surface.
\item
\cadre{$n=6$, $b=(\zeta_6)^{-1}$}: {$F_4=\lambda x^4$ \it and $F_6=\mu x^3y^3$,}\\
for some $\lambda,\mu \in \K$. We also get a singular surface.
\item 
\cadre{$n>6$}: {\it both $F_6$ \it and $F_4$ are monomials}.\\
 Then, either $y^2$ or $x^2$ divides $F_6$, which implies that $F_4$ is a multiple of $x^4$ or $y^4$. Up to an exchange of coordinates, we may suppose then that $F_4=x^4$ and $F_6=xy^5$ or $F_6=y^6$.

In the first case, the equation of the surface is $w^2=z^3+x^4z+xy^5$ whose group of automorphisms is $\pi^{-1}(\gamma_{10})\cong \Z{20}$, generated by $\DiaG{\im}{1}{\zeta_{10}}{-1}$.
 So $G$ is the entire group $\pi^{-1}(\gamma_{10})$ and we get $\nump{1.B10}$.

In the second case, the equation of the surface is $w^2=z^3+x^4z+y^6$, whose group of automorphisms is $\pi^{-1}(\gamma_{12})\cong \Z{2}\times\Z{12}$, generated by $\sigma$ and $g=\DiaG{\im}{1}{\zeta_{12}}{-1}$. Let us describe the eigenvalues of $g$ on $\Pic{S}$. Note that $g^3=\DiaG{-\im}{1}{\im}{-1}=\DiaG{1}{-\im}{1}{1}$ and has eigenvalues $<(1)^2,(-1)^3,(\im)^2,(-\im)^2>$ (see above the description of the case $\nump{1.B4.1}$). Then, $1$ is an eigenvalue of $g$ of multiplicity $2$ and $-1$ of mulitiplicity $\geq 1$. So $\rkPic{S}^{g}=2$ and $\rkPic{S}^{\sigma g}\geq 2$. The group $G$ must be the whole group $\pi^{-1}(\gamma_{12})$ and we get $\nump{1.B12}$.

In fact, using the fact that $g^2$ has eigenvalues $<1^3,(\omega)^1,(\omega^2)^1,(-\omega)^2,(-\omega^2)^2>$, the eigenvalues of $g$ are $<1^2,-1,-\omega,-\omega^2,\zeta_{12},\zeta_{12}^5,\zeta_{12}^7,\zeta_{12}^{11}>$, so the trace on $\Pic{S}$ is $2$. This confirms the Lefschetz formula, as the automorphism fixes $4$ points on $S$ ($v$, $(\pm 1:0:1:0)$ and $(0:1:0:0)$).\proofend
\end{itemize}
\end{proof}
\chapter{Conjugation between cases}
\label{Chap:ConjBetweenTheCases}
\markboths{Chapter \thechapter. Conjugation between cases}{Chapter \thechapter. Conjugation between cases}
\ChapterBegin{In this chapter, we give the possible birational conjugations between the cases obtained in previous chapters.}
In the previous chapters, we enumerated all possibilities for finite abelian subgroups of the Cremona group. We now have to decide whether these represent different conjugacy classes of the Cremona group. First of all, the main distinction between the minimal pairs $(G,S)$ is given by $\rkPic{S}^G$ (see Proposition \refd{Prp:TwoCases}), which can be either $1$ or $2$. But, as we explained earlier, although the two cases are isomorphically distinct, they are not birationally distinct.

In fact, the existence of a representation as a finite group of automorphisms of a conic bundle is equivalent to the existence of an invariant pencil of rational curves.\\
 This is the main conjugacy invariant (see Section \refd{Sec:PencilRatCurves}). Let us state the following important result:
\begin{Prp}\fauxtitre
\label{Prp:ConjBetweenRankPic}
\index{Del Pezzo surfaces!automorphisms}
\index{Del Pezzo surfaces!conic bundles on the surfaces}\index{Conic bundles!on Del Pezzo surfaces}
Let $G\subset \Aut(S)$ be a finite abelian group of automorphisms of a rational surface $S$, and suppose that $\rkPic{S}^G=1$. Then, up to birational conjugation, one and only one of the following occurs:
\begin{itemize}
\item[\upshape 1.]
$G$ preserves a pencil of rational curves of $S$.
\item[\upshape 2.]
$S=\Pn$ and $G=V_9\cong (\Z{3})^2$ is generated by  $(x:y:z)\mapsto (x:\omega y:\omega^2 z)$ and $(x:y:z) \mapsto (y:z:x)$.
\item[\upshape 3.]
$S$ is a Del Pezzo surface of degree $4$ and the pair $(G,S)$ is isomorphic to $\nump{4.222}$, $\nump{4.2222}$ or $\nump{4.42}$.
\item[\upshape 4.]
$S$ is a Del Pezzo surface of degree $3$.
\item[\upshape 5.]
$S$ is a Del Pezzo surface of degree $2$.
\item[\upshape 6.]
$S$ is a Del Pezzo surface of degree $1$.
\end{itemize}
\end{Prp}
\begin{proof}\fauxtitred\upshape
Observe first of all that $S$ is a Del Pezzo surface (Lemma \refd{Lem:rkDelPezzo}). 
Let us now prove that one of the cases listed above occurs.
\begin{itemize}
\item
If the degree of the surface is $\geq 5$, the group $G$ is birationally conjugate to a subgroup of $\Aut(\Pn)$ or $\Aut(\mathbb{P}^1\times\mathbb{P}^1)$ (Proposition \refd{Prp:LargeDegree}). \begin{itemize}
\item
In the first case, either the group is diagonalisable, and then leaves invariant the lines through some of the fixed points, or the group is conjugate to $V_9$, defined earlier (Proposition \refd{Prp:PGL3Cr}).
 \item
 In the second case, the group preserves a pencil of rational curves (Proposition \refd{Prp:AutP1Bir}).
 \end{itemize}
 \item
 If the degree of the surface is $4$, by Proposition \refd{Prp:DelP4} either $G$ preserves a pencil of rational curves or $(G,S)$ is isomorphic to one of the pairs   $\nump{4.222}$, $\nump{4.2222}$, $\nump{4.42}$.
 \end{itemize}
 
We now prove that all these possibilities are distinct, by showing that case $i$ is not birationally conjugate to case $j$, for any $j<i$.
\begin{itemize}
\item[$i=2$]
Suppose that $V_9$ preserves a pencil of rational curves. By Proposition \refd{Prp:CBdistCases}, it would be conjugate to a subgroup of $\Aut(\mathbb{P}^1\times\mathbb{P}^1)$. We use then Proposition \refd{Prp:AutP1Bir} and see that the group would be conjugate to the group generated by $(x,y) \mapsto (\omega x,y)$ and $(x,y) \mapsto (x,\omega y)$, that fixes four points of $\mathbb{P}^1\times\mathbb{P}^1$. But this is not possible, as $V_9$ fixes no point of\hspace{0.2 cm}$\Pn$. 
\item[$i=3$]
Note that none of the groups of pairs $\nump{4.222}$, $\nump{4.2222}$ or $\nump{4.42}$ is isomorphic to $(\Z{3})^2$. Cases $2$ and $3$ are then birationally distinct. Furthermore, none of these four groups leaves invariant a pencil of rational curves (see Proposition \refd{Prp:DelP4}).
\item[$i=4$]
We use Proposition \refd{Prp:Cubics} to find the possibilities for a pair $(G,S)$, when $S$ is a Del Pezzo surface of degree $3$. Except for the last pair of the list ( $\nump{3.6.2}$), every group contains an element of the type $\DiaG{\omega}{1}{1}{1}$, which has order $3$ and fixes an elliptic curve.\index{Curves!elliptic!birational maps that fix an elliptic curve|idi} This is not the case for the previous cases 1, 2, and 3. 

There remains to study the pair  $\nump{3.6.2}$. The surface has equation $wx^2+w^3+y^3+z^3+\lambda wyz$ in $\mathbb{P}^3$ and the group is cyclic of order $6$, generated by $\alpha=\DiaG{1}{-1}{\omega}{\omega^2}$. The involution $\alpha^3$ fixes the elliptic curve of equation $x=0$,\index{Curves!elliptic!birational maps that fix an elliptic curve|idi} $w^3+y^3+z^3+\lambda wyz$. Note that this group is isomorphic neither to $V_9$, nor to any group of the case $j=3.$ It remains then to show that this group does not leave invariant a pencil of rational curves. Suppose this false, which implies that the group is birationally conjugate to a group $\tilde{G}=<\tilde{\alpha}>\subset \Aut(\tilde{S},\pi)$ of some conic bundle. Then, $\tilde{\alpha}^3$ acts trivially on the fibration, as it fixes an elliptic curve,\index{Curves!elliptic!birational maps that fix an elliptic curve|idi} and twists $4$ singular fibres. The action of $\tilde{\alpha}$ on the fibration therefore has order $3$ (as no twisting de Jonqui\`eres involution has a root which acts trivially on the fibration, see Lemma \refd{Lem:NoRoot}).
But this is not possible, since (Proposition \refd{Prp:DirectProduct}), the genus of the curve fixed by the involution (here $1$) must be equal to $3 \pmod{4}$.
\item[$i=5$]
If $S$ is a Del Pezzo surface of degree $2$, then $G$ contains either the Geiser involution of the surface or an element of the type $\DiaG{1}{1}{1}{\im}$, which has order $4$ and fixes an elliptic curve\index{Curves!elliptic!birational maps that fix an elliptic curve|idi} (see Proposition \refd{Prp:DelP2}). The description of the fixed point sets shows that neither of these two elements belongs to any group of the previous cases.
\item[$i=6$]
If $S$ is a Del Pezzo surface of degree $1$, then by Proposition \refd{Prp:DelP1}, $G$ contains one of the following elements:
\begin{itemize}
\item
the Bertini involution of the surface;
\item
an element of order $3$ that fixes a curve of genus $2$ (called $\rho$);
\item
an element of order $5$ or $6$ that fixes an elliptic curve\index{Curves!elliptic!birational maps that fix an elliptic curve|idi} (see $\nump{1.5}$ and $\nump{1.6}$);
\end{itemize}
The description of the fixed point set shows that none of these three elements belongs to any group of the previous cases.
\proofend
\end{itemize}
\end{proof}

\bigskip
\bigskip

\index{Conic bundles!automorphisms}
Let us then study the groups that preserve a pencil of rational curves. These are birationally conjugate to a subgroup of automorphisms of a conic bundle, and the action on the fibration gives rise to an exact sequence (introduced in Section \refd{Sec:Exactsequence})
\begin{equation}
1 \rightarrow G' \rightarrow G \stackrel{\overline{\pi}}{\rightarrow} \overline{\pi}(G) \rightarrow 1.
\tag{\refd{eq:ExactSeqCB}}\end{equation}
In Chapter \refd{Chap:FiniteAbConicBundle}, we classified these groups, depending on the structures of $G$ and $G'$ (which determine that of $\overline{\pi}(G)\cong G/G'$). It is now interesting to decide if some group $G$ may be birationally conjugate to a group $H$, with $G'\not\cong H'$. Note that any element  of the group that fixes a hyperelliptic curve of positive genus always acts trivially on the fibration. We prove a more precise result:\index{Curves!hyperelliptic!birational maps that fix a hyperelliptic curve}
\begin{Prp}\titreProp{Possible changes on the exact sequence}\\
\label{Prp:HardChgExacrSeq}
Let $G\subset \Aut(S,\pi)$ be an finite abelian group of automorphisms of a conic bundle $(S,\pi)$.
\begin{itemize}
\item[\upshape 1.]\index{Curves!of positive genus!birational maps that fix a curve of positive genus}
If $G$ contains at least two non-trivial elements that fix a curve of positive genus, $G'\cong (\Z{2})^2$.
\item[\upshape 2.]
If $G$ contains exactly one non-trivial element that fixes a curve of positive genus, this element belongs to $G'$ and  $G'\cong \Z{2}$ or $G'\cong (\Z{2})^2$. 

If $G'\cong (\Z{2})^2$, then
 $G=G'$ if and only if $G$ is birationally conjugate to a group $H \subset \Aut(\tilde{S},\tilde{\pi})$ with $H'\cong \Z{2}$.
\item[\upshape 3.]
If $G$ contains no non-trivial element that fixes a curve of positive genus, it is birationally conjugate either to a subgroup of $\Aut(\mathbb{P}^1\times\mathbb{P}^1)$, or to the group of $\nump{Cs.24}$.
\end{itemize}
\end{Prp}
\begin{proof}\fauxtitred\upshape
Recall that any element of $G$ that fixes a curve of positive genus is a $(S,\pi)$-twisting de Jonqui\`eres involution, and belongs to $G'$ (Lemma \ref{Lem:CurveFixDeJ}). In this case (Proposition \refd{Prp:CBdistCases}), $G'\cong \Z{2}$ or $G'\cong (\Z{2})^2$.

\begin{enumerate}
\item
If two elements of the group fix a curve of positive genus, they generate a group isomorphic to $(\Z{2})^2$, which is then equal to $G'$.
\item
Suppose that $G$ contains only one element $a$, whose fixed curve is of positive genus, and assume that $G\cong (\Z{2})^2$. We denote by $s$ and $as$ the two other involutions of $G'$. We now prove the following assertion:
\begin{equation}\label{eq:as}\begin{array}{c}a\mbox{ \it twists }4\mbox{ \it singular fibres }F_1,...,F_4\mbox{\it , and up to a change in their numbering}\\
s\mbox{ \it twists the fibres }F_3,F_4\mbox{ \it and }sa\mbox{ \it twists the fibres }F_1,F_2.\end{array}\end{equation}
Indeed, we know that $s$ and $as$ do not fix a curve of positive genus, by hypothesis, so  each one twists at most two singular fibres (Lemma \refd{Lem:DeJI}). If $s$ (respectively $as$) twists no singular fibre, then $as$ (respectively $s$) twists four, which is not possible. We obtain then assertion (\refd{eq:as}).
\item[$2.G=G'$]
We suppose that $\overline{\pi}(G)=1$ (i.e.\ $G=G'$) and that the triple $(G,S,\pi)$ is minimal, which implies that the number of singular fibres of $\pi$ is exactly $4$. We now prove that the surface is a special Del Pezzo surface of degree $4$, whose group of automorphisms contains $(\mathbb{F}_2)^4 \rtimes \Z{2}$, and see that $G$ preserves another fibration, where $s$ does not act trivially.

Let $\eta:(S,\pi)\rightarrow (S',\pi')$ denote the $s$-equivariant birational morphism  of conic bundles that contracts one component in each of the two singular fibres not twisted by $s$ (i.e.\ in $F_1$ and $F_2$).  Note that  $\pi'$ has two singular fibres, both twisted by $s':=\eta s \eta^{-1}\in \Aut(S',\pi')$. This implies that $S'$ is the Del Pezzo surface of degree $6$ (see Lemma \refd{Lem:Degree567CB}) of equation $S'=\{ (x:y:z) \times (u:v:w)\ | \ ux=vy=wz\} \subset \Pn\times\Pn$.

Up to isomorphism, we may suppose that $\pi'$ is given as
\begin{center}
$\pi'((x:y:z) \times (u:v:w))=\left\{\begin{array}{ll}
(y:z)& \mbox{ if } (x:y:z) \not= (1:0:0),\\
(w:v)& \mbox{ if } (u:v:w) \not= (1:0:0),\end{array}\right.$\end{center}
 and find that $s'=\kappa_{\alpha,\beta}:(x:y:z) \times (u:v:w) \mapsto (u:\alpha w:\beta v)\times (x:\alpha^{-1} z:\beta^{-1} y)$, for some $\alpha,\beta \in \K^{*}$ (see Lemma \refd{Lem:KappaDP6}). 
As $s$, and hence $s'$, acts trivially on the fibration, then $\alpha=\beta$, and up to a change of variables we may suppose that $\alpha=1$. 

We denote by $\eta':S'\rightarrow \Pn$ the projection on the first factor that is the blow-up of $A_3=(0:0:1)$, $A_4=(0:1:0)$ and $A_5=(1:0:0)$ (we use the standard notation for the Del Pezzo surface of degree $6$, introduced in Section \refd{Subsec:DelPezzo6} but change the numbering to $1,2,3\rightarrow 5,4,3$). The fibres of $\pi'$ correspond then to the lines of\hspace{0.2 cm}$\Pn$ passing through $A_5=(1:0:0)$.

So the situation is as follows:
\begin{center}
$s'=\eta s\eta^{-1}=\kappa_{1,1}:(x:y:z) \times (u:v:w) \mapsto (u: w: v)\times (x: z:y)$.\\
\includegraphics[width=60.00mm]{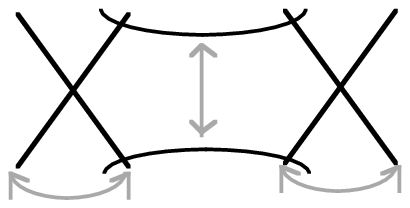}%
\drawat{-33.62mm}{25.47mm}{$D_{34}$}%
\drawat{-20.5mm}{19.05mm}{$E_3$}%
\drawat{-23mm}{11.81mm}{$D_{35}$}%
\drawat{-32.98mm}{5.3mm}{$E_5$}%
\drawat{-43mm}{11.81mm}{$D_{45}$}%
\drawat{-43mm}{19.05mm}{$E_4$}%
\drawat{-61.65mm}{3.66mm}{$\grey{s'}$}
\end{center}
We denote by $A_1',A_2'$ the points of $S'$ blown-up by $\eta$.
As these points are both fixed by $s'$, they do not belong to $E_5$ or $D_{34}$, which are permuted by $s'$, and neither belongs to a singular fibre, as $(S,\pi)$ is a conic bundle. So no blown-up point belongs to an exceptional section of $S$. Furthermore, we observe that the projection by $\eta'$ of the points of $S'$ fixed by $s'$ is the conic of\hspace{0.2 cm}$\Pn$ of equation $x^2=yz$.  The points  $A_1=\eta'(A_1')$ and $A_2=\eta'(A_2')$ of\hspace{0.2 cm}$\Pn$ belong then to the conic and so the points $\{A_1,A_2,A_3,A_4,A_5\}$ are in general position (no $3$ are collinear). This implies that $S$ is a Del Pezzo surface of degree $4$.

Recall that $\eta'\eta:S\rightarrow \Pn$ is the blow-up of $A_1,...,A_5$. As usual, we set $E_i=(\eta\eta')^{-1}(A_i)$ for $i=1,...,5$ and denote by $L$ the pull-back of a line of\hspace{0.2 cm}$\Pn$ by the morphism $\eta'\eta$.
The Picard group of $S$ is generated by $E_1,...,E_5$ and $L$, and the divisor of a fibre of $\pi$ is $L-E_5$ (the pull-back of lines of\hspace{0.2 cm}$\Pn$ by $A_5$, as we mentioned above). Denoting by $D_{15}$ (respectively $D_{25}$) the divisor $L-E_1-E_5$ (respectively $L-E_2-E_5$), the situation is then the following: 
 \begin{center}
 \includegraphics[width=100.00mm]{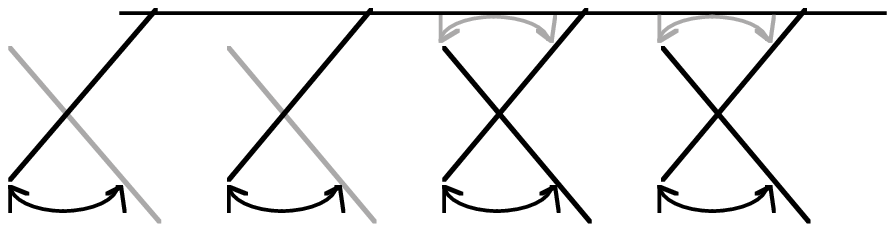}%
 \drawat{-86.25mm}{15mm}{$D_{15}$}%
\drawat{-63.6mm}{15mm}{$D_{25}$}%
\drawat{-41.5mm}{15mm}{$D_{35}$}%
\drawat{-19.05mm}{15mm}{$D_{45}$}%
\drawat{-84.5mm}{7.7mm}{$E_1$}%
\drawat{-62mm}{7.7mm}{$E_2$}%
\drawat{-40.14mm}{7.7mm}{$E_3$}%
\drawat{-17.43mm}{7.7mm}{$E_4$}%
\drawat{-4mm}{22mm}{$E_5$}%
\drawat{-101mm}{22mm}{$\grey{s}$}%
\drawat{-101mm}{3mm}{$a$}
.\end{center}
Note that both $a$ and $s$ fix the divisor $L-E_5$, since $a,s \in \Aut(S,\pi)$. We deduce the actions of $s$ and $a$ on the Picard group from their action on the singular fibres:
\begin{itemize}
\item
$a(E_5)$ is an exceptional divisor that intersects $E_1$, $E_2$, $E_3$ and $E_4$: so is $C=2L-E_1-E_2-E_3-E_4-E_5$, the pull-back by $p$ of the conic passing through the blown-up points, so $a(L)=a(L-E_5)+a(E_5)=3L-E_1-E_2-E_3-E_4-2E_5$. 
\item
$s(E_5)$ is an exceptional divisor that intersects $D_{15}$, $D_{25}$, $E_3$ and $E_4$. It is thus $D_{34}$,  so $s(L)=s(L-E_5)+s(E_5)=2L-E_3-E_4-E_5$. 
\end{itemize}
The matrices of the actions on $\Pic{S}$ of $a$ and $s$, with respect to the basis $E_1,...,E_5,L$ are then respectively:
\begin{center}
$\left(\begin{array}{rrrrrr}
-1 & 0  & 0 & 0 & -1 & -1 \\
0  & -1 & 0 & 0 & -1 & -1 \\
0  & 0  & -1 & 0  & -1 &  -1\\
0  & 0  & 0  & -1 & -1 & -1 \\
-1 & -1 & -1 & -1 & -1 & -2\\
1  & 1 &  1 &  1 &  2 & 3\end{array}\right)$ \ 
$\left(\begin{array}{rrrrrr}
1 & 0  & 0 & 0 & 0 & 0 \\
0  & 1 & 0 & 0 & 0 & 0 \\
0  & 0  & -1 & 0  & -1 &  -1\\
0  & 0  & 0  & -1 & -1 & -1 \\
0 & 0 & -1 & -1 & 0 & -1\\
0  & 0 &  1 &  1 &  1 & 2\end{array}\right).$
\end{center}
We deduce from these matrices the actions of $s$, $a$ and $sa$ on the $16$ exceptional divisors:
\begin{center}
$\begin{array}{rrccl}
s:&\h{1}&  (E_5\ D_{34})(C\ D_{12})&\h{0.5}(D_{13}\ D_{14})(D_{23}\ D_{24})&\h{0.5}(E_3\ D_{35})(E_4\ D_{45})\\
a:&\h{1}(E_1\ D_{15})(E_2\ D_{25})&\h{0.5}  (E_5\ C)(D_{12}\ D_{34})&\h{0.5}(D_{13}\ D_{24})(D_{14}\ D_{23})&\h{0.5}(E_3\ D_{35})(E_4\ D_{45})\\
sa:&\h{1}(E_1\ D_{15})(E_2\ D_{25})&\h{0.5}(E_5\ D_{12})(C\ D_{34})&\h{0.5}(D_{13}\ D_{23})(D_{14}\ D_{24}). \end{array}$
\end{center}

 As $G$ leaves invariant the divisor $L-E_5$, it also leaves invariant the divisor $-K_S-L-E_5$, which is the fibre of a conic bundle $\pi'':S\rightarrow \mathbb{P}^1$, whose singular fibres are $\{E_5,C\}$, $\{D_{12},D_{34}\}$, $\{D_{13},D_{24}\}$ and $\{D_{14},D_{23}\}$. From the above computation, we see that the situation on the conic bundle $(S,\pi'')$  is as follows:
 \begin{center}
 \includegraphics[width=100.00mm]{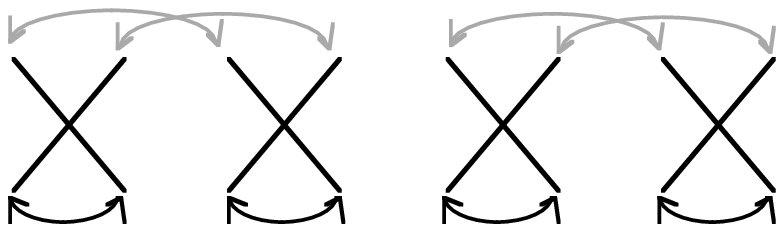}%
\drawat{-84.5mm}{18mm}{$D_{13}$}%
\drawat{-59mm}{18mm}{$D_{23}$}%
\drawat{-33mm}{18mm}{$D_{12}$}%
\drawat{-7.5mm}{18mm}{$E_5$}%
\drawat{-84.5mm}{10mm}{$D_{24}$}%
\drawat{-59mm}{10mm}{$D_{14}$}%
\drawat{-33mm}{10mm}{$D_{34}$}%
\drawat{-7mm}{10mm}{$C$}%
\drawat{-101mm}{25mm}{$\grey{s}$}%
\drawat{-101mm}{4mm}{$a$}%
.\end{center}
 We see then that $s$ does not act trivially on this fibration. More precisely, the exact sequence associated to this conic bundle structure is
 \begin{center}$ 1 \rightarrow <a> \rightarrow G \rightarrow \Z{2} \rightarrow 1$.\end{center}
 \item[$2.G\not=G'$]
We suppose as before that $G$ contains only one element $a$ that fixes a curve of positive genus and that $G'$ is generated by $a$ and $s$, satisfying assertion $(\refd{eq:as})$. But here, we assume that $G'\subsetneqq G$. We now prove that $G$ is not birationally conjugate to a group $G_1 \subset \Aut(S_1,\pi_1)$  with $G_1'\cong \Z{2}$.

First, we prove that any element of 
$G$ that does not act trivially on the fibration is a square root of an involution of $G'$:
 \begin{itemize}
 \item
 Let $h \in G$, with $\overline{\pi}(h)\not=1$ and suppose that $<h>\cap\ G'=\{1\}$. This implies that the group $H$ generated by $h$ and $s$ is equal to $<s>\times <h>\cong \Z{2}\times <h>$ and that the exact sequence associated to this group (which is $1\rightarrow <s>\rightarrow \Z{2}\times<h>\rightarrow <h>\rightarrow 1$) splits. In this case (see Proposition \refd{Prp:DirectProduct})  the genus of the curve fixed by $s$ must be equal to $n-1 \pmod{n}$, where $n$ is the order of $h$. This is not possible, since $s$ fixes a rational curve. Then, $h$ is a root of some involution of $G'$. 
 Note then that $h$ acts on the sets of singular fibres twisted by $s$ (respectively $sa$), as it commutes with $s$ (respectively with $sa$). So the action of $h$ on the fibration has order $2$.
 \end{itemize}
 Then, one of the following situations occurs:
 \begin{enumerate}
 \item
 $G\cong \Z{2}\times\Z{4}$, generated by $s$ and a square root of $a$.
 \item
 $G\cong \Z{2}\times\Z{4}$, generated by $a$ and a square root of $s$ (or of $sa$).
 \item
 $G\cong(\Z{4})^2$, generated by a square root of $a$ and a square root of $s$.
 \end{enumerate}
  In each case, we show that $G$ is not birationally conjugate to a group $G_1 \subset \Aut(S_1,\pi_1)$
  with $G_1'\cong \Z{2}$. If this the case, we denote by $a'$ the conjugate of $a$, and see that $a'$ generates $G_1'$. The group $\overline{\pi_1}(G_1)$ is then isomorphic to $G_1 /G_1'\cong G/<a>$. This gives in the above cases respectively $(\Z{2})^2$, $\Z{4}$ and $\Z{2}\times \Z{4}$.  
 
 The last case is clearly not possible, as $\overline{\pi_1}(G_1)\subset \PGLn{2}$. In the second case, the exact sequence associated to $G_1 \subset \Aut(S_1,\pi_1)$ splits. Using once again Proposition \refd{Prp:DirectProduct}, the genus of the curve fixed by $a$ (here $4$) must be equal to $3 \pmod{4}$. This case is thus also impossible.
 
 It remains to study the case where  $G\cong \Z{2}\times\Z{4}$, generated by $s$ and a square root of $a$, that we denote by $b$. In this case, $\overline{\pi_1}(G_1)$ would be isomorphic to $(\Z{2})^2$ and then $G_1$ fixes no point of $S_1$, as it does not leave invariant any fibre, so $G$ fixes no point of $\tilde{S}$ (Section \refd{Sec:ExistenceFixedPoints}). In this case, none of the two fibres invariant by $b$  singular (otherwise the singular point of each fibre would be fixed by the entire group $G$). Then $b$ twists no fibre of $\tilde{\pi}$. Assuming that $(G,\tilde{S},\tilde{\pi})$ is minimal, each singular fibre is twisted by an element of the group, so the number of singular fibres of $\tilde{\pi}$ is $4$ and we may assume that the surface $S$ is the Del Pezzo surface of degree $4$, exactly as in the previous case ($2.G=G'$).

Let us give some information on the elements $a$ and $s$ and on the surface $S$. 
Note that the action of $a$ and $s$ on the exceptional pairs (pairs of fibres of conic bundle $\{L-E_i,-K_S-L+E_i\}$, see Section \refd{Sec:DP4}) is as follows:
\begin{itemize}
\item
$a$ permutes $L-E_i$ and $-K_S-L+E_i$ for $i=1,...,4$. It fixes both divisors $L-E_5$ and $-K_S-L+E_5$.
\item
The action of $s$ on the exceptional pairs is $(1\ 2)(3\ 4)$ and $s$ permutes the divisors in the following way: $L-E_1\leftrightarrow -K_S-L+E_2,L-E_2\leftrightarrow-K_S-L+E_1,L-E_3\leftrightarrow L-E_4,-K_S-L+E_3\leftrightarrow -K_S-L+E_4$.
\end{itemize}
As the action of $s$ on the exceptional pairs is $(1\ 2)(3\ 4)$, $\Aut(S)$ contains a subgroup $(\mathbb{F}_2)^4\rtimes \Z{2}$ (see Proposition \refd{Prp:ClassDelPezzoQSpec}), and the representations of $a$ and $s$ in this subgroup are $((1,1,1,1,0),id)$, $((1,1,0,0,0),(1\ 2)(3\ 4))$.

 Recall that $\Aut(S) = (\mathbb{F}_2)^4\rtimes H_S$ (Proposition \refd{Prp:AutSF24S5}) and that $H_S$ contains $(1\ 2)(3\ 4)$. Note that $H_S$ is isomorphic either to $\Z{2}$, to $\Z{4}$, or to a diehedral group of order $6$ or $10$  (see Proposition \refd{Prp:ClassDelPezzoQSpec}). So no other element of order $2$ of $H_S$ commutes with $(1\ 2)(3\ 4)$.
 As the representation of $a=b^2$ in  $(\mathbb{F}_2)^4\rtimes H_S$ is $((1,1,1,1,0),id)$, the element $b$ must be equal to $((b_1,b_2,b_3,b_4,b_5),(1\ 2)(3\ 4))$, for some $b_1,...,b_5 \in \mathbb{F}_2$, with $\sum_{i=1}^5 b_i=0$.
 
As $b$ commutes with $s$, whose representation in $(\mathbb{F}_2)^4\rtimes H_S$ is $((1,1,0,0,0),(1\ 2)(3\ 4))$, we find that $b_1=b_2$ and $b_3=b_4$. So $b^2= ((b_1+b_2,b_1+b_2,b_3+b_4,b_3+b4,0),id)$ is the identity, which is impossible.
\item
Suppose now that no element fixes a curve of positive genus. Using once again Propositon \refd{Prp:CBdistCases}, one of the following occurs:
\begin{itemize}
\item
 the group is birationally conjugate to a subgroup of $\Aut(\mathbb{P}^1\times\mathbb{P}^1,\pi_1)$,
 \item
 the group is birationally conjugate either to $\nump{Cs.24}$ or $\nump{P1s.24}$,
 \item
 $G'$ contains an $(S,\pi)$- twisting de Jonqui\`eres involution and $G' \cong \Z{2}$ or $G'\cong (\Z{2})^2$.\end{itemize}
As $\nump{P1s.24}$ is a subgroup of $\Aut(\mathbb{P}^1\times \mathbb{P}^1)$, we need only consider the last case. 

First, we prove this assertion: 
\begin{equation}\label{eq:asshroot}\mbox{\it if }h\in G\mbox{\it\ and\ }\overline{\pi}(h)\not=1\mbox{\it, then }h\mbox{\it\ is a root of some involution of }G'.\end{equation}
Suppose that $<h>\cap\ G'=\{1\}$. This implies that the group $H$ generated by $h$ and one of the $(S,\pi)$-twisting involutions of $G$ is isomorphic to $\Z{2}\times <h>$ and that the exact sequence associated to this group (which is $1\rightarrow \Z{2}\rightarrow \Z{2}\times<h>\rightarrow <h>\rightarrow 1$) splits. Using Proposition \refd{Prp:DirectProduct}, the genus of the curve fixed by the $(S,\pi)$-twisting involution (here $0$) must be equal to $n-1 \pmod{n}$. We get a contradiction, which proves assertion (\refd{eq:asshroot}).

Observe that if the number of singular fibres is $\leq 3$, then the group is conjugate to a subgroup of $\Aut(\mathbb{P}^1\times\mathbb{P}^1)$ (see Proposition \refd{Prp:LargeDegree}). We assume then that $(G,S,\pi)$ is minimal and that the number of singular fibres is at least $4$. Recall that any singular fibre is twisted by an element of $G$ (Lemma \refd{Lem:MinTripl}).

Let us describe the possible elements of $G$:
\begin{itemize}
\item
All elements of $G'$  are involutions, since $G'\cong \Z{2}$ or $G'\cong (\Z{2})^2$. These may be either $(S,\pi)$-twisting involutions, or involutions that do not twist which any singular fibre.
\item
An element which is a root of a $(S,\pi)$-twisting de Jonqui\`eres involution that fixes a rational curve is in fact a square root and twists one singular fibre (see Proposition \refd{Prp:rootsDeJ}, with $k=1$, so $n=2,l=3$).
\item
An element which is a root of an involution of $G'$ that does not twist a singular fibre may twist either $0$ or $2$ fibres (see Proposition \refd{Prp:TwocasesKappaTwist}).
\end{itemize}

Three situations are possible, depending on the number of $(S,\pi)$-twisting involutions (which is positive, as we assumed above):
\begin{itemize}
\item
{\it The group $G'$ contains exactly one $(S,\pi)$-twisting involution, denoted by $s$.}\\
Since the composition of a twisting and a non-twisting involution gives a twisting one, no other involution belongs to $G'$, so $G'=<s>$. Suppose that $G$ is generated by some root $h$ of $s$. The number of singular fibres is then $3$, as the number of singular fibres twisted by $s$ and $h$ are respectively $2$ and $1$. We must then assume that two different square roots of $h$ belong to $G$. Denoting these roots by $h_1$ and $h_2$, we see that $h_1h_2\not=s$, so $\overline{\pi}(h_1h_2)\not=1$, and $(h_1h_2)^2=id$, which contradicts (\refd{eq:asshroot}).
\item
{\it The group $G'$ contains exactly two $(S,\pi)$-twisting involutions.}\\
The group $G'$ is then isomorphic to $(\Z{2})^2$.
Let us denote by $s$ one of the twisting involutions and by $a$ the non-twisting involution. We see that $sa$ and $s$ twist the same $2$ singular fibres, as $a$ does not twist any fibre. 

Since the number of singular fibres is at least $4$, some other singular fibres are twisted by elements of $G$, which may be roots of $a$ or $s$ or $sa$. So one of the following two situations arises:
\begin{itemize}
\item
Another singular fibre is twisted by a root $h$ of $a$. We now prove that this is not possible.  As $s$ and $h$ commute, $h$ permutes the two singular fibres twisted by $s$, so its action on the fibration is of order $2$ and $h^2=a$. Supposing that $G$ is generated by $h$ and $s$, the number of singular fibres is $4$ and the situation is as follows:
 \begin{center}
 \includegraphics[width=100.00mm]{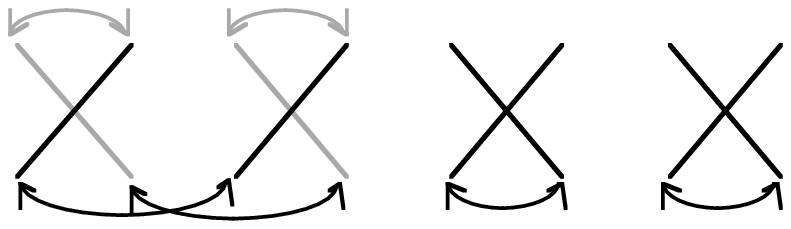}%
\drawat{-101mm}{25mm}{$\grey{s}$}%
\drawat{-101mm}{4mm}{$h$}%
.\end{center}
We denote by $\eta:(S,\pi)\rightarrow (S',\pi')$ the $h$-equivariant birational morphism  of conic bundles which consists of the contraction of one component in every singular fibre twisted by $s$ (in grey on the diagram). Note that  $\pi'$ has two singular fibres and $\eta h \eta^{-1}\in \Aut(S',\pi')$ twists the two singular fibres. This implies that $S'$ is the Del Pezzo surface of degree $6$ (see Lemma \refd{Lem:Degree567CB}) and 
\begin{center}
$h'=\eta h\eta^{-1}=\kappa_{\alpha,\beta}:(x:y:z) \times (u:v:w) \mapsto (u:\alpha w:\beta v)\times (x:\alpha^{-1} z:\beta^{-1} y)$,\end{center} for some $\alpha,\beta \in \K^{*}$  (see Lemma \refd{Lem:KappaDP6}). As the action of $h$ on the fibration has order $2$, we deduce that $\alpha=\beta$. Using the notation of Section \refd{SubSec:DelPezzo6Kappa}, the situation is the following:
 \begin{center}
 \includegraphics[width=90.00mm]{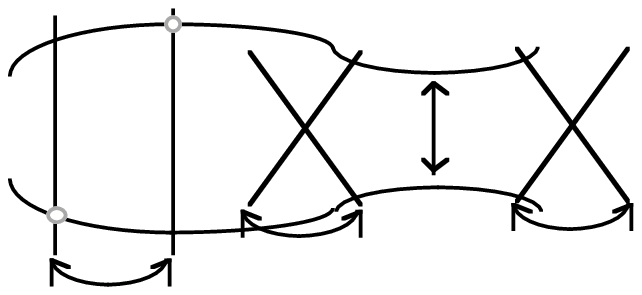}%
\drawat{-89.55mm}{3.83mm}{$h'$}%
\drawat{-43.5mm}{18mm}{$D_{12}$}%
\drawat{-43.5mm}{26.1mm}{$E_2$}%
\drawat{-24.5mm}{18mm}{$D_{13}$}%
\drawat{-21.5mm}{25.3mm}{$E_3$}%
\drawat{-34mm}{31.6mm}{$D_{23}$}%
\drawat{-33.8mm}{11.7mm}{$E_1$}%
.\end{center}
 
 Note that the set of points of $S'$ fixed by \begin{center}$(h')^2=(x:y:z) \times (u:v:w) \mapsto (-x:y:z)\times (-u:v:w)$\end{center} is the union of the two singular points of the singular fibres \begin{center}$(0:1:0)\times (0:0:1)$ and $(0:0:1)\times (0:1:0)$\end{center} and the sections \begin{center}$\begin{array}{l}E_1=\{(1:0:0) \times (0:a:b)\ | \ (a:b) \in \mathbb{P}^1\},\\
 D_{23}=\{(0:a:b) \times (1:0:0)\ | \ (a:b) \in \mathbb{P}^1\}.\end{array}$\end{center}
 
 Let us denote by $\tilde{D_{23}}$, $\tilde{E_1}$, $\tilde{E_2}$, $\tilde{D_{12}}$ the strict pull-backs by $\eta$ of respectively $D_{23}$, $E_1$, $E_2$ and $D_{12}$. 
 
 As $s$ commutes with $h^2$, it must leave invariant the set $\{\tilde{D_{23}},\tilde{E_1}\}$. As it twists two fibres, it must invert $D_{23}$ and $\tilde{E_1}$. This contradicts the fact that $s$ leaves invariant both $\tilde{E_2}$ and $\tilde{D_{12}}$.
\item
Another singular fibre is twisted by a root of $s$ or $sa$. 
As we mentioned above, a root of a $(S,\pi)$-twisting de Jonqui\`eres involution is a square root and twists only one singular fibre. So there must exist two different elements $h_1$ and $h_2$ that are roots of $(S,\pi)$-twisting de Jonqui\`eres involutions, each one twisting a different fibre. We may assume that either  $(h_1)^2=s$ and $(h_2)^2=sa$ or $(h_1)^2=(h_2)^2=s$:
\begin{itemize}
\item
If $(h_1)^2=s$ and $(h_2)^2=sa$, $h=h_1h_2$ is such that $h^2=a$, and we come back to the precedent case.
\item
If $(h_1)^2=(h_2)^2=s$, then $(h_1h_2)^2=1$, so $h_1h_2$ belongs to $G'$, using assertion (\refd{eq:asshroot}). This is not possible, as $h_1h_2$ does not twist the same fibres as $s$ and $sa$.
\end{itemize}
\end{itemize}
\item
{\it The group $G$ contains three $(S,\pi)$-twisting de Jonqui\`eres involutions. }

Note also that two different $(S,\pi)$-twisting involutions do not twist singular fibres which are all distinct, as the product of the two involutions would give a $(S,\pi)$-twisting involution that twists $4$ singular fibres and then fixes an elliptic curve. So two different involutions twist one common singular fibre.\index{Curves!elliptic!birational maps that fix an elliptic curve|idi}

As the number of singular fibres is at least $4$, at least one other singular fibre is twisted by elements of $G$, which is then a square root of one of the twisting involutions. Let us denote the root by $h$ and by $s$ a $(S,\pi)$-twisting involution that is not $h^2$.

Note that $h$ exchanges the two singular fibres twisted by $h^2$. One of these is twisted by $s$ and the other one is not, so $h$ and $s$ do not commute.\proofend
\end{itemize}
\end{enumerate}

\end{proof}

\chapter{The classification and other main results}
\label{Chap:List}
\section{\h{1}Classification of finite cyclic subgroups of the Cremona group}
\pagetitless{The classification and other main results}{Classification of finite cyclic subgroups of the Cremona group}
\begin{flushleft}{\bf THEOREM \refThmCyclicGroup: \sc Classification of finite cyclic groups}\\ \itshape
\index{Del Pezzo surfaces!automorphisms}\index{Conic bundles!automorphisms}
Let $G \subset \CrP$ be a finite cyclic group of order $n$. It is conjugate, in the Cremona group, to the group generated by $\alpha$, where $\alpha$ is one (and only one) of the following:\\
(we use once again the notations $\zeta_m=e^{2\im \pi/m}$, $\omega=\zeta_3$):
\begin{itemize}
\item
{\bf roots of de Jonqui\`eres involutions (see Section \refd{Sec:ExtCyclDeJ})}\\
{\bf \emph{n}$\mathbf{=2}$\emph{m}$\mathbf{>2}$}, $\alpha^m$ is a de Jonqui\`eres involution that fixes a hyperelliptic curve $\Gamma$ of genus $k-1>0${\upshape :}\index{Curves!hyperelliptic!birational maps that fix a hyperelliptic curve|idb}
\begin{itemize}
\item[\nump{C.ro.m}]
{\bf \emph{m} is odd:}\\
$\alpha:(x_1:x_2)\times (y_1:y_2) \dasharrow (x_1:\zeta_m x_2) \times(y_2 \prod_{i=1}^k(x_1-b_i x_2):y_1 \prod_{i=1}^k(x_1-a_i x_2))$\\ 
where $a_1,...,a_k,b_1,...,b_k \in \K^{*}$ are all distinct, and the sets $\{a_1,...,a_k\}$ and $\{b_1,...,b_k\}$ are both invariant by multiplication by $\zeta_m$.
\item[\nump{C.re.m}]
{\bf \emph{m} is even:}\\
$\alpha$ is an automorphism of some conic bundle $(S,\pi)$ with $l$ singular fibres. The action of $\alpha$ on the ramification set of $\pi:\Gamma\rightarrow \mathbb{P}^1$ has $r$ orbits, each consisting of $m$ points (so $rm=2k$). 
Furthermore, one of the following holds:
\begin{itemize}
\item[-]
$l=2k$, $r$ is even, $\alpha$ acts on $\Gamma$ with $4$ fixed points;
\item[-]
$l=2k+1$, $r$ is odd, $\alpha$ acts on $\Gamma$ with $2$ fixed points;
\item[-]
$l=2k+2$, $r$ is even, $\alpha$ fixes no point of $\Gamma$.
\end{itemize}
\end{itemize}
\item
{\bf Linear automorphisms (also described in \cite{bib:BeB})}
\begin{itemize}
\item[$\nump{0.n}$]
$\alpha=\Diag{1}{1}{\zeta_n}$ is a linear automorphism of\hspace{0.2 cm}$\Pn$.
\end{itemize}
\item
{\bf \emph{n}$\mathbf{=2}$: Involutions (also described in \cite{bib:BaB})}
\begin{itemize}
\item[$\nump{C.2}$]
$\alpha$ is a de Jonqui\`eres involution acting on some conic bundle, fixing a hyperelliptic curve of positive genus.\index{Curves!hyperelliptic!birational maps that fix a hyperelliptic curve|idb}\index{Curves!of positive genus!birational maps that fix a curve of positive genus}\\
The conjugacy classes are parametrised by non-rational hyperelliptic curves.\footnote{\label{footnote1}In this statement, when we say that some conjugacy classes are parametrised by some curves, we mean that  the map $NFC$ (see Definition \refd{Def:IsoFixedPoints}) that associates to the group the normalisation of its set of fixed points induces a bijection between the conjugacy classes of  groups of type $\nump{x.y}$ to isomorphism classes of the given curves.}
\item[$\nump{2.G}$]
$\alpha$ is a Geiser involution acting on some Del Pezzo surface of degree $2$. \\
The conjugacy classes are parametrised by smooth quartic plane curves.
\item[$\nump{1.B}$]\index{Involutions!Bertini}\index{Involutions!Geiser}
$\alpha$ is a Bertini involution acting on some Del Pezzo surface of degree $1$. \\
The conjugacy classes are parametrised by smooth curves of genus $4$ lying on a quadric cone.
\end{itemize}
\item
{\bf \emph{n}$\mathbf{=3}$ (also described in \cite{bib:deF})}
\begin{itemize}
\item[$\nump{3.3}$]
$\alpha=\DiaG{\omega}{1}{1}{1}$ acts on the cubic surface $w^3+L_3(x,y,z)$ in $\mathbb{P}^3$. \\
The conjugacy classes are parametrised by elliptic curves.\index{Curves!elliptic!birational maps that fix an elliptic curve}
\item[$\nump{1.\rho}$]
$\alpha=\DiaG{1}{1}{1}{\omega}$ acts on the surface $w^2=z^3+L_6(x,y)$ in $\mathbb{P}(3,1,1,2)$. \\
The conjugacy classes are parametrised by curves of genus $2$.
\end{itemize}
\item
{\bf \emph{n}$\mathbf{=4}$}
\begin{itemize}
\item[$\nump{2.4}$]
$\alpha=\DiaG{1}{1}{1}{\im}$ acts on the surface $w^2=L_4(x,y)+z^4$ in $\mathbb{P}(2,1,1,1)$. \\
The conjugacy classes are parametrised by elliptic curves, $\alpha^2$ is a de Jonqui\`eres involution.\index{Curves!elliptic!birational maps that fix an elliptic curve}
\item[$\nump{1.B2.2}$]
$\alpha=\DiaG{\im}{1}{-1}{-1}$ acts on the surface $w^2=z^3+zL_2(x^2,y^2)+xyL_2'(x^2,y^2)$ in $\mathbb{P}(3,1,1,2)$.  $\alpha^2$ is the Bertini involution of the surface.\\
\end{itemize}
\item
{\bf \emph{n}$\mathbf{=5}$ (also described in \cite{bib:deF} and \cite{bib:BeB})}
\begin{itemize}
\item[$\nump{1.5}$]
$\alpha=\DiaG{1}{1}{\zeta_5}{1}$ acts on the surface $w^2=z^3+\lambda x^4z+x(\mu x^5+y^5)$ in $\mathbb{P}(3,1,1,2)$. \\
The conjugacy classes are parametrised by elliptic curves.\index{Curves!elliptic!birational maps that fix an elliptic curve}
\end{itemize}
\item
{\bf \emph{n}$\mathbf{=6}$}\vspace{0.1cm}\\
\begin{tabular}{|llll|}
\hline
{\it name} &  {\it description} & {\it equation of}& {\it in the } \\
{\it of $\mathit{(G,S)}$} & {\it of $\mathit{\alpha}$}  & {\it the surface} & {\it space } \\
\hline 
${\nump{3.6.1}}$ & $\DiaG{\omega}{1}{1}{-1}$& $w^3+x^3+y^3+xz^2+\lambda yz^2$& $\mathbb{P}^3\hts$  \\
${\nump{3.6.2}}$ & $\DiaG{1}{-1}{\omega}{\omega^2}$& $wx^2+w^3+y^3+z^3+\lambda wyz$& $\mathbb{P}^3\hts$  \\
${\nump{2.G3.1}}$ & $\DiaG{-1}{1}{1}{\omega}$& $w^2=L_4(x,y)+z^3L_1(x,y)$& $\mathbb{P}(2,1,1,1)\hts$  \\
${\nump{2.G3.2}}$ & $\DiaG{-1}{1}{\omega}{\omega^2}$& $w^2=x(x^3+y^3+z^3)+yzL_1(x^2,yz)$& $\mathbb{P}(2,1,1,1)\hts$  \\
${\nump{2.6}}$ & $\DiaG{-1}{\omega}{1}{-1}$& $w^2=x^3y+y^4+z^4+\lambda y^2z^2$& $\mathbb{P}(2,1,1,1)\hts$  \\
$\nump{1.\sigma\rho}$ & $\DiaG{-1}{1}{1}{\omega}$ & $w^2=z^3+L_6(x,y)$ & $\mathbb{P}(3,1,1,2)\hts$  \\
$\nump{1.\rho 2}$ & $\DiaG{1}{1}{-1}{\omega}$ & $w^2=z^3+L_3(x^2,y^2)$ & $\mathbb{P}(3,1,1,2)\hts$  \\
${\nump{1.B3.1}}$ & $\DiaG{-1}{1}{\omega}{1}$ & $w^2=z^3+xL_1(x^3,y^3)z+L_2(x^3,y^3)$ & $\mathbb{P}(3,1,1,2)\hts$  \\
${\nump{1.B3.2}}$ & $\DiaG{-1}{1}{\omega}{\omega}$ & $w^2=z^3+\lambda x^2 y^2z+L_2(x^3,y^3)$ & $\mathbb{P}(3,1,1,2)\hts$  \\
${\nump{1.6}}$  & $\DiaG{1}{1}{-\omega}{1}$ & $w^2=z^3+\lambda x^4z+\mu x^6+y^6$ & $\mathbb{P}(3,1,1,2)\hts$  \\
\hline
\end{tabular}

The fixed curves of $\alpha^2$ and $\alpha^3$ are given in the following table:\vspace{0.1cm}\\
\begin{tabular}{|lll|ll|}
\hline
{\it name} &  \multicolumn{2}{l}{\h{2.5}\it curve fixed } &  \multicolumn{2}{l|}{\it curve fixed }\\
{\it of \h{0.5}$\mathit{(G,S)}$} &\multicolumn{2}{l}{ \h{2.5}\it by $\mathit{\alpha^3}$} &\multicolumn{2}{l|}{\it by $\mathit{\alpha^2}$}  \\
\hline 
${\nump{3.6.1}}$ & \h{3}elliptic &\h{1.5} $w^3+x^3+y^3$& elliptic\index{Curves!elliptic!birational maps that fix an elliptic curve}  & \h{1.5}$w^3+x^3+y^3\hts$ \\
${\nump{3.6.2}}$ & \h{3}elliptic & \h{1.5}$wx^2+w^3+z^3\hts$&  &\h{1.5} \\
${\nump{2.G3.1}}$ & \h{3}quartic & \h{1.5}$L_4(x,y)+z^3L_1(x,y)$& elliptic  & \h{1.5}$w^2=L_4(x,y)\hts$ \\
${\nump{2.G3.2}}$ & \h{3}quartic & \h{1.5}$x(x^3\h{0.7}+\h{0.7}y^3\h{0.7}+\h{0.7}z^3)\h{0.7}+\h{0.7}yzL_1(x^2,yz)\hts$& &  \\
${\nump{2.6}}$ &\h{3}  &\h{1.5} & elliptic & \h{1.5}$w^2=y^4+z^4+\lambda y^2z^2\hts$ \\
$\nump{1.\sigma\rho}$ & \h{3}genus \h{0.5}$4$ & \h{1.5}$z^3+L_6(x,y)$ & genus \h{0.5}$2$ & \h{1.5}$w^2=L_6(x,y)\hts$ \\
$\nump{1.\rho 2}$ & \h{3}elliptic  & \h{1.5}$w^2=z^3+L_3(x^2)$ & genus \h{0.5}$2$ & \h{1.5}$w^2= L_3(x^2,y^2)\hts$ \\
${\nump{1.B3.1}}$ & \h{3}genus \h{0.5}$4$ & \h{1.5}$z^3\h{0.7}+\h{0.7}xL_1(x^3,y^3)z\h{0.7}+\h{0.7}L_2(x^3,y^3)$ & elliptic & \h{1.5}$w^2=z^3\h{0.7}+\h{0.7}L_1(x^4)z\h{0.7}+\h{0.7}L_2(x^3)\hts$\h{2}   \\
${\nump{1.B3.2}}$ & \h{3}genus \h{0.5}$4$ & \h{1.5}$z^3+\lambda x^2 y^2z+L_2(x^3,y^3)\hts$ & &\h{1.5} \\
${\nump{1.6}}$  & \h{3}elliptic & \h{1.5}$w^2=z^3+\lambda x^4z+\mu x^6$ & elliptic & \h{1.5}$w^2=z^3+x^4z+\lambda x^6\hts$ \\
\hline
\end{tabular}
\item
{\bf \emph{n}$\mathbf{=8}$}
\begin{itemize}
\item[$\nump{1.B4.2}$]
$\alpha=\DiaG{\zeta_8}{1}{\im}{-\im}$ acts on the surface $w^2=\lambda x^2y^2z +xy(x^4+y^4)$ in $\mathbb{P}(3,1,1,2)$.\\ $\alpha^4$ is the Bertini involution of the surface.\end{itemize}
\item
{\bf \emph{n}$\mathbf{=9}$}
\begin{itemize}
\item[${\nump{3.9}}$]
$\alpha=\DiaG{\zeta_9}{1}{\omega}{\omega^2}$ acts on the surface $w^3+xz^2+x^2y+y^2z$ in $\mathbb{P}^3$.\\ $\alpha^3$ fixes an elliptic curve.\index{Curves!elliptic!birational maps that fix an elliptic curve}\end{itemize}
\item
{\bf \emph{n}$\mathbf{=10}$}
\begin{itemize}
\item[${\nump{1.B5}}$]
$\alpha=\DiaG{-1\h{0.5}}{\h{0.5}1\h{0.5}}{\h{0.5}\zeta_5\h{0.5}}{\h{0.5}1}$ acts on the surface $w^2=z^3\h{0.5}+\h{0.5}\lambda x^4z\h{0.5}+\h{0.5}x(\mu x^5\h{0.5}+\h{0.5}y^5)$ in $\mathbb{P}(3,1,1,2)$.\\ $\alpha^5$ is the Bertini involution of the surface and $\alpha^2$ fixes an elliptic curve.\index{Curves!elliptic!birational maps that fix an elliptic curve}\end{itemize}
\item
{\bf \emph{n}$\mathbf{=12}$}
\begin{itemize}
\item[${\nump{3.12}}$]
$\alpha=\DiaG{\omega}{1}{-1}{\im}$ acts on the surface $w^3+x^3+yz^2+y^2x$ in $\mathbb{P}^3$.\\ Both $\alpha^6$ and $\alpha^4$ fix elliptic curves.\index{Curves!elliptic!birational maps that fix an elliptic curve}
\item[${\nump{2.12}}$]
$\alpha=\DiaG{1}{\omega}{1}{\im}$ acts on the surface $w^2=x^3y+y^4+z^4$ in $\mathbb{P}(2,1,1,1)$.\\ Both $\alpha^4$ and $\alpha^3$ fix elliptic curves.
\item[${\nump{1.\sigma\rho 2.2}}$]
$\alpha=\DiaG{\im}{1}{-1}{-\omega}$ acts on the surface $w^2=z^3+xyL_2(x^2,y^2)$ in $\mathbb{P}(3,1,1,2)$.\\ $\alpha^6$ is the Bertini involution of the surface and $\alpha^4$ fixes an elliptic curve.
\end{itemize}
\item
{\bf \emph{n}$\mathbf{=14}$}
\begin{itemize}
\item[${\nump{2.G7}}$]
$\alpha=\DiaG{-1\h{0.5}}{\h{0.5}\zeta_7\h{0.5}}{\h{0.5}(\zeta_7)^4\h{0.5}}{\h{0.5}(\zeta_7)^2}$ acts on the surface $w^2=\h{0.5}x^3y\h{0.5}+\h{0.5}y^3z\h{0.5}+\h{0.5}xz^3$ in $\mathbb{P}(2,1,1,1)$.\\ $\alpha^7$ is the Geiser involution of the surface.
\end{itemize}
\item
{\bf \emph{n}$\mathbf{=15}$}
\begin{itemize}
\item[${\nump{1.\rho 5}}$]
$\alpha=\DiaG{1}{1}{\zeta_5}{\omega}$ acts on the surface $w^2=z^3+x(x^5+y^5)$ in $\mathbb{P}(3,1,1,2)$.\\ Both $\alpha^3$ and $\alpha^5$ fix elliptic curves.
\end{itemize}
\item
{\bf \emph{n}$\mathbf{=18}$}
\begin{itemize}
\item[${\nump{2.G9}}$]
$\alpha=\DiaG{-1}{(\zeta_9)^6}{1}{\zeta_9}$ acts on the surface $w^2=x^3y+y^4+xz^3$ in $\mathbb{P}(2,1,1,1)$.\\ $\alpha^9$ is the Geiser involution of the surface.
\end{itemize}
\item
{\bf \emph{n}$\mathbf{=20}$}
\begin{itemize}
\item[${\nump{1.B10}}$]
$\alpha=\DiaG{\im}{1}{\zeta_{10}}{-1}$ acts on the surface $w^2=z^3+x^4z+xy^5$ in $\mathbb{P}(3,1,1,2)$.\\ $\alpha^{10}$ is the Bertini involution of the surface.
\end{itemize}
\item
{\bf \emph{n}$\mathbf{=24}$}
\begin{itemize}
\item[${\nump{1.\sigma\rho 4}}$]
$\alpha=\DiaG{\zeta_{8}}{1}{\im}{-\im\omega}$ acts on the surface $w^2=z^3+xy(x^4+y^4)$ in $\mathbb{P}(3,1,1,2)$.\\ $\alpha^{12}$ is the Bertini involution of the surface and $\alpha^8$ fixes an elliptic curve.\index{Curves!elliptic!birational maps that fix an elliptic curve}
\end{itemize}
\item
{\bf \emph{n}$\mathbf{=30}$}
\begin{itemize}
\item[${\nump{1.\sigma\rho 5}}$]
$\alpha=\DiaG{-1}{1}{\zeta_5}{\omega}$ acts on the surface $w^2=z^3+x(x^5+y^5)$ in $\mathbb{P}(3,1,1,2)$.\\ $\alpha^{15}$ is the Bertini involution of the surface and both $\alpha^{6}$ and $\alpha^5$ fix elliptic curves.
\end{itemize}
\end{itemize} 
\end{flushleft}
\begin{proof}\fauxtitred\upshape
As $G$ is finite, it is birationally conjugate to a group of automorphisms of a rational surface $S$. We may thus assume that $(G,S)$ is minimal and then one of the following situations occurs (Proposition \refd{Prp:TwoCases}):
\begin{enumerate}
\item
The surface $S$ has a conic bundle structure invariant by $G$, and $\rkPic{S}^G=2$, i.e.\ the fixed part of the Picard group is generated by the canonical divisor and the divisor class of a fibre.
\item
$\rkPic{S}^{G}=1$, i.e.\ the fixed part of the Picard group is generated over $\mathbb{Q}$ by the canonical divisor.
\end{enumerate}
\begin{itemize}
\item {$\mathbf{rk\ Pic(}${\bf \emph{S}}$\mathbf )^{\mbox{\scriptsize \bf \emph{G}}}\mathbf{=2}$.}\\
We use Proposition \refd{Prp:CBdistCases} and obtain one of the following cases:
\begin{itemize}
\item
$(S,\pi)=(\mathbb{P}^1\times \mathbb{P}^1,\pi_1)$;
\item
some twisting $(S,\pi)$-de Jonqui\`eres involution belongs to $G$.
\end{itemize}
The first case gives a linear automorphism (Proposition \refd{Prp:AutP1Bir}). 

In the second case the group $G$, generated by $\alpha$, is isomorphic to $G \cong \Z{2m}$ and $\alpha^m$ is a twisting $(S,\pi)$-de Jonqui\`eres involution, that fixes a hyperelliptic curve of genus $k-1$;\index{Curves!hyperelliptic!birational maps that fix a hyperelliptic curve} the restriction of $\pi$ to $\Gamma$ presents $\Gamma$ as the double covering of $\mathbb{P}^1$ ramified over $2k$ points (see Lemma \refd{Lem:DeJI}). The action of $G$ on the conic bundle gives an exact sequence (Section \refd{Sec:Exactsequence})
\begin{equation}
1 \rightarrow G' \rightarrow G \stackrel{\overline{\pi}}{\rightarrow} \overline{\pi}(G) \rightarrow 1
\tag{\refd{eq:ExactSeqCB}}\end{equation}
and $G'$ is generated by $\alpha^m$ (using once again Proposition \refd{Prp:CBdistCases}, and the fact that $G$ is cyclic). Then, $\overline{\pi}(G) \cong \Z{m}$; we are in the case of the extension of a cyclic group by a de Jonqui\`eres involution (described in Section \refd{Sec:ExtCyclDeJ}).
\begin{itemize}
\item 
If $m=1$, $G$ is generated by a de Jonqui\`eres involution. If the involution does not fix a curve of positive genus, we get a linear automorphism of\hspace{0.2 cm}$\Pn$ ($\nump{0.n}$, for $n=2$). Otherwise, we get the case $\nump{C.2}$. Proposition \refd{Prp:classDeJICremona} gives the parametrisation by the fixed curve. 
\item
If $m>1$ is odd, then (\refd{eq:ExactSeqCB}) splits, and using Lemma \refd{Lem:splitsequence}, we get $k>1$ and the case mentioned in the statement of the theorem.
\item
If $m$ is even, then (\refd{eq:ExactSeqCB}) does not split, and using Proposition \refd{Prp:rootsDeJ}, we get the statement of the theorem.

It remains to prove that we can remove the case $k=1$. As $rm=2k$, we have $m=2$ and $r=1$. The number of singular fibres is then $3$. Using Proposition \refd{Prp:LargeDegree}, the group is conjugate to a linear group of automorphisms of\hspace{0.2 cm}$\Pn$.
\end{itemize}
\item {$\mathbf{rk\ Pic(}${\bf \emph{S}}$\mathbf )^{\mbox{\scriptsize \bf \emph{G}}}\mathbf{=1}$.}\\
Using Lemma \refd{Lem:rkDelPezzo}, we see that the surface is a Del Pezzo surface. If its degree is $\geq 5$, the group is conjugate to a subgroup of $\Aut(\Pn)$ (see Proposition \refd{Prp:LargeDegree}). We enumerate the possible cases of automorphisms of Del Pezzo surfaces of degree $1$, $2$, $3$ and $4$ (given respectively in Propositions \refd{Prp:DelP1}, \refd{Prp:DelP2}, \refd{Prp:Cubics} and \refd{Prp:DelP4}) and give the cyclic ones.
\end{itemize}
Whenever an element of the groups in the list fixes some curve of positive genus, we said so and wrote down the nature of the curve.  Looking at them, we find that two groups of different cases represent different conjugacy classes.

It remains to prove the parametrisations by the fixed curve stated in the proposition. First if two groups are birationally conjugate, they have isomorphic sets of fixed points (see Section \refd{Sec:CurveFixedPoints}). In cases $\nump{2.G}$, $\nump{1.B}$, $\nump{3.3}$, $\nump{1,\rho}$, $\nump{2.4}$ and $\nump{1.5}$, the groups are cyclic of order $n$ and the surface where the group acts is a $n$-covering over some other surface, ramified over the curve of fixed points and some finite number of points.

The isomorphism class of the curve then induces the isomorphism class of the surface and we remark that two groups of the same surface of cases $\nump{2.G}$, $\nump{1.B}$, ..., $\nump{1.5}$ are conjugate by an isomorphism of the surface.
 The case of $\nump{C.2}$ follows from Proposition \refd{Prp:classDeJICremona}.
\proofend
\end{proof}

\break

\section{Classification of finite abelian non-cyclic subgroups of the Cremona group}
\pagetitless{The classification and other main results}{Classification of finite abelian non-cyclic subgroups of the Cremona group}
\begin{flushleft}{\bf THEOREM \refThmNonCyclicGroup: \sc Classification of finite non-cyclic abelian groups}\\ \itshape
\index{Del Pezzo surfaces!automorphisms}\index{Conic bundles!automorphisms}
A finite non-cyclic abelian subgroup $G$ of $\CrP$ is conjugate, in the Cremona group, to one (and only one) of the following (we use once again the notation $\zeta_m=e^{2\im \pi/m}$, $\omega=\zeta_3$):
\begin{flushleft}{\Large \sc Automorphisms of conic bundles:}\end{flushleft}
\begin{itemize}
\item
{\bf Automorphisms of $(\mathbb{P}^1\times\mathbb{P}^1)$:}
\begin{flushleft}\begin{tabular}{|lll|}
\hline
{\it name} & {\it structure} & {\it generators} \\
\hline 
${\nump{0.mn}}$ & $ \Z{n}\times\Z{m}$&$(x,y) \mapsto (\zeta_n x,y)$, 
$(x,y) \mapsto (x,\zeta_m y)\hts$\\
${\nump{P1.22n}}$ & $ \Z{2}\times\Z{2n} $&$(x,y) \mapsto (x^{-1},y)$, 
$(x,y) \mapsto (-x,\zeta_{2n} y)\hts$\\
${\nump{P1.222n}}$ & $ (\Z{2})^2\times\Z{2n} $&$(x,y) \mapsto (\pm x^{\pm 1},y)$, 
$(x,y) \mapsto (x,\zeta_{2n} y)\hts$\\
${\nump{P1.22.1}}$ & $ (\Z{2})^2 $&$(x,y) \mapsto (\pm x^{\pm 1},y)\hts$\\
${\nump{P1.222}}$ & $ (\Z{2})^3 $&$(x,y) \mapsto (\pm x,\pm y)$, 
$(x,y) \mapsto (x^{-1},y)\hts$\\
${\nump{P1.2222}}$ & $(\Z{2})^4$&$(x,y) \mapsto (\pm x^{\pm 1},\pm y^{\pm 1})\hts$\\
${\nump{P1s.24}}$ & $\Z{2}\times \Z{4}$&  $(x,y) \mapsto (x^{-1},y^{-1})$, $(x,y) \mapsto (-y,x)\hts$\\
${\nump{P1s.222}}$ & $(\Z{2})^3$& $(x,y) \mapsto (-x,-y)$, $(x,y) \mapsto (x^{-1},y^{-1})\hts$, \\
 & & and $(x,y) \mapsto (y,x)$ \\
\hline
\end{tabular}\end{flushleft}
\item
{\bf A special group that does not contain a twisting de Jonqui\`eres involution:}
\begin{flushleft}\begin{tabular}{|lll|}
\hline
{\it name} & {\it structure} & {\it generators} \\
\hline
${\nump{Cs.24}}$ & $\Z{2}\times\Z{4}$ & $(x:y:z)\dasharrow (yz:xy:-yz)$ \\
& & and $(x:y:z)\dasharrow (yz(y-z):xz(y+z):xy(y+z)).$ \\
\hline
\end{tabular}\end{flushleft}
\item
{\bf Groups containing exactly one twisting de Jonqui\`eres involution, that fixes a hyperelliptic curve of genus $k-1$, $G'\cong \Z{2}$:}\index{Curves!hyperelliptic!birational maps that fix a hyperelliptic curve}
\begin{itemize}
\item[$\nump{C.2,\h{0.5}2n}$]
 $k$ is divisible by $2n$,  $n\geq 1$, ${G \cong \Z{2}\times\Z{2n}}$ is generated by 
 \begin{center}$\begin{array}{rcccl}
(x_1:x_2)\times (y_1:y_2)\h{3}& \dasharrow& \h{4}(x_1:x_2) &\h{3}\times\h{3} & (y_2 \prod_{i=1}^k(x_1-b_i x_2):y_1 \prod_{i=1}^k(x_1-a_i x_2)),\\
(x_1:x_2)\times (y_1:y_2)\h{3}& \mapsto& \h{4}(x_1:\zeta_{2n} x_2)&\h{3}\times\h{3} & (y_1:y_2),\end{array}$\end{center}
where $a_1,...,a_k,b_1,...,b_k \in \K^{*}$ are all distinct, and the sets $\{a_1,...,a_k\}$ and $\{b_1,...,b_k\}$ are both invariant by multiplication by $\zeta_{2n}$.
\item[$\nump{C.2,\h{0.5}22}$]
 $k\geq 4$ is divisible by $4$, ${G \cong (\Z{2})^3}$ is generated by
\begin{center}
$\begin{array}{rcccl}
(x_1:x_2)\times (y_1:y_2)& \dasharrow &(x_1:x_2) &\times&(y_2 \prod_{i=1}^{k/4}P(b_i):y_1 \prod_{i=1}^{k/4}P(a_i)),\\
(x_1:x_2)\times (y_1:y_2)& \mapsto &(x_1:- x_2)&\times& (y_1:y_2),\\
(x_1:x_2)\times (y_1:y_2)& \mapsto &(x_2: x_2)&\times& (y_1:y_2),\end{array}$\end{center}
where $P(a)=(x_1-ax_2)(x_1+ax_2)(x_1-a^{-1}x_2)(x_1+a^{-1}x_2) \in \K[x_1,x_2]$ and $a_1,...,a_{k/4},b_1,...,b_{k/4} \in \K\backslash\{0,\pm 1\}$ are all distinct.
\item[$\nump{C.24}$]
$k\geq 2$ is divisible by $2$, $G \cong (\Z{2}\times\Z{4})$ is generated by $\alpha,\rho$, elements of order respectively $4$ and $2$:
\begin{center}
$\begin{array}{rclcccc}
 \alpha:(x,y)& \dasharrow &(\h{3}&-x&\h{3},\h{3}&\nu(x)\cdot \frac{y-\nu(-x)}{y-\nu(x)}&\h{4}),\\
 \alpha^2:(x,y)& \dasharrow &(\h{3}&x&\h{3},\h{3}&\frac{\nu(x)\nu(-x)}{y}&\h{4}),\\
\rho:(x,y)&\mapsto &(\h{3}&x^{-1}&\h{3},\h{3}&y&\h{4}),\end{array}$\end{center}
where $\nu(x) \in \K(x)^{*}$ is invariant by the action of $x \mapsto x^{-1}$, $\nu(x)\not=\nu(-x)$.
\end{itemize}

\item
{\bf Groups containing at least two twisting de Jonqui\`eres involutions, and at least one of them fixes a curve of positive genus, $G'\cong (\Z{2})^2$:}\index{Curves!of positive genus!birational maps that fix a curve of positive genus}
\begin{itemize}
\item[$\nump{C.22}$]
 ${G \cong (\Z{2})^2}$ generated by two twisting de Jonqui\`eres involutions, given by
\begin{center}
$\begin{array}{rclcccc}
(x,y)& \h{3}\dasharrow\h{3} &(&\h{3}x&\h{3},\h{3}&g(x)/y&\h{4}),\\
(x,y)& \h{3}\dasharrow\h{3} &(&\h{3}x&\h{3},\h{3}&(h(x)y-g(x))/(y-h(x))&\h{4}),\end{array}$\end{center}
for some $g(x),h(x) \in \K(x) \backslash \{0\}$, and where each of the two involutions fixes a hyperelliptic curve of positive genus.\index{Curves!hyperelliptic!birational maps that fix a hyperelliptic curve}\index{Curves!of positive genus!birational maps that fix a curve of positive genus}
\item[$\nump{C.22,\h{0.5}n}$]
 ${G \cong (\Z{2})^2\times \Z{n}}$, $n>1$, generated by
\begin{center}
$\begin{array}{rclcccc}
(x,y)& \h{3}\dasharrow\h{3} &(&\h{3}x&\h{3},\h{3}&g(x)/y&\h{4}),\\
(x,y)& \h{3}\dasharrow\h{3} &(&\h{3}x&\h{3},\h{3}&(h(x)y-g(x))/(y-h(x))&\h{4}),\\
(x,y)& \h{3}\dasharrow\h{3} &(&\h{3}\zeta_n x&\h{3},\h{3}&y&\h{4}),\end{array}$\end{center}
for some $g(x),h(x) \in \K(x) \backslash \{0\}$ invariant by $\zeta_n$. The genus of the hyperelliptic curve fixed by any involution of $G'$ is equal to $-1 \pmod{n}$.\index{Curves!hyperelliptic!birational maps that fix a hyperelliptic curve}
\item[$\nump{C.22,\h{0.5}22}$]
 ${G \cong (\Z{2})^4}$ generated by
\begin{center}
$\begin{array}{rclcccc}
(x,y)& \h{3}\dasharrow\h{3} &(&\h{3}x&\h{3},\h{3}&g(x)/y&\h{4}),\\
(x,y)& \h{3}\dasharrow\h{3} &(&\h{3}x&\h{3},\h{3}&(h(x)y-g(x))/(y-h(x))&\h{4}),\\
(x,y)& \h{3}\dasharrow\h{3} &(&\h{3}\pm x^{\pm 1}&\h{3},\h{3}&y&\h{4}),\end{array}$\end{center}
for some $g(x),h(x) \in \K(x) \backslash \{0\}$ invariant by $x \mapsto \pm x^{\pm 1}$, the genus of the hyperelliptic curve fixed by any involution of $G'$ is equal to $3 \pmod{4}$. (These groups are already described in \cite{bib:Be2}.)\index{Curves!hyperelliptic!birational maps that fix a hyperelliptic curve}
\item[$\nump{C.2n2}$]
 ${G \cong (\Z{2n})^2\times \Z{2}}$ is generated by the elements $\alpha$ and $\beta$ that have order respectively $2n$ and $2$:
\begin{center}
$\begin{array}{rrclcccc}
\alpha:&(x,y)& \h{3}\dasharrow\h{3} &(&\h{3}\zeta_n x&\h{3},\h{3}&(a(x)y+b(x))/(c(x)y+d(x))&\h{4}),\\
\alpha^n:&(x,y)& \h{3}\dasharrow\h{3} &(&\h{3}x&\h{3},\h{3}&g(x)/y&\h{4}),\\
\beta:& (x,y)& \h{3}\dasharrow\h{3} &(&\h{3}x&\h{3},\h{3}&(h(x)y-g(x))/(y-h(x))&\h{4}),\end{array}$\end{center}
for some $a,b,c,d,g,h \in \K(x) \backslash \{0\}$ such that $\alpha$ and $\beta$ commute.
\item[$\nump{C.422}$]
 ${G \cong (\Z{4})\times (\Z{2})^2}$ is generated by the elements $\alpha$, $\beta$ and $\gamma$, that have order respectively $4$, $2$ and $2$:
\begin{center}
$\begin{array}{rrclcccc}
\alpha:&(x,y)& \h{3}\dasharrow\h{3} &(&\h{3}x^{-1}&\h{3},\h{3}&(a(x)y+b(x))/(c(x)y+d(x))&\h{4}),\\
\alpha^2:&(x,y)& \h{3}\dasharrow\h{3} &(&\h{3}x&\h{3},\h{3}&g(x)/y&\h{4}),\\
\beta:& (x,y)& \h{3}\dasharrow\h{3} &(&\h{3}x&\h{3},\h{3}&(h(x)y-g(x))/(y-h(x))&\h{4}),\\
\gamma:& (x,y)& \h{3}\dasharrow\h{3} &(&\h{3}-x&\h{3},\h{3}&y&\h{4}),\end{array}$\end{center}
for some $a,b,c,d,g,h \in \K(x) \backslash \{0\}$ such that $\alpha$ and $\beta$ commute,  and $g,h$ invariant by $x \mapsto -x$.
 \item[$\nump{C.44}$]
 ${G \cong (\Z{4})^2}$ is generated by the elements $\alpha$, $\beta$ of order $4$:
\begin{center}
$\begin{array}{rrclcccc}
\alpha:&(x,y)& \h{3}\dasharrow\h{3} &(&\h{3}-x&\h{3},\h{3}&(a(x)y+b(x))/(c(x)y+d(x))&\h{4}),\\
\alpha^2:&(x,y)& \h{3}\dasharrow\h{3} &(&\h{3}x&\h{3},\h{3}&g(x)/y&\h{4}),\\
\beta:& (x,y)& \h{3}\dasharrow\h{3} &(&\h{3}x^{-1}&\h{3},\h{3}&(p(x)y+q(x))/(r(x)y+s(x))&\h{4}),\\
\beta^2:& (x,y)& \h{3}\dasharrow\h{3} &(&\h{3}x&\h{3},\h{3}&(h(x)y-g(x))/(y-h(x))&\h{4}),\end{array}$\end{center}
for some $a,b,c,d,g,h,p,q,r,s \in \K(x) \backslash \{0\}$ such that $\alpha$ and $\beta$ commute.
\end{itemize}
\end{itemize}
\begin{flushleft}{\Large \sc Groups that do not leave invariant a pencil of rational curves:}\end{flushleft}
\begin{itemize}
\item
{\bf Automorphisms of $\bf \Pn$}
\begin{itemize}
\item[$\nump{0.V9}$]
$G=V_9\cong (\Z{3})^2$, generated by $(x:y:z) \mapsto (x:\omega y:\omega^2 z)$ and $(x:y:z) \mapsto (y:z:x)$.
\end{itemize}
\item
{\bf Automorphisms of Del Pezzo surfaces of degree $4$}\\
The surface $S\subset \mathbb{P}^4$ is the blow-up of $(1:0:0)$, $(0:1:0)$, $(0:0:1)$, $(1:1:1)$, $(a:b:c) \in \Pn$, of equations
\begin{center}$\begin{array}{cccccc}cx_1^2&&-ax_3^2&-(a-c)x_4^2&-ac(a-c)x_5^2&=0,\\
&cx_2^2&-bx_3^2&+(c-b)x_4^2&-bc(c-b)x_5^2&=0.\end{array}$\end{center}
\begin{itemize}
\item[$\nump{4.222}$]
$G\cong (\Z{2})^3$, given by $G=\{(x_1:x_2:x_3:x_4:x_5) \mapsto (\pm x_1:\pm x_2:\pm x_3:x_4:x_5)\}$.
\item[$\nump{4.2222}$]
$G\cong (\Z{2})^4$, given by $G=\{(x_1:x_2:x_3:x_4:x_5) \mapsto (\pm x_1:\pm x_2:\pm x_3:\pm x_4:x_5)\}$.
\item[$\nump{4.42}$]
$(a:b:c)=(1:\xi:1+\xi)$, for any $\xi\in \K\backslash \{0,\pm 1\}$\\
$G\cong \Z{4} \times \Z{2}$,  generated by $(x_1:x_2:x_3:x_4:x_5) \mapsto (-x_2:x_1:x_4:x_3: - x_5)$
 and $(x_1:x_2:x_3:x_4:x_5) \mapsto (x_1:x_2:x_3:x_4:-x_5)$.
\end{itemize}

\item
{\bf Automorphisms of cubic surfaces \vspace{0.1cm}}\\
\begin{tabular}{|llll|}
\hline
{\it name} & {\it structure} & {\it generators} & {\it equation of}\\
{\it of $\mathit (G,S)$} & {\it of $\mathit G$} & {\it of $\mathit G$} & $\mathit F$\\
\hline 
${\nump{3.33.1}}$ & $(\Z{3})^2$& $\DiaG{\omega}{1}{1}{1}$, $\DiaG{1}{1}{1}{\omega}$& $w^3+x^3+y^3+z^3\hts$ \\ 
${\nump{3.33.2}}$ & $(\Z{3})^2$& $\DiaG{\omega}{1}{1}{1}$, $\DiaG{1}{1}{\omega}{\omega^2}$& $w^3+x^3+y^3+z^3+\lambda xyz\hts$ \\
${\nump{3.36}}$ & $\Z{3}\times\Z{6}$& $\DiaG{\omega}{1}{1}{1}$, $\DiaG{1}{1}{-1}{\omega}$& $w^3+x^3+xy^2+z^3\hts$ \\ 
${\nump{3.333}}$ & $(\Z{3}) ^3$&$\DiaG{\omega}{1}{1}{1}$, $\DiaG{1}{\omega}{1}{1}$& $w^3+x^3+y^3+z^3\hts$\\
& &\mbox{and } $\DiaG{1}{1}{\omega}{1}$ & \\
\hline
\end{tabular}
\item
{\bf Automorphisms of Del Pezzo surfaces of degree $2$}\\
The surface is given in the weighted projective space $\mathbb{P}(2,1,1,1)$ by the equation $w^2=F(x,y,z)$, where $F$ is a non-singular form of degree $4$.
\begin{itemize}
\item
groups containing the Geiser involution $\sigma$:\\
$G=<\sigma>\times H$, where $H \subset
\PGLn{3}$ is an abelian group of automorphisms of the quartic curve $\Gamma \subset \Pn$ defined by the equation $F(x,y,z)=0$:\vspace{0.1cm}\\
\begin{tabular}{|llll|}
\hline
{\it name} & {\it structure} & {\it generators} & {\it equation of}\\
{\it of $\mathit (G,S)$} & {\it of $\mathit G$} & {\it of $\mathit G$} & $\mathit F$\\
\hline 
${\nump{2.G2}}$ & $(\Z{2})^2$&$\DiaG{\pm 1}{1}{1}{-1}$& $L_4(x,y)+L_2(x,y)z^2+z^4\hts$ \\ 
${\nump{2.G4.1}}$ & $\Z{2}\times\Z{4}$&$\DiaG{\pm 1}{1}{1}{\im}$& $L_4(x,y)+z^4\hts$ \\
${\nump{2.G4.2}}$ & $\Z{2}\times\Z{4}$&$\DiaG{\pm 1}{1}{-1}{\im}$& $x^3y+y^4+z^4+xyL_1(xy,z^2)\hts$ \\ 
${\nump{2.G6}}$ & $\Z{2}\times\Z{6}$&$\DiaG{\pm 1}{\omega}{1}{-1}$& $x^3y+y^4+z^4+\lambda y^2z^2\hts$ \\
${\nump{2.G8}}$ & $\Z{2}\times\Z{8}$&$\DiaG{\pm 1}{\zeta_8}{-\zeta_8}{1}$& $x^3y+xy^3+z^4\hts$ \\ 
${\nump{2.G12}}$ & $\Z{2}\times\Z{12}$&$\DiaG{\pm 1}{\omega}{1}{\im}$& $x^3y+y^4+z^4\hts$ \\ 
${\nump{2.G22}}$ & $(\Z{2})^3$&$\DiaG{\pm 1}{1}{\pm 1}{\pm 1}$ & $L_2(x^2,y^2,z^2)\hts$ \\ 
${\nump{2.G24}}$ & $(\Z{2})^2\times \Z{4}$&$\DiaG{\pm 1}{1}{\pm 1}{\im}$ & $x^4+y^4+z^4+\lambda x^2y^2\hts$ \\ 
${\nump{2.G44}}$ & $\Z{2}\times (\Z{4})^2$&$\DiaG{\pm 1}{1}{1}{\im}$& $x^4+y^4+z^4\hts$\\
& & and $\DiaG{1}{1}{\im}{1}$& \\
\hline
\end{tabular}
\item
groups not containing the Geiser involution:\vspace{0.1cm}\\
\begin{tabular}{|llll|}
\hline
{\it name} & {\it structure} & {\it generators} & {\it equation of}\\
{\it of $\mathit (G,S)$} & {\it of $\mathit G$} & {\it of $\mathit G$} & $\mathit F$\\
\hline 
${\nump{2.24.1}}$ & $\Z{2}\times \Z{4}$& $\DiaG{1}{1}{1}{\im}$, $\DiaG{1}{1}{-1}{1}$& 
$x^4+y^4+z^4+\lambda x^2y^2\hts$ \\ 
${\nump{2.24.2}}$ & $\Z{2}\times \Z{4}$& $\DiaG{1}{1}{1}{\im}$, $\DiaG{-1}{1}{-1}{1}$& 
$x^4+y^4+z^4+\lambda x^2y^2\hts$ \\
${\nump{2.44.1}}$ & $(\Z{4})^2$& $\DiaG{1}{1}{1}{\im}$, $\DiaG{1}{1}{\im}{1}$& 
$x^4+y^4+z^4\hts$ \\ 
${\nump{2.44.2}}$ & $(\Z{4})^2$& $\DiaG{1}{1}{1}{\im}$, $\DiaG{-1}{1}{\im}{1}$& 
$x^4+y^4+z^4\hts$\\
\hline
\end{tabular}
\end{itemize}
\item {\bf Automorphisms of Del Pezzo surfaces of degree $1$}\vspace{0.1cm}\\
\begin{tabular}{|lllll|}
\hline
{\it name} & {\it structure} & {\it generators} & {\it equation of}& {\it equation of}\\
{\it of $\mathit (G,S)$} & {\it of $\mathit G$} & {\it of $\mathit G$} & $\mathit F_4$& $\mathit F_6$\\
\hline 
${\nump{1.B2.1}}$ & $(\Z{2})^2$&$\sigma$, $\DiaG{1}{1}{-1}{1}$& $L_2(x^2,y^2)$ & $L_3(x^2,y^2)\hts$\\ 
${\nump{1.\sigma\rho 2.1}}$ & $\Z{6}\times \Z{2}$& $\sigma\rho$, $\DiaG{1}{1}{-1}{1}$& $0$ & $L_3(x^2,y^2)\hts$ \\ 
${\nump{1.\sigma\rho 3}}$ & $\Z{6}\times\Z{3}$& $\sigma\rho$, $\DiaG{1}{1}{\omega}{1}$& $0$ &  $L_2(x^3,y^3)\hts$ \\ 
${\nump{1.\rho 3}}$ & $ (\Z{3})^2$& $\rho$, $\DiaG{1}{1}{\omega}{1}$& $0$ & $L_2(x^3,y^3)\hts$ \\
${\nump{1.B4.1}}$ & $\Z{2}\times\Z{4}$&$\sigma$, $\DiaG{1}{1}{\im}{1}$& $L_1(x^4,y^4)$ &  $x^2L_1'(x^4+y^4)\hts$ \\ 
${\nump{1.B6.1}}$ & $\Z{2}\times\Z{6}$&$\sigma$, $\DiaG{1}{1}{-\omega}{1}$& $\lambda x^4$ & $\mu x^6+y^6\hts$\\
${\nump{1.\sigma\rho 6}}$ & $ (\Z{6})^2$& $\sigma\rho$, $\DiaG{1}{1}{-\omega}{1}$&$0$& $x^6+y^6\hts$ \\ 
${\nump{1.\rho 6}}$ & $\Z{3}\times \Z{6}$& $\rho$, $\DiaG{1}{1}{-\omega}{1}$&$0$& $x^6+y^6\hts$\\ 
${\nump{1.B6.2}}$ & $\Z{2}\times\Z{6}$&$\sigma$, $\DiaG{1}{1}{-\omega}{\omega}$& $\lambda x^2y^2$ & $x^6+y^6\hts$ \\
${\nump{1.B12}}$ & $\Z{2}\times\Z{12}$&$\sigma, \DiaG{\im}{1}{\zeta_{12}}{-1}$& $x^4$ & $y^6\hts$\\
\hline
\end{tabular}\\
where $L_i$ denotes a form of degree $i$, $\lambda,\mu \in \K$ denote parameters, and the equations are such that $S$ is non-singular.
\end{itemize}
\end{flushleft}
\begin{remark}\fauxtitred
We prove that two groups of different types represent different conjugacy classes of finite abelian subgroups of the Cremona group. However, we do not make precise the conjugacy class within the same type, when this depends on parameters.  We did it however for some cases in previous chapters. Note that the action on the sets of fixed points is often sufficient.
\end{remark}
\begin{proof}\fauxtitred\upshape
Let us first prove that $G$ is birationally conjugate to one of the examples of the list above.

We distinguish two cases, depending on the existence or not of an invariant pencil of rational curves.
\begin{itemize}
\item
{\bf \emph{G} preserves a pencil of rational curves:}\\
Up to birational conjugation, $G\subset \Aut(S,\pi)$ acts on a conic bundle $(S,\pi)$.
We use Proposition \refd{Prp:CBdistCases} and get one of the following cases:
\begin{itemize}
\item
Some twisting $(S,\pi)$-de Jonqui\`eres involution belongs to $G$;
\item
$(S,\pi)=(\mathbb{P}^1\times \mathbb{P}^1,\pi_1)$;
\item
the group $G$ is isomorphic to $\Z{2}\times\Z{4}$ and the pair $(G,S)$ is either $\nump{Cs.24}$ or $\nump{P1s.24}$.
\end{itemize}
The second case gives groups of the list above, using Proposition \refd{Prp:AutP1Bir}. The two possibilities appear in the list. In remains to study the first case.

Note that we may suppose that at least one involution of $G$ fixes a curve of positive genus, using Proposition \refd{Prp:HardChgExacrSeq}.
We will use the exact sequence defined in Section \refd{Sec:Exactsequence}:
\begin{equation}
1 \rightarrow G' \rightarrow G \stackrel{\overline{\pi}}{\rightarrow} \overline{\pi}(G) \rightarrow 1.
\tag{\refd{eq:ExactSeqCB}}\end{equation}
Using once again Proposition \refd{Prp:CBdistCases}, the group $G'$ contains one $(S,\pi)$-twisting de Jonqui\`eres involution and is either cyclic of order $2$ or isomorphic to $(\Z{2})^2$, generated by two twisting de Jonqui\`eres involutions. We study the two cases separately:

$G'\cong \Z{2}$: as $G$ is not cyclic, then $\overline{\pi}(G)\not=1$, and three possibilities arise:
\begin{itemize}
\item
$G\cong \Z{2}\times\Z{2n}$ and  $\overline{\pi}(G)\cong \Z{2n}$: $\nump{C.2,\h{0.5}2n}$ (Proposition \refd{Prp:DirectProduct}).
\item
$G\cong (\Z{2})^3$ and  $\overline{\pi}(G)\cong (\Z{2})^2$: $\nump{C.2,\h{0.5}22}$ (Proposition \refd{Prp:DirectProductC2C2C2}).
\item
$G\cong \Z{2}\times \Z{4}$ and  $\overline{\pi}(G)\cong (\Z{2})^2$: $\nump{C.24}$ (Proposition \refd{Prp:DoesnotsplitC2-C2C2}).
\end{itemize}

$G'\cong (\Z{2})^2$: the possible cases are:
\begin{itemize}
\item
$\overline{\pi}(G)=1$: \begin{itemize}
\item
If at least two of the involutions of $G$ fix (each one) a curve of positive genus, we get $\nump{C.22}$, (Proposition \refd{Prp:C2C2PGL2Ct}).
\item
If exactly one of the involutions fixes a curve of positive genus, then $G$ is conjugate to a group $H$ acting on another conic bundle, with $H'\cong \Z{2}$ (Proposition \refd{Prp:HardChgExacrSeq}).
\end{itemize}
\item
$\overline{\pi}(G)\not=1$ and (\refd{eq:ExactSeqCB}) splits:  $\nump{C.22,\h{0.5}n}$,  $\nump{C.22,\h{0.5}22}$, (Proposition \refd{Prp:SplitC2C2});
\item
$\overline{\pi}(G)\not=1$ and (\refd{eq:ExactSeqCB}) does not split:
 $\nump{C.2n2}$,  $\nump{C.422}$,  $\nump{C.44}$, (Proposition \refd{Prp:NotSplitC2C2}).
 \end{itemize}
\item
{\bf \emph{G} does not preserve a pencil of rational curves:}\\
By Proposition \refd{Prp:ConjBetweenRankPic}, the group is conjugate to one of the groups stated above.
\end{itemize}

It remains to show that the different cases of the proposition represent distinct conjugacy classes of the Cremona group. Firstly, Proposition \refd{Prp:ConjBetweenRankPic} shows that neither $\nump{0.V9}$ nor any of the groups that appear after it in the list  leaves invariant a pencil of rational curves. As the groups that appear before $\nump{0.V9}$ keep invariant a pencil of rational curves, this establishes a distinction between the birational conjugacy classes.
(The existence of a pencil of rational curves is an invariant of birational conjugation, see Section \refd{Sec:PencilRatCurves}.)

Let us look at the groups that preserve a pencil of rational curves:

\item
{\bf Automorphisms of conic bundles}:\\ 
We enumerate the cases given in the list and prove that each one is not birationally conjugate to the previous cases.
\begin{itemize}
\item
We proved in Proposition \refd{Prp:AutP1Bir} that the subgroups of $\Aut(\mathbb{P}^1\times\mathbb{P}^1)$ given in the list are not birationally conjugate.
\item
Lemma \refd{Lem:Cs24notconjugate} shows that $\nump{Cs.24}$ is not birationally conjugate to a subgroup of $\Aut(\mathbb{P}^1\times\mathbb{P}^1)$.
\item
 Each group that follows contains an involution that fixes a curve of positive genus. So none of these groups is birationally conjugate to $\nump{Cs.24}$ or a subgroup of $\Aut(\mathbb{P}^1\times\mathbb{P}^1)$.
 
 It remains to show that the $9$ cases  $\nump{C.2,\h{0.5}2n}$, ..., $\nump{C.44}$ are distinct. As these are classified by the structure of $G'$, $G$ (and $\overline{\pi}(G)$, which depends on the other two), we need only to show that we cannot conjugate one of the groups with $G'\cong (\Z{2})^2$ to a group $H$ with $H'\cong \Z{2}$.
 
 Proposition \refd{Prp:HardChgExacrSeq} shows that such a phenomenon is possible only if  $G=G'\cong (\Z{2})^2$ and $G$ contains exactly one involution that fixes a curve of positive genus. But this case does not appear in the list (as it is conjugate to a group $H$ with $H'\cong \Z{2}$).
 \end{itemize}

 Secondly,  using once again Proposition \refd{Prp:ConjBetweenRankPic}, the group $\nump{0.V9}$ and the four types of Del Pezzo surfaces of degree $1,2,3,4$ are distinct. It remains to show that within one such type, the cases are birationally distinct.
\begin{itemize}
\item
{\bf Del Pezzo surfaces of degree $\bf 4$}: the groups of different cases are not isomorphic.
\item
{\bf Del Pezzo surfaces of degree $\bf 3$}: only $\nump{3.33.1}$ and $\nump{3.33.2}$ contain isomorphic groups (isomorphic to $(\Z{3})^2$). The first group contains $4$ elements of order $3$ that fix an elliptic curve ($\DiaG{\omega^{\pm 1}}{1}{1}{1}$, $\DiaG{1}{1}{1}{\omega^{\pm 1}}$) and the second one contains only two such elements.
\item
{\bf Del Pezzo surfaces of degree $\bf 2$}: the only pairs that contain isomorphic groups are distinguished by the set of fixed points of their elements.\\
Firstly, we study pairs with groups isomorphic to $\Z{2}\times\Z{4}$:
\begin{itemize}
\item
$\nump{2.G4.1}$: the group contains the Geiser involution, that fixes a curve of genus $2$, and elements of order $4$ of type $\DiaG{1}{1}{1}{\im}$, that fix an elliptic curve. 
\item
$\nump{2.G4.2}$: contains the Geiser involution, and every element of order $4$ fixes a finite number of points.
\item
$\nump{2.24.1}$: every involution of the group fixes an elliptic curve. \index{Curves!elliptic!birational maps that fix an elliptic curve|idi}
\item
$\nump{2.24.2}$: one involution ($\DiaG{1}{1}{1}{-1}$) fixes an elliptic curve, the other two fix only a finite number of points.
\end{itemize}
Secondly, we study pairs with groups isomorphic to $(\Z{4})^2$:
\begin{itemize}
\item
$\nump{2.44.1}$: the group contains four elements of order $4$ that fix an elliptic curve ($\DiaG{1}{1}{1}{\pm \im}$, $\DiaG{1}{1}{\pm \im}{1}$).
\item
$\nump{2.44.2}$: contains only two elements of order $4$ that fix an elliptic curve. The others fix only a finite number of points.
\end{itemize}
\item
{\bf Del Pezzo surfaces of degree $\bf 1$}: the only isomorphism class that appears several times is  $\Z{2}\times \Z{6}$. The pairs are distinguished by the set of fixed points of elements of the groups:
\begin{itemize}
\item
$\nump{1.\sigma\rho 2.1}$: the group contains $\rho=\DiaG{1}{1}{1}{\omega}$, that fixes a curve of genus $2$.
\item
$\nump{1.B6.1}$: contains $\DiaG{1}{1}{-\omega}{1}$, an element of order $6$ that fixes an elliptic curve. No element of order $3$ fixes a curve of genus $2$.
\item
$\nump{1.B6.2}$: every element of order $3$ or $6$ fixes only a finite number of points.
\proofend\end{itemize}
\end{itemize}
\end{proof}
\section{Other main results}
\pagetitless{The classification and other main results}{Other main results}
Using our classification, we are able to prove some new results, already stated in Section \refd{Sec:Results}.

\begin{flushleft}{\bf THEOREM 1: \sc Non-linear birational maps of large order}\\ \itshape
\begin{itemize}\itshape
\item
For any integer $n\geq 1$, there are infinitely many conjugacy classes of birational maps of the plane of order $2n$, that are non-conjugate to a linear automorphism. 
\item
If $n>15$, a birational map of order $2n$ is an $n$-th root of a de Jonqui\`eres involution and preserves a pencil of rational curves.
\item
If a birational map is of finite odd order and is not conjugate to a linear automorphism of the plane, then its order is $3,5,9$ or $15$. In particular,  any birational map of the plane of odd order $>15$ is conjugate to a linear automorphism of the plane.
\end{itemize}\end{flushleft}
\begin{proof}\fauxtitred\upshape
\begin{itemize}
\item
The first assertion follows from the existence of $n$-th roots of de Jonqui\`eres involutions for any integer $n$. 
\begin{itemize}
\item
The case of odd $n$ is given by elements of the form 
\begin{center}
$\begin{array}{rcccl}
(x_1: x_2)\times (y_1:y_2) \dasharrow (x_1:\zeta_n x_2) \times(y_2 \prod_{i=1}^k(x_1-b_i x_2):y_1 \prod_{i=1}^k(x_1-a_i x_2))\\\end{array}$\end{center}
where $a_1,...,a_k,b_1,...,b_k \in \K^{*}$ are all distinct, $\zeta_n=e^{2\im \pi/n}$ and the sets $\{a_1,...,a_k\}$ and $\{b_1,...,b_k\}$ are both invariant by multiplication by $\zeta_n$.
\item
The case $n$ even is given by Proposition \refd{Prp:ExistRoot}. 
\end{itemize}
In both cases, if $k>1$, the root is not conjugate to a linear automorphism of the plane, as the de Jonqui\`eres involution fixes a curve of positive genus.
\item
The second and third assertions follow directly from the classification of finite cyclic groups (Theorem  \refThmCyclicGroup). Note that a linear automorphism of even order is also a root of a de Jonqui\`eres involution.\proofend
\end{itemize}
\end{proof}
\begin{remark}\fauxtitred
Note that Proposition \refd{Prp:ExistRoot} ensures the existence of $n$-th roots of de Jonqui\`eres involutions, for $n$ even, but not a simple expression as in the odd case.\end{remark}

\bigskip
\bigskip

\noindent{\bf THEOREM 2: \sc Roots of linear automorphisms}\\ \itshape
Any birational map which is a root of a non-trivial linear automorphism of finite order of the plane  is conjugate to a linear automorphism of the plane. \upshape
\begin{proof}\fauxtitred\upshape
This follows directly from the classification of finite cyclic groups (Theorem  \refThmCyclicGroup). Indeed, we see that any cyclic group which is not conjugate to a subgroup of $\Aut(\Pn)$ contains an element that fixes a curve of positive genus. \proofend
\end{proof}

\bigskip
\bigskip

\begin{flushleft}{\bf THEOREM 3: \sc Groups which fix a curve of positive genus}\\ \itshape
Let $G$ be a finite abelian group which fixes some curve of positive genus. Then $G$ is cyclic, of order $2$, $3$, $4$, $5$ or $6$, and all these cases occur. If the curve has genus $>1$, the order is $2$ or $3$.\index{Curves!elliptic!birational maps that fix an elliptic curve|idi}\index{Curves!of positive genus!birational maps that fix a curve of positive genus}\end{flushleft}
\begin{proof}\fauxtitred\upshape
Using  the classification of finite abelian non-cyclic groups (Theorem 
\refThmNonCyclicGroup), we see that none of them fixes a curve of positive genus. Hence, the group $G$ must be cyclic of order $n$.

Using the classification of finite cyclic groups (Theorem  \refThmCyclicGroup), we get
\begin{itemize}
\item
$n=2$: de Jonqui\`eres, Geiser and Bertini involutions.\index{Involutions!Bertini|idi}\index{Involutions!Geiser|idi}
\index{Del Pezzo surfaces!automorphisms}
\item
$n=3:$ Automorphisms of cubic surfaces, that fix an elliptic curve. Automorphisms of \\
$\hphantom{n=3:}$ Del Pezzo surfaces of degree $1$, that fix curves of genus $2$.
\item
$n=4$: Automorphisms of Del Pezzo surfaces of degree $2$ that fix elliptic curves.
\item
$n=5$: Automorphisms of Del Pezzo surfaces of degree $1$ that fix elliptic curves.
\item
$n=6$: Automorphisms of Del Pezzo surfaces of degree $1$ that fix elliptic curves.
 \proofend
\end{itemize}
\end{proof}

\bigskip
\bigskip

\begin{flushleft}{\bf THEOREM 4: \sc Cyclic groups whose non-trivial elements do not fix a curve of positive genus}\\ \itshape
\index{Curves!of positive genus!birational maps that do not fix a curve of positive genus}
Let $G$ be a finite cyclic subgroup of the Cremona group. The following conditions are equivalent:
\begin{itemize}
\item
If $g \in G$, $g\not=1$, then $g$ does not fix a curve of positive genus.
\item
$G$ is birationally conjugate to a subgroup of $\Aut(\Pn)$.
\item
$G$ is birationally conjugate to a subgroup of $\Aut(\mathbb{P}^1\times\mathbb{P}^1)$.
\end{itemize}
\end{flushleft}
\begin{proof}\fauxtitred\upshape
Let $n$ be the order of the group.
\begin{itemize}
\item
The second assertion is equivalent to the fact that $G$ is birationally conjugate to a diagonal cyclic group generated by $(x:y:z)\mapsto (\zeta_n x:y:z)$ (see Proposition \refd{Prp:PGL3Cr}).
\item
The third assertion is equivalent to the fact that $G$ is birationally conjugate to a diagonal cyclic group generated by $(x,y) \mapsto (\zeta_n x,y)$ (see Proposition \refd{Prp:AutP1Bir}).
\end{itemize}
So the second and third assertions are equivalent and imply the first one. It remains to prove that the first one implies one of the others.  Since $G$ is cyclic, it is birationally conjugate to a subgroup of the list of Theorem \refThmCyclicGroup. The only case that contains no element which fixes a curve of positive genus is a group of diagonal automorphisms of\hspace{0.2 cm}$\Pn$.\proofend
\end{proof}
\bigskip
\bigskip

\begin{flushleft}{\bf THEOREM 5: \sc Abelian groups whose non-trivial elements do not fix a curve of positive genus}\\ \itshape
\index{Curves!of positive genus!birational maps that do not fix a curve of positive genus}
Let $G$ be a finite abelian subgroup of the Cremona group. The following conditions are equivalent:
\begin{itemize}
\item
If $g \in G$, $g\not=1$, then $g$ does not fix a curve of positive genus.
\item
$G$ is birationally conjugate to a subgroup of $\Aut(\Pn)$, or to a subgroup of $\Aut(\mathbb{P}^1\times\mathbb{P}^1)$ or to the group isomorphic to $ \Z{2}\times\Z{4}$, generated by the two elements $\begin{array}{lll}(x:y:z)&\dasharrow&(yz:xy:-xz),\\
(x:y:z)&\dasharrow &( yz(y-z):xz(y+z):xy(y+z)).\end{array}$
\end{itemize}
Moreover, this last group is conjugate neither to a subgroup of $\Aut(\Pn)$, nor to a subgroup of $\Aut(\mathbb{P}^1\times\mathbb{P}^1)$.\end{flushleft}
\begin{proof}\fauxtitred\upshape
Note first of all that the second assertion implies the first. Indeed, no non-trivial element of finite order of $\Aut(\Pn)$, or $\Aut(\mathbb{P}^1\times\mathbb{P}^1)$ fixes a curve of positive genus (see Theorem \refThmCyclicGenusCurve), and the group isomorphic to $\Z{2}\times\Z{4}$ of the statement is the group of $\nump{Cs.24}$; none of its elements fixes a curve of positive genus (see Lemma \refd{Lem:Cs24Properties}).

Let us now prove that the first assertion implies the second. 
\begin{itemize}
\item
 If $G$ is cyclic, this follows from Theorem \refThmCyclicGenusCurve.
 \item
 If $G$ is non-cyclic, it is birationally conjugate to a subgroup of the list of Theorem \refThmNonCyclicGroup. The cases that contain no element which fixes a curve of positive genus are $\nump{0.V9}$, which is a subgroup of $\Aut(\Pn)$, $\nump{Cs.24}$, which is the group mentioned in this theorem, and subgroups of $\Aut(\mathbb{P}^1\times\mathbb{P}^1)$.
 \end{itemize}
 The last assertion was proved in Lemma \refd{Lem:Cs24Properties}.\proofend
\end{proof}

\bigskip

\begin{flushleft}{\bf THEOREM 6: \sc Isomorphy classes of finite abelian groups}\\ \itshape
The \emph{isomorphism} classes of finite abelian subgroups of the Cremona group are the following:
\begin{itemize}
\item
$\Z{m}\times\Z{n}$, for any integers $m,n \geq 1;$
\item
$\Z{2n}\times(\Z{2})^2$, for any integer $n\geq 1;$
\item
$(\Z{4})^2\times\Z{2};$
\item
$(\Z{3})^3;$
\item
$(\Z{2})^4$.
\end{itemize}
\end{flushleft} \upshape
\begin{proof}\fauxtitred\upshape
This follows directly from Theorems  \refThmCyclicGroup\ and \refThmNonCyclicGroup.
\proofend
\end{proof}
\chapter{Appendix}
We give in this chapter some auxiliary results.
\section{Appendix 1 - Automorphisms of curves}
\pagetitless{Appendix}{Appendix 1 - Automorphisms of curves}
We give some properties of the group of automorphisms of curves that appear in this paper. Most of these results are classical.
\subsection*{11.1.1\hspace{0.4 cm}Automorphisms of $\mathbb{P}^1$}

We find it useful to write an element of $\Aut(\mathbb{P}^1)=\PGL(2,\K)$ in either of the two forms
\begin{center}
$(x_1:x_2) \mapsto (a x_1 +bx_2:c x_1 + d x_2)$\hspace{0.5 cm} or \hspace{0.5 cm}
$x \mapsto \frac{a x +b}{c x +d}$, for some $\left( \begin{array}{cc} a & b \\ c & d\end{array}\right) \in \PGL(2,\K)$.
\end{center}
Let us recall some classical easy results on $\Aut(\mathbb{P}^1)$:
\begin{Lem}\fauxtitre
\label{Lem:AutP1}
Let $G \subset \Aut(\mathbb{P}^1)$ be a finite abelian group. Then:
\begin{itemize}
\item
$G$ is either cyclic or isomorphic to $(\Z{2})^2$.
\item
Up to conjugation, either $G$ is generated by $(x_1:x_2) \mapsto (\alpha x_1:x_2)$, for some $\alpha \in \K^{*}$ (cyclic case), or $G$ is the group $x\mapsto \pm x^{\pm 1}$.
\item
$G\not=\{1\}$ is cyclic if and only if it fixes exactly two points of $\mathbb{P}^1$. If it does so, all elements of $G\backslash \{1\}$ fix the same two points.
\end{itemize}
\end{Lem}
\begin{proof}\fauxtitred\upshape
All these assertions are well known and can be proved directly using the description of automorphisms of $\mathbb{P}^1$.\proofend\end{proof}

\bigskip
\begin{remark}\fauxtitred
Note that this lemma is also true for every algebraically closed field; by scalar extension to the algebraic closure, we see that assertion $1$ is still true for any field {\upshape (}like $\K(x)${\upshape )}, although the others are not.
 For example, the involution {\tiny $\left(\begin{array}{cc}0 & g\\ 1 & 0 \end{array}\right)$} is diagonalizable in $\PGL(2,\K(x))$ if and only if $g$ is a square in $\K(x)$ (see Proposition \refd{Prp:InvDeJPgl2} ).
\end{remark}

\subsection*{11.1.2\hspace{0.4 cm}Automorphisms of elliptic curves}
We recall that an elliptic curve is an irreducible curve of genus $1$. It has three standard classical forms, given in the following lemma:\index{Curves!elliptic!definition|idb}
\begin{Lem}\fauxtitred
\begin{itemize}
\item
Any elliptic curve is isomorphic to the double covering of $\mathbb{P}^1$ ramified over $4$ points (it is a hyperelliptic curve). This gives the \defn{Weierstrass form} \begin{center}$W_{\lambda}=\{y_1^2=x_1x_2(x_1-x_2)(x_1-\lambda x_2)\}\subset \mathbb{P}(1,1,2)$.\end{center}
The linear equivalence of the set of ramification points determines the isomorphism class of the curve.
\item
Any elliptic curve is isomorphic to some smooth cubic curve in $\Pn$. This gives the \defn{Hesse form}
\begin{equation}\label{eqn:HesseForm}H_{\lambda}=\{(x:y:z) \in \Pn \ |\ x^3+y^3+z^3+\lambda xyz=0\}. \end{equation}
The parameter $\lambda$ is defined up to an action of the Hesse group {\upshape (}generated by permutations of the coordinates, multiplication by a $3$-rd root of unity, and the automorphism\\ $(x:y:z) \mapsto (x+y+z:x+\omega y+\omega^2 z:x+\omega^2 y +\omega z)${\upshape )}, and determines the curve up to isomorphism.
\item
Any elliptic curve is isomorphic to $\K /\Lambda$, for some lattice $\Lambda \subset \K$. The lattice, up to homothety, determine the isomorphism classes of the curve.
\end{itemize}\proofend
\end{Lem}
Note that the last form gives a natural group structure on the elliptic curve. \index{Curves!elliptic!group structure|idb}This is isomorphic to the structure given by the classical cubic law using one reference point and the geometry of a smooth cubic.
We now give the group of automorphisms of an elliptic curve:\index{Curves!elliptic!automorphism group|idb}
\begin{Lem}\fauxtitre
Let $\Gamma$ be an elliptic curve with a group structure. We denote by $H_{\Gamma}$ the group of group automorphisms of $\Gamma$.
 Then, 
$\Aut(\Gamma)=\Gamma \rtimes H_{\Gamma}$.

In particular, writing $\Gamma=\K /_{\mathbb{Z} + \lambda\mathbb{Z}}$ we have 
\begin{itemize}
\item
$\Aut(\K /_{\mathbb{Z} + \im\mathbb{Z}})=\K /_{\mathbb{Z} + \im\mathbb{Z}} \rtimes \Z{4}$, where $\Z{4}$ is generated by $x\mapsto \im x$.
\item
$\Aut(\K /_{\mathbb{Z} + \omega\mathbb{Z}})=\K /_{\mathbb{Z} + \omega\mathbb{Z}} \rtimes \Z{6}$, where $\Z{6}$ is generated by $x \mapsto -\omega x$, $\omega=e^{2\im \pi/3}$.
\item
$\Aut(\K /_{\mathbb{Z} + \lambda\mathbb{Z}})=\K /_{\mathbb{Z} + \lambda\mathbb{Z}} \rtimes \Z{2}$, where $\Z{2}$ is generated by $x\mapsto -x$ if the curve is not isomorphic to one of the two given above.
\end{itemize}
\end{Lem}

\subsection*{11.1.3\hspace{0.4 cm}Automorphisms of hyperelliptic curves}
\label{Sec:HyperellipticCurves}
Recall that an irrreducible curve is \defn{hyperelliptic} if it is isomorphic to some double covering of $\mathbb{P}^1$, ramified over $2k$ points.\index{Curves!hyperelliptic!definition|idi}

Note that the isomorphism class of a hyperelliptic curve depends only on the linear equivalence of the set of ramification points.
\index{Curves!hyperelliptic!isomorphism class|idi}

We describe in the following proposition the group of automorphisms of the curve which are compatible with the double covering (i.e.\ which permute the fibres). If $k>2$ it is the entire group of automorphisms of the curve; if $k=1,2$, we refer respectively to Sections 11.1.1 and 11.1.2 for a more precise description of the entire group of automorphisms of the curve.

\begin{Prp}\titreProp{automorphisms of hyperelliptic curves}\\
\index{Curves!hyperelliptic!automorphism group|idb}
\label{Prp:Authyperellic}
Let $\Gamma$ be a hyperelliptic curve of genus $k-1$, with a double covering $g^1_2:\Gamma\rightarrow \mathbb{P}^1$ ramified over a set $I$ of $2k$ points of $\mathbb{P}^1$.

Let $\Aut(\Gamma,g^1_2)$ be the group of automorphisms of $\Gamma$ which leave the covering invariant.
\begin{enumerate}
\item[\upshape 1.]
There is an exact sequence
\begin{center}$\begin{array}{cccccccc}
1& \rightarrow& \Z{2}=<\sigma> &\rightarrow &\Aut(\Gamma,g^1_2)& \stackrel{\rho}{\rightarrow} &\PGL(2,\K)_I\rightarrow &1,\end{array}$\end{center}
where $\sigma$ is the automorphism that exchanges the two points of the general fibres of $g^1_2$ and $\PGL(2,\K)_I$ is the subgroup of automorphisms of $\mathbb{P}^1$ which leave invariant the set $I \subset \mathbb{P}^1$.
\item[\upshape 2.]
The involution $\sigma$ commutes with all the elements of $\Aut(\Gamma,g^1_2)$.
\item[\upshape 3.]
If $k$ is even, the elements $\alpha$ and $\sigma\alpha$ are not conjugate in $\Aut(\Gamma,g^1_2)$, for any element $\alpha \in \Aut(\Gamma,g^1_2)$.\end{enumerate}
\end{Prp}
\begin{proof}\fauxtitred\upshape
The curve $\Gamma$ can be defined in the weighted projective plane $\mathbb{P}(k,1,1)$ by the equation $w^2=g(x_1,x_2)$, where $g(x_1,x_2)$ is the form of degree $2k$ whose set of roots is $I$.
As any element of $\Aut(\Gamma,g^1_2)$ acts on the basis and leaves invariant the set of ramification points, the morphism $\rho:\Aut(\Gamma,g^1_2) \rightarrow \PGL(2,\K)_I$ is induced by the projection of $\mathbb{P}(k,1,1)$ on the two last coordinates. We see then that the kernel of $\rho$ is generated by the automorphism $\sigma:(w:x_1:x_2) \mapsto (-w:x_1:x_2)$. As $g(x_1,x_2)$ is an eigenvector of any automorphism $\beta \in \PGL(2,\K)_I$, the automorphisms $(w:x_1:x_2) \mapsto (\pm \lambda w: \beta(x_1,x_2))$ are sent by $\rho$ on $\beta$, where $\lambda^2$ is the eigenvalue. We thus get assertion $1$.

Assertion $2$ follows directly from the exact sequence, because the group $<\sigma>$ is of order $2$ and is normal in $G$. 

Let us prove the last assertion. We suppose that the automorphisms $\alpha,\sigma\alpha \in \Aut(\Gamma,g^1_2)$ are conjugate by some automorphism $\nu \in \Aut(\Gamma,g^1_2)$. We write the elements explicitly:
\begin{center}$\begin{array}{rllccccl}
\alpha:&(w:x_1:x_2)& \mapsto (&\h{3}\lambda w&\h{2}:\h{2}& a(x_1,x_2)&\h{3}),\\
\sigma\alpha:&(w:x_1:x_2)& \mapsto (&\h{3}-\lambda w&\h{2}:\h{2}& a(x_1,x_2)&\h{3}),\\
\nu:&(w:x_1:x_2)& \mapsto (&\h{3}\mu w&\h{2}:\h{2}& v(x_1,x_2)&\h{3}),\\
\nu\alpha:&(w:x_1:x_2)& \mapsto (&\h{3}\lambda\mu w&\h{2}:\h{2}& v(a(x_1,x_2))&\h{3}),\\
\sigma\alpha\nu:&(w:x_1:x_2)& \mapsto (&\h{2}-\lambda\mu w&\h{2}:\h{2}& a(v(x_1,x_2))&\h{3})&\h{3},
\end{array}$\end{center}
where $\lambda,\mu \in \K^{*}$, and $a,v \in \PGL(2,\K)_I$.
From the relation $\nu\alpha=\sigma\alpha\nu$, we see that $a$ and $v$ commute in $\PGL(2,\K)$, but in $GL(2,\K)$ we have $a\circ v=\zeta v\circ a$, where $\zeta^{k}=-1$. Using the determinant, we see that $\zeta=\pm1$, which is impossible if $k$ is even.
\proofend
\end{proof}

\section{Appendix 2 - Results on sums}
\pagetitless{Appendix}{Appendix 2 - Results on sums}
\label{Sec:Arithmeticalresults}
We need some results on sums. The following one is the most important:
\begin{Lem}\fauxtitre
\label{Lem:Tecaicarr}
Let $a_1,...,a_k$ be  real numbers; then $$\big(\sum_{i=1}^k a_i\big)^2\leq k\sum_{i=1}^k {a_i}^2$$ and equality holds if and only if all the $a_i$'s are equal.
\end{Lem}
\begin{proof}\fauxtitred\upshape
We prove this by induction on $k$. The case $k=1$ is clear; for $k>1$ we have
\begin{eqnarray*}
\big(\sum_{i=1}^k a_i\big)^2\h{0.5}-\h{0.5}k\h{0.3}\sum_{i=1}^k {a_i}^2\h{2.5}&=\h{2}&\big(a_k+\sum_{i=1}^{k-1} a_i\big)^2-k\sum_{i=1}^{k-1} {a_i}^2-k{a_k}^2\\
&=\h{2}&a_k^2+\big(\sum_{i=1}^{k-1} a_i\big)^2+2a_k\sum_{i=1}^{k-1} a_i-k\sum_{i=1}^{k-1} {a_i}^2-k{a_k}^2\\
&=\h{2}&\big(\sum_{i=1}^{k-1} a_i\big)^2-(k-1)\sum_{i=1}^{k-1} {a_i}^2-\sum_{i=1}^{k-1} {a_i}^2-(k-1){a_k}^2+2a_k\sum_{i=1}^{k-1} a_i\\
&=\h{2}&\Big(\big(\sum_{i=1}^{k-1} a_i\big)^2-(k-1)\sum_{i=1}^{k-1} {a_i}^2\Big)-\sum_{i=1}^{k-1}(a_i-a_k)^2\leq 0.
\end{eqnarray*}
We then get the result by using the induction hypothesis.\proofend
\end{proof}

\bigskip

\begin{remark}\fauxtitred
The inequality of Lemma \refd{Lem:Tecaicarr} is a particular case of Cauchy's inequality
\begin{equation}\label{eq:CauchySchwartz}
\big(\sum_{i=1}^ka_ib_i\big)^2\leq\big(\sum_{i=1}^ka_i^2\big)\big(\sum_{i=1}^kb_i^2\big)
\end{equation}
{\upshape (}take $b_i=1$ for all $i${\upshape )}, itself a consequence of Lagrange's identity
\begin{equation}
\big(\sum_{i=1}^ka_ib_i\big)^2=\big(\sum_{i=1}^ka_i^2\big)\big(\sum_{i=1}^kb_i^2\big)\ -\sum_{1\leq i<j\leq k}(a_ib_j-a_jb_i)^2.
\end{equation}
This identity, if not immediately obvious, can be proved by induction along the lines of the proof of Lemma \refd{Lem:Tecaicarr}. {\upshape (}Equality holds in (\refd{eq:CauchySchwartz}) if and only if $(a_1,...,a_k)$ and $(b_1,...,b_k)$ are proportional.{\upshape )} (See {\upshape \cite{bib:HaLiPo}}, Theorem 7, page 16.)
\end{remark}

\break

\begin{Lem}\fauxtitre
\label{Lem:ArCondSum}
Let $m$, $(s_i)_{i=1}^4$ be integers, $m\geq 1$ and $s_i\geq 0$ for all $i$. If
\begin{itemize}
\item[\upshape 1.]
$\displaystyle\sum_{i=1}^4 s_i^2=m^2-1$,
\item[\upshape 2.]
$\displaystyle\sum_{i=1}^4 s_i=2(m-1)$,
\item[\upshape 3.]
$\displaystyle s_i+s_j \leq m$, for $1\leq i,j\leq 4$, $i\not=j$,
\end{itemize}
then $m=1$ and $s_1=s_2=s_3=s_4=0$.
\end{Lem}\begin{proof}\fauxtitred\upshape
We may assume that $s_1 \geq s_2 \geq s_3 \geq s_4$. The hypotheses imply that $s_1<m$, since $s_1\geq 0$ and $\displaystyle s_1^2\leq \sum_{i=1}^4 s_i^2=m^2-1<m^2$.

\bigskip

We distinguish two cases: $0\leq s_1\leq \frac{m}{2}$ and $\frac{m}{2}<s_1<m$.
\begin{enumerate}
\item
If $0\leq s_1\leq \frac{m}{2}$, then $0\leq s_i\leq \frac{m}{2}$ for all $i$, and 
$$m^2-1=\sum_{i=1}^4 s_i^2\leq \frac{m}{2}\cdot \sum_{i=1}^4 s_i=\frac{m}{2}\cdot 2(m-1)=m^2-m,$$
whence $m\leq 1$. Therefore $m=1$;\ \ and then $a_i=0$ for all $i$, because of the first condition.
\item
We now show that $\frac{m}{2}<s_1<m$ cannot occur, because of the first and third conditions. Define the integer $k$ by $k=2s_1-m$; we then have $s_1=\frac{1}{2}(m+k)$ and $0<k<m$. By the third condition, $s_i\leq \frac{1}{2}(m-k)$ for $i\geq 2$. Hence
$$m^2-1=s_1^2+ \sum_{i=2}^4 s_i^2\leq \bigg(\frac{m+k}{2}\bigg)^{\!2}+3\bigg(\frac{m-k}{2}\bigg)^{\!2}=m^2+k^2-mk,$$
whence $k(m-k)\leq 1$, so that $k=1$ and $m=2$. But then $s_1=\frac{3}{2}$, not an integer.\proofend
\end{enumerate}
\end{proof}
\vspace{5cm}
\bigskip
\break

\section{Appendix 3 - Morphisms $S\rightarrow \mathbb{P}^1$}
\pagetitless{Appendix}{Appendix 3 - morphisms $S\rightarrow \mathbb{P}^1$}
\begin{Lem}\fauxtitre
\label{Lem:Negativity}
Let $S$ be a surface and $\pi:S\rightarrow \mathbb{P}^1$ be a morphism, with one fibre isomorphic to $\mathbb{P}^1$.

Let $F=\sum_{i=1}^m a_i F_i$ be some fibre, with $a_i>0$, $F_i$ an irreducible curve, for $i=1,...,m$. Then, one of the following holds:
\begin{enumerate}
\item[\upshape 1.]
$m=1$, $F=F_1\cong \mathbb{P}^1$.
\item[\upshape 2.]
$m=2$, $F=F_1+F_2$, with $F_1\cong F_2 \cong \mathbb{P}^1$, $F_1\cdot F_2=1$, $F_1^2=F_2^2=-1$.
\item[\upshape 3.]
$m>2${\upshape :} \begin{itemize}
\item $F_i^2<0$ and $a_i=1$ for $i=1,...,m$.
\item
Two of the $F_i's$ that have self-intersection $-1$ do not intersect.
\item
One of the $F_i's$ has self-intersection $\leq -2$;
\end{itemize}
\end{enumerate}
\end{Lem}
\begin{proof}\fauxtitred\upshape
The result follows from a Noether-Enriques theorem (a recent proof of which can be found in \cite{bib:Be1}, Theorem III.4). This theorem implies that there exists a birational map $\varphi:S\dasharrow \mathbb{P}^1\times\mathbb{P}^1$ such that the following diagram commutes:

\begin{center}
\hspace{0 cm}\xymatrix{ S\ar@{-->}^{\varphi}[rr] \ar_{\pi}[rd] & & \mathbb{P}^1\times\mathbb{P}^1 \ar_{\pi_1}[ld] \\
& \mathbb{P}^1, & }
\end{center}
where $\pi_1$ denotes the projection on the first factor. Resolving the indeterminacy of $\varphi$, this is a composition of blow-ups of points and then blow-down of curves. As the diagram commutes, every curve contracted is contained in some fibre.

Any fibre of $S$ is then obtained from $\mathbb{P}^1$, by blowing-up some points , and then blowing-down irreducible component of the fibres, of self-intersection $-1$. 

This implies that the fibre $F$ is connected, that $a_i=1$ for $i=1,...,m$ and that we may number the $F_i's$ so that for $i\not= j$, $F_i \cdot F_j=1$  if and only if $i-j=\pm 1$. Then, every $F_i$ is smooth and is isomorphic to $\mathbb{P}^1$.

The case $m=1$ is then clear.

Using $F^2=0$ and $F \cdot K_S=-2$ (adjunction formula) we obtain the case $m=2$.

 Suppose then that $m>3$. Remark first of all that $F_i^2<0$ for $i=1,...,m$. Indeed, $F_i^2=F_i \cdot (F-\sum_{j\not=i} a_j F_j)<0$ as $F\cdot F_i=0$ ($F$ is equivalent to the divisor of any other general fibre) and $F_i\cdot F_j>0$ for some $j$, as the singular fibre is connected. 
Suppose then that $F_i$ intersects $F_{i+1}$, with $F_i^2=F_{i+1}^2=-1$. Contracting $F_i$, the self-intersection of the image of $F_{i+1}$ is $0$, so the image of $F$ becomes a fibre equal to the image of $F_{i+1}$, so $m=2$. The last assertion follows from the previous one. \proofend
\end{proof}

\chapter*{Résumé de la thèse}
\addcontentsline{toc}{chapter}{Résumé de la thèse}
\ChapterBegin{Le but de ce travail est de classifier à conjugaison près les sous-groupes abéliens finis du groupe de Cremona du plan. Cette classification induit beaucoup de résultats et contribue à la compréhension des propriétés de ce groupe classique, étudié depuis déjà plus d'un siècle.}

\bigskip

Nous passons ici en revue chacun des $10$ chapitres de la thèse, et résumons ceux-ci. Pour plus de cohérence, la numérotation des propositions, définitions, lemmes,... est ici la même que dans la partie principale de la thèse.
\section*{Chapitre 1 - Introduction}
\titrepageresume{1}
Tout d'abord, quelques notions basiques sont données, afin d'introduire le sujet au lecteur, qu'il soit peu spécialiste du domaine, ou même très peu habitué au langage mathématique. La première section est donc élémentaire. La notion de conjugaison, fondamentale dans le sujet, est expliquée et la recherche de classes de conjugaison de groupes finis dans divers groupes de transformations géométriques est décrite au travers d'exemples concrets.

On introduit ensuite le groupe de Cremona  $\CrP$, qui est le groupe des transformations birationnelles du plan $\Pn$, ou de manière équivalente, le groupe des $\K$-automorphismes du corps $\K(X,Y)$.

L'étude des sous-groupes finis du groupe de Cremona a fait l'objet de nombreux articles de la fin du $\mbox{XIX}^e$ siècle à nos jours et n'est toujours pas terminée. Un historique des résultats est donné au chapitre 
\refd{Chap:Introduction}. Il y est expliqué comment le présent travail apporte sa pierre à l'édifice, en donnant une classification des sous-groupes abéliens finis. Celle-ci, quoique pas complètement précise (les classes de conjugaison à l'intérieur de certaines familles à un ou plusieurs paramètres ne sont pas données), est la plus complète à ce jour sur le sujet. De plus, on indique que plusieurs résultats intéressants découlent de cette classification, répondant à des questions soulevées par les anciens travaux, et généralisant également des anciens résultats. Citons-en quelques uns:

\begin{itemize}
\item
Tout d'abord, une question naturelle est de savoir combien il y a de classes de conjugaisons d'éléments d'un ordre donné. La réponse pour les éléments d'ordre $2$ fut donnée par L. Bayle et A. Beauville \cite{bib:BaB} et pour les éléments d'ordre premier par T. de Fernex \cite{bib:deF}.

Nous généralisant ceci en démontrant au théorème \refThmExistenceNonLinear\ l'existence d'une infinité de classes de conjugaison d'éléments d'ordre $2n$, pour $n\geq 1$ quelconque, et l'existence par contre d'une seule classe de conjugaison d'éléments d'ordre $2n+1$, pour tout $n\geq 7$.
\item
Ensuite, il fut démontré par G. Castelnuovo qu'un élément d'ordre fini qui fixe une courbe de genre $\geq 2$ est d'ordre $2$, $3$, $4$ ou $6$.

Ceci soulève donc directement la question de savoir si tous les ordres sont possibles, et de se poser la même question pour des courbe de genre $1$ (le cas des courbes de genre $0$ étant facile, tout les ordres sont possibles). De plus, on peut également se demander si un groupe abélien fini non cyclique fixe une courbe de genre positif.

Nous prouvons au théorème  \refThmCastelnuovoGeneralis\ que si un groupe abélien fini fixe une courbe de genre strictement positif, alors le groupe est cyclique d'ordre $2$, $3$, $4$, $5$ ou $6$ et que tous les cas se présentent. Si le genre de la courbe est $\geq 2$, le groupe est d'ordre $2$ ou $3$ (et les deux cas se présentent).
\item
Il est facile de vérifier que les automorphismes non triviaux de $\Pn$ ou $\mathbb{P}^1\times\mathbb{P}^1$ ne fixent pas de courbe de genre strictement positif. Une question naturelle est de savoir si un sous-groupe fini du groupe de Cremona dont tous les éléments non triviaux ne fixent pas de courbe de genre strictement positif est conjugué à un groupes d'automorphismes d'une de ces deux surfaces.

Le résultat est vrai pour les groupes cycliques d'ordre premier $p$ (démontré pour  $p=2$ dans \cite{bib:BaB}, pour $p=3$ et $p\geq 7$ dans \cite{bib:deF} et pour $p=5$ dans \cite{bib:BeB}). Nous généralisons le résultat aux groupes cycliques d'ordre fini quelconque au théorème\refThmCyclicGenusCurve. Dans le cas des groupes abéliens finis, nous démontrons qu'il n'y a qu'une seule exception (théorème \refThmAbelGenusCurve), à savoir le groupe appelé {\it Cs.24}, engendré par les deux transformations suivantes:
\begin{center}$\begin{array}{lll}
(x:y:z)&\dasharrow &(yz:xy:-xz) \\
(x:y:z)&\dasharrow &( yz(y-z):xz(y+z):xy(y+z)).\end{array}$\end{center}
Ce groupe spécial, isomorphe à $\Z{2}\times\Z{4}$, n'est conjugué ni à un sous-groupe d'automorphismes de $\Pn$ ni à un sous-groupe d'automorphismes $\mathbb{P}^1\times\mathbb{P}^1$, bien qu'aucun de ses éléments non triviaux ne fixe une courbe de genre strictement positif.
\item
Finalement, les classifications précédentes de sous-groupes finis du groupe de Cremona ne donnent pas les classes d'isomorphismes de tels groupes, cette question étant pourtant naturellement plus simple que celle des classes de conjugaisons. Nous donnons les classes d'isomorphismes de sous-groupes abéliens finis du groupe de Cremona au théorème \refThmIsomorphyclasses. Celles-ci sont:
\begin{itemize}
\item
$\Z{m}\times\Z{n}$, pour tous entiers $m,n \geq 1;$
\item
$\Z{2n}\times(\Z{2})^2$, pour tout entier $n\geq 1;$
\item
$(\Z{4})^2\times\Z{2};$
\item
$(\Z{3})^3;$
\item
$(\Z{2})^4$.
\end{itemize}
\end{itemize}

Après avoir explicité les résultats attendus, on entre directement dans le vif du sujet.

\section*{Chapitre 2 - Groupes finis d'automorphismes de surfaces rationnelles}
\titrepageresume{2}
L'approche moderne  (utilisée dans $\cite{bib:BaB}$, $\cite{bib:deF}$, $\cite{bib:Be2}$ et $\cite{bib:Dol}$) est de voir les sous-groupes finis du groupe de Cremona comme des groupes d'automorphismes de surfaces rationnelles projectives et lisses, et de supposer l'action minimale. Le lecteur notera que souvent, pour abréger, nous parlerons uniquement de surface rationnelle, voire juste de surface, pour une surface rationnelle projective et lisse. Toutes les surfaces traitées ici ont ces trois propriétés.

Précisons maintenant ces notions.

\bigskip
\noindent {\bf Définition \refd{Def:Gsurfaces}}\fauxtitred
{\it \begin{itemize}
\item
On note \defn{$(G,S)$} une paire dans laquelle $S$ est une surface rationnelle et $G\subset \Aut(S)$ est un groupe d'automorphismes de la surface. Une paire $(G,S)$ est également appelée classiquement une \defn{$G$-surface}.
\item
Soit  $(G,S)$ une $G$-surface. On dit qu'une application birationnelle $\varphi: S \rightarrow S'$ est \defn{$G$-équivariante} si  la $G$-action sur $S'$ induite par $\varphi$ est birégulière. L'application birationnelle $\varphi$ est appelée \defn{application birationnelle de $G$-surfaces.}
\item
On dit qu'une paire $(G,S)$ est \defn{minimale}  si tout  morphisme birationnel $G$-équivariant   $\varphi:S\rightarrow S'$ est un isomorphisme.
\end{itemize}}

On se doit ici d'observer la chose suivante: une paire $(G,S)$ représente une classe de conjugaison de sous-groupes de $\CrP$. En effet, en choisissant une application birationnelle $\varphi:S\dasharrow \Pn$, on obtient un sous-groupe $\varphi G \varphi^{-1}$ de $\CrP$. Le choix de $\varphi$ ne change pas la classe de conjugaison du groupe et tout élément de la classe de conjugaison peut être obtenu de cette manière. La réciproque est vraie, si $G$ est fini:

\bigskip
\noindent {\bf Proposition \refd{Prp:GS}}\fauxtitre 
{\it Soit $G \subset \CrP$ un groupe fini. Il existe une surface rationnelle projective et lisse $S$ et une application birationnelle $\varphi: S\dasharrow\Pn$ telles que $\varphi^{-1}G \varphi$ agit birégulièrement sur $S$.}

\bigskip

De cette observation suit que notre problème se ramène à la recherche de telles paires, à conjugaison birationnelle près. On peut supposer que les paires sont minimales et utiliser alors la proposition suivante, due à Yu.\ Manin dans le cas abélien (voir \cite{bib:Man}) et V.A.\ Iskovskikh dans le cas général (voir \cite{bib:Isk3}). 
\vspace{3 cm}

\break

\noindent {\bf Proposition \refd{Prp:TwoCases}}\fauxtitre 
{\it Soit $S$ une surface (rationnelle, projective et lisse) et $G \subset \Aut(S)$ un sous-groupe fini d'automorphismes de $S$.
Si la paire $(G,S)$ est minimale, alors une \emph{et une seule}  des assertions suivantes se produit:
\begin{enumerate}
\item[\upshape 1.]
La surface $S$ a une structure de fibré en coniques invariante par $G$, et $\rkPic{S}^G=2$, i.e.\ la partie fixe du groupe de Picard est engendrée par le diviseur canonique et le diviseur d'une fibre.
\item[\upshape 2.]
La partie du groupe de Picard fixée par le groupe $G$ est engendrée par le diviseur canonique, i.e.\ $\rkPic{S}^{G}=1$.\proofend
\end{enumerate}}

Il est intéressant de remarquer ici que les deux cas sont distincts, mais pas birationnellement distincts (voir la section \refd{Sec:PencilRatCurves}).

Cette dichotomie organise nos recherches.  Nous classifions toutes les paires minimales dans les deux cas. Ensuite, il faut encore éliminer les doublons, à conjugaison birationnelle près. Cette question est considérée dans les chapitres \refd{Chap:ConjBetweenTheCases} et \refd{Chap:List}.

\section*{Chapitre 3 - Invariants de conjugaison}
\titrepagesresume{3,4}
Une classe de conjugaison de sous-groupes finis du groupe de Cremona a une infinité de représentations comme groupe de transformations birationnelles d'une surface. 

On exprime dans ce chapitre quelques invariants de ces représentations, qui permettent dans certains cas de prouver que deux groupes isomorphes agissant birationnellement sur deux surfaces ne sont pas birationnellement conjugués (i.e.\ représentent des classes de conjugaison distinctes de sous-groupes du groupe de Cremona).

On note $G \subset \Cr{S}$ un groupe abélien fini de transformations birationnelles d'une surface $S$. Les invariants utilisés dans ce travail sont les suivants (la numérotation est celle des sections où les invariants sont traités):
\begin{itemize}
\item[\refd{Sec:PencilRatCurves}]
Existence de pinceaux de courbes rationnelles laissés invariants.
\item[\refd{Sec:PairInvariantPencils}] Existence de deux pinceaux $G$-invariants de courbes rationnelles sur la surface $S$, avec une intersection libre donnée.
\item[\refd{Sec:ExistenceFixedPoints}] Existence de points fixes, lorsque $G \subset Aut(S)$.
\item[\refd{Sec:CurveFixedPoints}] L'ensemble des courbes non-rationnelles fixées par les éléments non triviaux de $G$.
\item[\refd{Sec:ActionFixedCurves}] L'action du groupe sur les courbes non-rationnelles fixées par les éléments non-triviaux de $G$.
\end{itemize}

\section*{Chapitre 4 - Surfaces de Del Pezzo}
\titrepagesresume{3,4}
Les surfaces de Del Pezzo (surfaces rationnelles dont le diviseur anti-canonique est ample) jouent un rôle prépondérant dans l'étude des sous-groupes finis de transformations birationnelles du plan. 

En effet, si $G \subset \Aut(S)$ est un groupe fini d'automorphismes de $S$ et si  $\rkPic{S}^G=1$ (cas $2$ de la proposition \refd{Prp:TwoCases}), alors $S$ est une surface de Del Pezzo. De plus, même dans le cas $1$ de la proposition \refd{Prp:TwoCases}, des surfaces de Del Pezzo apparaissent.

\break

Citons la caractérisation suivante:

\bigskip
\noindent {\bf Proposition \refd{Prp:ConditionsDelPezzo}}\fauxtitre 
{\it Soit $S$ une surface rationnelle projective et lisse. Les conditions suivantes sont équivalentes:
\begin{enumerate}
\item[\upshape 1.]
$S$ est une surface de Del Pezzo;
\item[\upshape 2.]
$S \cong \mathbb{P}^1\times\mathbb{P}^1$ ou $S\cong \Pn$ ou $S$ est l'éclatement de $1\leq r \leq 8$ points en position générale de $\Pn$ (i.e.\  $3$ ne sont pas alignés, $6$ ne sont pas sur la même conique, et il n'existe pas de cubique passant par tous les points, avec multiplicité $2$ à l'un des points);
\item[\upshape 3.]
$K_S^2\geq 1$ et toute courbe irréductible de $S$ a une self-intersection $\geq -1$;
\item[\upshape 4.]
$C\cdot (-K_S)>0$ pour tout diviseur effectif $C$.
\end{enumerate}}

On décrit également le groupe de Picard des surfaces de Del Pezzo, le nombre de diviseurs exceptionnels et les structures de fibré en coniques que l'on peut mettre sur de telles surfaces.

\section*{Chapitre 5 - Automorphismes de fibrés en coniques}
\titrepageresume{5}
Dans ce chapitre, quelques outils concernant les groupes d'automorphismes de fibrés en coniques sont introduits.

Rappelons qu'un fibré en coniques est une surface rationnelle projective et lisse $S$, équipée d'un morphisme $\pi:S\rightarrow \mathbb{P}^1$ dont la fibre générique est isomorphe à $\mathbb{P}^1$. Il y a un nombre fini de fibres singulières consistant en l'union transverse de deux courbes rationnelles de self-intersection $-1$.

Deux approches sont possibles:
\begin{itemize}
\item
On peut étudier géométriquement les groupes d'automorphismes de fibrés en coniques.\\
Dans ce point de vue, on a l'avantage d'avoir des automorphismes (biréguliers) et de pouvoir utiliser l'action du groupe sur les fibres singulières et leurs composantes.
\item
Birationnellement, les groupes d'automorphismes de fibrés en coniques sont conjugués a des sous-groupes du groupe $\CrPp$ des transformations birationnelles de $\mathbb{P}^1\times\mathbb{P}^1$ qui laissent invariante la première projection.\\

Ce deuxième point de vue est plus algébrique. 
Explicitement, $\CrPp$ est le groupe des transformations birationnelles de la forme :
$$(x,y)\dasharrow \left(\frac{ax+b}{cx+d},\frac{\alpha(x)y+\beta(x)}{\gamma(x)y+\delta(x)}\right),$$
où $a,b,c,d \in \K$, $\alpha,\beta,\gamma,\delta \in \K(x)$, et $(ad-bc)(\alpha\delta-\beta\gamma)\not=0$.
Ce groupe, appelé le \defn{groupe de de Jonqui\`eres}, est isomorphe à $\PGL(2,\K(x))\rtimes \PGL(2,\K)$. On pourra donc utiliser la structure du groupe.
\end{itemize}

Dans les deux cas, on a un groupe $G$ qui agit birationnellement sur un fibré en coniques $(S,\pi)$. Ceci induit un homomorphisme naturel
$\overline{\pi}: G
\rightarrow \Aut(\mathbb{P}^1)= \PGLn{2}$ qui satisfait $\overline{\pi}(g)\pi=\pi g$, pour tout $g \in G$.

On utilisera la suite exacte
\begin{equation}
1 \rightarrow G' \rightarrow G \stackrel{\overline{\pi}}{\rightarrow} \overline{\pi}(G) \rightarrow 1
\tag{\refd{eq:ExactSeqCB}}\end{equation} afin de réduire les possibilités pour la structure de $G$. 
En effet, si $G$ est abélien et fini, alors $G'$ et $\overline{\pi}(G)$ le sont également. Les sous-groupes abéliens finis de $\PGLn{2}$ et $\PGL(2,\K(x))$ sont cycliques ou isomorphes à $(\Z{2})^2$ (voir le lemme \refd{Lem:AutP1}). On examine toutes les possibilités au Chapitre  \refd{Chap:FiniteAbConicBundle}.

Citons finalement un dernier outil simple mais important:

\bigskip
\noindent {\bf Lemme \refd{Lem:MinTripl}}\fauxtitre 
{\it Soit $G \subset \Aut(S,\pi)$ un groupe fini d'automorphismes du fibré en coniques $(S,\pi)$.
Les conditions suivantes sont équivalentes:
\begin{enumerate}
\item[\upshape 1.]
Le triplet $(G,S,\pi)$ est minimal.
\item[\upshape 2.]
Pour toute fibre singulière $\{F_1,F_2\}$ de $\pi$, 
 il existe $g \in G$ tel que $g(F_1)=F_2$ (ce qui implique que $g(F_2)=F_1$).
\end{enumerate}}

En supposant donc le triplet $(G,S,\pi)$ minimal, on obtient pour chaque fibre singulière un élément qui \defn{tord} celle-ci, i.e.\ qui permute les deux composantes de la fibre.

\section*{Chapitre 6 - Groupes abéliens finis d'automorphismes de surfaces de grand degré}
\titrepagesresume{6,7}
On s'intéresse dans ce chapitre aux surfaces ayant un degré $\geq 5$ (le degré étant le carré du diviseur canonique). Parmi celles-ci, les deux surfaces les plus simples et les plus étudiées sont $\Pn$ est $\mathbb{P}^1\times\mathbb{P}^1$. On donne la classification complète des sous-groupes abéliens finis de $\Aut(\Pn)$ et $\Aut(\mathbb{P}^1\times\mathbb{P}^1)$, à conjugaison birationnelle près (voir respectivement les propositions \refd{Prp:PGL3Cr} et \refd{Prp:AutP1Bir}). 

De plus, on obtient la proposition importante suivante:

\bigskip
\noindent {\bf Proposition \refd{Prp:LargeDegree}}\fauxtitre 
{\it Soit $G \subset \Aut(S)$ un groupe abélien fini d'automorphismes d'une surface rationnelle $S$ et supposons que $K_S^2\geq 5$. Alors $G$ est birationnellement conjugué à un sous-groupe de $\Aut(\mathbb{P}^1\times\mathbb{P}^1)$ ou $\Aut(\Pn)$.}

\section*{Chapitre 7 - Groupes abéliens finis d'automorphismes de fibrés en coniques}\titrepagesresume{6,7}
Comme cité précédemment, un groupe agit de manière minimale sur un fibré en coniques $(S,\pi)$ si et seulement si chaque fibre est tordue par un des éléments du groupe. Il y a deux exemples fondamentaux d'éléments qui tordent des fibres singulières:
\begin{itemize}
\item
On note $S_6$ la surface de Del Pezzo de degré $6$ donnée dans $\Pn \times \Pn$ sous la forme
\begin{center}
$S_6=\{ (x:y:z) \times (u:v:w)\ | \ ux=vy=wz\}$
\end{center}
(par commodité de notations, une élément de $\Pn\times\Pn$ sera noté $P\times Q$ au lieu de $(P,Q)$) et soit $(S_6,\pi)$ le fibré en coniques induit par 
\begin{center}
$\pi((x:y:z) \times (u:v:w))=\left\{\begin{array}{ll}
(y:z)& \mbox{ si } (x:y:z) \not= (1:0:0),\\
(w:v)& \mbox{ si } (u:v:w) \not= (1:0:0).\end{array}\right.$
\end{center}
Il y a exactement deux fibres singulières, qui sont
\begin{center}
$\begin{array}{lll}
\pi^{-1}(1:0)=\{(0:1:0) \times (a:0:b) \ | \ (a:b) \in \mathbb{P}^1\} \cup \{(a:b:0) \times (0:0:1) \ | \ (a:b) \in \mathbb{P}^1\}\\
\pi^{-1}(0:1)=\{(0:0:1) \times (a:b:0) \ | \ (a:b) \in \mathbb{P}^1\} \cup \{(a:0:b) \times (0:1:0) \ | \ (a:b) \in \mathbb{P}^1\}
\end{array}$\end{center}
(Une description précise de la surface et de ses trois structures de fibré en coniques est donnée à la section  \refd{Subsec:DelPezzo6}.)
Pour tout $\alpha,\beta \in \K^{*}$, on définit $\kappa_{\alpha,\beta}$ comme l'automorphisme de $S_6$ suivant:
\begin{center}
$\kappa_{\alpha,\beta}:(x:y:z) \times (u:v:w) \mapsto (u:\alpha w:\beta v)\times (x:\alpha^{-1} z:\beta^{-1} y)$.\end{center}
L'automorphisme $\kappa_{\alpha,\beta}$ tord les deux fibres singulières de $\pi$. De plus, $\kappa_{\alpha,\beta}$ est une involution si et seulement si son action sur la fibration est triviale.
\item
\titrepageresume{7}
Le deuxième exemple correspond aux involutions de de Jonqui\`eres qui tordent un nombre pair de fibres singulières. Le lemme suivant les décrit:

\bigskip
\noindent {\bf Lemme \refd{Lem:DeJI}}\titreProp{involutions de de Jonqui\`eres qui tordent un fibré en coniques}\\
{\it Soit $g \in \Aut(S,\pi)$ un automorphisme d'ordre fini tel que  
$g$ tord au moins une fibre singulière de $\pi$.

Alors les assertions suivantes sont équivalentes:
\begin{enumerate}
\item[\upshape 1.]
$g$ est une involution;
\item[\upshape 2.]
$\overline{\pi}(g)=1$, i.e.\ $g$ agit trivialement sur la fibration;
\item[\upshape 3.]
L'ensemble des points de $S$ fixés par $g$ est une courbe hyperelliptique de genre $k-1$ (revêtement double de $\mathbb{P}^1$ par $\pi$, ramifié au-dessus de $2k$ points), plus éventuellement un nombre fini de points fixes isolés, qui sont les points singuliers des fibres singulières qui ne sont pas tordues par $g$.
\end{enumerate}
De plus, si les trois conditions sont satisfaites, le nombre de fibres singulières de $\pi$ tordues par $g$ est $2k$.}

\bigskip
\noindent {\bf Définition \refd{Def:DeJI}}\titreProp{involutions de de Jonqui\`eres qui tordent le fibré $(S,\pi)$}\\
{\it Soit $(S,\pi)$ un fibré en coniques. On dit qu'un automorphisme  $g \in \Aut(S,\pi)$ qui satisfait les conditions du lemme \refd{Lem:DeJI} est une \defn{involution de de Jonqui\`eres qui tord le fibré en coniques $(S,\pi)$}.}
\end{itemize}

On remarque alors que les deux exemples représentent presque tous les cas d'automorphismes qui tordent un fibré en coniques:

\bigskip
\noindent {\bf Proposition \refd{Prp:TwocasesKappaTwist}}\fauxtitre 
{\it Soit $(S,\pi)$ un fibré en coniques et soit $g \in \Aut(S,\pi)$ un automorphisme d'ordre fini qui tord au moins une fibre singulière.
\begin{enumerate}
\item[\upshape 1.]
Si $\overline{\pi}(g)=1$, i.e.\ si l'action de $g$ sur la base est triviale, alors $g$ est une involution qui tord le fibré en coniques $(S,\pi)$.
\item[\upshape 2.]
Si le groupe cyclique  $G$  engendré par $g$ ne contient pas d'involution qui tord $(S,\pi)$, il existe un morphisme birationnel de fibrés en coniques $g$-équivariant  $\eta:S\rightarrow S_6$ où $S_6$ est la surface de Del Pezzo de degré $6$ et $\eta g\eta^{-1}=\kappa_{\alpha,\beta}$ pour certains $\alpha,\beta \in \K^{*}$.\\
De plus on a:
\begin{enumerate}
\item[\upshape (a)]
$\overline{\pi}(g)\not=1$;
\item[\upshape (b)]
$g$ tord exactement deux fibres singulières;
\item[\upshape (c)]
l'ordre de $g$ est $4n$, pour un certain entier $n$;
\item[\upshape (d)]
la suite exacte induite par  $\overline{\pi}$ est:
$1 \rightarrow G' \rightarrow G \stackrel{\overline{\pi}}{\rightarrow} \Z{2n} \rightarrow 1$,\\
où $G' \cong \Z{2}$.
\end{enumerate}
\end{enumerate}}

On donne ensuite le plus important exemple de la thèse. C'est le seul sous-groupe abélien fini qui n'est pas conjugué à un sous-groupe d'automorphismes de $\Pn$ ou $\mathbb{P}^1\times\mathbb{P}^1$, bien qu'aucun de ses éléments non triviaux ne fixe une courbe de genre positif (voir  le théorème \refThmAbelGenusCurve).

Le groupe, appelé ${Cs.24}$, agit sur une surface de degré $4$, qui n'est pas du type Del Pezzo. Elle s'obtient en éclatant deux points particuliers de la surface de Del Pezzo de degré $6$.

Le groupe $Cs.24$ peut être vu dans $\CrP$ comme le groupe engendré par les applications birationnelles $(x:y:z)\dasharrow (yz:xy:-xz)$ et $(x:y:z)\dasharrow ( yz(y-z):xz(y+z):xy(y+z))$.
On décrit le groupe précisément, et on démontre que celui-ci n'est pas birationnellement conjugué à un sous-groupe d'automorphismes de $\Pn$ ou $\mathbb{P}^1\times\mathbb{P}^1$.
 
 De plus, on prouve la proposition suivante:

\bigskip
\noindent {\bf Proposition \refd{Prp:CBdistCases}}\fauxtitre 
{\it Soit $(S,\pi)$ un fibré en coniques et soit $G \subset \Aut(S,\pi)$ un groupe abélien fini tel que 
le triplet $(G,S,\pi)$ soit minimal.

Alors, à conjugaison birationnelle près, une des assertions suivantes se produit:
\begin{itemize}
\item
Une involution qui tord le fibré  $(S,\pi)$ appartient à $G$;
\item
le groupe $G$ est un sous-groupe d'automorphismes de $\mathbb{P}^1\times\mathbb{P}^1$;
\item
le groupe $G$ est ${Cs.24}$.
\end{itemize}}

Du fait que les groupes d'automorphismes de  $\mathbb{P}^1\times\mathbb{P}^1$ ont été classifiés à la section \refd{SubSec:AutP1}, il ne reste donc plus qu'à étudier les groupes d'automorphismes de fibrés en coniques qui contiennent une involution qui tord le fibré. 

Tout d'abord, on décrit le cas où $G$ est engendré par une telle involution. On montre que l'involution est conjuguée à une transformation birationnelle
\begin{center}$\sigma_{a,b}:(x_1:x_2)\times (y_1:y_2) \dasharrow (x_1:x_2) \times(y_2 \prod_{i=1}^k(x_1-b_i x_2):y_1 \prod_{i=1}^k(x_1-a_i x_2))$\end{center}de $\mathbb{P}^1\times\mathbb{P}^1$, pour un certain entier $k\geq 1$ et $a_1,...,a_k,b_1,...,b_k \in \K$ tous distincts.

On étudie également les classes de conjugaison de tels éléments, dans le groupe de Jonqui\`eres et le groupe de Cremona, prouvant de ce fait un résultat classique dont la première preuve complète est donnée dans \cite{bib:BaB}.

Puis, on étudie le cas où $G$ contient une involution qui tord le fibré, mais n'est pas engendré par celle-ci.
On peut voir que dans ce cas, $G'=G\cap \ker \overline{\pi}$ contient cette involution et est soit cyclique d'ordre $2$, soit isomorphe à $(\Z{2})^2$. On utilise la suite exacte
\begin{equation}
1 \rightarrow G' \rightarrow G \stackrel{\overline{\pi}}{\rightarrow} \overline{\pi}(G) \rightarrow 1
\tag{\refd{eq:ExactSeqCB}}\end{equation}
introduite précédemment, et en énumérant les cas pour $\overline{\pi}(G)$ et $G'$, on trouve en général deux possibilités pour $G$, suivant que la suite exacte est scindée ou non.

Dans le cas où la suite est scindée, on utilise quelques résultats de cohomologie des groupes de \cite{bib:Be2} et on obtient une forme ``standard''. Dans le cas où la suite n'est pas scindée, le groupe $G$ contient une racine paire d'une involution de de Jonqui\`eres. On prouve l'existence de telles racines, pour des ordres arbitrairement grands.

Voici quelques résultats et exemples qui illustrent cette situation.

Tout d'abord, deux des nombreux résultats dans le cas où la suite est scindée:

\bigskip
\noindent {\bf Proposition \refd{Prp:DirectProductC2C2C2}}\fauxtitre 
{\it Soit $G \subset \CrPp $ un sous-groupe fini abélien du groupe de de Jonqui\`eres tel que
\begin{itemize}
\item
La suite exacte (\refd{eq:ExactSeqCB}) est scindée.
\item
$\overline{\pi}(G)\cong (\Z{2})^2$.
\item
$G'\cong \Z{2}$ est engendré par une involution de de Jonqui\`eres dont l'ensemble des points fixes est une courbe hyperelliptique de genre $k-1$.
\end{itemize}
Alors $k\geq 4$ est divisible par $4$, et $G\cong (\Z{2})^3$ est conjugué dans le groupe de Cremona à un groupe du type ${\num{C.2,\h{0.5}22}}$, engendré par
\begin{center}
$\begin{array}{rrcccl}
&\sigma:(x_1:x_2)\times (y_1:y_2)& \dasharrow &(x_1:x_2) &\times&(y_2 \prod_{i=1}^{k/4}P(b_i):y_1 \prod_{i=1}^{k/4}P(a_i))\\
&\rho_1:(x_1:x_2)\times (y_1:y_2)& \mapsto &(x_1:- x_2)&\times& (y_1:y_2)\\
&\rho_2:(x_1:x_2)\times (y_1:y_2)& \mapsto &(x_2: x_1)&\times& (y_1:y_2),\end{array}$\end{center}
où $P(a)=(x_1-ax_2)(x_1+ax_2)(x_1-a^{-1}x_2)(x_1+a^{-1}x_2) \in \K[x_1,x_2]$ et $a_1,...,a_{k/4},b_1,...,b_{k/4} \in \K\backslash\{0,\pm 1\}$ sont tous distincts.

De plus, deux tels groupes sont conjugués dans $\CrPp$ (et donc dans $\CrP$) si et seulement si ils \emph{ont la même action sur les ensembles de points fixes de leurs éléments} 
(voir la définition \refd{Def:actionFixedPoints}). 
Explicitement, la courbe $y_1^2 \prod_{i=1}^{k/4}P(a_i)= y_2^2 \prod_{i=1}^{k/4}P(b_i)$ fixée par $\sigma$, et l'action de $\rho_1$ et $\rho_2$ (correspondant à $(x_1:x_2)\ \mapsto (x_2: x_1)$  et $(x_1:x_2)\ \mapsto (x_1: -x_2)$) sur celle-ci déterminent la classe de conjugaison du groupe.}

\bigskip
\noindent {\bf Proposition \refd{Prp:SplitC2C2}}\fauxtitre 
{\it  Soit $G \subset \CrPp$ un sous-groupe abélien fini du groupe de de Jonqui\`eres tel que $G'\cong (\Z{2})^2$, $\overline{\pi}(G)\not=1$ et que la suite exacte (\refd{eq:ExactSeqCB}) est scindée.
 
Alors, à conjugaison près dans le groupe de de Jonqui\`eres, $G$ est l'un des groupes suivants:
 \begin{itemize}
\item
 $\overline{\pi}(G)\cong \Z{n}$, $n>1$:  
 ${G \cong (\Z{2})^2\times \Z{n}}$ est engendré par
\begin{center}
$\begin{array}{rrcccl}
\num{C.22,\h{0.5}n}&(x,y)& \dasharrow &(x,g(x)/y)\\
&(x,y)& \dasharrow &(x,(h(x)y-g(x))/(y-h(x)))\\
&(x,y)& \dasharrow & (\zeta_n x,y),\end{array}$\end{center}
pour certains $g(x),h(x) \in \K(x) \backslash \{0\}$ invariants par l'action $x\mapsto \zeta_n x$, avec $\zeta_n=e^{2\im\pi/n}$.  Le genre de la courbe hyperelliptique fixée par n'importe quelle involution de $G'$ est congru à\hspace{0.1 cm} $-1 \pmod{n}$.
\item
 $\overline{\pi}(G)\cong (\Z{2})^2$:  
 ${G \cong (\Z{2})^4}$ est engendré par
\begin{center}
$\begin{array}{rrcccl}
\num{C.22,\h{0.5}22}&(x,y)& \dasharrow &(x,g(x)/y)\\
&(x,y)& \dasharrow &(x,(h(x)y-g(x))/(y-h(x)))\\
&(x,y)& \dasharrow & (\pm x^{\pm 1},y),\end{array}$\end{center}
pour certains $g(x),h(x) \in \K(x) \backslash \{0\}$ invariants par $x \mapsto \pm x^{\pm 1}$.  Le genre de la courbe hyperelliptique fixée par n'importe quelle involution de $G'$ est congru à $3 \pmod{4}$.
\end{itemize}}

Dans le cas où la suite n'est pas scindée, le groupe contient des racines d'ordre pair d'involutions de de Jonqui\`eres, dont nous décrivons la géométrie:

\bigskip
\noindent {\bf Proposition \refd{Prp:rootsDeJ}}\titreProp{racines d'ordre pair d'involutions de de Jonqui\`eres}\\
{\it 
Soit $(S,\pi)$ un fibré en coniques, soit $\alpha \in \Aut(S,\pi)$ et soit $G$ le groupe engendré par $\alpha$. On suppose que $\alpha^n$ est une involution qui tord le fibré $(S,\pi)$, pour un certain entier pair $n\geq 2$ et que le triplet $(G,S,\pi)$ est minimal. Soit $\Gamma$ la courbe hyperelliptique lisse fixée par $\alpha^n$, qui est un revêtement double de $\mathbb{P}^1$ par la projection $\pi$, ramifié au-dessus de $2k$ points. Alors:

\begin{enumerate}
\item[\upshape 1.]
L'action de $\alpha$ sur la base est d'ordre $n$.
\item[\upshape 2.]
L'action de $\alpha$ sur l'ensemble de ramification n'a pas de points fixes. En particulier, $n$ divise $2k$.
\item[\upshape 3.]
Notons $r$ le nombre d'orbites de l'action de $\alpha$ sur les points de ramification (de telle sorte que $rn=2k$) et par $l$ le nombre de fibres singulières de $\pi$. Alors une des situations suivantes se produit:
\begin{enumerate}
\item[\upshape (a)]
$l=2k$, $r$ est pair, l'action de $\alpha$ sur $\Gamma$ a $4$ points fixes;
\item[\upshape (b)]
$l=2k+1$, $r$ est impair, l'action de $\alpha$ sur $\Gamma$ a $2$ points fixes;
\item[\upshape (c)]
$l=2k+2$, $r$ est pair, l'action de $\alpha$ sur $\Gamma$ n'a pas de points fixes.
\end{enumerate}
\item[\upshape 4.]
L'involution $\alpha^n$ tord $2k$ des $l$ fibres singulières de $\pi$, et $\alpha$ tord les $l-2k$ autres.
\end{enumerate}}

Cette proposition décrit la géométrie de l'action sur le fibré en coniques et implique que l'action de $\alpha$ sur la courbe $\Gamma$ ne fixe aucun point de ramification. Réciproquement, étant donnée une telle courbe et une telle action sur elle, la proposition suivante prouve l'existence d'une racine $\alpha$ de l'involution de de Jonqui\`eres correspondant à $\Gamma$:

\bigskip
\noindent {\bf Proposition \refd{Prp:ExistRoot}}\titreProp{existence des racines d'involutions de de Jonqui\`eres}\\
{\it 
Soit $\beta \in \PGL(2,\K)$ un automorphisme de $\mathbb{P}^1$ d'ordre pair $n\geq 2$ et soit $\{p_1,...,p_{rn}\}$ un ensemble de $rn$ points distincts de $\mathbb{P}^1$ invariant par $\beta$.
On suppose qu'aucun point $p_i$ n'est fixé par $\beta$.

Alors, pour $l=rn+1$ si $r$ est impair (respectivement pour $l=rn$ et $l=rn+2$ si $r$ est pair), il existe un fibré en coniques $(S,\pi)$, qui est l'éclatement de $l$ points de $\mathbb{P}^1\times\mathbb{P}^1$, et $\alpha \in \Aut(S,\pi)$ tel que
\begin{itemize}
\item
le triplet $(<\alpha>,S,\pi)$ est minimal;
\item
$\alpha^n$ est une involution qui tord le fibré $(S,\pi)$, dont l'ensemble des points fixes est l'union disjointe de $l-rn$ points et d'une courbe  hyperelliptique lisse $\Gamma$, telle que $\pi:\Gamma \rightarrow \mathbb{P}^1$ est un revêtement double ramifié au-dessus de $p_1,...,p_{rn}$.
\item
L'action de $\overline{\alpha}$ sur $\mathbb{P}^1$ correspond à $\beta$.
\end{itemize}}

\bigskip

Donnons encore quelques exemples qui illustrent quelques possibilités de groupes:

\bigskip
\noindent {\bf Exemple \refd{Exa:Root2oddorder}}\fauxtitre 
{\it 
Soit $n$ un entier impair et soit $g \in \K(x)$ une fonction rationnelle. Alors l'application birationnelle
\begin{center}$\alpha:(x,y)\dasharrow (\zeta_{2n}\cdot x, -g(x^n)\cdot \frac{y+g(-x^n)}{y+g(x^n)})$\end{center} (où $\zeta_{2n}=e^{2\im\pi/2n}$) est une racine d'ordre $2n$ de l'involution de de Jonqui\`eres \begin{center}$\alpha^{2n}:(x,y)\dasharrow (x, \frac{g(x^n)\cdot g(-x^n)}{y})$\end{center} dont l'ensemble de points fixes est $y^2=g(x^n)\cdot g(-x^n)$. 
}

\bigskip

\bigskip
\noindent {\bf Exemple \refd{Exa:4throotDeJonq}}\fauxtitre 
{\it 
Une racine $4$-ième d'une involution de de Jonqui\`eres:
\begin{center}
$\begin{array}{llcccl}
\alpha&:(x,y) \dasharrow (&\im x&,&\frac{(x+1)((\sqrt{2}-1)-x)y+(x^4-1)}{y+(x+1)((\sqrt{2}-1)-x))}&)\\
\alpha^2&:(x,y) \dasharrow (&-x&,&\frac{-\im(x+1)(x-\im)y+x^4-1}{y-\im(x+1)(x-\im)}&)\\
\alpha^4&:(x,y) \dasharrow (&x&,&\frac{x^4-1}{y}&)\\
\alpha^8&:(x,y) \dasharrow (&x&,&y&)
\end{array}$
\end{center}}

\bigskip
\noindent {\bf Exemple \refd{Exa:C4C2ellratdoesnotsplit}}\fauxtitre 
{\it 
Soit $G\cong \Z{4}\times\Z{2}$ le groupe engendré par $\alpha$, $\beta$, éléments d'ordre respectivement $4$ et $2$:
\begin{center}
$\begin{array}{llcccl}
\alpha&:(x,y) \dasharrow (&-x&,&(1+\im x)\frac{y-(1-\im x)}{y-(1+\im x)}&)\\
\alpha^2&:(x,y) \dasharrow (&x&,&\frac{x^2+1}{y}&)\\
\beta&:(x,y) \dasharrow (&x&,&\frac{y(1+x)-(1+x^2)}{y-(1+x)}&)\\
\alpha^4=\beta^2&:(x,y) \dasharrow (&x&,&y&)
\end{array}$
\end{center}
Les éléments $\alpha^2$ et $\beta$ sont des involutions de de Jonqui\`eres agissant trivialement sur la fibration (qui tordent chacun $2$ fibres d'un fibré en coniques), et l'action de $\alpha$ sur la fibration est cyclique d'ordre $2$.

Les éléments $\alpha^2$ et $\beta$ fixent tous deux une courbe rationnelle (ramifiée respectivement au-dessus de $\{\pm \im\}$ et $\{0,\infty\}$), et $\alpha^2\beta$ fixe une courbe elliptique (ramifiée au-dessus de $\{0,\pm \im, \infty\}$).
}

\section*{Chapitre 8 - Groupes abéliens finis d'automorphismes de surfaces de Del Pezzo de degré $\leq 4$}
\titrepageresume{8}
On décrit dans ce chapitre les paires $(G,S)$, où  $S$ est une surface rationnelle et  $G \subset \Aut(S)$ est un groupe abélien fini tel que $\rkPic{S}^G=1$ (cas $2$ de la  proposition \refd{Prp:TwoCases}). Comme il fut prouvé au lemme \refd{Lem:rkDelPezzo}, la surface $S$ est du type Del Pezzo. Les groupes abéliens finis d'automorphismes de surfaces de degré $\geq 5$ ayant déjà été classés au chapitre \refd{Chap:LargeDegreeSurfaces} (voir en particulier la proposition 
\refd{Prp:LargeDegree}), nous considérerons uniquement les cas restants, à savoir les surfaces de Del Pezzo de degré 1, 2, 3 ou 4.

Dans le cas des surfaces de degré 4, en utilisant l'action sur les fibrés en coniques, nous prouvons le théorème classique (voir par exemple \cite{bib:Ho1}) qui donne la structure des groupes d'automorphismes de telles surfaces: c'est un produit semi-direct $(\Z{2})^4\rtimes H$. Nous prouvons également que $H$ correspond aux automorphismes du plan  qui laissent invariant l'ensemble des $5$ points éclatés pour obtenir la surface. De plus, nous donnons la géométrie de tous les éléments du groupe d'automorphismes.

Le plongement classique de $S$ dans $\mathbb{P}^4$ comme intersection de deux quadriques est également donné de manière explicite, directement à partir des courbes elliptiques fixes par les éléments de $(\Z{2})^4$.
Nous donnons ensuite les groupes abéliens possibles, ainsi que quelques résultats concernant leurs classes de conjugaison birationnelle, généralisant notamment un théorème de \cite{bib:Be2}.

Citons un résultat récapitulatif:

\bigskip
\noindent {\bf Proposition \refd{Prp:DelP4}}\fauxtitre 
{\it 
Soit $G \subset \Aut(S)$ un groupe abélien fini d'automorphismes d'une surface de Del Pezzo $S$ de degré $4$, telle que $\rkPic{S}^G=1$. Alors:
\begin{enumerate}
\item[\upshape 1.]
À isomorphisme près, 
$S\subset \mathbb{P}^4$ est l'éclatement de $(1:0:0)$, $(0:1:0)$, $(0:0:1)$, $(1:1:1)$, $(a:b:c) \in \Pn$, et est donné par les équations
\begin{center}$\begin{array}{cccccc}cx_1^2&&-ax_3^2&+(a-c)x_4^2&-ac(a-c)x_5^2&=0\\
&cx_2^2&-bx_3^2&+(c-b)x_4^2&-bc(c-b)x_5^2&=0.\end{array}$\end{center}
Une des circonstances suivantes se produit:
\begin{itemize}
\item[$\nump{4.22}$]
$G=\{(x_1:x_2:x_3:x_4:x_5) \mapsto (\pm x_1:\pm x_2:x_3:x_4:x_5)\}\cong (\Z{2})^2$.
\item[$\nump{4.222}$]
$G=\{(x_1:x_2:x_3:x_4:x_5) \mapsto (\pm x_1:\pm x_2:\pm x_3:x_4:x_5)\}\cong (\Z{2})^3$.
\item[$\nump{4.2222}$]
$G=\{(x_1:x_2:x_3:x_4:x_5) \mapsto (\pm x_1:\pm x_2:\pm x_3:\pm x_4:x_5)\}\cong (\Z{2})^4$.
\item[$\nump{4.42}$]
$(a:b:c)=(1:\xi:1+\xi)$, pour un certain $\xi\in \K\backslash \{0,\pm 1\}$\\
$G\cong \Z{4} \times \Z{2}$,  engendré par $(x_1:x_2:x_3:x_4:x_5) \mapsto (-x_2:x_1:x_4:x_3: - x_5)$
 et $(x_1:x_2:x_3:x_4:x_5) \mapsto (x_1:x_2:x_3:x_4:-x_5)$.
\item[$\nump{4.223}$]
$(a:b:c)=(1:\zeta_6:1+\zeta_6)$, où $\zeta_6=e^{2\im\pi/6}$ \\
$G\cong (\Z{2})^2\times\Z{3}$, engendré par
$(x_1:x_2:x_3:x_4:x_5) \mapsto (x_1:x_2:\pm x_3:\pm x_4:x_5)$
et $(x_1:x_2:x_3:x_4:x_5) \mapsto ((\zeta_6-1)x_5:x_1:(\zeta_6-1)x_3:-\zeta_6 x_4:-\zeta_6 x_2)$
\item[$\nump{4.8}$]
$(a:b:c)=(1:\im:1+\im)$ \\
$G\cong \Z{8}$ engendré par $(x_1:x_2:x_3:x_4:x_5) \mapsto (x_4:x_3:-x_1:x_2:\im x_5)$
\end{itemize}
\item[\upshape 2.]
$G$ laisse invariant un pinceau de courbes rationnelles de $S$ si et seulement si la paire $(G,S)$ est conjuguée à $\nump{4.22}$ ou $\nump{4.223}$. \proofend
\end{enumerate}}
\begin{remarque}\fauxtitred
Chaque surface de Del Pezzo de degré $4$ possède des groupes d'automorphismes du type $\nump{4.22}$, $\nump{4.222}$ et $\nump{4.2222}$. Certaines surfaces possèdent de plus des groupes du type $\nump{4.42}$, et parmi celles-ci, deux différentes possèdent chacune également un groupe de la liste ci-dessus, du type respectivement $\nump{4.223}$ ou $\nump{4.8}$.
\end{remarque}

\bigskip

Les surfaces de Del Pezzo de degré respectivement $3$, $2$ et $1$ sont obtenues comme hypersurfaces dans les espaces projectifs $\mathbb{P}^3$, $\mathbb{P}(2,1,1,1)$ et $\mathbb{P}(3,1,1,2)$; celles de degré $1$ et $2$ sont des revêtements doubles de respectivement $\Pn$ et $\mathbb{P}(1,1,2)$ et donc possèdent une involution associée au revêtement. Ces involutions sont appelées involutions de Geiser, respectivement de Bertini. En général, les groupes d'automorphismes de surfaces de Del Pezzo de degré $3$, $2$ et $1$ sont respectivement $\{0\}$, $\Z{2}$ et $\Z{2}$ (voir par exemple \cite{bib:Koi}), mais il y a tout de même plusieurs dizaines de surfaces différentes avec des groupes d'automorphismes plus grands. On décrit les sous-groupes abéliens finis minimaux dans tous les cas, en utilisant plusieurs outils, notamment l'action sur les diviseurs exceptionnels et la formule de Lefschetz.

Ce qui donne les résultats suivants:

\bigskip
\noindent {\bf Proposition \refd{Prp:Cubics}} \titreProp{automorphismes de surfaces cubiques}\\
{\it 
Soit $G$ un groupe abélien fini d'automorphismes d'une cubique non-singulière $S=V(F)\subset\mathbb{P}^3$ (surface de Del Pezzo de degré $3$), telle que $\rkPic{S}^{G}=1$.\\
Alors, à isomorphisme près, une des circonstances suivantes se produit:\vspace{0.1 cm}\\
\begin{tabular}{|llll|}
\hline
{\it nom} & {\it structure} & {\it générateurs} & {\it équation}\\
{\it de $\mathit (G,S)$} & {\it de $\mathit G$} & {\it de  $\mathit G$} & de $\mathit F$\\
\hline 
${\num{3.3}}$ & $\Z{3}$& $\DiaG{\omega}{1}{1}{1}$& $w^3+L_3(x,y,z) \hts$ \\
${\num{3.6.1}}$ & $\Z{6} $& $\DiaG{\omega}{1}{1}{-1}$& $w^3+x^3+y^3+xz^2+\lambda yz^2 \hts$ \\
${\num{3.33.1}}$ & $(\Z{3})^2$& $\DiaG{\omega}{1}{1}{1}$, $\DiaG{1}{1}{1}{\omega}$& $w^3+x^3+y^3+z^3\hts$ \\
${\num{3.9}}$ & $\Z{9}$& $\DiaG{\zeta_9}{1}{\omega}{\omega^2}$& $w^3+xz^2+x^2y+y^2z\hts$  \\
${\num{3.33.2}}$ & $(\Z{3})^2$& $\DiaG{\omega}{1}{1}{1}$, $\DiaG{1}{1}{\omega}{\omega^2}$& $w^3+x^3+y^3+z^3+\lambda xyz\hts$  \\
${\num{3.12}}$ & $\Z{12}$& $\DiaG{\omega}{1}{-1}{\im}$& $w^3+x^3+yz^2+y^2x \hts$ \\
${\num{3.36}}$ & $\Z{3}\times\Z{6}$& $\DiaG{\omega}{1}{1}{1}$, $\DiaG{1}{1}{-1}{\omega}$& $w^3+x^3+xy^2+z^3 \hts$ \\
${\num{3.333}}$ & $(\Z{3}) ^3$&$\DiaG{\omega}{1}{1}{1}$, $\DiaG{1}{\omega}{1}{1}$& $w^3+x^3+y^3+z^3\hts$ \\
& &\mbox{et } $\DiaG{1}{1}{\omega}{1}$ &  \\
${\num{3.6.2}}$ & $\Z{6}$ &$\DiaG{1}{-1}{\omega}{\omega^2}$ & $wx^2+w^3+y^3+z^3+\lambda wyz\hts$\\
\hline
\end{tabular}\vspace{0.1 cm}\\
où $L_3$ est une forme homogène non-singulière de degré~3, et $\lambda \in \K$ est un paramètre tel que la surface soit lisse, $\omega=e^{2\im\pi/3},\zeta_9=e^{2\im\pi/9}$. 

De plus, tous les cas ci-dessus sont des paires minimales $(G,S)$ avec $\rkPic{S}^G=1$.}

\bigskip

\noindent {\bf Proposition \refd{Prp:DelP2}}\titreProp{automorphismes de surfaces de Del Pezzo de degré $2$}\\
{\it Soit $S$ une surface de Del Pezzo de degré $2$ donnée par l'équation $w^2=F(x,y,z)$ et $G$ 
un sous-groupe abélien fini de $\Aut(S)$ tel que $\rkPic{S}^{G}=1$.\\
Alors, à isomorphisme près, une des circonstances suivantes se produit:
\begin{itemize}
\item
groupes contenant l'involution de Geiser $\sigma$:\\
$G=<\sigma>\times H$, où $H \subset
\PGLn{3}$ est un groupe abélien fini d'automorphismes de la quartique lisse $\Gamma \subset \Pn$ définie par l'équation $F(x,y,z)=0$:\\
\begin{tabular}{|llll|}
\hline
{\it nom} & {\it structure} & {\it générateurs} & {\it équation de}\\
{\it de $\mathit (G,S)$} & {\it de $\mathit G$} & {\it de $\mathit G$} & $\mathit F$\\
\hline 
${\num{2.G}}$ & $\Z{2}$& $\sigma$& $L_4(x,y,z)\hts$ \\
${\num{2.G2}}$ & $(\Z{2})^2$&$\DiaG{\pm 1}{1}{1}{-1}$& $L_4(x,y)+L_2(x,y)z^2+z^4\hts$ \\
${\num{2.G3.1}}$ & $\Z{6}$&$\DiaG{-1}{1}{1}{\omega}$& $L_4(x,y)+z^3L_1(x,y)\hts$ \\
${\num{2.G3.2}}$ & $\Z{6}$&$\DiaG{-1}{1}{\omega}{\omega^2}$& $x(x^3+y^3+z^3)+yzL_1(x^2,yz)\hts$ \\
${\num{2.G4.1}}$ & $\Z{2}\times\Z{4}$&$\DiaG{\pm 1}{1}{1}{\im}$& $L_4(x,y)+z^4\hts$ \\
${\num{2.G4.2}}$ & $\Z{2}\times\Z{4}$&$\DiaG{\pm 1}{1}{-1}{\im}$& $x^3y+y^4+z^4+xyL_1(xy,z^2)\hts$ \\
${\num{2.G6}}$ & $\Z{2}\times\Z{6}$&$\DiaG{\pm 1}{\omega}{1}{-1}$& $x^3y+y^4+z^4+\lambda y^2z^2\hts$ \\
${\num{2.G7}}$ & $\Z{14}$&$\DiaG{-1}{\zeta_7}{(\zeta_7)^4}{(\zeta_7)^2}$& $x^3y+y^3z+xz^3\hts$ \\
${\num{2.G8}}$ & $\Z{2}\times\Z{8}$&$\DiaG{\pm 1}{\zeta_8}{-\zeta_8}{1}$& $x^3y+xy^3+z^4\hts$ \\
${\num{2.G9}}$ & $\Z{18}$&$\DiaG{-1}{(\zeta_9)^6}{1}{\zeta_9}$& $x^3y+y^4+xz^3\hts$\\
${\num{2.G12}}$ & $\Z{2}\times\Z{12}$&$\DiaG{\pm 1}{\omega}{1}{\im}$& $x^3y+y^4+z^4\hts$ \\
${\num{2.G22}}$ & $(\Z{2})^3$&$\DiaG{\pm 1}{1}{\pm 1}{\pm 1}$ & $L_2(x^2,y^2,z^2)\hts$ \\
${\num{2.G24}}$ & $(\Z{2})^2\times \Z{4}$&$\DiaG{\pm 1}{1}{\pm 1}{\im}$ & $x^4+y^4+z^4+\lambda x^2y^2\hts$ \\
${\num{2.G44}}$ & $\Z{2}\times (\Z{4})^2$&$\DiaG{\pm 1}{1}{1}{\im}$& $x^4+y^4+z^4\hts$\\
& & and $\DiaG{1}{1}{\im}{1}$& \\
\hline
\end{tabular}
\item
groupes ne contenant pas l'involution de Geiser:\\
\begin{tabular}{|llll|}
\hline
{\it nom} & {\it structure} & {\it générateurs} & {\it équation de}\\
{\it de $\mathit (G,S)$} & {\it de $\mathit G$} & {\it de $\mathit G$} & $\mathit F$\\
\hline 
${\num{2.4}}$ & $\Z{4}$& $\DiaG{1}{1}{1}{\im}$& $L_4(x,y)+z^4\hts$ \\
${\num{2.6}}$ & $\Z{6}$& $\DiaG{-1}{\omega}{1}{-1}$& $x^3y+y^4+z^4+\lambda y^2z^2\hts$ \\
${\num{2.12}}$ & $\Z{12}$& $\DiaG{1}{\omega}{1}{\im}$& $x^3y+y^4+z^4\hts$ \\
${\num{2.24.1}}$ & $\Z{2}\times \Z{4}$& $\DiaG{1}{1}{1}{\im}$, $\DiaG{1}{1}{-1}{1}$& 
$x^4+y^4+z^4+\lambda x^2y^2\hts$ \\
${\num{2.24.2}}$ & $\Z{2}\times \Z{4}$& $\DiaG{1}{1}{1}{\im}$, $\DiaG{-1}{1}{-1}{1}$& 
$x^4+y^4+z^4+\lambda x^2y^2\hts$ \\
${\num{2.44.1}}$ & $(\Z{4})^2$& $\DiaG{1}{1}{1}{\im}$, $\DiaG{1}{1}{\im}{1}$& 
$x^4+y^4+z^4\hts$ \\
${\num{2.44.2}}$ & $(\Z{4})^2$& $\DiaG{1}{1}{1}{\im}$, $\DiaG{-1}{1}{\im}{1}$& 
$x^4+y^4+z^4\hts$\\
\hline
\end{tabular}
\end{itemize}
où $\lambda,\mu \in \K$ sont des paramètres et $L_i$ des formes homogènes de degré $i$, tels que les surfaces soient lisses, et $\zeta_n=e^{2\im\pi/n},\omega=\zeta_3$. 

De plus, tous les cas ci-dessus sont des paires minimales $(G,S)$ avec $\rkPic{S}^G=1$.}

\bigskip

\break
\noindent {\bf Proposition \refd{Prp:DelP1}}\titreProp{automorphismes de surfaces de Del Pezzo de degré $1$}\\
{\it 
Soit $S$ une surface de Del Pezzo de degré $1$ donnée par l'équation 
\begin{center}$w^2=z^3+F_4(x,y)z+F_6(x,y)$ \end{center}
et soit $G$ un sous-groupe abélien fini de $\Aut(S)$ tel que $\rkPic{S}^{G}=1$.\\
Alors, à isomorphisme près, une des circonstances suivantes se produit:\\
\begin{tabular}{|lllll|}
\hline
{\it nom} & {\it structure} & {\it générateurs} & {\it équation de}& {\it équation de}\\
{\it de $\mathit (G,S)$} & {\it de $\mathit G$} & {\it de $\mathit G$} & $\mathit F_4$& $\mathit F_6$\\
\hline 
${\num{1.B}}$ & $\Z{2}$& $\sigma$& $L_4(x,y)$ & $L_6(x,y)\hts$ \\
${\num{1.\rho}}$ & $\Z{3}$& $\rho$& $0$ &$L_6(x,y)\hts$ \\
${\num{1.\sigma\rho}}$ & $\Z{6}$& $\sigma\rho$& $0$ & $L_6(x,y)\hts$ \\
${\num{1.B2.1}}$ & $(\Z{2})^2$&$\sigma$, $\DiaG{1}{1}{-1}{1}$& $L_2(x^2,y^2)$ & $L_3(x^2,y^2)\hts$ \\
${\num{1.\sigma\rho 2.1}}$ & $\Z{6}\times \Z{2}$& $\sigma\rho$, $\DiaG{1}{1}{-1}{1}$& $0$ & $L_3(x^2,y^2)\hts$ \\
${\num{1.\rho 2}}$ & $\Z{6}$& $\DiaG{1}{1}{-1}{\omega}$& $0$ & $L_3(x^2,y^2)\hts$ \\
${\num{1.B2.2}}$ & $\Z{4}$& $\DiaG{\im}{1}{-1}{-1}$& $L_2(x^2,y^2)$ & $xyL_2'(x^2,y^2)\hts$ \\
${\num{1.\sigma\rho 2.2}}$ & $\Z{12}$& $\DiaG{\im}{1}{-1}{-\omega}$& $0$ & $xyL_2(x^2,y^2)\hts$ \\
${\num{1.B3.1}}$ & $\Z{6}$&$\DiaG{-1}{1}{\omega}{1}$& $xL_1(x^3,y^3)$ & $L_2(x^3,y^3)\hts$ \\
${\num{1.\sigma\rho 3}}$ & $\Z{6}\times\Z{3}$& $\sigma\rho$, $\DiaG{1}{1}{\omega}{1}$& $0$ &  $L_2(x^3,y^3)\hts$ \\
${\num{1.\rho 3}}$ & $ (\Z{3})^2$& $\rho$, $\DiaG{1}{1}{\omega}{1}$& $0$ & $L_2(x^3,y^3)\hts$ \\
${\num{1.B3.2}}$ & $\Z{6}$& $\DiaG{-1}{1}{\omega}{\omega}$& $\lambda x^2 y^2$ & $L_2(x^3,y^3)\hts$ \\
${\num{1.B4.1}}$ & $\Z{2}\times\Z{4}$&$\sigma$, $\DiaG{1}{1}{\im}{1}$& $L_1(x^4,y^4)$ &  $x^2L_1'(x^4+y^4)\hts$ \\
${\num{1.B4.2}}$ & $\Z{8}$&$\DiaG{\zeta_8}{1}{\im}{-\im}$& $\lambda x^2y^2$ & $xy(x^4+y^4)\hts$ \\
${\num{1.\sigma\rho 4}}$ & $\Z{24}$& $\DiaG{\zeta_{8}}{1}{\im}{-\im\omega}$& $0$& $xy(x^4+y^4)\hts$ \\
${\num{1.5}}$ & $\Z{5}$&$\DiaG{1}{1}{\zeta_5}{1}$& $\lambda x^4$ & $x(\mu x^5+y^5)\hts$ \\
${\num{1.B5}}$ & $\Z{10}$&$\DiaG{-1}{1}{\zeta_5}{1}$& $\lambda x^4$ & $x(\mu x^5+y^5)\hts$ \\
${\num{1.\sigma\rho 5}}$ & $\Z{30}$& $\DiaG{-1}{1}{\zeta_5}{\omega}$ & $0$ & $x(x^5+y^5)\hts$ \\
${\num{1.\rho 5}}$ & $\Z{15}$& $\DiaG{1}{1}{\zeta_5}{\omega}$& $0$ & $x(x^5+y^5)\hts$ \\
${\num{1.6}}$ & $\Z{6}$&$\DiaG{1}{1}{-\omega}{1}$& $\lambda x^4$ & $\mu x^6+y^6\hts$ \\
${\num{1.B6.1}}$ & $\Z{2}\times\Z{6}$&$\sigma$, $\DiaG{1}{1}{-\omega}{1}$& $\lambda x^4$ & $\mu x^6+y^6\hts$ \\
${\num{1.\sigma\rho 6}}$ & $ (\Z{6})^2$& $\sigma\rho$, $\DiaG{1}{1}{-\omega}{1}$&$0$& $x^6+y^6\hts$ \\
${\num{1.\rho 6}}$ & $\Z{3}\times \Z{6}$& $\rho$, $\DiaG{1}{1}{-\omega}{1}$&$0$& $x^6+y^6\hts$ \\
${\num{1.B6.2}}$ & $\Z{2}\times\Z{6}$&$\sigma$, $\DiaG{1}{1}{-\omega}{\omega}$& $\lambda x^2y^2$ & $x^6+y^6\hts$ \\
${\num{1.B10}}$ & $\Z{20}$&$\DiaG{\im}{1}{\zeta_{10}}{-1}$& $x^4$& $xy^5\hts$ \\
${\num{1.B12}}$ & $\Z{2}\times\Z{12}$&$\sigma, \DiaG{\im}{1}{\zeta_{12}}{-1}$& $x^4$ & $y^6\hts$ \\
\hline
\end{tabular}\vspace{0.1 cm}\\
o\`u $L_i$ est une forme homogène de degré $i$, $\lambda,\mu \in \K$ sont des paramètres et les équations sont telles que $S$ est lisse.

De plus, tous les cas ci-dessus sont des paires minimales $(G,S)$ avec $\rkPic{S}^G=1$.}
\section*{Chapitre 9 - Conjugaison entre les cas}
\titrepageresume{9}
Aux chapitres précédents, nous avons énuméré toutes les possibilités pour les sous-groupes abéliens finis du groupe de Cremona. Nous devons maintenant décider si les groupes trouvés représentent des classes de conjugaison différentes ou non. 

Dans le cas de la proposition \refd{Prp:TwoCases} où $\rkPic{S}^G=1$, on classifie les groupes ainsi:

\bigskip
\noindent {\bf Proposition \refd{Prp:ConjBetweenRankPic}}\fauxtitre 
{\it 
Soit $G\subset \Aut(S)$ un groupe abélien fini d'automorphismes d'une surface rationnelle $S$, et supposons que $\rkPic{S}^G=1$. Alors, à conjugaison birationnelle près, une et une seule des situations suivantes se produit:
\begin{enumerate}
\item[\upshape 1.]
$G$ préserve un pinceau de courbes rationnelles de $S$.
\item[\upshape 2.]
$S=\Pn$ et $G=V_9\cong (\Z{3})^2$ est engendré par  $(x:y:z)\mapsto (x:\omega y:\omega^2 z)$ et $(x:y:z) \mapsto (y:z:x)$.
\item[\upshape 3.]
$S$ est une surface de Del Pezzo de degré $4$ et la paire $(G,S)$ est isomorphe à $\nump{4.222}$, $\nump{4.2222}$, $\nump{4.42}$ ou $\nump{4.8}$
\item[\upshape 4.]
$S$ est une surface de Del Pezzo de degré $3$.
\item[\upshape 5.]
$S$ est une surface de Del Pezzo de degré $2$.
\item[\upshape 6.]
$S$ est une surface de Del Pezzo de degré $1$.
\end{enumerate}}
\begin{remarque}\fauxtitred
L'assertion $1$ est équivalente à ce que le groupe soit birationnellement conjugué à un groupe d'automorphismes de fibré en coniques (cas $2$ de la proposition $\refd{Prp:TwoCases}$). La proposition $\refd{Prp:ConjBetweenRankPic}$ prouve donc non seulement que les cas $2.$ à $6.$ sont birationnellement distincts, mais qu'ils sont aussi distincts du cas des automorphismes de fibrés en coniques, que nous étudions ensuite.
\end{remarque}

\bigskip

Dans le cas de la proposition \refd{Prp:TwoCases} où $\rkPic{S}^G=2$, on a classé les groupes à partir de leur présentation par la suite exacte (\refd{eq:ExactSeqCB}) (en fonction de la structure de $G$ et de $G'$). Nous montrons dans quelle mesure le groupe $G$ peut possèder plusieurs représentations de type (\refd{eq:ExactSeqCB}).

\bigskip
\noindent {\bf Proposition \refd{Prp:HardChgExacrSeq}}\titreProp{Changements possibles de la suite exacte}\\
{\it Soit $G\subset \Aut(S,\pi)$ un groupe abélien fini d'automorphismes du fibré en coniques $(S,\pi)$.
\begin{enumerate}
\item[\upshape 1.]
Si $G$ contient au moins deux éléments non triviaux qui fixent une courbe de genre positif, alors $G'\cong (\Z{2})^2$.
\item[\upshape 2.]
Si $G$ contient exactement un élément non trivial qui fixe une courbe de genre positif, celui-ci appartient à $G'$ et  $G'\cong \Z{2}$ ou $G'\cong (\Z{2})^2$. 

Si $G'\cong (\Z{2})^2$, alors:
 $G=G'$ si et seulement si $G$ est birationnellement conjugué à un groupe $H \subset \Aut(\tilde{S},\tilde{\pi})$ avec $H'\cong \Z{2}$.
\item[\upshape 3.]
Si $G$ ne contient pas d'élément non trivial qui fixe une courbe de genre positif, il est birationnellement conjugué à un sous-groupe de $\Aut(\mathbb{P}^1\times\mathbb{P}^1)$, ou au groupe $\nump{Cs.24}$.
\end{enumerate}}

\section*{Chapitre 10 - La classification et autres résultats intéressants}
En utilisant les résultats de tous les chapitres précédents, on finit par donner une classification des sous-groupes abéliens finis du groupe de Cremona. Celle-ci étant très longue, nous préférons ne pas la mettre à nouveau ici.

Celle-ci est la liste la plus précise donnée jusqu'à présent. Les classes de conjugaison de certains cas avec paramètres ne sont pas totalement décrites, mais la liste permet toutefois de démontrer les conséquences citées au début du travail.
\section*{Chapitre 11 - Appendices}
\markboths{Résumé de la thèse}{Résumé de la thèse: Chapitres 10,11}
Ce chapitre regroupe quelques résultats annexes utilisés au cours de la thèse. On y retrouve notamment des résultats classiques sur les groupes d'automorphismes de courbes, ainsi que des résultats sur les sommes de nombres réels ou entiers. 
\clearpage
\addcontentsline{toc}{chapter}{Index \& Bibliography}
\printindex

\end{document}